\documentclass[leqno,openany,draft]{book}

\usepackage{makeidx}
\makeindex

\def\supp{\mathop{\rm supp}}
\def\re{\mathop{\rm Re}}
\def\im{\mathop{\rm Im}}
\def\hom{\mathop{\rm Hom}}
\def\id{\mathop{\rm id}}

\newtheorem{theorem}{Theorem}
\newtheorem{lemma}[theorem]{Lemma}
\newtheorem{proposition}[theorem]{Proposition}
\newtheorem{definition}[theorem]{Definition}
\newtheorem{corollary}[theorem]{Corollary}

\newcommand{\begintheorem}{\addtocounter{equation}{1}\begin{theorem}}
\newcommand{\beginlemma}{\addtocounter{equation}{1}\begin{lemma}}
\newcommand{\beginproposition}{\addtocounter{equation}{1}\begin{proposition}}
\newcommand{\begindefinition}{\addtocounter{equation}{1}\begin{definition}}
\newcommand{\begincorollary}{\addtocounter{equation}{1}\begin{corollary}}

\begin{document}

\frontmatter

\title{Topics in Fourier Analysis}

\author{Stephen Semmes \\
        Rice University}

\date{}

\maketitle

\chapter*{Preface}

        Of course, the most classical versions of Fourier analysis
deal with functions on the circle and the real line, as well as their
higher-dimensional counterparts.  At the same time, there are more
``fractal'' types of situations, such as infinite products of cyclic
groups, $p$-adic numbers, and solenoids.  An overview of some of the
relevant notions is given here, with examples like these especially
in mind.  In particular, although many cases have a lot of features
in common, there is also some nice variety in other ways.

\tableofcontents

\mainmatter

\chapter{Preliminaries}
\label{preliminaries}

\section{Real and complex numbers}
\label{real, complex numbers}

        Remember that the \emph{absolute value}\index{absolute value}
$|x|$ of a real number $x$ is defined to be equal to $x$ when $x \ge 0$
and to $-x$ when $x \le 0$.  Thus $|x| \ge 0$ for every $x$ in the real
line\index{real line} ${\bf R}$,\index{R@${\bf R}$} and $|x| = 0$ if and 
only if $x = 0$.  It is easy to see that
\begin{equation}
\label{|x + y| le |x| + |y|}
        |x + y| \le |x| + |y|
\end{equation}
and
\begin{equation}
\label{|x y| = |x| |y|}
        |x \, y| = |x| \, |y|
\end{equation}
for every $x, y \in {\bf R}$.

        A complex number\index{complex numbers} $z$ can be expressed as 
$x + y \, i$, where $x$ and $y$ are real numbers, and $i^2 = -1$.  
More precisely, $x$ and $y$ are known as the \emph{real}\index{real
  part} and \emph{imaginary} parts\index{imaginary part} of $z$, which
may be denoted $\re z$ and $\im z$, respectively.  The \emph{complex
  conjugate}\index{complex conjugates} of $z$ is defined by
\begin{equation}
\label{overline{z} = x - y i}
        \overline{z} = x - y \, i,
\end{equation}
and it is easy to see that
\begin{equation}
\label{overline{z + w} = overline{z} + overline{w}}
        \overline{z + w} = \overline{z} + \overline{w}
\end{equation}
and
\begin{equation}
\label{overline{z w} = overline{z} overline{w}}
        \overline{z \, w} = \overline{z} \, \overline{w}
\end{equation}
for every $z$, $w$ in the complex plane ${\bf C}$.\index{C@${\bf C}$}

        The \emph{modulus}\index{modulus} of a complex number
$z = x + y \, i$, $x, y \in {\bf R}$, is defined by
\begin{equation}
\label{|z| = (x^2 + y^2)^{1/2}}
        |z| = (x^2 + y^2)^{1/2}.
\end{equation}
Note that this reduces to the absolute value of $z$ when $y = 0$, and
that $|z|$ is the same as the standard Euclidean norm of $(x, y) \in
{\bf R}^2$ for every $z \in {\bf C}$.  In particular, it is well known
that
\begin{equation}
\label{|z + w| le |z| + |w|}
        |z + w| \le |z| + |w|
\end{equation}
for every $z, w \in {\bf C}$.  One can also check that
\begin{equation}
\label{|z|^2 = z overline{z}}
        |z|^2 = z \, \overline{z}
\end{equation}
for every $z \in {\bf C}$, which implies that
\begin{equation}
\label{|z w| = |z| |w|}
        |z \, w| = |z| \, |w|
\end{equation}
for every $z, w \in {\bf C}$, by (\ref{overline{z w} = overline{z}
  overline{w}}).

\section{The complex exponential function}
\label{complex exponential function}

        The complex exponential function\index{complex exponential 
function}\index{exponential function} is defined for each $z \in {\bf C}$
by
\begin{equation}
\label{exp z = sum_{j = 0}^infty frac{z^j}{j!}}
        \exp z = \sum_{j = 0}^\infty \frac{z^j}{j!},
\end{equation}
where $z^j$ is interpreted as being equal to $1$ for each $z$ when $j
= 0$.  As usual, $j!$ is ``$j$ factorial'', the product of the integer
from $1$ to $j$, and interpreted as being equal to $1$ when $j = 0$.
It is easy to see that this series converges absolutely for each $z
\in {\bf C}$, using the ratio test or by comparison with geometric
series.  This implies that the partial sums converge uniformly on
bounded subsets of ${\bf C}$, and hence that $\exp z$ is a continuous
function on ${\bf C}$.

        If $w$ is another complex number, then
\begin{equation}
\label{(exp z) (exp w) = ...}
 (\exp z) \, (\exp w) = \Big(\sum_{j = 0}^\infty \frac{z^j}{j!}\Big) \,
                      \Big(\sum_{k = 0}^\infty \frac{w^k}{k!}\Big)
  = \sum_{n = 0}^\infty \Big(\sum_{j = 0}^n
                            \frac{z^j \, w^{n - j}}{j! \, (n - j)!}\Big),
\end{equation}
using the standard expression for the product of two infinite series
in terms of Cauchy products in the second step.  More precisely, it is
well known that the absolute convergence of the original series for
the exponential function implies the absolute convergence of the
Cauchy product series, and that the sum of the Cauchy product series
is the product of the original series.  The binomial theorem implies that
\begin{equation}
\label{(w + z)^n = sum_{j = 0}^n frac{n!}{j! (n - j)!} z^j w^{n - j}}
        (w + z)^n = \sum_{j = 0}^n \frac{n!}{j! \, (n - j)!} \, z^j \, w^{n - j},
\end{equation}
for every $z, w \in {\bf C}$, so that
\begin{equation}
\label{(exp z) (exp w) = exp (z + w)}
        (\exp z) \, (\exp w) = \exp (z + w).
\end{equation}

        If we take $w = -z$ in (\ref{(exp z) (exp w) = exp (z + w)}),
then we get that
\begin{equation}
\label{(exp z) (exp (-z)) = exp 0 = 1}
        (\exp z) \, (\exp (-z)) = \exp 0 = 1
\end{equation}
for every $z \in {\bf C}$.  In particular, $\exp z \ne 0$ for each $z
\in {\bf C}$.  Note that $\exp x$ is a real number when $x$ lies in
the real line ${\bf R}$.  It is easy to see that $\exp x \ge 1$ when
$x \ge 0$, and indeed that $\exp x$ is strictly increasing on the set
of nonnegative real numbers.  Since $\exp (-x) = 1/(\exp x)$, as in
(\ref{(exp z) (exp (-z)) = exp 0 = 1}), we get that $0 < \exp x \le 1$
when $x \le 0$, and that $\exp x$ is strictly increasing on all of
${\bf R}$.

        Observe that
\begin{equation}
\label{overline{exp z} = exp overline{z}}
        \overline{\exp z} = \exp \overline{z}
\end{equation}
for every $z \in {\bf C}$, because of (\ref{overline{z + w} =
  overline{z} + overline{w}}) and (\ref{overline{z w} = overline{z}
  overline{w}}).  Similarly,
\begin{equation}
\label{|exp z|^2 = (exp z) (exp overline{z}) = ... = exp (2 re z)}
        |\exp z|^2 = (\exp z) \, (\exp \overline{z}) = \exp (z + \overline{z})
                                                      = \exp (2 \, \re z),
\end{equation}
for every $z \in {\bf C}$, because of (\ref{overline{z w} =
  overline{z} overline{w}}) and (\ref{(exp z) (exp w) = exp (z + w)}).

\section{Metrics and norms}
\label{metrics, norms}

       A \emph{metric space}\index{metric spaces} is a set $M$ together
with a nonnegative real-valued function $d(x, y)$ defined for $x, y \in M$
such that $d(x, y) = 0$ if and only if $x = y$,
\begin{equation}
\label{d(x, y) = d(y, x)}
        d(x, y) = d(y, x)
\end{equation}
for every $x, y \in M$, and
\begin{equation}
\label{d(x, z) le d(x, y) + d(y, z)}
        d(x, z) \le d(x, y) + d(y, z)
\end{equation}
for every $x, y, z \in M$.  Similarly, a \emph{norm}\index{norms} on a
vector space $V$ over the real or complex numbers is a nonnegative
real-valued function $\|v\|$ on $V$ such that $\|v\| = 0$ if and only
if $v = 0$,
\begin{equation}
\label{||t v|| = |t| ||v||}
        \|t \, v\| = |t| \, \|v\|
\end{equation}
for every $v \in V$ and $t \in {\bf R}$ or ${\bf C}$, as appropriate, and
\begin{equation}
\label{||v + w|| le ||v|| + ||w||}
        \|v + w\| \le \|v\| + \|w\|
\end{equation}
for every $v, w \in V$.  Here $|t|$ denotes the absolute value of $t
\in {\bf R}$ in the real case, and the modulus of $t \in {\bf C}$ in
the complex case.  If $\|v\|$ is a norm on $V$, then it is easy to see
that
\begin{equation}
\label{d(v, w) = ||v - w||}
        d(v, w) = \|v - w\|
\end{equation}
defines a metric on $V$.

        Of course, the real line may be considered as a $1$-dimensional real
vector space, and the absolute value function defines a norm on ${\bf R}$.  
The complex plane may also be considered as a $1$-dimensional complex vector
space, and the modulus defines a norm on ${\bf C}$.  The corresponding metrics
are the standard Euclidean metrics on ${\bf R}$ and ${\bf C}$, respectively.

        Let $V$ be a real or complex vector space with a norm $\|v\|$
again.  An infinite series $\sum_{j = 1}^\infty v_j$ with terms $v_j
\in V$ is said to \emph{converge}\index{convergent series} in $V$ with
respect to a norm $\|v\|$ if the corresponding sequence of partial
sums $\sum_{j = 1}^n v_j$ converges in $V$ as $n \to \infty$ with
respect to the metric (\ref{d(v, w) = ||v - w||}) associated to the
norm, in which case the sum $\sum_{j = 1}^\infty v_j$ is defined to be
the limit of the partial sums.  Similarly, $\sum_{j = 1}^\infty v_j$
is said to converge \emph{absolutely}\index{absolute convergence} with
respect to $\|v\|$ if
\begin{equation}
\label{sum_{j = 1}^infty ||v_j||}
        \sum_{j = 1}^\infty \|v_j\|
\end{equation}
converges as an infinite series of nonnegative real numbers.  As in
the classical situation of real or complex numbers, one can check that
the sequence of partial sums $\sum_{j = 1}^n v_j$ of an absolutely
convrergent series in $V$ is a Cauchy sequence with respect to the
associated metric (\ref{d(v, w) = ||v - w||}), using the triangle
inequality.  If $V$ is complete as a metric space with respect to
(\ref{d(v, w) = ||v - w||}), so that every Cauchy sequence of elements
of $V$ converges to an element of $V$, then $V$ is said to be a
\emph{Banach space}\index{Banach spaces} with respect to the norm
$\|v\|$.  Thus an absolutely convergent series $\sum_{j = 1}^\infty
v_j$ in a Banach space $V$ converges in $V$, and it is easy to see
that
\begin{equation}
\label{||sum_{j = 1}^infty v_j|| le sum_{j = 1}^infty ||v_j||}
        \biggl\|\sum_{j = 1}^\infty v_j\biggr\| \le \sum_{j = 1}^\infty \|v_j\|,
\end{equation}
by standard arguments.

\section{Inner product spaces}
\label{inner product spaces}

        Let $V$ be a vector space over the real or complex numbers.
An \emph{inner product}\index{inner products} on $V$ is a real or
complex-valued function $\langle v, w \rangle$, as appropriate,
defined for $v, w \in V$, and satisfying the following properties.
First, $\langle v, w \rangle$ is a linear function of $v$ for each
$w \in W$.  Second,
\begin{equation}
\label{langle w, v rangle = langle v, w rangle}
        \langle w, v \rangle = \langle v, w \rangle
\end{equation}
for every $v, w \in V$ in the real case, and
\begin{equation}
\label{langle w, v rangle = overline{langle v, w rangle}}
        \langle w, v \rangle = \overline{\langle v, w \rangle}
\end{equation}
for every $v, w \in V$ in the complex case.  It follows that $\langle
v, w \rangle$ is a linear function of $w$ for each $v \in V$ in the
real case, and that $\langle v, w \rangle$ is conjugate-linear in $w$
in the complex case.  In the complex case, (\ref{langle w, v rangle =
  overline{langle v, w rangle}}) also implies that $\langle v, v
\rangle$ is a real number for every $v \in V$.  The third and last
condition is that
\begin{equation}
\label{langle v, v rangle > 0}
        \langle v, v \rangle > 0
\end{equation}
for every $v \in V$ with $v \ne 0$.  Of course, $\langle v, v \rangle
= 0$ when $v = 0$, because of the first condition about linearity.

        If $\langle v, w \rangle$ is an inner product on $V$, then
we put
\begin{equation}
\label{||v|| = langle v, v rangle^{1/2}}
        \|v\| = \langle v, v \rangle^{1/2}
\end{equation}
for every $v \in V$.  The \emph{Cauchy--Schwarz
  inequality}\index{Cauchy--Schwarz inequality} states that
\begin{equation}
\label{|langle v, w rangle| le ||v|| ||w||}
        |\langle v, w \rangle| \le \|v\| \, \|w\|
\end{equation}
for every $v, w \in V$.  As usual, this can be shown using the fact
that
\begin{equation}
\label{langle v + t w, v + t w rangle ge 0}
        \langle v + t \, w, v + t \, w \rangle \ge 0
\end{equation}
for every $t \in {\bf R}$ or ${\bf C}$, as appropriate.  Using the
Cauchy--Schwarz inequality, one can check that
\begin{equation}
 \|v + w\|^2 \le \|v\|^2 + 2 \, \|v\| \, \|w\| + \|w\|^2 = (\|v\| + \|w\|)^2
\end{equation}
for each $v, w \in V$, so that $\|\cdot \|$ satisfies the triangle
inequality (\ref{||v + w|| le ||v|| + ||w||}).  This implies that
$\|v\|$ defines a norm on $V$, since the positivity and homogeneity
conditions for a norm follow directly from the definition of an inner
product.  In particular, (\ref{d(v, w) = ||v - w||}) defines a metric
on $V$, as in the previous section.  If $V$ is complete as a metric
space with respect to (\ref{d(v, w) = ||v - w||}), then $V$ is said to
be a \emph{Hilbert space}.\index{Hilbert spaces}

        Let $n$ be a positive integer, and let ${\bf R}^n$ and ${\bf C}^n$
be the usual spaces of $n$-tuples of real and complex numbers, respectively.
Remember that these are vector spaces with respect to coordinatewise
addition and scalar multiplication.  The standard inner product on
${\bf R}^n$ is defined by
\begin{equation}
\label{langle v, w rangle = sum_{j = 1}^n v_j w_j}
        \langle v, w \rangle = \sum_{j = 1}^n v_j \, w_j,
\end{equation}
and the standard inner product on ${\bf C}^n$ is defined by
\begin{equation}
\label{langle v, w rangle = sum_{j = 1}^n v_j overline{w_j}}
        \langle v, w \rangle = \sum_{j = 1}^n v_j \, \overline{w_j}.
\end{equation}
It is easy to see that these do define inner products on ${\bf R}^n$
and ${\bf C}^n$, for which the corresponding norm
\begin{equation}
\label{||v|| = (sum_{j = 1}^n |v_j|^2)^{1/2}}
        \|v\| = \Big(\sum_{j = 1}^n |v_j|^2\Big)^{1/2}
\end{equation}
is the standard Euclidean norm.

\section{Orthogonal vectors}
\label{orthogonal vectors}

        Let $V$ be a real or complex vector space, and let 
$\langle v, w \rangle$ be an inner product on $V$.  A pair of
vectors $v, w \in V$ are said to be \emph{orthogonal}\index{orthogonal
vectors} if
\begin{equation}
\label{langle v, w rangle = 0}
        \langle v, w \rangle = 0,
\end{equation}
which may be also expressed symbolically by $v \perp w$.  
If $v$ is orthogonal to $w$ in $V$, then it is easy to see that
\begin{equation}
\label{||v + w||^2 = ||v||^2 + ||w||^2}
        \|v + w\|^2 = \|v\|^2 + \|w\|^2.
\end{equation}
More precisely, (\ref{||v + w||^2 = ||v||^2 + ||w||^2}) is equivalent
to the orthgonality of $v$ and $w$ in the real case, and in the
complex case (\ref{||v + w||^2 = ||v||^2 + ||w||^2}) holds if and only
if $\re \langle v, w \rangle = 0$.  Note that a complex vector space
$V$ may be considered as a real vector space by forgetting about
scalar multiplication by $i$, and that the real part of an inner
product on $V$ as a complex vector space is an inner product on $V$ as
a real inner product space, which determines the same norm on $V$.

        Suppose that $v_1, \ldots, v_n$ are finitely many 
orthonormal\index{orthonormal vectors} vectors in $V$, so that $v_j
\perp v_l$ when $j \ne l$ and $\|v_j\| = 1$ for each $j$.  Let $W$
be the linear subspace of $V$ spanned by $v_1, \ldots, v_n$, and put
\begin{equation}
\label{P_W(v) = sum_{j = 1}^n langle v, v_j rangle v_j}
        P_W(v) = \sum_{j = 1}^n \langle v, v_j \rangle \, v_j
\end{equation}
for each $v \in V$.  Thus $P_W$ defines a linear mapping from $V$ into
$W$, and it is easy to see that $P_W(v_l) = v_l$ for each $l = 1,
\ldots, n$, which implies that $P_W(w) = w$ for every $w \in W$.
Similarly,
\begin{equation}
\label{langle P_W(v), v_l rangle = langle v, v_l rangle}
        \langle P_W(v), v_l\rangle = \langle v, v_l \rangle,
\end{equation}
for each $v \in V$ and $1 \le l \le n$, so that $v - P_W(v)$ is
orthogonal to $v_l$ for each $l$.  This implies that $(v - P_W(v))
\perp w$ for every $v \in V$ and $w \in W$.

        Conversely, let $v \in V$ be given, and suppose that $u \in W$
has the property that $(v - u) \perp w$ for every $w \in W$.  Thus
$P_W(v) - u \in W$, and 
\begin{equation}
\label{P_W(v) - u = (v - u) - (v - P_W(v))}
        P_W(v) - u = (v - u) - (v - P_W(v))
\end{equation}
is also orthogonal to every $w \in W$.  If we apply this to $w =
P_W(v) - u$, then we get that $P_W(v) - u$ is orthogonal to itself,
and hence is equal to $0$.  This shows that $P_W(v)$ is uniquely
determined by the conditions that $P_W(v) \in W$ and $(v - P_W(v))
\perp w$ for every $w \in W$.  In particular, $P_W(v)$ depends only on
$W$, and not on the choice of orthonormal basis $v_1, \ldots, v_n$ of
$W$.

        Note that
\begin{equation}
\label{||v||^2 = ||P_W(v)||^2 + ||v - P_W(v)||^2}
        \|v\|^2 = \|P_W(v)\|^2 + \|v - P_W(v)\|^2
\end{equation}
for each $v \in V$, because $P_W(v) \in W$ and $v - P_W(v)$ is
orthogonal to every element of $W$, so that $P_W(v) \perp (v - P_W(v))$.
This implies that
\begin{equation}
\label{||v||^2 = sum_{j = 1}^n |langle v, v_j rangle|^2 + ||v - P_W(v)||^2}
        \|v\|^2 = \sum_{j = 1}^n |\langle v, v_j \rangle|^2 + \|v - P_W(v)\|^2,
\end{equation}
and hence that
\begin{equation}
\label{sum_{j = 1}^n |langle v, v_j rangle|^2 le ||v||^2}
        \sum_{j = 1}^n |\langle v, v_j \rangle|^2 \le \|v\|^2.
\end{equation}
If $w \in W$, then $(v - P_W(v))$ is also orthogonal to $P_W(v) - w
\in W$, so that
\begin{eqnarray}
\label{||v - w||^2 = ... = ||v - P_W(v)||^2 + ||P_W(v) - w||^2}
        \|v - w\|^2 & = & \|(v - P_W(v)) - (P_W(v) - w)\|^2  \\
                    & = & \|v - P_W(v)\|^2 + \|P_W(v) - w\|^2. \nonumber
\end{eqnarray}
It follows that
\begin{equation}
\label{||v - w|| ge ||v - P_W(v)||}
        \|v - w\| \ge \|v - P_W(v)\|
\end{equation}
for each $w \in W$, and that equality holds in (\ref{||v - w|| ge ||v
  - P_W(v)||}) if and only if $w = P_W(v)$.  Thus $P_W(v)$ minimizes
the distance to $v$ among elements of $W$, and is uniquely determined
by this property.

\section{Orthogonal sequences}
\label{orthogonal sequences}

        Let $V$ be a real or complex vector space with an inner product
$\langle v, w \rangle$ again, and let $v_1, v_2, v_3, \ldots$ be a sequence
of pairwise-orthogonal vectors in $V$.  Note that
\begin{equation}
\label{||sum_{j = 1}^n v_j||^2 = sum_{j = 1}^n ||v_j||^2}
        \biggl\|\sum_{j = 1}^n v_j\biggr\|^2 = \sum_{j = 1}^n \|v_j\|^2
\end{equation}
for each positive integer $n$, as in (\ref{||v + w||^2 = ||v||^2 +
  ||w||^2}).  Thus the partial sums $\sum_{j = 1}^n v_j$ of $\sum_{j =
  1}^\infty v_j$ have bounded norm in $V$ if and only if
\begin{equation}
\label{sum_{j = 1}^infty ||v_j||^2}
        \sum_{j = 1}^\infty \|v_j\|^2
\end{equation}
converges as an infinite series of nonnegative real numbers.  In this
case, one can check that the partial sums of $\sum_{j = 1}^\infty v_j$
form a Cauchy sequence in $V$ with respect to the metric associated to
the norm.  If $V$ is complete, then it follows that $\sum_{j =
  1}^\infty v_j$ converges in $V$, and one can also check that
\begin{equation}
\label{||sum_{j = 1}^infty v_j||^2 = sum_{j = 1}^infty ||v_j||^2}
 \biggl\|\sum_{j = 1}^\infty v_j\biggr\|^2 = \sum_{j = 1}^\infty \|v_j\|^2.
\end{equation}

        Now let $v_1, v_2, v_3, \ldots$ be an infinite sequence of
orthonormal vectors in $V$, and let $W_n$ be the linear span of
$v_1, \ldots, v_n$ for each positive integer $n$.  Put
\begin{equation}
\label{P_n(v) = P_{W_n}(v) = sum_{j = 1}^n langle v, v_j rangle v_j}
        P_n(v) = P_{W_n}(v) = \sum_{j = 1}^n \langle v, v_j \rangle \, v_j
\end{equation}
for each $n$, so that
\begin{equation}
\label{||v||^2 = sum_{j = 1}^n |langle v, v_j rangle|^2 + ||v - P_n(v)||^2}
 \|v\|^2 = \sum_{j = 1}^n |\langle v, v_j \rangle|^2 + \|v - P_n(v)\|^2
\end{equation}
for each $n$, as in (\ref{||v||^2 = sum_{j = 1}^n |langle v, v_j
  rangle|^2 + ||v - P_W(v)||^2}).  Thus (\ref{sum_{j = 1}^n |langle v,
  v_j rangle|^2 le ||v||^2}) holds for each $n$, which implies that
$\sum_{j = 1}^\infty |\langle v, v_j \rangle|^2$ converges as an
infinite series of nonnegative real numbers, and satisfies
\begin{equation}
\label{sum_{j = 1}^infty |langle v, v_j rangle|^2 le ||v||^2}
        \sum_{j = 1}^\infty |\langle v, v_j \rangle|^2 \le \|v\|^2.
\end{equation}
If $V$ is complete, then it follows that
\begin{equation}
\label{sum_{j = 1}^infty langle v, v_j rangle v_j}
        \sum_{j = 1}^\infty \langle v, v_j \rangle \, v_j
\end{equation}
converges in $V$, as before.

        Observe that $\bigcup_{n = 1}^\infty W_n$ is a linear subspace
of $V$, which is the linear span of the $v_j$'s.  Let $W$ be the closure
of $\bigcup_{n = 1}^\infty W_n$ in $V$ with respect to the metric associated 
to the norm, which is a closed linear subspace of $V$.  Equivalently,
$v \in V$ is an element of $W$ if and only if
\begin{equation}
\label{lim_{n to infty} ||v - P_n(v)|| = 0}
        \lim_{n \to \infty} \|v - P_n(v)\| = 0,
\end{equation}
because $P_n(v)$ minimizes the distance from $v$ to $W_n$ for each
$n$, as in (\ref{||v - w|| ge ||v - P_W(v)||}), and $W_n \subseteq
W_{n + 1}$ for each $n$.  This says exactly that (\ref{sum_{j =
    1}^infty langle v, v_j rangle v_j}) converges and is equal to $v$
when $v \in W$, which works whether or not $V$ is complete.  

        Suppose that (\ref{sum_{j = 1}^infty langle v, v_j rangle v_j})
converges for some $v \in V$, and let the sum be denoted $P(v)$.
Thus $P(v) \in W$, and
\begin{equation}
\label{||P(v)||^2 = sum_{j = 1}^infty |langle v, v_j rangle|^2 le ||v||^2}
        \|P(v)\|^2 = \sum_{j = 1}^\infty |\langle v, v_j \rangle|^2 \le \|v\|^2.
\end{equation}
If $l$ is any positive integer, then
\begin{equation}
\label{langle P(v), v_l rangle = ... = langle v, v_l rangle}
 \langle P(v), v_l \rangle = \lim_{n \to \infty} \langle P_n(v), v_l \rangle
                            = \langle v, v_l \rangle,
\end{equation}
as in (\ref{langle P_W(v), v_l rangle = langle v, v_l rangle}).  This
implies that $(v - P(v)) \perp v_l$ for each $l$, and hence that $(v -
P(v)) \perp w$ for every $w \in \bigcup_{n = 1}^\infty W_n$, and
indeed for every $w \in W$.  As in the previous section, $P(v)$ is
uniquely determined by the conditions that $P(v) \in W$ and $(v -
P(v)) \perp w$ for every $w \in W$.  These conditions also imply that
the analogues of (\ref{||v - w||^2 = ... = ||v - P_W(v)||^2 + ||P_W(v)
  - w||^2}) and (\ref{||v - w|| ge ||v - P_W(v)||}) with $P_W(v)$
replaced by $P(v)$ hold for each $w \in W$, as before.  Hence $P(v)$
minimizes the distance from $v$ to $W$, and is also uniquely
determined by this property.  If $V$ is complete, then $P(v)$ is
defined for every $v \in V$, and defines a linear mapping from $V$
onto $W$.

\section{Minimizing distances}
\label{minimizing distances}

        Let $V$ be a real or complex vector space with an inner
product $\langle v, w \rangle$.  The \emph{parallelogram 
law}\index{parallelogram law} states that
\begin{equation}
\label{parallelogram law}
 \biggl\|\frac{x + y}{2}\biggr\|^2 + \biggl\|\frac{x - y}{2}\biggr\|^2
                           = \frac{\|x\|^2}{2} + \frac{\|y\|^2}{2}
\end{equation}
for every $x, y \in V$.  This is easy to check, by expanding the norms
in terms of the inner products, and observing that the cross terms
cancel.  Let $W$ be a linear subspace of $V$, let $v$ be an element of
$W$, and put
\begin{equation}
\label{r = inf {||v - w|| : w in W}}
        r = \inf \{\|v - w\| : w \in W\}.
\end{equation}
Also let $z_j$ be an element of $W$ such that
\begin{equation}
\label{||v - z_j|| < r + 1/j}
        \|v - z_j\| < r + 1/j
\end{equation}
for each positive integer $j$.  Applying the parallelogram law to
$x = v - z_j$ and $y = v - z_k$, we get that
\begin{equation}
\label{||v - frac{z_j + z_k}{2}||^2 + ||frac{z_j - z_k}{2}||^2 = ...}
        \biggl\|v - \frac{z_j + z_k}{2}\biggr\|^2 + 
                                        \biggl\|\frac{z_j - z_k}{2}\biggr\|^2
                    = \frac{\|v - z_j\|^2}{2} + \frac{\|v - z_k\|^2}{2}
\end{equation}
for every $j, k \ge 1$.  Because $(z_j + z_k)/2 \in W$, we have that
\begin{equation}
\label{||v - frac{z_j + z_k}{2}|| ge r}
        \biggl\|v - \frac{z_j + z_k}{2}\biggr\| \ge r
\end{equation}
for each $j$ and $k$, and hence
\begin{equation}
\label{r^2 + ||frac{z_j - z_k}{2}||^2 < (r + 1/j)^2/2 + (r + 1/k)^2/2}
        r^2 + \biggl\|\frac{z_j - z_k}{2}\biggr\|^2 
                         < \frac{(r + 1/j)^2}{2} + \frac{(r + 1/k)^2}{2},
\end{equation}
using also (\ref{||v - z_j|| < r + 1/j}).  This implies that
\begin{equation}
\label{||frac{z_j - z_k}{2}||^2 < r (1/j + 1/k) + (1/2) (1/j^2 + 1/k^2) to 0}
        \biggl\|\frac{z_j - z_k}{2}\biggr\|^2 
                   < r \, (1/j + 1/k) + (1/2) \, (1/j^2 + 1/k^2) \to 0
\end{equation}
as $j, k \to \infty$, so that $\{z_j\}_{j = 1}^\infty$ is a Cauchy sequence
in $V$.

        Suppose now that $V$ is complete and that $W$ is a closed
linear subspace of $V$, so that $\{z_j\}_{j = 1}^\infty$ converges to an
element $z$ of $W$.  By construction,
\begin{equation}
\label{||v - z|| = r}
        \|v - z\| = r.
\end{equation}
If $w \in W$ and $t \in {\bf R}$ or ${\bf C}$, as appropriate, then
$z - t \, w \in W$ too, and hence
\begin{equation}
\label{||v - z||^2 le ||v - z + t w||^2}
        \|v - z\|^2 \le \|v - z + t \, w\|^2,
\end{equation}
by the definition (\ref{r = inf {||v - w|| : w in W}}) of $r$.
In the real case, this implies that
\begin{equation}
\label{||v - z||^2 le ||v - z||^2 + 2 t langle v - z, w rangle + t^2 ||w||^}
        \|v - z\|^2 \le \|v - z\|^2 + 2 \, t \, \langle v - z, w \rangle
                                                           + t^2 \, \|w\|^2,
\end{equation}
and in the complex case, we get that
\begin{equation}
\label{||v - z||^2 le ||v - z||^2 + 2 re t (v - z, w) + |t|^2 ||w||^2}
        \|v - z\|^2 \le \|v - z\|^2 + 2 \, \re t \, \langle v - z, w \rangle
                                                     + |t|^2 \, \|w\|^2.
\end{equation}
In both cases, the minimum of the right side is attained at $t = 0$,
and one can use this to show that
\begin{equation}
\label{langle v - z, w rangle = 0}
        \langle v - z, w \rangle = 0
\end{equation}
for every $w \in W$.

        As in Section \ref{orthogonal vectors}, $z$ is uniquely determined
by the conditions that $z \in W$ and $(v - z) \perp w$ for every $w \in W$.
If we put $P_W(v) = z$, then it is easy to see that $P_W$ defines a
linear mapping from $V$ into $W$, because these conditions
characterizing $P_W(v)$ are linear in $v$.  Note that $P_W(v)$ is the
same as in in Section \ref{orthogonal vectors} when $W$ is
finite-dimensional, and that $P_W(v)$ is the same as $P(v)$ in Section
\ref{orthogonal sequences} when $W$ is the closure of the linear span
of an orthonormal sequence in $V$.  If $W$ is infinite-dimensional and
separable, in the sense that $W$ has a countable dense subset, then
one can use the Gram--Schmidt process to get an orthonormal sequence
of vectors in $W$ whose linear span is dense in $W$.  As in the previous
situations, $P_W(v) = v$ when $v \in W$, and
\begin{equation}
\label{||P_W(v)|| le ||v||}
        \|P_W(v)\| \le \|v\|
\end{equation}
for every $v \in V$.

        The \emph{orthogonal complement}\index{orthogonal complement}
$W^\perp$\index{W^perp@$W^\perp$} of $W$ is defined by
\begin{equation}
\label{W^perp = {y in V : y perp w for every w in W}}
        W^\perp = \{y \in V : y \perp w \hbox{ for every } w \in W\},
\end{equation}
which is automatically a closed linear subspace of $V$.  The previous
discussion implies that every element of $V$ has a unique
representation as a sum of elements of $W$ and $W^\perp$ when $V$ is
complete and $W$ is a closed linear subspace of $V$.

\section{Summable functions}
\label{summable functions}

        Let $E$ be a nonempty set, and let $f$ be a nonnegative real-valued
function on $E$.  If $A$ is a nonempty finite subset of $E$, then the
sum
\begin{equation}
\label{sum_{x in A} f(x)}
        \sum_{x \in A} f(x)
\end{equation}
of $f$ over the elements of $A$ can be defined in the usual way.  The sum
\begin{equation}
\label{sum_{x in E} f(x)}
        \sum_{x \in E} f(x)
\end{equation}
is defined as the supremum of the finite subsums (\ref{sum_{x in A}
  f(x)}) over all nonempty finite subsets $A$ of $E$.  More precisely,
$f$ is said to be \emph{summable}\index{summable functions} on $E$ if
the finite subsums (\ref{sum_{x in A} f(x)}) have a finite upper bound
in ${\bf R}$, and otherwise (\ref{sum_{x in E} f(x)}) is interpreted as
being $+\infty$.

        If $a$ is a nonnegative real number, then it is easy to see that
\begin{equation}
\label{sum_{x in E} a f(x) = a sum_{x in E} f(x)}
        \sum_{x \in E} a \, f(x) = a \, \sum_{x \in E} f(x),
\end{equation}
with the convention that $0 \cdot (+\infty) = 0$ when $a = 0$
and $f$ is not summable on $E$.  Similarly, if $g$ is another nonnegative
real-valued function on $E$, then one can check that
\begin{equation}
\label{sum_{x in E} (f(x) + g(x)) = sum_{x in E} f(x) + sum_{x in E} g(x)}
        \sum_{x \in E} (f(x) + g(x)) = \sum_{x \in E} f(x) + \sum_{x \in E} g(x),
\end{equation}
with the usual conventions that $b + (+\infty) = (+\infty) + b = +\infty$
for every real number $b$ and $(+\infty) + (+\infty) = +\infty$.  In
particular, if $f$ and $g$ are both summable functions on $E$, then
$f + g$ is summable too.

        Let $\epsilon > 0$ be given, and put
\begin{equation}
\label{E(f, epsilon) = {x in E : f(x) ge epsilon}}
        E(f, \epsilon) = \{x \in E : f(x) \ge \epsilon\}.
\end{equation}
If $f$ is summable on $E$, then $E(f, \epsilon)$ has only finitely
many elements, and in fact
\begin{equation}
\label{epsilon (num E(f, epsilon)) le sum_{x in E} f(x)}
        \epsilon \, (\# E(f, \epsilon)) \le \sum_{x \in E} f(x),
\end{equation}
where $\# A$ is the number of elements of a finite set $A$.  
It follows that
\begin{equation}
\label{{x in E : f(x) > 0} = bigcup_{n = 1}^infty E(f, 1/n)}
        \{x \in E : f(x) > 0\} = \bigcup_{n = 1}^\infty E(f, 1/n)
\end{equation}
has only finitely or countably many elements when $f$ is summable on
$E$.

        Now let $f$ be a real or complex-valued function on $E$,
and let us say that $f$ is \emph{summable}\index{summable functions}
on $E$ if $|f(x)|$ is summable on $E$.  Let 
$\ell^1(E)$\index{l^1(E)@$\ell^1(E)$} be the space
of all summable functions on $E$, which may also be denoted
$\ell^1(E, {\bf R})$ or $\ell^1(E, {\bf C})$, to indicate whether
real or complex valued functions are being considered.  It is easy
to see that $\ell^1(E)$ is a vector space with respect to pointwise
addition and scalar multiplication, and that
\begin{equation}
\label{||f||_1 = sum_{x in E} |f(x)|}
        \|f\|_1 = \sum_{x \in E} |f(x)|
\end{equation}
defines a norm on $\ell^1(E)$.

        There are a couple of equivalent ways in which to define the sum
(\ref{sum_{x in E} f(x)}) of a summable real or complex-valued function
$f$ on $E$.  One way is to express $f$ as a linear combination of
nonnegative real-valued summable functions, and apply the previous
definition to those.  Another way is to use the fact that $f(x) \ne 0$
for only finitely or countably many $x \in E$, by applying
(\ref{{x in E : f(x) > 0} = bigcup_{n = 1}^infty E(f, 1/n)}) to $|f(x)|$.
This permits the sum (\ref{sum_{x in E} f(x)}) to be reduced to either
a finite sum or an absolutely convergent infinite series.  This also
uses the fact that absolutely convergent series are invariant under
rearrangements, so that the definition of (\ref{sum_{x in E} f(x)})
does not depend on the way that the $x \in E$ with $f(x) \ne 0$ are
listed in a sequence.

        In both of these approaches, one can check that the sum
(\ref{sum_{x in E} f(x)}) is linear in $f$, and satisfies
\begin{equation}
\label{|sum_{x in E} f(x)| le sum_{x in E} |f(x)|}
        \Bigl|\sum_{x \in E} f(x)\Bigr| \le \sum_{x \in E} |f(x)|,
\end{equation}
where the right side is defined as before.  Of course, the definition
of the sum (\ref{sum_{x in E} f(x)}) is trivial when $f(x) = 0$ for
all but finitely many $x \in E$, and one can show that the sum
(\ref{sum_{x in E} f(x)}) for summable real or complex-valued
functions on $E$ is characterized by these properties.  This is
because the functions $f$ on $E$ with $f(x) = 0$ for all but finitely
many $x \in E$ form a dense linear subspace of $\ell^1(E)$ with
respect to the metric associated to the $\ell^1$ norm $\|f\|_1$,
and the sum (\ref{sum_{x in E} f(x)}) defines a uniformly continuous
function on $\ell^1(E)$, since it is linear in $f$ and satisfies
(\ref{|sum_{x in E} f(x)| le sum_{x in E} |f(x)|}).  One can also
use these properties to give another approach to the definition of
the sum (\ref{sum_{x in E} f(x)}), because any uniformly continuous
mapping from a dense subset of a metric space $M$ into a complete metric
space $N$ can be extended to a uniformly continuous mapping from $M$ into $N$.

        Suppose that $\{f_j\}_{j = 1}^\infty$ is a Cauchy sequence in
$\ell^1(E)$ with respect to the $\ell^1$ metric.  This implies that
$\{f_j(x)\}_{j = 1}^\infty$ is a Cauchy sequence in ${\bf R}$ or ${\bf C}$,
as appropriate, for each $x \in E$.  Because the real and complex numbers
are complete with respect to their standard metrics, it follows that
$\{f_j(x)\}_{j = 1}^\infty$ converges for each $x \in E$.  If the
limit is denoted $f(x)$, then one can check that $f(x)$ is also
summable on $E$, and that $\{f_j\}_{j = 1}^\infty$ converges to $f$
with respect to the $\ell^1$ norm, using the Cauchy condition for
$\{f_j\}_{j = 1}^\infty$ with respect to the $\ell^1$ norm.  It follows that
$\ell^1(E)$ is complete, and hence a Banach space.

\section{$p$-Summability}
\label{p-summability}

        Let $E$ be a nonempty set, and let $p$ be a positive real number.
A real or complex-valued function $f$ on $E$ is said to be 
\emph{$p$-summable}\index{p-summable functions@$p$-summable functions}
if $|f(x)|^p$ is a summable function on $E$, in which case we put
\begin{equation}
\label{||f||_p = (sum_{x in E} |f(x)|^p)^{1/p}}
        \|f\|_p = \Big(\sum_{x \in E} |f(x)|^p\Big)^{1/p}.
\end{equation}
Let $\ell^p(E)$\index{l^p(E)@$\ell^p(E)$} be the space of $p$-summable
functions on $E$, which may also be denoted $\ell^p(E, {\bf R})$ or
$\ell^p(E, {\bf C})$, to indicate whether real or complex-valued
functions are being considered.  As a substitute for $p = \infty$, let
$\ell^\infty(E)$ be the space of bounded real or complex-valued
functions on $E$, which may also be denoted $\ell^\infty(E, {\bf R})$
or $\ell^\infty(E, {\bf C})$ to indicate whether real or
complex-valued functions are being used.  
If $f \in \ell^\infty(E)$, then we put
\begin{equation}
\label{||f||_infty = sup_{x in E} |f(x)|}
        \|f\|_\infty = \sup_{x \in E} |f(x)|.
\end{equation}

        If $f$ is $p$-summable on $E$ for any $p > 0$, then 
$|f(x)| \le \|f\|_p$ for each $x \in E$, so that $f$ is bounded on $E$
and
\begin{equation}
\label{||f||_infty le ||f||_p}
        \|f\|_\infty \le \|f\|_p.
\end{equation}
Similarly, if $p < q < \infty$, then $f$ is $q$-summable and
\begin{equation}
\label{||f||_q le ||f||_p}
        \|f\|_q \le \|f\|_p.
\end{equation}
To see this, observe that
\begin{equation}
\label{|f(x)|^q le ||f||_infty^{q - p} |f(x)|^p le ||f||_p^{q - p} |f(x)|^p}
 |f(x)|^q \le \|f\|_\infty^{q - p} \, |f(x)|^p \le \|f\|_p^{q - p} \, |f(x)|^p
\end{equation}
for each $x \in E$, so that
\begin{equation}
\label{||f||_q^q = ... le ||f||_p^{q - p} ||f||_p^p = ||f||_p^q}
 \|f\|_q^q = \sum_{x \in E} |f(x)|^q \le \|f\|_p^{q - p} \sum_{x \in E} |f(x)|^p
                                    \le \|f\|_p^{q - p} \, \|f\|_p^p = \|f\|_p^q.
\end{equation}

        If $p$ is a positive real number, $f$ is a $p$-summable function
on $E$, and $a$ is a real or complex number, as appropriate, then it is
easy to see that $a \, f(x)$ is also $p$-summable on $E$, and that
\begin{equation}
\label{||a f||_p = |a| ||f||_p}
        \|a \, f\|_p = |a| \, \|f\|_p.
\end{equation}
Similarly, if $f$ and $g$ are $p$-summable functions on $E$, then one can
check that $f + g$ is also $p$-summable, because
\begin{eqnarray}
\label{|f(x) + g(x)|^p le ... le 2^p (|f(x)|^p + |g(x)|^p)}
 |f(x) + g(x)|^p & \le & (|f(x)| + |g(x)|)^p \le (2 \max(|f(x)|, |g(x)|))^p \\
 & = & 2^p \max(|f(x)|^p, |g(x)|^p) \le 2^p \, (|f(x)|^p + |g(x)|^p) \nonumber
\end{eqnarray}
for every $x \in E$.  Thus $\ell^p(E)$ is a vector space with respect to
pointwise addition and scalar multiplication when $0 < p < \infty$,
which can also be verified directly from the definitions when $p = \infty$.

        As in the previous section, $\|f\|_1$ defines a norm on $\ell^1(E)$,
and it is easy to see that $\|f\|_\infty$ defines a norm on $\ell^\infty(E)$.
It is well known that $\|f\|_p$ also defines a norm on $\ell^p(E)$ when
$1 < p < \infty$.  One way to show this is to check that the unit ball
in $\ell^p(E)$ is convex when $p \ge 1$, using the convexity of the
function $r^p$ on the set of nonnegative real numbers.  If $0 < p < 1$
and $E$ has at least two elements, then the unit ball in $\ell^p(E)$ is
not convex, and $\|f\|_p$ does not define a norm on $\ell^p(E)$.
There is however an alternative to this, as follows.

        If $0 < p \le 1$ and $a$, $b$ are nonnegative real numbers, then
\begin{equation}
\label{(a + b)^p le a^p + b^p}
        (a + b)^p \le a^p + b^p.
\end{equation}
This can be derived from (\ref{||f||_q le ||f||_p}) with $q = 1$, by
taking $E$ to be a set with exactly two elements.  If $f$ and $g$ are 
$p$-summable functions on any set $E$ and $0 < p \le 1$, then we get that
\begin{equation}
\label{|f(x) + g(x)|^p le (|f(x)| + |g(x)|)^p le |f(x)|^p + |g(x)|^p}
        |f(x) + g(x)|^p \le (|f(x)| + |g(x)|)^p \le |f(x)|^p + |g(x)|^p
\end{equation}
for each $x \in E$, and hence
\begin{equation}
\label{||f + g||_p^p le ||f||_p^p + ||g||_p^p}
        \|f + g\|_p^p \le \|f\|_p^p + \|g\|_p^p.
\end{equation}
This implies that
\begin{equation}
\label{d_p(f, g) = ||f - g||_p^p}
        d_p(f, g) = \|f - g\|_p^p
\end{equation}
defines a metric on $\ell^p(E)$ when $0 < p \le 1$, in place of
(\ref{d(v, w) = ||v - w||}).

        Suppose that $\{f_j\}_{j = 1}^\infty$ is a Cauchy sequence
in $\ell^p(E)$ for some $p$, $0 < p \le \infty$, with respect to the
metric $\|f - g\|_p$ associated to the $\ell^p$ norm when $p \ge 1$,
and with respect to the metric (\ref{d_p(f, g) = ||f - g||_p^p})
when $0 < p \le 1$.  In both cases, $\{f_j(x)\}_{j = 1}^\infty$ is a 
Cauchy sequence in ${\bf R}$ or ${\bf C}$ for each $x \in E$, as
appropriate, and hence converges to a real or complex number $f(x)$.
Using the Cauchy condition with respect to the $\ell^p$ metric,
one can show that $f \in \ell^p(E)$ too, and that $\{f_j\}_{j = 1}^\infty$
converges to $f$ with respect to the $\ell^p$ metric.  Thus $\ell^p(E)$
is complete with respect to the $\ell^p$ metric for each $p > 0$,
and hence $\ell^p(E)$ is a Banach space when $p \ge 1$.

        A real or complex-valued function $f$ on $E$ is said to
\emph{vanish at infinity}\index{vanishing at infinity} if
\begin{equation}
\label{E(f, epsilon) = {x in E : |f(x)| ge epsilon}}
        E(f, \epsilon) = \{x \in E : |f(x)| \ge \epsilon\}
\end{equation}
is a finite set for each $\epsilon > 0$.  Let
$c_0(E)$\index{c_0(E)@$c_0(E)$} be the space of functions on $E$ that
vanish at infinity, which may also be denoted $c_0(E, {\bf R})$ or
$c_0(E, {\bf C})$, to indicate whether real or complex-valued
functions are being used.  It is well known and not difficult to check
that $c_0(E)$ is a closed linear subspace of $\ell^\infty(E)$.
If $f$ is a $p$-summable function on $E$ for some $p$, $0 < p < \infty$,
then it is easy to see that $f \in c_0(E)$, with
\begin{equation}
\label{epsilon^p (num E(f, epsilon))^p le ||f||_p^p}
        \epsilon^p \, (\# E(f, \epsilon))^p \le \|f\|_p^p
\end{equation}
for each $\epsilon > 0$.  Note that $f(x) = 0$ for all but finitely or
countably many elements of $E$ when $f$ is any function that vanishes
at infinity on $E$, as in the previous section.

\section{Square-summability}
\label{square-summability}

        Let $E$ be a nonempty set again, and let us restrict our attention
now to $p = 2$.  If $a$, $b$ are nonnegative real numbers, then
\begin{equation}
\label{a b le max(a^2, b^2) le a^2 + b^2}
        a \, b \le \max(a^2, b^2) \le a^2 + b^2,
\end{equation}
and in fact $2 \, a \, b \le a^2 + b^2$, since $(a - b)^2 \ge 0$.  If
$f$, $g$ are $2$-summable functions on $E$, then we can apply either
of these inequalities to $a = |f(x)|$ and $b = |g(x)|$, to conclude
that the product $|f(x)| \, |g(x)|$ is summable on $E$.  Put
\begin{equation}
\label{langle f, g rangle = sum_{x in E} f(x) g(x)}
        \langle f, g \rangle = \sum_{x \in E} f(x) \, g(x)
\end{equation}
in the real case, and
\begin{equation}
\label{langle f, g rangle = sum_{x in E} f(x) overline{g(x)}}
        \langle f, g \rangle = \sum_{x \in E} f(x) \, \overline{g(x)}
\end{equation}
in the complex case.  It is easy to see that these define inner products
on $\ell^2(E, {\bf R})$ and $\ell^2(E, {\bf C})$, respectively, for which
the corresponding norms are equal to the $\ell^2$ norm discussed in
the previous section.

        Now let $V$ be a real or complex vector space with an inner
product $\langle v, w \rangle_V$, and let $\|v\|_V$ be the correspponding
norm on $V$.  Also let $A$ be a nonempty set, and suppose that for each
$\alpha \in A$ we have a vector $v_\alpha \in V$ such that $\|v_\alpha\|_V = 1$
and $v_\alpha \perp v_\beta$ when $\alpha \ne \beta$.  If $v \in V$ and
$\alpha_1, \ldots, \alpha_n$ are finitely many distinct elements of $A$, then
\begin{equation}
\label{sum_{j = 1}^n |langle v, v_{alpha_j} rangle_V|^2 le ||v||_V^2}
        \sum_{j = 1}^n |\langle v, v_{\alpha_j} \rangle_V|^2 \le \|v\|_V^2,
\end{equation}
as in (\ref{sum_{j = 1}^n |langle v, v_j rangle|^2 le ||v||^2}).  Thus
\begin{equation}
\label{f_v(alpha) = langle v, v_alpha rangle_V}
        f_v(\alpha) = \langle v, v_\alpha \rangle_V
\end{equation}
is a $2$-summable function on $A$, with
\begin{equation}
\label{sum_{alpha in A} |f_v(alpha)|^2 le ||v||_V^2}
        \sum_{\alpha \in A} |f_v(\alpha)|^2 \le \|v\|_V^2.
\end{equation}
Note that the mapping from $v \in V$ to $f_v \in \ell^2(A)$ is linear.

        Let $f$ be any $2$-summable function on $A$, which is real or 
complex-valued depending on whether $V$ is real or complex.  If $V$ is
complete, then we would like to define
\begin{equation}
\label{T(f) = sum_{alpha in A} f(alpha) v_alpha}
        T(f) = \sum_{\alpha \in A} f(\alpha) \, v_\alpha
\end{equation}
as an element of $V$.  This reduces to an ordinary finite sum when
$f(\alpha) = 0$ for all but finitely many $\alpha \in A$, and
otherwise $f(\alpha) = 0$ for all but countably many $\alpha \in A$,
because $f$ is $2$-summable on $A$.  In the latter case, the sum may
be considered as an infinite series, as in Section \ref{orthogonal
sequences}, and one can also show that the sum does not depend on the
way that the terms are listed, as for absolutely convergent series.
Of course, if $A$ has only finitely or countably many elements,
then one can use the same listing of elements of $A$ for every
$f \in \ell^2(A)$.

        It is easy to see that the mapping from $f \in \ell^2(A)$ to 
$T(f) \in V$ is linear.  As in (\ref{||sum_{j = 1}^infty v_j||^2 = 
sum_{j = 1}^infty ||v_j||^2}),
\begin{equation}
\label{||T(f)||_V = ||f||_{ell^2(A)}}
        \|T(f)\|_V = \|f\|_{\ell^2(A)}
\end{equation}
for every $f \in \ell^2(A)$, where $\|f\|_{\ell^2(A)} = \|f\|_2$ is
the $\ell^2$ norm of $f$ on $A$.  Similarly,
\begin{equation}
\label{langle T(f), T(g) rangle_V = langle f, g rangle_{ell^2(A)}}
        \langle T(f), T(g) \rangle_V = \langle f, g \rangle_{\ell^2(A)}
\end{equation}
for every $f, g \in \ell^2(A)$, where $\langle f, g
\rangle_{\ell^2(A)}$ is the $\ell^2$ inner product on $A$ defined
earlier.  The mapping $f \mapsto T(f)$ may be characterized as the
unique linear mapping from $\ell^2(A)$ into $V$ that satisfies
(\ref{||T(f)||_V = ||f||_{ell^2(A)}}) and agrees with (\ref{T(f) =
  sum_{alpha in A} f(alpha) v_alpha}) when $f(\alpha) = 0$ for all but
finitely many $\alpha \in A$, since the $f \in \ell^2(A)$ with
$f(\alpha) = 0$ for all but finitely many $\alpha \in A$ are dense in
$\ell^2(A)$.  Alternatively, one can start with $T(f)$ as in
(\ref{T(f) = sum_{alpha in A} f(alpha) v_alpha}) when $f(\alpha) = 0$
for all but finitely many $\alpha \in A$, and show that this has a
unique extension to a linear mapping from $\ell^2(A)$ into $V$ that
satisfies (\ref{||T(f)||_V = ||f||_{ell^2(A)}}) when $V$ is complete.

        Let $W$ be the set of vectors in $V$ of the form $T(f)$ for some
$f \in \ell^2(A)$.  Thus $W$ is a linear subspace of $V$, $W$ contains
the linear span of the $v_\alpha$'s in $V$, and $W$ is contained in
the closure of the linear span of the $v_\alpha$'s in $V$.  Using the
completeness of $\ell^2(A)$ and (\ref{||T(f)||_V = ||f||_{ell^2(A)}}),
one can show that $W$ is complete with respect to the restriction of
the norm $\|\cdot \|_V$ to $W$, which implies that $W$ is a closed set
in $V$.  Hence $W$ is equal to the closure of the linear span of the
$v_\alpha$'s in $V$.  In particular, $W = V$ when the linear span of
the $v_\alpha$'s is dense in $V$.

        If $v \in V$, then we can  define $f_v \in \ell^2(A)$ as in
(\ref{f_v(alpha) = langle v, v_alpha rangle_V}), and then consider
$T(f_v) \in V$ as in (\ref{T(f) = sum_{alpha in A} f(alpha) v_alpha}).
It is easy to see that $T(f_v)$ is the same as the orthogonal
projection of $v$ onto $W$, as in Sections \ref{orthogonal vectors}
and \ref{orthogonal sequences}.  In particular, $v = T(f_v)$ for every
$v \in W$, and for every $v \in V$ when the linear span of the
$v_\alpha$'s is dense in $V$.  Similarly, if $f \in \ell^2(A)$, then
\begin{equation}
        \langle T(f), v_\beta \rangle_V = f(\beta)
\end{equation}
for each $\beta \in A$, so that $f = f_v$ with $v = T(f)$.

        An orthonormal family of vectors $\{v_\alpha\}_{\alpha \in A}$
in a Hilbert space $V$ is said to be an \emph{orthonormal
  basis}\index{orthonormal bases} for $V$ if the linear span of the
$v_\alpha$'s is dense in $V$.  Thus the mapping $T$ defines an
isometric linear isomorphism from $\ell^2(A)$ onto $V$ in this case.
If $V$ is separable, then one can use the Gram--Schmidt process to get
an orthonormal basis in $V$ with only finitely or countably many elements.

\section{Bounded linear functionals}
\label{bounded linear functionals}

        Let $V$ be a real or complex vector space with a norm $\|v\|$.
Remember that a \emph{linear functional}\index{linear functionals}
on $V$ is a linear mapping from $V$ into the real or complex numbers,
as appropriate.  A linear functional $\lambda$ on $V$ is said to be
\emph{bounded}\index{bounded linear functionals} if
\begin{equation}
\label{|lambda(v)| le C ||v||}
        |\lambda(v)| \le C \, \|v\|
\end{equation}
for some nonnegative real number $C$ and every $v \in V$.  This implies that
\begin{equation}
        |\lambda(v) - \lambda(w)| = |\lambda(v - w)| \le C \, \|v - w\|
\end{equation}
for every $v, w \in V$, so that $\lambda$ is uniformly continuous with
respect to the metric on $V$ associated to the norm $\|\cdot \|$.
Conversely, if a linear functional $\lambda$ on $V$ is continuous at
$0$, then there is a $\delta > 0$ such that $|\lambda(v)| < 1$ for
every $v \in V$ with $\|v\| < \delta$, and one can check that
(\ref{|lambda(v)| le C ||v||}) holds with $C = 1/\delta$.

        It is easy to see that the dual space\index{dual spaces}
$V^*$\index{V^*@$V^*$} of bounded linear functionals on $V$ is also 
a real or complex vector space, as appropriate, with respect to
pointwise addition and scalar multiplication.  
If $\lambda \in V^*$, then put
\begin{equation}
\label{||lambda||_* = sup {|lambda(v)| : v in V, ||v|| le 1}}
        \|\lambda\|_* = \sup \{|\lambda(v)| : v \in V, \ \|v\| \le 1\},
\end{equation}
which is the same as the smallest $C \ge 0$ for which
(\ref{|lambda(v)| le C ||v||}) holds.  One can check that this defines
a norm on $V^*$, known as the \emph{dual norm}\index{dual norms}
associated to the norm $\|\cdot \|$ on $V$.  It is well known that
$V^*$ is automatically complete with respect to the dual norm.  
As usual, one can start by showing that a Cauchy sequence in $V^*$
converges pointwise to a linear functional on $V$, and then use the
Cauchy condition with respect to the dual norm to show that the limit
is a bounded linear functional on $V$, and that the sequence converges
to this limit with respect to the dual norm.

        Suppose that $\langle v, w \rangle$ is an inner product on $V$,
and that $\|v\|$ is the norm on $V$ corresponding to this inner product.
If $w \in V$, then
\begin{equation}
\label{lambda_w(v) = langle v, w rangle}
        \lambda_w(v) = \langle v, w \rangle
\end{equation}
defines a linear functional on $V$, which is bounded by the
Cauchy--Schwarz inequality.  More precisely, the Cauchy--Schwarz
inequality implies that the dual norm of $\lambda_w$ is less than or
equal to $\|w\|$, and one can check that the dual norm of $\lambda_w$
is equal to $w$, since
\begin{equation}
\label{lambda_w(w) = ||w||^2}
        \lambda_w(w) = \|w\|^2.
\end{equation}
Conversely, if $V$ is a Hilbert space, then one can show that every
bounded linear functional $\lambda$ on $V$ is of the form $\lambda_w$
for some $w \in W$.  One way to do this is to look at vectors in $V$
that are orthogonal to the kernel of $\lambda$, using orthogonal
projections.

        Now let $E$ be a nonempty set, and suppose that 
$1 \le p, q \le \infty$ are conjugate exponents, in the sense that
\begin{equation}
\label{frac{1}{p} + frac{1}{q} = 1}
        \frac{1}{p} + \frac{1}{q} = 1,
\end{equation}
where $1/\infty = 0$, as usual.  If $f \in \ell^p(E)$ and $g \in
\ell^q(E)$, then \emph{H\"older's inequality}\index{Holder's
  inequality@H\"older's inequality} states that $f(x) \, g(x)$ is a
summable function on $E$, and that
\begin{equation}
\label{|sum_{x in E} f(x) g(x)| le ||f||_p ||g||_q}
        \Bigl|\sum_{x \in E} f(x) \, g(x) \Bigr| \le \|f\|_p \, \|g\|_q.
\end{equation}
This is very simple when $p = 1$ and $q = \infty$ or $q = 1$ and $p =
\infty$, and it follows from the Cauchy--Schwarz inequality when
$p = q = 2$.  This implies that
\begin{equation}
        \lambda_g(f) = \sum_{x \in E} f(x) \, g(x)
\end{equation}
defines a bounded linear functional on $\ell^p(E)$ when $g \in
\ell^q(E)$, with dual norm less than or equal to $\|g\|_q$.  One can
check that the dual norm of $\lambda_g$ on $\ell^p(E)$ is actually
equal to $\|g\|_q$, by considering suitable choices of $f \in \ell^p(E)$.

        Conversely, if $\lambda$ is a bounded linear functional on $\ell^p(E)$
and $1 \le p < \infty$, then one can show that $\lambda = \lambda_g$ for some
$g \in \ell^q(E)$.  More precisely, for each $y \in E$, let $\delta_y(x)$
be the function on $E$ equal to $1$ when $x = y$ and to $0$ when $x \ne y$.
Thus $\delta_y \in \ell^p(E)$, and one can define $g$ on $E$ by
\begin{equation}
\label{g(y) = lambda(delta_y)}
        g(y) = \lambda(\delta_y).
\end{equation}
Using the boundedness of $\lambda$ on $\ell^p(E)$, one can show that
$g \in \ell^q(E)$, and that $\lambda_g(f) = \lambda(f)$ for every
$f \in \ell^p(E)$.  This also uses the fact that the functions $f$
on $E$ such that $f(x) = 0$ for all but finitely many $x \in E$
are dense in $\ell^p(E)$, which works when $p < \infty$ and not when
$p = \infty$.

        As an alternative for $p = \infty$, one can consider bounded
linear functionals on $c_0(E)$, with respect to the $\ell^\infty$ norm
on $c_0(E)$.  If $g \in \ell^1(E)$, then $\lambda_g$ defines a bounded
linear functional on $\ell^\infty(E)$ as before, and hence its restriction
to $c_0(E)$ is a bounded linear functional on $c_0(E)$.  One can also
check that the dual norm of the restriction of $\lambda_g$ to $c_0(E)$
is equal to $\|g\|_1$.  Conversely, one can show that every bounded
linear functional on $c_0(E)$ is of this form, by the same type of argument
as for $\ell^p(E)$ when $1 \le p < \infty$.  This uses the fact that
the functions $f$ on $E$ with $f(x) = 0$ for all but finitely many $x
\in E$ are dense in $c_0(E)$.

        Of course, $\ell^2(E)$ is a Hilbert space, and the
bounded linear functionals on $\ell^2(E)$ can also be described as in
that case.  The two descriptions are equivalent, even if they
are presented in slightly different ways in the complex case.

        If $V$ is any real or complex vector space with a norm $\|\cdot \|$,
$v \in V$, and $v \ne 0$, then the theorem of Hahn and Banach\index{Hahn--Banach
theorem} implies that there is a bounded linear functional $\lambda$ on $V$
such that $\lambda(v) = \|v\|$ and $\|\lambda\|_* = 1$.  In many situations,
this can be verified more directly, or at least approximately so.

\section{Bounded linear mappings}
\label{bounded linear mappings}

        Let $V$ and $W$ be vector spaces, both real or both complex,
and equipped with norms $\|v\|_V$ and $\|w\|_W$, respectively.
A linear mapping $T$ from $V$ into $W$ is said to be 
\emph{bounded}\index{bounded linear mappings} if
\begin{equation}
\label{||T(v)||_W le C ||v||_V}
        \|T(v)\|_W \le C \, \|v\|_V
\end{equation}
for some nonnegative real number $C$ and every $v \in V$.  Thus a
bounded linear functional on $V$ is the same as a bounded linear
mapping from $V$ into the real or complex numbers, as appropriate.  As
for bounded linear functionals, it is easy to see that a bounded
linear mapping from $V$ into $W$ is uniformly continuous with respect
to the corresponding metrics, and conversely that a linear mapping
from $V$ into $W$ is bounded if it is continuous at $0$.  The space of
bounded linear mappings from $V$ into $W$ is denoted $\mathcal{BL}(V,
W)$,\index{BL(V, W)@$\mathcal{BL}(V, W)$} and is a vector space with
respect to pointwise addition and scalar multiplication.

        If $T$ is a bounded linear mapping from $V$ into $W$, then put
\begin{equation}
\label{||T||_{op} = sup {||T(v)||_W : v in V, ||v||_V le 1}}
        \|T\|_{op} = \sup \{\|T(v)\|_W : v \in V, \ \|v\|_V \le 1\},
\end{equation}
which is the same as the smallest $C \ge 0$ for which (\ref{||T(v)||_W
  le C ||v||_V}) holds.  It is easy to see that this defines a norm on
$\mathcal{BL}(V, W)$, known as the \emph{operator norm}\index{operator
  norms} associated to the given norms on $V$ and $W$.  As before, one
can show that $\mathcal{BL}(V, W)$ is complete with respect to the
operator norm when $W$ is complete.  Suppose now that $V_1$, $V_2$,
and $V_3$ are vector spaces, all real or all complex, and equipped
with norms $\|\cdot \|_1$, $\|\cdot \|_2$, and $\|\cdot \|_3$,
respectively.  If $T_1 : V_1 \to V_2$ and $T_2 : V_2 \to V_3$ are
bounded linear mappings, then one can check that their composition
$T_2 \circ T_1$ is a bounded linear mapping from $V_1$ into $V_3$,
and that
\begin{equation}
\label{||T_2 circ T_1||_{op, 13} le ||T_1||_{op, 12} ||T_2||_{op, 23}}
        \|T_2 \circ T_1\|_{op, 13} \le \|T_1\|_{op, 12} \, \|T_2\|_{op, 23},
\end{equation}
where the subscripts indicate the norms and spaces being used.

        Let us restrict our attention for the rest of this section to
the case where $V$ is a real or complex Hilbert space, with inner
product $\langle v, w \rangle$ and corresponding norm $\|v\|$.  
If $T$ is a bounded linear mapping from $V$ into itself, then the
operator norm of $T$ can also be given by
\begin{equation}
\label{||T||_{op} = ..., 2}
        \|T\|_{op} = \sup \{|\langle T(v), w \rangle| : v, w \in V, \
                                                    \|v\|, \|w\| \le 1\}.
\end{equation}
Indeed, this expression for the operator norm is clearly less than or
equal to the one in (\ref{||T||_{op} = sup {||T(v)||_W : v in V,
    ||v||_V le 1}}), because of the Cauchy--Schwarz inequality.  To
get the opposite inequality, one can choose $w$ in (\ref{||T||_{op} =
  ..., 2}) so that $\langle T(v), w \rangle = \|T(v)\|$.

        The \emph{adjoint}\index{adjoints of linear mappings}
$T^*$\index{T^*@$T^*$} of a bounded linear mapping $T$ from $V$ into itself 
is defined as follows.  It is easy to see that
\begin{equation}
\label{mu_w(v) = langle T(v), w rangle}
        \mu_w(v) = \langle T(v), w \rangle
\end{equation}
is a bounded linear functional on $V$ for each $w \in V$, since $T$
is bounded on $V$.  As in the previous section, there is an element
$T^*(w)$ of $V$ such that
\begin{equation}
\label{langle T(v), w rangle = langle v, T^*(w) rangle}
        \langle T(v), w \rangle = \langle v, T^*(w) \rangle
\end{equation}
for every $v \in V$.  One can also check that $T^*(w)$ is uniquely
determined by $w$, and that $T^*$ is linear as a mapping from $V$ into
itself.  Moreover,
\begin{equation}
\label{||T^*(w)|| le ||T||_{op} ||w||}
        \|T^*(w)\| \le \|T\|_{op} \, \|w\|
\end{equation}
for every $w \in V$, because $\|T^*(w)\|$ is equal to the dual norm of
$\mu_w$, which is less than or equal to $\|T\|_{op} \, \|w\|$.  This
implies that $T^*$ is a bounded linear mapping on $V$, with operator
norm less than or equal to the operator norm of $T$.  In fact,
\begin{equation}
        \|T^*\|_{op} = \|T\|_{op},
\end{equation}
by (\ref{||T||_{op} = ..., 2}) and (\ref{langle T(v), w rangle =
  langle v, T^*(w) rangle}).

        Observe that $(T_1 + T_2)^* = T_1^* + T_2^*$ for any two bounded
linear operators $T_1$ and $T_2$ on $V$.  Similarly, if $T$ is a bounded
linear operator on $V$, then $(a \, T)^* = a \, T^*$ for every $a \in {\bf R}$
in the real case, and $(a \, T)^* = \overline{a} \, T^*$ for every
$a \in {\bf C}$ in the complex case.  If $T_1$ and $T_2$ are bounded linear
operators on $V$ again, then their composition $T_2 \circ T_1$ is also
a bounded linear operator on $V$, and
\begin{equation}
\label{(T_2 circ T_1)^* = T_1^* circ T_2^*}
        (T_2 \circ T_1)^* = T_1^* \circ T_2^*.
\end{equation}
The adjoint of the identity operator $I$ on $V$ is equal to itself.
If $T$ is any bounded linear operator on $V$, then one can check that
\begin{equation}
\label{(T^*)^* = T}
        (T^*)^* = T.
\end{equation}

        A bounded linear operator $T$ on $V$ is said to be 
\emph{invertible}\index{invertible linear operators} if $T$
is a one-to-one mapping of $V$ onto itself for which the inverse
mapping $T^{-1}$ is also bounded on $V$. In this case, $T^*$ is 
invertible on $V$ too, and
\begin{equation}
\label{(T^*)^{-1} = (T^{-1})^*}
        (T^*)^{-1} = (T^{-1})^*,
\end{equation}
since one can apply (\ref{(T_2 circ T_1)^* = T_1^* circ T_2^*})
to $T^{-1} \circ T = T \circ T^{-1} = I$.

        A bounded linear operator $T$ on $V$ is said to be 
\emph{self-adjoint}\index{self-adjoint linear operators} if $T^* = T$,
which is equivalent to asking that
\begin{equation}
\label{langle T(v), w rangle = langle v, T(w) rangle}
        \langle T(v), w \rangle = \langle v, T(w) \rangle
\end{equation}
for every $v, w \in V$.  The sum of two bounded self-adjoint linear
operators on $V$ is also self-adjoint, as is a real number times a
bounded self-adjoint linear operator on $V$.  Thus the bounded
self-adjoint linear operators on $V$ form a vector space over the real
numbers in a natural way, and it is important to use the real numbers
for this even when $V$ is a complex Hilbert space.  As a related
point, if $V$ is a complex Hilbert space and $T$ is a bounded
self-adjoint linear operator on $V$, then
\begin{equation}
\label{langle T(v), v rangle in {bf R}}
        \langle T(v), v \rangle \in {\bf R}
\end{equation}
for each $v \in V$.  Indeed, (\ref{langle T(v), w rangle = langle v,
  T(w) rangle}) and the definition of an inner product imply that
\begin{equation}
\label{langle T(v), v rangle = ... = overline{langle T(v), v rangle}}
        \langle T(v), v \rangle = \langle v, T(v) \rangle 
                                = \overline{\langle T(v), v \rangle}
\end{equation}
for every $v \in V$.

        A bounded self-adjoint linear operator $T$ on a real or complex
Hilbert space $V$ is said to be \emph{nonnegative}\index{nonnegative 
self-adjoint operators} if
\begin{equation}
\label{langle T(v), v rangle ge 0}
        \langle T(v), v \rangle \ge 0
\end{equation}
for every $v \in V$.  Note that the sum of two nonnegative bounded
self-adjoint linear operators on $V$ is nonnegative, as is a
nonnegative real number times a nonnegative bounded self-adjoint
linear operator on $V$.  If $T$ is any
bounded linear operator on $V$, then $T^* \circ T$ is self-adjoint,
because
\begin{equation}
\label{(T^* circ T)^* = T^* circ (T^*)^* = T^* circ T}
        (T^* \circ T)^* = T^* \circ (T^*)^* = T^* \circ T.
\end{equation}
Moreover, $T^* \circ T$ is nonnegative, since
\begin{equation}
\label{langle (T^* circ T)(v), v rangle = ... = ||T(v)||^2 ge 0}
        \langle (T^* \circ T)(v), v \rangle = \langle v, T^*(T(v)) \rangle
                               = \langle T(v), T(v) \rangle = \|T(v)\|^2 \ge 0
\end{equation}
for every $v \in V$.

        Let $W$ be a closed linear subspace of $V$, and let $P_W$
be the orthogonal projection of $V$ onto $W$.  Remember that $P_W$
is characterized by the conditions that $P_W(v) \in W$ and $v - P_W(v)
\in W^\perp$ for every $v \in V$, and that $P_W$ is a bounded linear
operator on $V$ with operator norm equal to $1$, except in the trivial
case where $W = \{0\}$ and $P_W = 0$.  Using this, one can check that
\begin{equation}
\label{langle P_W(v), w rangle = ... = langle v, P_W(w) rangle}
        \langle P_W(v), w \rangle = \langle P_W(v), P_W(w) \rangle
                                   = \langle v, P_W(w) \rangle
\end{equation}
for every $v, w \in V$, so that $P_W$ is self-adjoint on $V$.  In
particular,
\begin{equation}
\label{langle P_W(v), v rangle = ... = ||P_W(v)||^2 ge 0}
\langle P_W(v), v \rangle = \langle P_W(v), P_W(v) \rangle = \|P_W(v)\|^2 \ge 0
\end{equation}
for each $v \in V$, so that $P_W$ is nonnegative.  Note that $P_W
\circ P_W = P_W$, and hence $P_W = P_W^* \circ P_W$, because $P_W$ is
self-adjoint.

        If $T$ is any bounded linear operator on $V$, then
\begin{equation}
        \|T^* \circ T\|_{op} \le \|T^*\|_{op} \, \|T\|_{op} = \|T\|_{op}^2.
\end{equation}
Using (\ref{langle (T^* circ T)(v), v rangle = ... = ||T(v)||^2 ge 0})
and the Cauchy--Schwarz inequality, we also have that
\begin{equation}
\label{||T(v)||^2 le ||(T^* circ T)(v)|| ||v|| le ||T^* circ T||_{op} ||v||^}
 \|T(v)\|^2 \le \|(T^* \circ T)(v)\| \, \|v\| \le \|T^* \circ T\|_{op} \, \|v\|^2
\end{equation}
for every $v \in V$.  This implies that $\|T\|_{op}^2 \le \|T^* \circ T\|_{op}$,
and hence
\begin{equation}
        \|T^* \circ T\|_{op} = \|T\|_{op}^2,
\end{equation}
which is known as the \emph{$C^*$-identity}.\index{C^*-identity@$C^*$-identity}

        A linear mapping $T$ of $V$ onto itself is said to be
\emph{unitary}\index{unitary operators} if
\begin{equation}
\label{langle T(v), T(w) rangle = langle v, w rangle}
        \langle T(v), T(w) \rangle = \langle v, w \rangle
\end{equation}
for every $v, w \in V$.  In the real case, one might say instead that
$T$ is an \emph{orthogonal transformation}.\index{orthogonal transformation}
If we take $v = w$ in (\ref{langle T(v), T(w) rangle = langle v, w rangle}),
then we get that
\begin{equation}
\label{||T(v)|| = ||v||}
        \|T(v)\| = \|v\|
\end{equation}
for every $v \in V$.  Conversely, it is well known that (\ref{||T(v)||
  = ||v||}) implies (\ref{langle T(v), T(w) rangle = langle v, w
  rangle}), because of polarization identities.  Note that
(\ref{||T(v)|| = ||v||}) implies that $T$ is a bounded linear operator
on $V$ with trivial kernel, and that the inverse operator $T^{-1}$ is
also bounded when $T$ maps $V$ onto itself.  It is easy to see that
(\ref{langle T(v), T(w) rangle = langle v, w rangle}) holds if and
only if $T^* \circ T = I$.  If $T$ maps $V$ onto itself, then this is
the same as saying that $T$ is invertible, with $T^{-1} = T^*$.

        Suppose now that $V$ is a complex Hilbert space, and let
$T$ be a bounded linear operator on $V$.  Observe that
\begin{equation}
\label{A = frac{T + T^*}{2} and B = frac{T - T^*}{2 i}}
        A = \frac{T + T^*}{2} \quad\hbox{and}\quad B = \frac{T - T^*}{2 \, i}
\end{equation}
are bounded self-adjoint linear operators on $V$, and that $T = A + i
\, B$.  If $T$ and $T^*$ commute with each other, which is to say that
$T \circ T^* = T^* \circ T$, then $T$ is said to be
\emph{normal}.\index{normal linear operators} Equivalently, $T$ is
normal if and only if $A$ and $B$ commute.  Thus self-adjoint linear
operators are automatically normal, and unitary operators are normal
as well, since an invertible operator automatically commutes with its
inverse.

\section{Double sums}
\label{double sums}

        Let $E_1$ and $E_2$ be nonempty sets, and consider their Cartesian
product $E = E_1 \times E_2$,  which is the set of all ordered pairs 
$(x, y)$ with $x \in E_1$ and $y \in E_2$.  If $f(x, y)$ is a nonnegative
real-valued function on $E$, then we can define the sums
\begin{equation}
\label{f_1(x) = sum_{y in E_2} f(x, y)}
        f_1(x) = \sum_{y \in E_2} f(x, y)
\end{equation}
for every $x \in E_1$ and
\begin{equation}
\label{f_2(y) = sum_{x in E_1} f(x, y)}
        f_2(y) = \sum_{x \in E_1} f(x, y)
\end{equation}
for every $y \in E_2$ as the suprema of the corresponding finite subsums,
as in Section \ref{summable functions}.  We can then define the iterated sums
\begin{equation}
\label{sum_{x in E_1} f_1(x) and sum_{y in E_2} f_2(y)}
        \sum_{x \in E_1} f_1(x) \quad\hbox{and}\quad \sum_{y \in E_2} f_2(y)
\end{equation}
also as in Section \ref{summable functions}, with the obvious convention
that these sums are infinite when any of their terms are infinite.
Under these conditions, one can show that the two sums in (\ref{sum_{x in E_1} 
f_1(x) and sum_{y in E_2} f_2(y)}) are equal to each other and to the double sum
\begin{equation}
\label{sum_{(x, y) in E} f(x, y)}
        \sum_{(x, y) \in E} f(x, y),
\end{equation}
also defined as in Section \ref{summable functions}.  In particular,
$f(x, y)$ is summable on $E$ if and only if the sums in (\ref{sum_{x
    in E_1} f_1(x) and sum_{y in E_2} f_2(y)}) are finite.

        Now let $f(x, y)$ be a real or complex-valued summable function on $E$,
so that $|f(x, y)|$ is a nonnegative real-valued summable function on $E$.
In this case, $f(x, y)$ is summable as a function of $y \in E_2$ for each
$x \in E_1$, and similarly $f(x, y)$ is summable as a function of $x \in E_1$
for each $y \in E_2$.  Thus $f_1(x)$ and $f_2(y)$ may be defined as in
(\ref{f_1(x) = sum_{y in E_2} f(x, y)}) and (\ref{f_2(y) = sum_{x in
    E_1} f(x, y)}), and we have that
\begin{equation}
\label{|f_1(x)| le sum_{y in E_2} |f(x, y)|}
        |f_1(x)| \le \sum_{y \in E_2} |f(x, y)|
\end{equation}
for each $x \in E_1$, and
\begin{equation}
\label{|f_2(y)| le sum_{x in E_1} |f(x, y)|}
        |f_2(y)| \le \sum_{x \in E_1} |f(x, y)|
\end{equation}
for each $y \in E_2$.  Because $f(x, y)$ is summable on $E$, it follows that
\begin{equation}
\label{sum_{x in E_1} |f_1(x)| le sum_{x in E_1} (sum_{y in E_2} |f(x, y)|)}
 \sum_{x \in E_1} |f_1(x)| \le \sum_{x \in E_1} \Big(\sum_{y \in E_2} |f(x, y)|\Big)
\end{equation}
and
\begin{equation}
\label{sum_{y in E_2} |f_2(y)| le sum_{y in E_2} (sum_{x in E_1} |f(x, y)|)}
 \sum_{y \in E_2} |f_2(y)| \le \sum_{y \in E_2} \Big(\sum_{x \in E_1} |f(x, y)|\Big)
\end{equation}
are finite, so that $f_1(x)$ and $f_2(y)$ are summable on $E_1$ and
$E_2$, respectively.  This implies that the iterated sums (\ref{sum_{x
    in E_1} f_1(x) and sum_{y in E_2} f_2(y)}) are defined in this
situation, and one can check that they are equal to each other and to
the double sum (\ref{sum_{(x, y) in E} f(x, y)}).  One way to do this
is to express $f(x, y)$ as a linear combination of nonnegative
real-valued summable functions on $E$, to reduce to the previous case.
Alternatively, one can approximate $f(x, y)$ by functions with finite
support on $E$.

\section{Continuous functions}
\label{continuous functions}

        Let $X$ be a topological space, and let $C(X)$\index{C(X)@$C(X)$}
be the space of continuous real or complex-valued functions on $X$.
As usual, this may also be denoted $C(X, {\bf R})$ or $C(X, {\bf C})$
to indicate whether real or complex-valued functions are being used,
and similarly for other spaces of functions on $X$.  Note that $C(X)$
is a vector space with respect to pointwise addition and scalar
multiplication.

        Let $C_b(X)$\index{C_b(X)@$C_b(X)$} be the linear subspace of $X$ 
consisting of bounded continuous functions on $X$.  Of course, every
continuous function on $X$ is bounded when $X$ is compact.  If $f \in
C_b(X)$, then the \emph{supremum norm}\index{supremum norm} of $f$ is
defined by
\begin{equation}
\label{||f||_{sup} = sup_{x in X} |f(x)|}
        \|f\|_{sup} = \sup_{x \in X} |f(x)|.
\end{equation}
It is easy to see that this defines a norm on $C_b(X)$, and it is well
known that $C_b(X)$ is complete with respect to this norm.  This uses the
fact that if $\{f_j\}_{j = 1}^\infty$ is a sequence of continuous functions
on $X$ that converges uniformly to a function $f$ on $X$, then $f$ is
continuous on $X$ too.

        As usual, the \emph{support}\index{support of a function} $\supp f$ 
of a function $f$ on $X$ is the closure of the set of $x \in X$ such that
$f(x) \ne 0$.  Let $C_{com}(X)$\index{C_{com}(X)@$C_{com}(X)$} be the space
of continuous functions on $X$ with compact support, which is a linear
subspace of $C_b(X)$.  Urysohn's lemma implies that there are plenty
of continuous functions on $X$ with compact support when $X$ is a
locally compact Hausdorff space.

        Suppose from now on in this section that $X$ is a locally compact 
Hausdorff space.  A continuous function $f$ on $X$ is said to
\emph{vanish at infinity}\index{vanishing at infinity} if for each
$\epsilon > 0$ there is a compact set $K \subseteq X$ such that
\begin{equation}
\label{|f(x)| < epsilon}
        |f(x)| < \epsilon
\end{equation}
for every $x \in X \backslash K$.  In particular, this implies that
$f$ is bounded on $X$, and one can check that the space
$C_0(X)$\index{C_0(X)@$C_0(X)$} of continuous functions on $X$ that
vanish at infinity is a closed linear subspace of $C_b(X)$ with
respect to the supremum norm.  Of course, $C_{com}(X) \subseteq
C_0(X)$, and one can also check that $C_0(X)$ is the closure of
$C_{com}(X)$ in $C_b(X)$.  More precisely, if $f \in C_0(X)$, then one
can approximate $f$ uniformly on $X$ by continuous functions on $X$
with compact support, by multiplying $f$ by suitable cut-off functions
obtained from Urysohn's lemma.

        A linear functional $\lambda$ on $C_{com}(X)$ is said to be
\emph{nonnegative}\index{nonnegative} if $\lambda(f)$ is a nonnegative
real number for every nonnegative real-valued continuous function $f$
with compact support on $X$.  In this case, the Riesz representation
theorem\index{Riesz representation theorem} implies that there is a
unique nonnegative Borel measure $\mu$ on $X$ with certain regularity
properties such that
\begin{equation}
\label{lambda(f) = int_X f d mu}
        \lambda(f) = \int_X f \, d\mu
\end{equation}
for every $f \in C_{com}(X)$.  In particular, $\mu(K) < \infty$ for every
compact set $K \subseteq X$ under these conditions.

        Now let $\lambda$ be a bounded linear functional on $C_0(X)$
with respect to the supremum norm.  Another version of the Riesz 
representation theorem\index{Riesz representation theorem} implies
that there is a unique real or complex Borel measure $\mu$ on $X$,
as appropriate, with certain regularity properties such that
(\ref{lambda(f) = int_X f d mu}) holds for every $f \in C_0(X)$.
Remember that there is a nonnegative measure $|\mu|$ associated to
any real or complex measure $\mu$, known as the \emph{total variation
measure}.\index{total variation measures}  The regularity of $\mu$
means that $|\mu|$ is regular as a nonnegative Borel measure on $X$,
and part of the theorem is that the dual norm of $\lambda$ on $C_0(X)$
is equal to $|\mu|(X)$.

        A set $A \subseteq X$ is said to be
\emph{$\sigma$-compact}\index{sigma-compact sets@$\sigma$-compact sets}
if there is a sequence $K_1, K_2, K_3, \ldots$ of compact subsets of $X$
such that $A = \bigcup_{l = 1}^\infty K_l$.  If every open set in $X$
is $\sigma$-compact, and if $\mu$ is a nonnegative Borel measure on
$X$ such that $\mu(K) < \infty$ for every compact set $K \subseteq X$,
then it is well known that $\mu$ satisfies the same regularity conditions
as in the Riesz representation theorem.  See Theorem 2.18 on p50 of
\cite{r2}, for instance.  In particular, this condition holds when there
is a countable base for the topology of $X$.  To see this, note that
a locally compact Hausdorff space $X$ is regular as a topological
space, which is to say that $X$ satisfies the third separation condition.
Together with local compactness, this implies that for every open set
$U$ in $X$ and every point $p \in U$ there is an open set $V(p)$ in $X$
such that $p \in V(p)$ and the closure $\overline{V(p)}$ of $V(p)$ in $X$
is a compact set contained in $U$.  If there is a countable base for
the topology of $X$, then it follows that $U$ can be expressed as the
union of only finitely or countably many $V(p)$'s, and hence as the
union of only finitely or countably many of their closures.

\section{Double integrals}
\label{double integrals}

        Let $X$ and $Y$ be locally compact Hausdorff topological
spaces, and let $\mu$ and $\nu$ be Borel measures on $X$ and $Y$,
respectively.  More precisely, $\mu$ and $\nu$ may be real or complex
measures, or nonnegative measures which are finite on compact sets.
In order to apply the usual construction of product measures on $X
\times Y$, one normally asks $\mu$ and $\nu$ to be $\sigma$-finite.

        However, there is a technical problem with the usual product
construction of measurable subsets of $X \times Y$, which is that
one would like open subsets of $X \times Y$ to be measurable, and
which would imply that Borel subsets of $X \times Y$ are measurable.
Of course, products of open subsets of $X$ and $Y$ are measurable in
$X \times Y$, and hence countable unions of products of open subsets
of $X$ and $Y$ are measurable in $X \times Y$.  If there are countable
bases for the topologies of $X$ and $Y$, then one can get a countable
base for the topology of $X \times Y$ by taking products of the basic
open sets in $X$ and $Y$, and it follows that every open set in $X \times Y$
is a countable union of products of open subsets of $X$ and $Y$.
If an open set $W$ in $X \times Y$ is $\sigma$-compact, then one can
also check that $W$ can be expressed as the union of countably many
products of open subsets of $X$ and $Y$.  One can also consider subclasses
of the Borel sets, for which continuous functions with compact support
or which vanish at infinity are still measurable.

        Alternatively, one can approach integration on $X \times Y$
in terms of linear functionals on spaces of continuous functions
on $X \times Y$, as in the previous section.  More precisely, given
nonnegative linear functionals $\lambda_X$ and $\lambda_Y$ on
$C_{com}(X)$ and $C_{com}(Y)$, respectively, one would like to define
a nonnegative linear functional $\lambda_{X \times Y}$ on $C_{com}(X \times Y)$.
Similarly, if $\lambda_X$ and $\lambda_Y$ are bounded linear functionals
on $C_0(Y)$ and $C_0(Y)$, respectively, then one would like to define a
bounded linear functional $\lambda_{X \times Y}$ on $C_0(X \times Y)$.
In both cases, one would like $\lambda_{X \times Y}$ to satisfy
\begin{equation}
\label{lambda_{X times Y}(f g) = lambda_X(f) lambda_Y(g)}
        \lambda_{X \times Y}(f \, g) = \lambda_X(f) \, \lambda_Y(g)
\end{equation}
when $f(x)$ and $g(y)$ are continuous functions of the appropriate type
on $X$ and $Y$, so that $f(x) \, g(y)$ is continuous on $X \times Y$.
One can show that $\lambda_{X \times Y}$ is uniquely determined by this
condition, by approximating continuous functions on $X \times Y$ by
finite sums of products of continuous functions on $X$ and $Y$.
If $h(x, y)$ is a continuous function on $X \times Y$ with compact
support or that vanishes at infinity, then one can define 
$\lambda_{X \times Y}(h)$ by first applying $\lambda_X$ to $h(x, y)$ as
a function of $x \in X$ for each $y \in Y$, and then apply $\lambda_Y$
to the resulting function of $y$.  It is not too difficult to check
that this has the desired properties.  One could also apply $\lambda_X$
and $\lambda_Y$ in the opposite order, and this would lead to the same
value of $\lambda_{X \times Y}(h)$, by the uniqueness argument mentioned
earlier.  As in the previous section, one could then use the appropriate
version of the Riesz representation theorem to get a Borel measure
on $X \times Y$ corresponding to $\lambda_{X \times Y}$.

        Now let $I$ be an infnite set, and suppose that $X_j$ is a compact
Hausdorff topological space for each $j \in I$.  Also let $X = \prod_{j \in I}
X_j$ be the Cartesian product of the $X_j$'s, which is a compact Hausdorff
space with respect to the product topology, by Tychonoff's theorem.
If $\mu_j$ is a probability measure on $X_j$ for each $j \in I$,
then there is a well known construction of a product probability measure
$\mu$ on $X$.  Alternatively, one can approach this in terms of nonnegative
linear functionals on $C(X_j)$, as follows.  Suppose that $\lambda_j$ is a 
nonnegative linear functional on $C(X_j)$ such that
\begin{equation}
\label{lambda_j({bf 1}_{X_j}) = 1}
        \lambda_j({\bf 1}_{X_j}) = 1
\end{equation}
for each $j \in I$, where ${\bf 1}_{X_j}$ is the constant function
equal to $1$ on $X_j$.  Thus $\lambda_j$ corresponds to a regular
Borel probability measure on $X_j$ for each $j \in I$, by the Riesz
representation theorem, and one would like to define a corresponding
product nonnegative linear functional $\lambda$ on $C(X)$.  If $f$ is
a continuous function on $X$ that depends on only finitely many
variables $x_j \in X_j$, then $\lambda(f)$ can be defined by applying
$\lambda_j$ to $f$ as a function of $x_j$ for those finitely many $j
\in I$, as before.  Otherwise, one can use compactness to show that
every continuous function $f$ on $X$ can be approximated uniformly by
continuous functions that depend on only finitely many variables, and
then use this to extend $\lambda$ to a nonnegative linear functional on
$C(X)$.

\section{Ultrametrics}
\label{ultrametrics}

        Let $(M, d(x, y))$ be a metric space.  The metric $d(x, y)$
on $M$ is said to be an \emph{ultrametric}\index{ultrametrics} if
\begin{equation}
\label{d(x, z) le max(d(x, y), d(y, z))}
        d(x, z) \le \max(d(x, y), d(y, z))
\end{equation}
for every $x, y, z \in M$.  Of course, this is stronger than the
ordinary triangle inequality (\ref{d(x, z) le d(x, y) + d(y, z)}) in
Section \ref{metrics, norms}.  As a basic example, the discrete
metric\index{discrete metric} on any set $M$ is defined by putting
$d(x, y)$ equal to $1$ when $x \ne y$ and equal to $0$ when $x = y$,
and is an ultrametric.

        Let $(M_1, d_1(x_1, y_1))$ and $(M_2, d_2(x_2, y_2))$ be metric spaces, 
and consider their Cartesian product.  As usual, it is easy to see that
\begin{equation}
\label{D((x_1, x_2), (y_1, y_2)) = d_1(x_1, y_1) + d_2(x_2, y_2)}
        D((x_1, x_2), (y_1, y_2)) = d_1(x_1, y_1) + d_2(x_2, y_2)
\end{equation}
and
\begin{equation}
\label{D'((x_1, x_2), (y_1, y_2)) = max(d_1(x_1, y_1), d_2(x_2, y_2))}
        D'((x_1, x_2), (y_1, y_2)) = \max(d_1(x_1, y_1), d_2(x_2, y_2))
\end{equation}
define metrics on $M_1 \times M_2$, for which the corresponding topologies
are the same as the product topology associated to the topologies on $M_1$
and $M_2$ by the metrics $d_1(x_1, y_1)$ and $d_2(x_2, y_2)$, respectively.
If $d_1(x_1, y_1)$ and $d_2(x_2, y_2)$ are ultrametrics on $M_1$ and $M_2$,
then (\ref{D'((x_1, x_2), (y_1, y_2)) = max(d_1(x_1, y_1), d_2(x_2, y_2))})
is also an ultrametric on $M_1 \times M_2$.  Note that a ball in $M_1
\times M_2$ with respect to (\ref{D'((x_1, x_2), (y_1, y_2)) =
  max(d_1(x_1, y_1), d_2(x_2, y_2))}) is the same as the Cartesian product
of balls in $M_1$ and $M_2$ with the same radius.

        Now let $(M_j, d_j(x_j, y_j))$, $j = 1, 2, 3, \ldots$, be a sequence
of metric spaces, and let $M = \prod_{j = 1}^\infty M_j$ be their Cartesian
product.  Thus $M$ consists of sequences $x = \{x_j\}_{j = 1}^\infty$, where
$x_j \in M_j$ for each $j$.  Also let $t = \{t_j\}_{j = 1}^\infty$ be a
sequence of positive real numbers that converges to $0$, and put
\begin{equation}
\label{d_j'(x_j, y_j) = min(d_j(x_j, y_j), t_j)}
        d_j'(x_j, y_j) = \min(d_j(x_j, y_j), t_j)
\end{equation}
for each $j$.  It is easy to see that $d_j'(x_j, y_j)$ is also a
metric on $M_j$ for each $j$, which determines the same topology on $M_j$
as $d_j(x_j, y_j)$.  Put
\begin{equation}
\label{d(x, y) = max_{j ge 1} d_j'(x_j, y_j)}
        d(x, y) = \max_{j \ge 1} d_j'(x_j, y_j)
\end{equation}
for each $x, y \in M$, which is obviously equal to $0$ when $x = y$.
If $x \ne y$, then $x_{j_0} \ne y_{j_0}$ for some $j_0 \ge 1$, so that
$d_{j_0}'(x_{j_0}, y_{j_0}) > 0$.  This implies that
\begin{equation}
\label{d_j'(x_j, y_j) le t_j le d_{j_0}'(x_{j_0}, y_{j_0})}
        d_j'(x_j, y_j) \le t_j \le d_{j_0}'(x_{j_0}, y_{j_0})
\end{equation}
for all but finitely many $j$, since $t_j \to 0$ as $j \to \infty$, so
that the maximum in (\ref{d(x, y) = max_{j ge 1} d_j'(x_j, y_j)})
always exists.

        One can check that (\ref{d(x, y) = max_{j ge 1}  d_j'(x_j, y_j)})
is a metric on $M$, for which the corresponding topology is the same
as the product topology associated to the topologies on the $M_j$'s
determined by the metrics $d_j(x_j, y_j)$.  More precisely, an open
ball in $M$ of radius $r > 0$ with respect to (\ref{d(x, y) = max_{j
    ge 1} d_j'(x_j, y_j)}) can be expressed as a product $\prod_{j =
  1}^\infty B_j$, where $B_j$ is an open ball in $M_j$ with radius $r$
when $r \le t_j$, and $B_j = M_j$ when $r > t_j$.  In particular, $B_j
= M_j$ for all but finitely many $j$, since $t_j \to 0$ as $j \to
\infty$, which implies that open balls in $M$ with respect to
(\ref{d(x, y) = max_{j ge 1} d_j'(x_j, y_j)}) are open sets with
respect to the product topology.  By taking $r$ sufficiently small,
one gets that $r \le t_j$ for any finite set of $j \ge 1$, so
that open balls in $M$ with respect to (\ref{d(x, y) = max_{j ge 1}
  d_j'(x_j, y_j)}) generate the product topology on $M$.  If $d_j(x_j,
y_j)$ is an ultrametric on $M_j$ for each $j$, then $d_j'(x_j, y_j)$
is an ultrametric on $M_j$ as well, and (\ref{d(x, y) = max_{j ge 1}
  d_j'(x_j, y_j)}) is an ultrametric on $M$ too.

        Let $(M, d(x, y))$ be an arbitrary metric space again.
The open ball\index{open balls} in $M$ with center $x \in M$ and
radius $r > 0$ is defined as usual by
\begin{equation}
\label{B(x, r) = {y in M : d(x, y) < r}}
        B(x, r) = \{y \in M : d(x, y) < r\},
\end{equation}
and similarly the closed ball\index{closed balls} with center $x \in
M$ and radius $r \ge 0$ is given by
\begin{equation}
\label{overline{B}(x, r) = {y in M : d(x, y) le r}}
        \overline{B}(x, r) = \{y \in M : d(x, y) \le r\}.
\end{equation}
Let us also put
\begin{equation}
\label{V(x, r) = {y in M : d(x, y) > r}}
        V(x, r) = \{y \in M : d(x, y) > r\}
\end{equation}
for each $x \in M$ and $r \ge 0$, which is the same as the complement
of $\overline{B}(x, r)$.  It is well known that $B(x, r)$ and $V(x,
r)$ are open sets in $M$, and that $\overline{B}(x, r)$ is a closed
set.  More precisely, if $z \in B(z, r)$, then one can check that
\begin{equation}
\label{B(z, t) subseteq B(x, r)}
        B(z, t) \subseteq B(x, r)
\end{equation}
with $t = r - d(x, z) > 0$, using the triangle inequality.  Similarly,
if $z \in V(x, r)$, then we have that
\begin{equation}
\label{B(z, t) subseteq V(x, r)}
        B(z, t) \subseteq V(x, r)
\end{equation}
with $t = d(x, z) - r > 0$.  One can check directly that
$\overline{B}(x, r)$ is a closed set, in the sense that it contains
all of its limit points, or derive this from the fact that $V(x, r)$
is an open set.

        Suppose now that $d(x, y)$ is an ultrametric on $M$.
In this case, it is easy to see that (\ref{B(z, t) subseteq B(x, r)})
holds for every $z \in B(x, r)$ with $t = r$.  Similarly, 
\begin{equation}
\label{overline{B}(z, r) subseteq overline{B}(x, r)}
        \overline{B}(z, r) \subseteq \overline{B}(x, r)
\end{equation}
for every $z \in \overline{B}(x, r)$.  Observe that
\begin{equation}
\label{d(x, z) le d(x, y)}
        d(x, z) \le d(x, y)
\end{equation}
for every $y \in M$ with $d(y, z) < d(x, z)$, since $d(x, y) < d(x,
z)$ would imply that $d(x, z) < d(x, z)$.  This shows that (\ref{B(z,
  t) subseteq V(x, r)}) holds with $t = d(x, z)$.  

        Put
\begin{equation}
\label{W(x, r) = {y in M : d(x, y) ge r}}
        W(x, r) = \{y \in M : d(x, y) \ge r\}
\end{equation}
for each $x \in M$ and $r \ge 0$, which is the same as the complement
of $B(x, r)$ in $M$.  Using (\ref{d(x, z) le d(x, y)}) again,
we get that
\begin{equation}
        B(z, d(x, z)) \subseteq W(x, r)
\end{equation}
for every $z \in W(x, r)$.  In particular, $W(x, r)$ is an open set in
an ultrametric space, so that $B(x, r)$ is both open and closed.  It
follows from (\ref{overline{B}(z, r) subseteq overline{B}(x, r)}) that
$\overline{B}(x, r)$ is an open set in an ultrametric space, and hence
is both open and closed as well.  Thus ultrametric spaces are totally
disconnected, in the sense that they do not contain any connected
subsets with more than one element.

\chapter{Fourier series}
\label{fourier series}

\section{Basic notions}
\label{basic notions}

        Let ${\bf T}$\index{T@${\bf T}$} be the unit circle\index{unit circle}
in the complex plane, which is the set of $z \in {\bf C}$ with $|z| =
1$.  Also let $f$ be a complex-valued integrable function on ${\bf T}$
with respect to arc-length measure $|dz|$ on ${\bf T}$, which corresponds 
to Lebesgue measure on an interval in the real line when ${\bf T}$ is
parameterized by arc length.  The $n$th 
\emph{Fourier coefficient}\index{Fourier coefficients} of $f$ is defined by
\begin{equation}
\label{widehat{f}(n) = frac{1}{2 pi} int_{bf T} f(z) overline{z}^n |dz|}
 \widehat{f}(n) = \frac{1}{2 \pi} \int_{\bf T} f(z) \, \overline{z}^n \, |dz|
\end{equation}
for each integer $n$.\index{f(n)@$\widehat{f}(n)$} The corresponding
\emph{Fourier series}\index{Fourier series} is given by
\begin{equation}
\label{sum_{n = -infty}^infty widehat{f}(n) z^n}
        \sum_{n = -\infty}^\infty \widehat{f}(n) \, z^n,
\end{equation}
where for the moment this is a formal series in $z \in {\bf T}$.  
Note that
\begin{equation}
\label{|widehat{f}(n)| le frac{1}{2 pi} int_{bf T} |f(z)| |dz|}
        |\widehat{f}(n)| \le \frac{1}{2 \pi} \int_{\bf T} |f(z)| \, |dz|
\end{equation}
for each $n \in {\bf Z}$, where ${\bf Z}$\index{Z@${\bf Z}$} denotes
the integers.

        It is well known that
\begin{equation}
\label{int_{bf T} z^n |dz| = 0}
        \int_{\bf T} z^n \, |dz| = 0
\end{equation}
for each $n \in {\bf Z}$ with $n \ne 0$.  One way to see this is to
use the fact that $\exp (i \, t)$ parameterizes ${\bf T}$ by arc
length for $0 \le t \le 2 \pi$.  Of course, $|\exp (i \, t)| = 1$ for
each $t \in {\bf R}$, as in (\ref{|exp z|^2 = (exp z) (exp
  overline{z}) = ... = exp (2 re z)}), and the derivative of $\exp (i
\, t)$ is equal to $i$ times $\exp (i \, t)$, as one can see by 
differentiating the power series for the exponential function term by term.
Thus the modulus of the derivative of $\exp (i \, t)$ is also equal to $1$
for every $t \in {\bf R}$, so that $\exp (i \, t)$ goes around the
unit circle at unit speed.  This permits (\ref{int_{bf T} z^n |dz| = 0}) to
be reduecd to
\begin{equation}
\label{int_0^{2 pi} exp(i n t) dt = 0}
        \int_0^{2 \pi} \exp(i \, n \, t) \, dt = 0
\end{equation}
when $n \ne 0$, which can be derived from the fundamental theorem of calculus.

        If $f$ and $g$ are complex-valued square-integrable functions on 
${\bf T}$, then put
\begin{equation}
\label{langle f, g rangle = frac{1}{2 pi} int_{bf T} f(z) overline{g(z)} |dz|}
        \langle f, g \rangle = \frac{1}{2 \pi} \int_{\bf T} f(z) \, 
                                               \overline{g(z)} \, |dz|.
\end{equation}
This defines an inner product on $L^2({\bf T})$, with the corresponding norm
\begin{equation}
\label{(frac{1}{2 pi} int_{bf T} |f(z)|^2 |dz|)^{1/2}}
        \Big(\frac{1}{2 \pi} \int_{\bf T} |f(z)|^2 \, |dz|\Big)^{1/2}.
\end{equation}
It is well known that $L^2({\bf T})$ is complete with respect to this
norm, and is thus a Hilbert space.  Put
\begin{equation}
\label{e_n(z) = z^n}
        e_n(z) = z^n
\end{equation}
for each $n \in {\bf Z}$ and $z \in {\bf T}$, and observe that the
$e_n$'s form an orthonormal collection of functions in $L^2({\bf T})$
with respect to the inner product (\ref{langle f, g rangle = frac{1}{2
    pi} int_{bf T} f(z) overline{g(z)} |dz|}), because of
(\ref{int_{bf T} z^n |dz| = 0}).  The $n$th Fourier coefficient of $f
\in L^2({\bf T})$ can be expressed as
\begin{equation}
\label{widehat{f}(n) = langle f, e_n rangle}
        \widehat{f}(n) = \langle f, e_n \rangle
\end{equation}
for each $n \in {\bf Z}$, and we have that
\begin{equation}
\label{sum_n |widehat{f}(n)|^2 le frac{1}{2 pi} int_{bf T} |f(z)|^2 |dz|}
        \sum_{n = -\infty}^\infty |\widehat{f}(n)|^2
                          \le \frac{1}{2 \pi} \int_{\bf T} |f(z)|^2 \, |dz|,
\end{equation}
as in (\ref{sum_{j = 1}^n |langle v, v_j rangle|^2 le ||v||^2}),
(\ref{sum_{j = 1}^infty |langle v, v_j rangle|^2 le ||v||^2}), and
(\ref{sum_{alpha in A} |f_v(alpha)|^2 le ||v||_V^2}) in Sections
\ref{orthogonal vectors}, \ref{orthogonal sequences}, and
\ref{square-summability}, respectively.

        In particular, the Fourier series (\ref{sum_{n = -infty}^infty 
widehat{f}(n) z^n}) converges in $L^2({\bf T})$ when $f \in L^2({\bf T})$,
as in Sections \ref{orthogonal sequences} and \ref{square-summability}.
Using the Stone--Weierstrass theorem, one can show that the linear span
of the $e_n$'s is dense in the space $C({\bf T})$ of continuous complex-valued
functions on ${\bf T}$ with respect to the supremum norm, and hence is also
dense in $L^2({\bf T})$.  This implies that the Fourier series of $f$
converges to $f$ in $L^2({\bf T})$, and we shall see another proof of this
later on.  Note that the convergence of the sum in (\ref{sum_n 
|widehat{f}(n)|^2 le frac{1}{2 pi} int_{bf T} |f(z)|^2 |dz|}) implies that
\begin{equation}
\label{lim_{|n| to infty} |widehat{f}(n)| = 0}
        \lim_{|n| \to \infty} |\widehat{f}(n)| = 0
\end{equation}
for every $f \in L^2({\bf T})$.  This also holds when $f \in L^1({\bf
  T})$, as one can show using the fact that $L^2({\bf T})$ is dense in
$L^1({\bf T})$, and the simple estimate (\ref{|widehat{f}(n)| le
  frac{1}{2 pi} int_{bf T} |f(z)| |dz|}).

\section{Abel sums}
\label{abel sums}

        Let us say that an infinite series $\sum_{n = 0}^\infty a_n$ of
complex numbers is \emph{admissible}\index{admissible series} if
$\sum_{n = 0}^\infty |a_n| \, r^n$ converges for every $r \in {\bf R}$
with $0 \le r < 1$.  Of course, the convergence of $\sum_{n = 0}^\infty a_n \,
r^n$ implies that $\{a_n \, r^n\}_{n = 0}^\infty$ converges to $0$,
and hence is bounded.  Conversely, if $\{a_n \, t^n\}_{n = 0}^\infty$ is
bounded for some $t \in {\bf R}$ with $0 < t < 1$, then $\sum_{n = 0}^\infty
|a_n| \, r^n$ converges for every $r \in {\bf R}$ with $0 \le r < t$,
by comparison with the convergent geometric series $\sum_{n = 0}^\infty
(r/t)^n$.  Thus $\sum_{n = 0}^\infty$ is admissible if and only if
$\{a_n \, t^n\}_{n = 0}^\infty$ is bounded for every $t \in {\bf R}$
with $0 \le t < 1$.  In particular, $\sum_{n = 0}^\infty a_n$ is admissible
when $\{a_n\}_{n = 0}^\infty$ is bounded.

        Suppose that $\sum_{n = 0}^\infty a_n$ is admissible, and put
\begin{equation}
\label{A(r) = sum_{n = 0}^infty a_n r^n}
        A(r) = \sum_{n = 0}^\infty a_n \, r^n
\end{equation}
for each $r \in {\bf R}$ with $0 \le r < 1$.  If the limit
\begin{equation}
\label{lim_{r to 1-} A(r)}
        \lim_{r \to 1-} A(r)
\end{equation}
exists, then $\sum_{n = 0}^\infty a_n$ is said to be \emph{Abel
  summable}.\index{Abel summability} If $\sum_{n = 0}^\infty |a_n|$
converges, then $\sum_{n = 0}^\infty a_n$ is obviously admissible, and
one can check that $\sum_{n = 0}^\infty a_n$ is Abel summable, with
Abel sum (\ref{lim_{r to 1-} A(r)}) equal to the usual sum $\sum_{n =
  0}^\infty a_n$.  One way to do this is to use the analogue of
Lebesgue's dominated convergence theorem for sums.  Alternatively, if
$\sum_{n = 0}^\infty |a_n|$ converges, then the partial sums $\sum_{n
  = 0}^N a_n \, r^n$ converge uniformly on $[0, 1]$ as $N \to \infty$
to $\sum_{n = 0}^\infty a_n \, r^n$, by Weierstrass' $M$-test.  This
implies that $\sum_{n = 0}^\infty a_n \, r^n$ defines a continuous
function of $r$ on $[0, 1]$ in this case, so that (\ref{lim_{r to 1-}
  A(r)}) exists and is equal to $\sum_{n = 0}^\infty a_n$.  If $a_n$
is a nonnegative real number for each $n$, then it is easy to see that
the Abel sums $A(r)$ are uniformly bounded for $0 \le r < 1$ if and
only if $\sum_{n = 0}^\infty a_n$ converges.

        If $\sum_{n = 0}^\infty a_n$ is a convergent series of
complex numbers, then $\lim_{n \to \infty} a_n = 0$, and hence
$\sum_{n = 0}^\infty a_n$ is admissible.  It is well known that
$\sum_{n = 0}^\infty a_n$ is also Abel summable in this case, with
(\ref{lim_{r to 1-} A(r)}) equal to $\sum_{n = 0}^\infty a_n$.
To see this, let
\begin{equation}
\label{s_n = sum_{j = 0}^n a_j}
        s_n = \sum_{j = 0}^n a_j
\end{equation}
be the partial sums of $\sum_{j = 0}^\infty a_j$, and put $s_{-1}= 0$
for convenience.  Thus $a_n = s_n - s_{n - 1}$ for each $n \ge 0$, so that
\begin{eqnarray}
\label{A(r) = ... = (1 - r) sum_{n = 0}^infty s_n r^n}
        A(r) = \sum_{n = 0}^\infty s_n \, r^n - \sum_{n = 0}^\infty s_{n - 1} \, r^n
 & = & \sum_{n = 0}^\infty s_n \, r^n - \sum_{n = 0}^\infty s_n \, r^{n + 1} \\
 & = & (1 - r) \sum_{n = 0}^\infty s_n \, r^n                \nonumber
\end{eqnarray}
for every $r \in [0, 1)$.  This uses the fact that $s_{-1} = 0$ in the
second step, and the boundedness of the partial sums $s_n$ to get the
convergence of these series when $0 \le r < 1$.  Put $s = \lim_{n \to
  \infty} s_n = \sum_{j = 0}^\infty a_j$, and observe that
\begin{equation}
\label{A(r) - s = (1 - r) sum_{n = 0}^infty (s_n - s) r^n}
        A(r) - s = (1 - r) \sum_{n = 0}^\infty (s_n - s) \, r^n
\end{equation}
for each $r \in [0, 1)$, since $(1 - r) \sum_{n = 0}^\infty r^n = 1$.
One can show that (\ref{A(r) - s = (1 - r) sum_{n = 0}^infty (s_n - s)
  r^n}) tends to $0$ as $r \to 1-$, because $\{s_n\}_{n = 0}^\infty$
converges to $s$, and $(1 - r) \, s_n \, r^n \to 0$ as $r \to 1$ for
each $n$.  If $a_j = z^j$ for some $z \in {\bf C}$ with $|z| = 1$ and
$z \ne 1$, then $\sum_{j = 0}^\infty a_j$ does not converge, but
$A(r) = (1 - r \, z)^{-1}$ for each $r \in [0, 1)$, which tends to
$(1 - z)^{-1}$ as $r \to 1-$.

        Suppose now that $\sum_{n = -\infty}^\infty a_n$ is a doubly-infinite
series of complex numbers.  The preceding discussion can be applied to
each of the ordinary infinite series
\begin{equation}
\label{sum_{n = 0}^infty a_n and sum_{n = 1}^infty a_{-n}}
        \sum_{n = 0}^\infty a_n \quad\hbox{and}\quad \sum_{n = 1}^\infty a_{-n},
\end{equation}
or to the series
\begin{equation}
\label{a_0 + sum_{n = 1}^infty (a_n + a_{-n})}
      a_0 + \sum_{n = 1}^\infty (a_n + a_{-n}).
\end{equation}
If both of the series in (\ref{sum_{n = 0}^infty a_n and sum_{n =
    1}^infty a_{-n}}) are admissible, then (\ref{a_0 + sum_{n =
    1}^infty (a_n + a_{-n})}) is also admissible, and the Abel sums
for (\ref{a_0 + sum_{n = 1}^infty (a_n + a_{-n})}) are the same as the
sum of the Abel sums for the series in (\ref{sum_{n = 0}^infty a_n and
  sum_{n = 1}^infty a_{-n}}), which can be expressed as
\begin{equation}
\label{A(r) = sum_{n = -infty}^infty a_n r^{|n|}}
        A(r) = \sum_{n = -\infty}^\infty a_n \, r^{|n|}.
\end{equation}
Note that the Fourier series (\ref{sum_{n = -infty}^infty
  widehat{f}(n) z^n}) of an integrable function $f$ on ${\bf T}$
satisfies these conditions for each $z \in {\bf T}$, since the Fourier
coefficients are bounded, as in (\ref{|widehat{f}(n)| le frac{1}{2 pi}
  int_{bf T} |f(z)| |dz|}).

\section{The Poisson kernel}
\label{poisson kernel}

        Let $f$ be an integrable complex-valued function on the unit
circle ${\bf T}$, and consider the Abel sums (\ref{A(r) = sum_{n =
    -infty}^infty a_n r^{|n|}}) corresponding to the Fourier series
(\ref{sum_{n = -infty}^infty widehat{f}(n) z^n}) of $f$ for each
$z \in {\bf T}$.  This is the same as
\begin{eqnarray}
\label{sum_{n = -infty}^infty widehat{f}(n) r^{|n|} z^n = ...}
        \sum_{n = -\infty}^\infty \widehat{f}(n) \, r^{|n|} \, z^n
              & = & \sum_{n = 0}^\infty \widehat{f}(n) \, r^n \, z^n 
                       + \sum_{n = 1}^\infty \widehat{f}(-n) \, r^n \, z^{-n} \\
 & = & \sum_{n = 0}^\infty \widehat{f}(n) \, r^n \, z^n 
    + \sum_{n = 1}^\infty \widehat{f}(-n) \, r^n \, \overline{z}^n \nonumber
\end{eqnarray}
for each $r \in [0, 1)$ and $z \in {\bf T}$.  Alternatively, if we put
  $\zeta = r \, z$, then $|\zeta| < 1$, and (\ref{sum_{n =
      -infty}^infty widehat{f}(n) r^{|n|} z^n = ...}) is equal to
\begin{equation}
\label{sum widehat{f}(n) zeta^n + sum widehat{f}(-n) overline{zeta}^n}
        \sum_{n = 0}^\infty \widehat{f}(n) \, \zeta^n
             + \sum_{n = 1}^\infty \widehat{f}(-n) \, \overline{\zeta}^n.
\end{equation}
Note that the first sum in (\ref{sum widehat{f}(n) zeta^n + sum
  widehat{f}(-n) overline{zeta}^n}) defines a holomorphic function on
the open unit disk\index{unit disk}
\begin{equation}
\label{U = {zeta in {bf C} : |zeta| < 1}}
        U = \{\zeta \in {\bf C} : |\zeta| < 1\},
\end{equation}
the second sum in (\ref{sum widehat{f}(n) zeta^n + sum widehat{f}(-n)
  overline{zeta}^n}) is the complex-conjugate of a holomorphic
function on $U$, and so the sum of these two function is harmonic on
$U$.

        The \emph{Poisson kernel}\index{Poisson kernel} $P(\zeta, w)$
is defined for $\zeta \in U$ and $w \in {\bf T}$ by
\begin{equation}
\label{P(zeta, w) = sum zeta^n overline{w}^n + sum overline{zeta}^n w^n}
        P(\zeta, w) = \sum_{n = 0}^\infty \zeta^n \, \overline{w}^n
                          + \sum_{n = 1}^\infty \overline{\zeta}^n \, w^n.
\end{equation}
Observe that the partial sums of these series converge uniformly over
$|\zeta| \le \rho$ and $w \in {\bf T}$ for each $\rho < 1$, by
Weierstrass' $M$-test, and in particular that $P(\zeta, w)$ is
continuous on $U \times {\bf T}$.  By construction, (\ref{sum
  widehat{f}(n) zeta^n + sum widehat{f}(-n) overline{zeta}^n}) is
equal to
\begin{equation}
\label{frac{1}{2 pi} int_{bf T} P(zeta, w) f(w) |dw|}
        \frac{1}{2 \pi} \int_{\bf T} P(\zeta, w) \, f(w) \, |dw|
\end{equation}
for every $\zeta \in U$, using considerations of uniform convergence
to interchange the order of summation and integration.  It is easy to
see that
\begin{equation}
\label{frac{1}{2 pi} int_{bf T} P(zeta, w) |dw| = 1}
        \frac{1}{2 \pi} \int_{\bf T} P(\zeta, w) \, |dw| = 1
\end{equation}
for every $\zeta \in U$, using (\ref{int_{bf T} z^n |dz| = 0}),
or by applying the previous remarks to the constant function $f$
equal to $1$ on ${\bf T}$.

        In order to compute the Poisson kernel, let us re-express it as
\begin{equation}
\label{P(zeta, w) = ... = 2 re sum_{n = 0}^infty (zeta overline{w})^n - 1}
 P(\zeta, w) = \sum_{n = 0}^\infty \zeta^n \, \overline{w}^n
     + \overline{\Big(\sum_{n = 0}^\infty \zeta^n \, \overline{w}^n\Big)} - 1
             = 2 \, \re \sum_{n = 0}^\infty (\zeta \, \overline{w})^n - 1.
\end{equation}
Of course, we can sum the geometric series, to get that
\begin{equation}
\label{sum_{n = 0}^infty (zeta overline{w})^n = ...}
        \sum_{n = 0}^\infty (\zeta \, \overline{w})^n 
                      = \frac{1}{1 - \zeta \, \overline{w}}
 = \frac{1}{(1 - \zeta \, \overline{w})} \,
                \frac{(1 - \overline{\zeta} \, w)}{(1 - \overline{\zeta} \, w)}
 = \frac{1 - \overline{\zeta} \, w}{|1 - \zeta \, \overline{w}|^2}
\end{equation}
for each $\zeta \in U$ and $w \in {\bf T}$.  This implies that
\begin{equation}
\label{P(zeta, w) = ... = frac{1 - |zeta|^2}{|1 - zeta overline{w}|^2}}
 P(\zeta, w) = \frac{2 - 2 \re (\overline{\zeta} \, w) - 
                 |1 - \zeta \, \overline{w}|^2}{|1 - \zeta \, \overline{w}|^2}
             = \frac{1 - |\zeta|^2}{|1 - \zeta \, \overline{w}|^2}
\end{equation}
for every $\zeta \in U$ and $w \in {\bf T}$, since
\begin{equation}
\label{|1 - zeta overline{w}|^2 = ... = 1 - 2 re (overline{zeta} w) + |zeta|^2}
 |1 - \zeta \, \overline{w}|^2 
           = (1 - \zeta \, \overline{w})(1 - \overline{\zeta} \, w)
           = 1 - 2 \, \re (\overline{\zeta} \, w) + |\zeta|^2
\end{equation}
when $|w| = 1$.  In particular, $P(\zeta, w) \ge 0$ for every $\zeta
\in U$ and $w \in {\bf T}$.

        Because $|w| = 1$, the Poisson kernel can also be given by
\begin{equation}
\label{P(zeta, w) = frac{1 - |zeta|^2}{|zeta - w|^2}}
        P(\zeta, w) = \frac{1 - |\zeta|^2}{|\zeta - w|^2}
\end{equation}
for every $\zeta \in U$ and $w \in {\bf T}$.  Thus, for each $\eta > 0$,
we have that
\begin{equation}
\label{P(zeta, w) le eta^{-2} (1 - |zeta|^2)}
        P(\zeta, w) \le \eta^{-2} \, (1 - |\zeta|^2)
\end{equation}
for every $\zeta \in U$ and $w \in {\bf T}$ such that $|\zeta - w| \ge \eta$.
This implies that
\begin{equation}
\label{P(zeta, w) to 0}
        P(\zeta, w) \to 0
\end{equation}
uniformly as $|\zeta| \to 1$ on the set where $|\zeta - w| \ge \eta$,
for each $\eta > 0$.

\section{Continuous functions}
\label{continuous functions, 2}

        Let $f$ be a continuous complex-valued function on the unit
circle, and define a complex-valued function $u$ on the closed unit
disk\index{unit disk} 
\begin{equation}
        \overline{U} = \{\zeta \in {\bf C} : |\zeta| \le 1\}
\end{equation}
by putting $u(\zeta) = f(\zeta)$ when $|\zeta| = 1$ and
$u(\zeta)$ equal to (\ref{sum widehat{f}(n) zeta^n + sum
  widehat{f}(-n) overline{zeta}^n}) when $|\zeta| < 1$, which is the
same as (\ref{frac{1}{2 pi} int_{bf T} P(zeta, w) f(w) |dw|}).
Thus the restrictions of $u$ to the open unit disk $U$ and to the unit
circle ${\bf T}$ are continuous, and one can show that $u$ is actually
continuous on the closed unit disk $\overline{U}$.  This means that
for each $z \in {\bf T}$,
\begin{equation}
\label{lim_{zeta in U atop zeta to z} u(zeta) = f(z)}
        \lim_{\zeta \in U \atop \zeta \to z} u(\zeta) = f(z),
\end{equation}
where the limit as $\zeta \to z$ is taken only over $\zeta \in U$.
This can be derived from the properties of the Poisson kernel
discussed in the previous section, which imply that $u(\zeta)$ is
basically an average of $f$ on ${\bf T}$ that is concentrated near $z$
as $\zeta \in U$ approaches $z$.  Note that $u$ is uniformly continuous
on $\overline{U}$, since continuous functions on compact metric spaces
are always uniformly continuous, and indeed one can also get uniformity
of the limit in (\ref{lim_{zeta in U atop zeta to z} u(zeta) = f(z)})
from the uniform continuity of $f$ on ${\bf T}$ by the same argument.

        Alternatively, put
\begin{equation}
\label{f_r(z) = u(r z)}
        f_r(z) = u(r \, z)
\end{equation}
for each $z \in {\bf T}$ and $r \in [0, 1)$, which is the same as
(\ref{sum_{n = -infty}^infty widehat{f}(n) r^{|n|} z^n = ...}).
This can also be expressed as
\begin{equation}
\label{f_r(z) = frac{1}{2 pi} int_{bf T} P_r(z, w) f(w) |dw|}
        f_r(z) = \frac{1}{2 \pi} \int_{\bf T} P_r(z, w) \, f(w) \, |dw|,
\end{equation}
where
\begin{eqnarray}
\label{P_r(z, w) = P(r z, w) = ... = frac{1 - r^2}{|1 - r z overline{w}|^2}}
 P_r(z, w) = P(r \, z, w) & = & \sum_{n = 0}^\infty r^n \, z^n \, \overline{w}^n
                       + \sum_{n = 1}^\infty r^n \, \overline{z}^n \, w^n \\
  & = & \frac{1 - r^2}{|1 - r \, z \, \overline{w}|^2}         \nonumber
\end{eqnarray}
is another version of the Poisson kernel.\index{Poisson kernel} Under
these conditions,
\begin{equation}
\label{lim_{r to 1-} f_r(z) = f(z)}
        \lim_{r \to 1-} f_r(z) = f(z)
\end{equation}
uniformly over $z \in {\bf T}$.  This follows from the uniform continuity
of $u$ on the closed unit disk, and it can also be derived from the uniform
continuity of $f$ on ${\bf T}$, using the same type of argument as for
(\ref{lim_{zeta in U atop zeta to z} u(zeta) = f(z)}).

        Note that the partial sums of the series in 
(\ref{sum_{n = -infty}^infty widehat{f}(n) r^{|n|} z^n = ...})
converge uniformly over $z \in {\bf T}$ for each $r \in [0, 1)$.
This follows from Weierstrass' $M$-test and the boundedness of the
Fourier coefficients of $f$, as in (\ref{|widehat{f}(n)| le frac{1}{2
    pi} int_{bf T} |f(z)| |dz|}).  In particular, $f_r$ can be approximated
uniformly on ${\bf T}$ by finite linear combinations of the functions
$z^n$ with $n \in {\bf Z}$ for each $r \in [0, 1)$.  This implies that every
continuous function $f$ on ${\bf T}$ can be approximated uniformly on
${\bf T}$ by finite linear combinations of the $z^n$'s, because of the
uniform convergence in (\ref{lim_{r to 1-} f_r(z) = f(z)}).

        As in Section \ref{basic notions}, the linear span of the $z^n$'s
is also dense in $L^2({\bf T})$, because continuous functions are dense
in $L^2({\bf T})$.  Thus the $z^n$'s with $n \in {\bf Z}$ form an
orthonormal basis for $L^2({\bf T})$.  This implies that the Fourier
series of $f \in L^2({\bf T})$ converges in $L^2({\bf T})$, and hence that
\begin{equation}
\label{sum |widehat{f}(n)|^2 = frac{1}{2 pi} int_{bf T} |f(z)|^2 |dz|}
        \sum_{n = -\infty}^\infty |\widehat{f}(n)|^2
                        = \frac{1}{2 \pi} \int_{\bf T} |f(z)|^2 \, |dz|.
\end{equation}

\section{Integrable functions}
\label{integrable functions}

        Let $f$ be a complex-valued integrable function on the unit
circle, and for each $r \in [0, 1)$, let $f_r$ be the function on the
unit circle which is given by (\ref{sum_{n = -infty}^infty
  widehat{f}(n) r^{|n|} z^n = ...}), or equivalently (\ref{f_r(z) =
  frac{1}{2 pi} int_{bf T} P_r(z, w) f(w) |dw|}).  Thus
\begin{equation}
\label{|f_r(z)| le frac{1}{2 pi} int_{bf T} P_r(z, w) |f(w)| |dw|}
 |f_r(z)| \le \frac{1}{2 \pi} \int_{\bf T} P_r(z, w) \, |f(w)| \, |dw|
\end{equation}
for every $z \in {\bf T}$ and $r \in [0, 1)$, since $P_r(z, w) \ge 0$,
as before.  It follows that
\begin{eqnarray}
\label{frac{1}{2 pi} int_{bf T} |f_r(z)| |dz| le ...}
 \frac{1}{2 \pi} \int_{\bf T} |f_r(z)| \, |dz| & \le &  \frac{1}{(2 \pi)^2} 
                \int_{\bf T} \int_{\bf T} P_r(z, w) \, |f(w)| \, |dw| \, |dz| \\
 & = & \frac{1}{(2 \pi)^2} \int_{\bf T} \Big(\int_{\bf T} P_r(z, w) \, |dz|\Big)
                                                \, |f(w)| \, |dw|, \nonumber
\end{eqnarray}
using Fubini's theorem in the second step.  Observe that
\begin{equation}
\label{P_r(z, w) = P_r(w, z)}
        P_r(z, w) = P_r(w, z)
\end{equation}
for every $z, w \in {\bf T}$ and $r \in [0, 1)$, so that
\begin{equation}
\label{frac{1}{2 pi} int_{bf T} P_r(z, w) |dz| = 1}
        \frac{1}{2 \pi} \int_{\bf T} P_r(z, w) \, |dz| = 1
\end{equation}
for every $w \in {\bf T}$ and $r \in [0, 1)$, by (\ref{frac{1}{2 pi}
    int_{bf T} P(zeta, w) |dw| = 1}).  Plugging this into
  (\ref{frac{1}{2 pi} int_{bf T} |f_r(z)| |dz| le ...}), we get that
\begin{equation}
\label{frac{1}{2 pi} int_{bf T} |f_r(z)| |dz| le ..., 2}
        \frac{1}{2 \pi} \int_{\bf T} |f_r(z)| \, |dz|
                          \le \frac{1}{2 \pi} \int_{\bf T} |f(w)| \, |dw|
\end{equation}
for every $r \in [0, 1)$.

        Using (\ref{frac{1}{2 pi} int_{bf T} |f_r(z)| |dz| le ..., 2}),
one can show that
\begin{equation}
\label{lim_{r to 1-} f_r = f}
        \lim_{r \to 1-} f_r = f
\end{equation}
for every $f \in L^1({\bf T})$, where the convergence takes place with
respect to the $L^1$ norm.  More precisely, if $f$ is a continuous
function on ${\bf T}$, then we already know that (\ref{lim_{r to 1-}
  f_r = f}) holds uniformly on ${\bf T}$, and hence with respect to
the $L^1$ norm.  If $f$ is an integrable function on ${\bf T}$, then
one can show that (\ref{lim_{r to 1-} f_r = f}) holds with respect to
the $L^1$ norm, by approximating $f$ by continuous functions on ${\bf
  T}$ with respect to the $L^1$ norm, and using (\ref{frac{1}{2 pi}
  int_{bf T} |f_r(z)| |dz| le ..., 2}) to estimate the errors.

        Suppose now that $f \in L^p({\bf T})$ for some $p$,
$1 < p < \infty$.  Because $|t|^p$ is a convex function on the real line
when $p \ge 1$, one can use (\ref{frac{1}{2 pi} int_{bf T} P(zeta, w)
  |dw| = 1}), (\ref{frac{1}{2 pi} int_{bf T} |f_r(z)| |dz| le ...}),
and Jensen's inequality to get that
\begin{equation}
\label{|f_r(z)|^p le frac{1}{2 pi} int_{bf T} P_r(z, w) |f(w)|^p |dw|}
 |f_r(z)|^p \le \frac{1}{2 \pi} \int_{\bf T} P_r(z, w) \, |f(w)|^p \, |dw|
\end{equation}
for every $z \in {\bf T}$ and $r \in [0, 1)$.  This implies that
\begin{equation}
\label{frac{1}{2 pi} int_{bf T} |f_r(z)|^p |dz| le ...}
        \frac{1}{2 \pi} \int_{\bf T} |f_r(z)|^p \, |dz|
          \le \frac{1}{2 \pi} \int_{\bf T} |f(w)|^p \, |dw|
\end{equation}
for every $r \in [0, 1)$, by integrating (\ref{|f_r(z)|^p le frac{1}{2
      pi} int_{bf T} P_r(z, w) |f(w)|^p |dw|}) over $z \in {\bf T}$
  and interchanging the order of integration, as before.  One can also
  show that (\ref{lim_{r to 1-} f_r = f}) holds with respect to the
  $L^p$ norm when $f \in L^p({\bf T})$ and $1 < p < \infty$, by
  approximating $f$ by continuous functions and using (\ref{frac{1}{2
      pi} int_{bf T} |f_r(z)|^p |dz| le ...}) to estimate the errors
  again.  If $p = 2$, then the series expansion (\ref{sum_{n =
      -infty}^infty widehat{f}(n) r^{|n|} z^n = ...}) for $f_r(z)$
  implies that
\begin{equation}
\label{frac{1}{2 pi} int_{bf T} |f_r(z)|^2 |dz| = sum r^{|n|} |widehat{f}(n)|^2}
        \frac{1}{2 \pi} \int_{\bf T} |f_r(z)|^2 \, |dz|
             = \sum_{n = -\infty}^\infty r^{|n|} \, |\widehat{f}(n)|^2
\end{equation}
for each $r \in [0, 1)$, because of the orthonormality of the $z^n$'s
  in $L^2({\bf T})$.  In particular, this implies (\ref{frac{1}{2 pi}
    int_{bf T} |f_r(z)|^p |dz| le ...}) in this case, because of
  (\ref{sum |widehat{f}(n)|^2 = frac{1}{2 pi} int_{bf T} |f(z)|^2
    |dz|}).  One can also use the series expansion (\ref{sum_{n =
      -infty}^infty widehat{f}(n) r^{|n|} z^n = ...}) and the
  orthonormality of the $z^n$'s to show that (\ref{lim_{r to 1-} f_r =
    f}) holds with respect to the $L^2$ norm.  This is analogous to
  the convergence of the Abel sums of an absolutely convergent series.

        If $f$ is a bounded measurable function on the unit circle,
then (\ref{|f_r(z)| le frac{1}{2 pi} int_{bf T} P_r(z, w) |f(w)|
  |dw|}) and (\ref{frac{1}{2 pi} int_{bf T} P(zeta, w) |dw| = 1})
imply that
\begin{equation}
\label{sup_{z in {bf T}} |f_r(z)| le ||f||_{L^infty({bf T})}}
        \sup_{z \in {\bf T}} |f_r(z)| \le \|f\|_{L^\infty({\bf T})}
\end{equation}
for every $r \in [0, 1)$, where $\|f\|_{L^\infty({\bf T})}$ is the
usual $L^\infty$ norm of $f$, which is the essential supremum of $|f|$
on ${\bf T}$.  However, if $f_r$ converges to $f$ as $r \to 1-$ with
respect to the $L^\infty$ norm, then $f$ has to be the same as a
continuous function almost everywhere on ${\bf T}$.  Remember that the
$L^\infty$ norm of a continuous function on ${\bf T}$ is equal to its
supremum norm, so that the convergence of a sequence of continuous
functions on ${\bf T}$ with respect to the $L^\infty$ norm implies that
the sequence is also a Cauchy sequence with respect to the supremum norm.
This implies that the sequence converges with respect to the supremum
norm on ${\bf T}$, and that the limit is a continuous function on ${\bf T}$.
Of course, the limit of the sequence with respect to the $L^\infty$ norm
is equal to the limit with respect to the supremum norm almost everywhere
on ${\bf T}$.

\section{Borel measures}
\label{borel measures}

        Let $\mu$ be a complex Borel measure on the unit circle.
The \emph{Fourier coefficients}\index{Fourier coefficients} of $\mu$
are defined by
\begin{equation}
\label{widehat{mu}(n) = int_{bf T} overline{z}^n d mu(z)}
        \widehat{\mu}(n) = \int_{\bf T} \overline{z}^n \, d\mu(z)
\end{equation}
for each $n \in {\bf Z}$, and the corresponding \emph{Fourier
  series}\index{Fourier series} is given by
\begin{equation}
\label{sum_{n = -infty}^infty widehat{mu}(n) z^n}
        \sum_{n = -\infty}^\infty \widehat{\mu}(n) \, z^n.
\end{equation}
If $f$ is an integrable function on ${\bf T}$, and if $\mu$ is the
Borel measure defined by
\begin{equation}
\label{mu(A) = frac{1}{2 pi} int_A f(z) |dz|}
        \mu(A) = \frac{1}{2 \pi} \int_A f(z) \, |dz|
\end{equation}
for each Borel set $A \subseteq {\bf T}$, then (\ref{widehat{mu}(n) =
  int_{bf T} overline{z}^n d mu(z)}) and (\ref{sum_{n = -infty}^infty
  widehat{mu}(n) z^n}) are the same as (\ref{widehat{f}(n) = frac{1}{2
    pi} int_{bf T} f(z) overline{z}^n |dz|}) and (\ref{sum_{n =
    -infty}^infty widehat{f}(n) z^n}) for $f$, respectively.
As in (\ref{|widehat{f}(n)| le frac{1}{2 pi} int_{bf T} |f(z)| |dz|}),
\begin{equation}
\label{|widehat{mu}(n)| le |mu|({bf T})}
        |\widehat{\mu}(n)| \le |\mu|({\bf T})
\end{equation}
for each $n \in {\bf Z}$, where $|\mu|$ denotes the total variation
measure on ${\bf T}$ associated to $\mu$.  However, it is not necessary
for $\widehat{\mu}(n)$ to tend to $0$ as $|n| \to \infty$, as one can
see by taking $\mu$ to be a Dirac mass at a point in ${\bf T}$.

        As before, put
\begin{equation}
\label{mu_r(z) = sum_{n = -infty}^infty widehat{mu}(n) r^{|n|} z^n}
        \mu_r(z) = \sum_{n = -\infty}^\infty \widehat{\mu}(n) \, r^{|n|} \, z^n
\end{equation}
for each $z \in {\bf T}$ and $r \in [0, 1)$, so that
\begin{equation}
\label{mu_r(z) = int_{bf T} P_r(z, w) d mu(w)}
        \mu_r(z) = \int_{\bf T} P_r(z, w) \, d\mu(w),
\end{equation}
as in Section \ref{continuous functions, 2}.  Thus
\begin{equation}
\label{|mu_r(z)| le int_{bf T} P_r(z, w) d|mu|(w)}
        |\mu_r(z)| \le \int_{\bf T} P_r(z, w) \, d|\mu|(w)
\end{equation}
for each $z \in {\bf T}$ and $r \in [0, 1)$, as in the previous
  section.  This implies that
\begin{equation}
\label{frac{1}{2 pi} int_{bf T} |mu_r(z)| |dz| le |mu|({bf T})}
        \frac{1}{2 \pi} \int_{\bf T} |\mu_r(z)| \, |dz| \le |\mu|({\bf T})
\end{equation}
for each $r \in [0, 1)$, by interchanging the order of integration and
using (\ref{frac{1}{2 pi} int_{bf T} P_r(z, w) |dz| = 1}), as in the
case of integrable functions.  Note that $\mu_r(z)$ is a nonnegative
real number for each $z \in {\bf T}$ and $r \in [0, 1)$ when $\mu$
is a nonnegative real-valued measure on ${\bf T}$, since the Poisson
kernel is real-valued and nonnegative.

        If $f$ is a continuous complex-valued function on ${\bf T}$,
then it is easy to see that
\begin{equation}
\label{frac{1}{2 pi} int_{bf T} mu_r(z) f(z) |dz| = int_{bf T} f_r(z) d mu(z)}
        \frac{1}{2 \pi} \int_{\bf T} \mu_r(z) \, f(z) \, |dz|
                                     = \int_{\bf T} f_r(z) \, d\mu(z)
\end{equation}
for every $r \in [0, 1)$, where $f_r(z)$ is as in (\ref{f_r(z) = u(r z)}).
We have already seen that $f_r \to f$ uniformly on ${\bf T}$ as $r \to 1-$,
which implies that
\begin{equation}
\label{lim_{r to 1-} int_{bf T} f_r(z) d mu(z) = int_{bf T} f(z) d mu(z)}
 \lim_{r \to 1-} \int_{\bf T} f_r(z) \, d\mu(z) = \int_{\bf T} f(z) \, d\mu(z).
\end{equation}
This shows that
\begin{equation}
\label{lim_{r to 1-} frac{1}{2 pi} int_T mu_r(z) f(z) |dz| = int_T f(z) d mu(z)}
        \lim_{r \to 1-} \frac{1}{2 \pi} \int_{\bf T} \mu_r(z) \, f(z) \, |dz|
                                              = \int_{\bf T} f(z) \, d\mu(z)
\end{equation}
for every continuous function $f$ on ${\bf T}$.

        Of course,
\begin{equation}
\label{lambda(f) = int_{bf T} f(z) d mu(z)}
        \lambda(f) = \int_{\bf T} f(z) \, d\mu(z)
\end{equation}
defines a linear functional on the vector space $C({\bf T})$ of continuous
complex-valued functions on ${\bf T}$.  More precisely, this is a bounded
linear functional on $C({\bf T})$ with respect to the supremum norm, because
\begin{equation}
\label{|lambda(f)| le int_T |f(z)| d |mu|(z) le (sup_{z in T} |f(z)|) |mu|(T)}
        |\lambda(f)| \le \int_{\bf T} |f(z)| \, d|\mu|(z)
                   \le \Big(\sup_{z \in {\bf T}} |f(z)|\Big) \, |\mu|({\bf T})
\end{equation}
for every $f \in C({\bf T})$.  This implies that the dual norm of
$\lambda$ with respect to the supremum norm on $C({\bf T})$ is less
than or equal to $|\mu|({\bf T})$, and one can show that the dual norm
of $\lambda$ is actually equal to $|\mu|({\bf T})$.  Conversely, a
version of the Riesz representation theorem states that every bounded
linear functional on $C({\bf T})$ is of this form, as in Section
\ref{continuous functions}.  Put
\begin{equation}
\label{lambda_r(f) = frac{1}{2 pi} int_{bf T} mu_r(z) f(z) |dz|}
        \lambda_r(f) = \frac{1}{2 \pi} \int_{\bf T} \mu_r(z) \, f(z) \, |dz|
\end{equation}
for each $r \in [0, 1)$, so that $\lambda_r$ is a bounded linear
functional on $C({\bf T})$ for each $r \in [0, 1)$.  Thus (\ref{lim_{r
      to 1-} frac{1}{2 pi} int_T mu_r(z) f(z) |dz| = int_T f(z) d
    mu(z)}) says exactly that $\lambda_r \to \lambda$ as $r \to 1-$
  with respect to the weak$^*$ topology on the dual of $C({\bf T})$.
Note that $\lambda_r$ does not normally converge to $\lambda$ as $r
\to 1-$ with respect to the dual norm associated to the supremum norm
on $C({\bf T})$, which is the same as saying that $(1/2 \pi) \,
\mu_r(z) \, |dz|$ does not normally converge to $\mu$ as $r \to 1-$
with respect to the total variation norm on the space of complex Borel
measures on ${\bf T}$.  This type of convergence would imply that
$\mu_r(z)$ converges in $L^1({\bf T})$ as $r \to 1-$, in which case
$\mu$ would be absolutely continuous, with density equal to the limit
of $\mu_r$ in $L^1({\bf T})$ as $r \to 1-$.

\section{Absolutely convergent series}
\label{absolutely convergent series}

        Let $\{a_j\}_{j = -\infty}^\infty$ be a doubly-infinite sequence
of complex numbers such that
\begin{equation}
\label{sum_{j = -infty}^infty |a_j| = ...}
        \sum_{j = -\infty}^\infty |a_j| = \sum_{j = 0}^\infty |a_j| 
                                         + \sum_{j = 1}^\infty |a_{-j}|
\end{equation}
converges, so that $\{a_j\}_{j = -\infty}^\infty$ corresponds exactly to
a summable function on ${\bf Z}$, as in Section \ref{summable functions}.
Put
\begin{equation}
\label{f(z) = sum_{j = -infty}^infty a_j z^j}
        f(z) = \sum_{j = -\infty}^\infty a_j \, z^j
\end{equation}
for each $z \in {\bf T}$, which may be considered as the sum of the
two absolutely convergent series $\sum_{j = 0}^\infty a_j \, z^j$ and
$\sum_{j = 1}^\infty a_{-j} \, z^{-j}$.  Weierstrass' $M$-test implies
that the partial sums of these series converge uniformly on ${\bf T}$,
and hence that $f$ is a continuous function on ${\bf T}$.  It is easy
to see that
\begin{equation}
\label{widehat{f}(j) = a_j}
        \widehat{f}(j) = a_j
\end{equation}
for each $j \in {\bf Z}$, using the uniform convergence of the partial
sums to interchange the order of summation and integration in the
definition (\ref{widehat{f}(n) = frac{1}{2 pi} int_{bf T} f(z)
  overline{z}^n |dz|}) of the Fourier coefficients of $f$.  Thus
(\ref{f(z) = sum_{j = -infty}^infty a_j z^j}) is the same as the
Fourier series of $f$ in this case.

        Conversely, suppose that $f$ is an integrable function on ${\bf T}$,
and that
\begin{equation}
\label{sum_{j = -infty}^infty |widehat{f}(j)| = ...}
 \sum_{j = -\infty}^\infty |\widehat{f}(j)| = \sum_{j = 0}^\infty |\widehat{f}(j)|
                                       + \sum_{j = 1}^\infty |\widehat{f}(-j)|
\end{equation}
converges, so that the Fourier series (\ref{sum_{n = -infty}^infty
  widehat{f}(n) z^n}) of $f$ converges absolutely for each $z \in {\bf T}$.
Under these conditions, we have seen in Section \ref{integrable functions}
that the Abel sums of the Fourier series of $f$ converge to $f$ with
respect to the $L^1$ norm, which implies that $f$ is equal to the function
defined by its Fourier series almost everywhere on ${\bf T}$.  Similarly,
if $\mu$ is a Borel measure on ${\bf T}$, then we have seen in the
previous section that the Abel sums of the Fourier series of $\mu$
converge to $\mu$ with respect to the weak$^*$ topology on the space
bounded linear functionals on $C({\bf T})$.  If the Fourier coefficients
of $\mu$ are absolutely summable, then it follows that $\mu$ is absolutely
continuous with respect to Lebesgue measure on ${\bf T}$, with density
equal to the function defined by the Fourier series of $\mu$.

        Let $\{a_j\}_{j = -\infty}^\infty$ be a doubly-infinite sequence
of complex numbers for which (\ref{sum_{j    = -infty}^infty |a_j| = ...})
converges, as before, and let $f(z)$ be the corresponding function 
defined on ${\bf T}$ as in (\ref{f(z) = sum_{j = -infty}^infty a_j z^j}).
Also let $\{b_k\}_{k = -\infty}^\infty$ be another doubly-infinite sequence
of complex numbers such that $\sum_{k = -\infty}^\infty |b_k|$ converges,
and put
\begin{equation}
\label{g(z) = sum_{k = -infty}^infty b_k z^k}
        g(z) = \sum_{k = -\infty}^\infty b_k \, z^k
\end{equation}
for each $z \in {\bf T}$.  The product of $f(z)$ and $g(z)$ can be expressed
formally as
\begin{equation}
\label{f(z) g(z) = sum_{n = -infty}^infty c_n z^n}
        f(z) \, g(z) = \sum_{n = -\infty}^\infty c_n \, z^n,
\end{equation}
where
\begin{equation}
\label{c_n = sum_{j = -infty}^infty a_j b_{n - j}}
        c_n = \sum_{j = -\infty}^\infty a_j \, b_{n - j}
\end{equation}
for each integer $n$.  More precisely, it is easy to see that the
series in (\ref{c_n = sum_{j = -infty}^infty a_j b_{n - j}}) converges
absolutely for each $n$, using the convergence of (\ref{sum_{j =
    -infty}^infty |a_j| = ...}) and the fact that the $b_k$'s are uniformly
bounded, because of the convergence of $\sum_{k = -\infty}^\infty |b_k|$.
We also have that
\begin{equation}
\label{|c_n| le sum_{j = -infty}^infty |a_j| |b_{n - j}|}
        |c_n| \le \sum_{j = -\infty}^\infty |a_j| \, |b_{n - j}|
\end{equation}
for each integer $n$, and hence
\begin{equation}
\label{sum_n |c_n| le sum_n (sum_j |a_j| |b_{n - j}|)}
        \sum_{n = -\infty}^\infty |c_n| \le
 \sum_{n = -\infty}^\infty \Big(\sum_{j = -\infty}^\infty |a_j| \, |b_{n - j}|\Big).
\end{equation}
Interchanging the order of summation on the right side of (\ref{sum_n
  |c_n| le sum_n (sum_j |a_j| |b_{n - j}|)}), we get that
\begin{equation}
\label{sum_n |c_n| le sum_j (sum_n |a_j| |b_{n - j}|)}
        \sum_{n = -\infty}^\infty |c_n| \le
 \sum_{j = -\infty}^\infty \Big(\sum_{n = -\infty}^\infty |a_j| \, |b_{n - j}|\Big).
\end{equation}
Of course,
\begin{equation}
\label{sum_n |a_j| |b_{n - j}| = |a_j| sum_n |b_{n -j }| = |a_j| sum_k |b_k|}
        \sum_{n = -\infty}^\infty |a_j| \, |b_{n - j}| 
                           = |a_j| \, \sum_{n = -\infty}^\infty |b_{n -j }| 
                           = |a_j| \sum_{k = -\infty}^\infty |b_k|
\end{equation}
for each $j$, and so we can substitute this into (\ref{sum_n |c_n| le
  sum_j (sum_n |a_j| |b_{n - j}|)}) to get that
\begin{equation}
\label{sum_n |c_n| le (sum_j |a_j|) (sum_k |b_k|)}
 \sum_{n = -\infty}^\infty |c_n| \le \Big(\sum_{j = -\infty}^\infty |a_j|\Big) \, 
                           \Big(\sum_{k = -\infty}^\infty |b_k|\Big).
\end{equation}
This shows that $\sum_{n = -\infty}^\infty |c_n|$ converges, so that
the series on the right side of (\ref{f(z) g(z) = sum_{n =
    -infty}^infty c_n z^n}) converges absolutely for every $z \in {\bf
  T}$.  Similarly, one can check that (\ref{f(z) g(z) = sum_{n =
    -infty}^infty c_n z^n}) holds for every $z \in {\bf T}$, by
interchanging the order of summation.

        Suppose now that $f(z), g(z) \in L^2({\bf T})$, so that their
product $f(z) \, g(z)$ is an integrable function on ${\bf T}$.
Let us check that
\begin{equation}
\label{widehat{(f g)}(n) = sum_j widehat{f}(j) widehat{g}(n - j)}
        \widehat{(f \, g)}(n)
              = \sum_{j = -\infty}^\infty \widehat{f}(j) \, \widehat{g}(n - j)
\end{equation}
for each integer $n$.  Remember that $\sum_{j = -\infty}^\infty
|\widehat{f}(j)|^2$ converges when $f \in L^2({\bf T})$, as in
(\ref{sum_n |widehat{f}(n)|^2 le frac{1}{2 pi} int_{bf T} |f(z)|^2
  |dz|}), and similarly for $g$.  In particular, this implies that the
sum on the right side of (\ref{widehat{(f g)}(n) = sum_j widehat{f}(j)
  widehat{g}(n - j)}) converges absolutely, as in Section
\ref{square-summability}.  If $g(z) = z^l$ for some integer $l$, then
$\widehat{g}(k)$ is equal to $1$ when $k = l$ and to $0$ otherwise, in
which case (\ref{widehat{(f g)}(n) = sum_j widehat{f}(j) widehat{g}(n
  - j)}) can be verified directly from the definitions.  If $g(z)$ is
a linear combination of $z^l$'s for finitely many integers $l$, then
(\ref{widehat{(f g)}(n) = sum_j widehat{f}(j) widehat{g}(n - j)})
follows from the previous case by linearity.  If $g$ is any $L^2$
function on ${\bf T}$, then we have seen that $g(z)$ can be
approximated by linear combinations of $z^l$'s in the $L^2$ norm, and
one can use this to derive (\ref{widehat{(f g)}(n) = sum_j
  widehat{f}(j) widehat{g}(n - j)}) from the preceding case.

\section{Holomorphic functions}
\label{holomorphic functions}

         Let $f$ be an integrable complex-valued function on the unit
circle, and let $u(\zeta)$ be the function on the open unit disk $U$
defined by (\ref{sum widehat{f}(n) zeta^n + sum widehat{f}(-n)
  overline{zeta}^n}).  Suppose that $\widehat{f}(n) = 0$ for each negative
integer $n$, so that
\begin{equation}
\label{u(zeta) = sum_{n = 0}^infty widehat{f}(n) zeta^n}
        u(\zeta) = \sum_{n = 0}^\infty \widehat{f}(n) \, \zeta^n
\end{equation}
is a holomorphic function on the unit disk.  Put $f_r(z) = u(r \, z)$
for each $z \in {\bf T}$ and $r \in [0, 1)$, as in (\ref{f_r(z) = u(r
    z)}), so that
\begin{equation}
\label{f_r(z) = sum_{n = 0}^infty widehat{f}(n) r^n z^n}
        f_r(z) = \sum_{n = 0}^\infty \widehat{f}(n) \, r^n \, z^n.
\end{equation}
Remember that $f_r \to f$ as $r \to 1-$ with respect to the $L^p$ norm
when $f$ is an $L^p$ function on ${\bf T}$ and $1 \le p < \infty$, and
that $f_r \to f$ uniformly on ${\bf T}$ as $r \to 1-$ when $f$ is
continuous on ${\bf T}$, as in Sections \ref{continuous functions, 2} and
\ref{integrable functions}.  It follows that $f$ can be approximated
by finite linear combinations of the $z^n$'s with $n \ge 0$ with
respect to the $L^p$ norm when $f$ is an $L^p$ function on ${\bf T}$
and $1 \le p < \infty$, and with respect to the supremum norm when $f$
is continuous on ${\bf T}$, since the partial sums of (\ref{f_r(z) =
  sum_{n = 0}^infty widehat{f}(n) r^n z^n}) converge to $f_r$ uniformly
on ${\bf T}$ by Weierstrass' $M$-test.

        Conversely, if $u(\zeta)$ is a holomorphic function on $U$, then
Cauchy's theorem implies that
\begin{equation}
\label{oint_{bf T} u(r z) z^n dz = 0}
        \oint_{\bf T} u(r \, z) \, z^n \, dz = 0
\end{equation}
for every $r \in [0, 1)$ and nonnegative integer $n$.  This implies that
\begin{equation}
\label{int_{bf T} u(r z) z^n |dz| = 0}
        \int_{\bf T} u(r \, z) \, z^n \, |dz| = 0
\end{equation}
for each $r \in [0, 1]$ and positive integer $n$, since $dz - i \, z
\, |dz|$ on the unit circle.  If $u$ is continuous on the closed unit
disk $\overline{U}$, and if $f$ is the restriction of $u$ to ${\bf T}$,
then it follows that
\begin{equation}
\label{int_{bf T} f(z) z^n |dz| = 0}
        \int_{\bf T} f(z) \, z^n \, |dz| = 0
\end{equation}
for every positive integer $n$, by taking the limit as $r \to 1-$ in
(\ref{int_{bf T} u(r z) z^n |dz| = 0}).  This also works when $u(r \,
z) \to f(z)$ as $r \to 1-$ with respect to the $L^1$ norm on ${\bf
  T}$, instead of asking that $u$ be continuous on $\overline{U}$.  Of
course, (\ref{int_{bf T} f(z) z^n |dz| = 0}) is the same as saying
that $\widehat{f}(n) = 0$ for each negative integer $n$.

        Suppose now that $f, g \in L^2({\bf T})$ satisfy $\widehat{f}(n)
= \widehat{g}(n) = 0$ for every negative integer $n$.  In this case,
(\ref{widehat{(f g)}(n) = sum_j widehat{f}(j) widehat{g}(n - j)}) implies
that
\begin{equation}
\label{widehat{(f g)}(n) = 0}
        \widehat{(f \, g)}(n) = 0
\end{equation}
when $n < 0$, and
\begin{equation}
\label{widehat{(f, g)}(n) = sum_{j = 0}^n widehat{f}(j) widehat{g}(n - j)}
 \widehat{(f, g)}(n) = \sum_{j = 0}^n \widehat{f}(j) \, \widehat{g}(n - j)
\end{equation}
when $n \ge 0$.  This also works when $f \in L^p({\bf T})$, $g \in
L^q({\bf T})$, and $1 \le p, q \le \infty$ are conjugate exponents,
which is to say that $1/p + 1/q = 1$.  Note that the product $f(z) \,
g(z)$ is an integrable function on ${\bf T}$, because of H\"older's
inequality.  As before, one can verify (\ref{widehat{(f g)}(n) = 0})
and (\ref{widehat{(f, g)}(n) = sum_{j = 0}^n widehat{f}(j)
  widehat{g}(n - j)}) under these conditions by approximating $f$ or
$g$ by finite linear combinations of $z^l$'s with $l \ge 0$.  More
precisely, it is better to approximate $f$ by finite linear
combinations of $z^l$'s when $p = 1$ and $q = \infty$, and similarly
to approximate $g$ by finite linear combinations of $z^l$'s when $q =
1$ and $p = \infty$.  This is a bit simpler than (\ref{widehat{(f
    g)}(n) = sum_j widehat{f}(j) widehat{g}(n - j)}), in that one does
not have to be concerned with convergence of an infinite sum of
products of Fourier coefficients in the present situation.  Once one
has (\ref{widehat{(f g)}(n) = 0}) and (\ref{widehat{(f, g)}(n) =
  sum_{j = 0}^n widehat{f}(j) widehat{g}(n - j)}), it follows that one
gets a holomorphic function on $U$ associated to $f \, g$ as in
(\ref{u(zeta) = sum_{n = 0}^infty widehat{f}(n) zeta^n}), and that
this function is equal to the product of the holomorphic functions on
$U$ associated to $f$ and $g$ in the same way.

        If $\mu$ is a complex Borel measure on ${\bf T}$ such that
$\widehat{\mu}(n) = 0$ for each negative integer $n$, then it can be shown
that $\mu$ is absolutely continuous on ${\bf T}$, so that $\mu$ is defined 
by an integrable function on ${\bf T}$.  This is a famous theorem of
F.\ and M.~Riesz.

\chapter{Topological groups}
\label{topological groups}

\section{Definitions and basic properties}
\label{definitions, basic properties}

        Let $G$ be a group, in which the group operations are
expressed multiplicatively.  Thus the product of $x, y \in G$
is denoted $x \, y$, $e$ is the identity element in $G$, and
$x^{-1}$ is the inverse of $x \in G$.  Suppose that $G$ is also
equipped with a topology, and that the group operations are continuous
with respect to this topology.  More precisely, this means that
multiplication is continuous as a mapping from $(x, y) \in G \times G$
to $x \, y \in G$, using the product topology on $G \times G$ associated
to the given topology on $G$, and that $x \mapsto x^{-1}$ is continuous
as a mapping from $G$ into itself.  In order for $G$ to be a \emph{topological
group},\index{topological groups} it is customary to ask in addition
that $\{e\}$ be a closed subset of $G$.

        The real line ${\bf R}$ is a commutative topological group, with 
respect to addition and the standard topology.  The unit circle ${\bf T}$
is a commutative topological group with respect to multiplication, using
the topology induced on ${\bf T}$ by the standard topology on the
complex plane.  Any group $G$ is a topological group with respect to
the discrete topology on $G$.  We shall be especially interested in
topological groups that are locally compact as topological
spaces,\index{locally compact groups} which includes the examples just
mentioned.

        If $G$ is a topological group and $a \in G$, then the corresponding
left translation
\begin{equation}
\label{x mapsto a x}
        x \mapsto a \, x
\end{equation}
is continuous as a mapping from $G$ into itself, by continuity of
multiplication in $G$.  Similarly, $x \mapsto a^{-1} \, x$ is also
continuous, and is the inverse mapping associated to (\ref{x mapsto a
  x}), so that (\ref{x mapsto a x}) is a homeomorphism from $G$ onto
itself for each $a \in G$.  The same argument shows that the right
translation mapping
\begin{equation}
\label{x mapsto x b}
        x \mapsto x \, b
\end{equation}
is a homeomorphism from $G$ onto itself for each $b \in G$.  It
follows that $\{a\}$ is a closed subset of $G$ for every $a \in G$.

        If $a, b \in G$ and $A, B \subseteq G$, then put
\begin{equation}
\label{a B = {a b : b in B}}
        a \, B = \{a \, b : b \in B\}
\end{equation}
and
\begin{equation}
\label{A b = {a b : a in A}}
        A \, b = \{a \, b : a \in A\}.
\end{equation}
Equivalently, $a \, B$ is the image of $B$ under the left translation
(\ref{x mapsto a x}), and $A \, b$ is the image of $A$ under the right
translation (\ref{x mapsto x b}).  In particular, $a \, B$ and $A \,
b$ have the same topological properties as $A$ and $B$, respectively,
such as being open, closed, compact, or connected.  Also put
\begin{equation}
\label{A B = {a b : a in A, b in B}}
        A \, B = \{a \, b : a \in A, \, b \in B\},
\end{equation}
which is the same as
\begin{equation}
\label{A B = bigcup_{a in A} a B = bigcup_{b in B} A b}
        A \, B = \bigcup_{a \in A} a \, B = \bigcup_{b \in B} A \, b.
\end{equation}
This implies that $A \, B$ is an open set in $G$ when $A$ or $B$ is
an open set.

        Note that $x \mapsto x^{-1}$ is a homeomorphism from $G$ onto
itself, since it is continuous and its own inverse.  Thus
\begin{equation}
\label{A^{-1} = {a^{-1} : a in A}}
        A^{-1} = \{a^{-1} : a \in A\}
\end{equation}
also has the same topological properties as $A$.

        If $W$ is an open set in $G$ such that $e \in W$, then there
are open subsets $U$, $V$ of $G$ such that $e \in U$, $e \in V$, and
\begin{equation}
\label{U V subseteq W}
        U \, V \subseteq W,
\end{equation}
because of continuity of multiplication at $e$.  Let $x$, $y$ be
distinct elements of $G$, so that the set $W$ of $w \in G$ with $w \ne
x^{-1} \, y$ is an open set in $G$ that contains $e$.  If $U$, $V$ are
as before, then (\ref{U V subseteq W}) says that $u \, v \ne x^{-1} \,
y$ for every $u \in U$ and $v \in V$, which implies that $x \, u
\ne y \, v^{-1}$ for every $u \in U$ and $v \in V$, or equivalently
\begin{equation}
        (x \, U) \cap (y \, V^{-1}) = \emptyset.
\end{equation}
This shows that $G$ is Hausdorff as a topological space, since $x \,
U$ and $y \, V^{-1}$ are disjoint open subsets of $G$ that contain $x$
and $y$, respectively.

        Now let $x$ be an element of $G$, and let $E$ be a closed set
in $G$ such that $x \not\in E$.  Thus $x^{-1} \, E$ is a closed set
in $G$ that does not contain $e$, so that its complement $W$ is an
open set that contains $e$.  If $U$ and $V$ are as in the previous
paragraph, then $u \, v \not\in x^{-1} \, E$ for every $u \in U$ and
$v \in V$, which is the same as saying that $x \, u \not\in E \, v^{-1}$
for every $u \in U$ and $v \in V$.  Equivalently,
\begin{equation}
\label{(x U) cap (E V^{-1}) = emptyset}
        (x \, U) \cap (E \, V^{-1}) = \emptyset,
\end{equation}
which implies that $G$ is regular as a topological space, because $x
\, U$ is an open set that contains $x$, and $E \, V^{-1}$ is an open
set that contains $E$.  Note that we could have applied the same
argument to $E \, x^{-1}$ instead of $x^{-1} \, E$, to get open
subsets $\widetilde{U}$, $\widetilde{V}$ of $G$ containing $e$ such that
\begin{equation}
\label{(widetilde{V} x) cap (widetilde{U}^{-1} E) = emptyset}
        (\widetilde{V} \, x) \cap (\widetilde{U}^{-1} \, E) = \emptyset.
\end{equation}
This would also imply that $G$ is regular, since $\widetilde{V} \, x$
and $\widetilde{U}^{-1} \, E$ are disjoint open sets that contain $x$
and $E$, respectively.  Alternatively, one could get open sets like
these by applying the previous argument to $E^{-1}$ and $x^{-1}$, and
then using the mapping $a \mapsto a^{-1}$.

        Let $E$ be a closed set in $G$ again, and let $K$ be a compact
set such that $K \cap E = \emptyset$.  If $x \in K$, then $x \not\in E$,
and hence there are open subsets $U(x)$, $V(x)$ of $G$ containing $e$ such that
\begin{equation}
\label{(x U(x)) cap (E V(x)^{-1}) = emptyset}
        (x \, U(x)) \cap (E \, V(x)^{-1}) = \emptyset,
\end{equation}
as in (\ref{(x U) cap (E V^{-1}) = emptyset}).  Because of continuity of
multiplication at $e$, for each $x \in K$ there is an open set
$U_1(x)$ in $G$ such that $e \in U_1(x)$ and
\begin{equation}
\label{U_1(x) U_1(x) subseteq U(x)}
        U_1(x) \, U_1(x) \subseteq U(x).
\end{equation}
Of course, $K$ is covered by the open sets $x \, U_1(x)$ with $x \in K$,
and hence there are finitely many elements $x_1, \ldots, x_n$ of $K$ such that
\begin{equation}
\label{K subseteq bigcup_{j = 1}^n x_j U_1(x_j)}
        K \subseteq \bigcup_{j = 1}^n x_j \, U_1(x_j),
\end{equation}
by the compactness of $K$.  Put
\begin{equation}
\label{U_1 = bigcap_{j = 1}^n U_1(x_j) and V = bigcap_{j = 1}^n V(x_j)}
        U_1 = \bigcap_{j = 1}^n U_1(x_j) \quad\hbox{and}\quad
                               V = \bigcap_{j = 1}^n V(x_j),
\end{equation}
so that $U_1$, $V$ are open subsets of $G$ containing $e$.  Thus
\begin{equation}
\label{(x_j U_1(x_j) U_1) cap (E V^{-1}) = emptyset}
        (x_j \, U_1(x_j) \, U_1) \cap (E \, V^{-1}) = \emptyset
\end{equation}
for $j = 1, \ldots, n$, by (\ref{U_1(x) U_1(x) subseteq U(x)}),
which implies that
\begin{equation}
\label{(K U_1) cap (E V^{-1}) = emptyset}
        (K \, U_1) \cap (E \, V^{-1}) = \emptyset,
\end{equation}
because of (\ref{K subseteq bigcup_{j = 1}^n x_j U_1(x_j)}).  As
before, one could just as well get that
\begin{equation}
\label{(widetilde{U}_1 K) cap (widetilde{V}^{-1} E) = emptyset}
        (\widetilde{U}_1 \, K) \cap (\widetilde{V}^{-1} \, E) = \emptyset
\end{equation}
for some open sets $\widetilde{U}_1$, $\widetilde{V}$ containing $e$
by an analogous argument, or by applying the previous argument to
$K^{-1}$ and $E^{-1}$ and using the mapping $a \mapsto a^{-1}$.

        If $A$ and $B$ are compact subsets of $G$, then $A \, B$
is also compact, because $A \times B$ is compact with respect to the
product topology, and the group operation is a continuous mapping that
sends $A \times B$ onto $A \, B$.  Of course, compact subsets of $G$
are closed, because $G$ is Hausdorff.  Suppose now that $A$ is compact
and $B$ is closed, and let us check that $A \, B$ is also closed.
If $x \in G$ is not in $A \, B$, then $x \ne a \, b$ for every
$a \in A$ and $b \in B$, so that $a^{-1} \, x \ne b$ for every $a \in A$
and $b \in B$, and hence $(A^{-1} \, x) \cap B = \emptyset$.
Thus we can apply the argument in the previous paragraph to
$K = A^{-1} \, x$ and $E = B$, to get an open set $U_1$ in $G$ such that
$e \in U_1$ and
\begin{equation}
\label{(A^{-1} x U_1) cap B = emptyset}
        (A^{-1} \, x \, U_1) \cap B = \emptyset.
\end{equation}
This is the same as saying that $a^{-1} \, x \, u \ne b$ for each $a
\in A$, $b \in B$, and $u \in U_1$, so that $x \, u \ne a \, b$ for
every $a \in A$, $b \in B$, and $u \in U_1$, and hence
\begin{equation}
\label{(x U_1) cap (A B) = emptyset}
        (x \, U_1) \cap (A \, B) = \emptyset.
\end{equation}
This implies that the complement of $A \, B$ is an open set, so that
$A \, B$ is closed, as desired.  Similarly, if $A$ is closed and $B$
is compact, then $A \, B$ is closed too.

        Let $E$ be any subset of $G$, and let $\overline{E}$ be the closure
of $E$ in $G$.  Thus $x \in \overline{E}$ if and only if every neighborhood
of $x$ in $G$ contains an element of $E$, which is the same as saying that
$(x \, U) \cap E \ne \emptyset$ for every open set $U$ in $G$ with $e \in U$,
and that $(U \, x) \cap E \ne \emptyset$ for every such $U$.  As usual,
$(x \, U) \cap E \ne \emptyset$ if and only if $x \in E \, U^{-1}$, and
similarly $(U \, x) \cap E \ne \emptyset$ if and only if $x \in U^{-1} \, E$.
Hence
\begin{equation}
\label{overline{E} = bigcap E V = bigcap V E}
        \overline{E} = \bigcap E \, V = \bigcap V \, E,
\end{equation}
where the intersection is taken over all open subsets $V$ of $G$
with $e \in V$.  In particular, $\overline{E}$ is contained in $E \, V$
and in $V \, E$ for every such $V$.

\section{Metrizability}
\label{metrizability}

        If $X$ is any topological space, then a simple necessary condition
for the existence of a metric on $X$ that determines the same topology
is that for each $p \in X$ there be a local base for the topology of $X$
at $p$ with only finitely or countably many elements.  Although this
necessary condition is not sufficient for arbitrary topological spaces,
it turns out to be sufficient for topological groups.  Of course, if
a topological group $G$ has a local base for its topology at the identity
$e$ with only finitely or countably many elements, then it has such a
local base at every point, because of continuity of translations.
More precisely, if there is a local base for the topology of $G$ at
$e$ with only finitely many elements, then $\{e\}$ is an open set in
$G$, and hence $G$ is equipped with the discrete topology.

        A metric $d(x, y)$ on a topological group $G$ is said to be invariant
under left translations\index{left-invariant metrics} if
\begin{equation}
\label{d(a x, a y) = d(x, y)}
        d(a \, x, a \, y) = d(x, y)
\end{equation}
for every $a, x, y \in G$.  Similarly, $d(x, y)$ is invariant under right
translations\index{right-invariant metrics} if
\begin{equation}
        d(x \, b, y \, b) = d(x, y)
\end{equation}
for every $b, x, y \in G$.  If there is a local base for the topology
of $G$ at $e$ with only finitely or countably many elements, then a
well-known theorem states that there is a left-invariant metric on $G$
that determines the same topology.  One could instead get a
right-invariant metric, by the same argument, or using the mapping
$x \mapsto x^{-1}$ to switch between the two.

        The standard metric on the real line is invariant under
translations, as is the restriction of the standard metric on the
complex plane to the unit circle as a group with respect to
multiplication.  The discrete metric on any group is invariant under
left and right translations.  If $G$ is a Lie group, then one can
first get a smooth Riemannian metric on $G$ that is invariant under
left translations, by choosing an inner product on the tangent space
at $e$ and extending it to the rest of $G$ using left translations.
If $G$ is also connected, then the corresponding Riemannian distance
function defines a left-invariant metric on $G$.  One can get
right-invariant Riemannian metrics and distance functions in the same
way, and the distance function will be invariant under both left and
right translations when the Riemannian metric is invariant under both
left and right translations.

        Let $G$ be a topological group, and let $\mathcal{B}(e)$ be
a local base for the topology of $G$ at $e$.  If $E$ is any subset of
$G$, then it is easy to see that
\begin{equation}
\label{overline{E} = bigcap_{V in B(e)} E V = bigcap_{V in B(e)} V E}
        \overline{E} = \bigcap_{V \in \mathcal{B}(e)} E \, V
                     = \bigcap_{V \in \mathcal{B}(e)} V \, E,
\end{equation}
by (\ref{overline{E} = bigcap E V = bigcap V E}).  In particular,
if $\mathcal{B}(e)$ is countable, then every closed set in $E$
can be expressed as the intersection of a sequence of open sets.
Equivalently, every open set in $G$ can be expressed as the union
of a sequence of closed sets.  This is a well-known property of
metric spaces, which can be verified directly in this way when
$\mathcal{B}(e)$ is countable.

        Similarly, if $A$ is a dense subset of $G$, then one can check that
the collection of open sets of the form $a \, U$ with $a \in A$ and
$U \in \mathcal{B}(e)$ is a base for the topology of $G$.  If $A$
and $\mathcal{B}(e)$ are countable, then it follows that $G$ has a countable
base for its topology.  This is another well-known property of metric spaces
that can be verified more directly when $\mathcal{B}(e)$ is countable.
Of course, the collection of open sets of the form $U \, a$ with $a \in A$
and $U \in \mathcal{B}(e)$ is also a base for the topology of $G$
under these conditions.

        If $G_1$ and $G_2$ are topological groups, then their Cartesian
product $G_1 \times G_2$ is also a topological group, with respect to
the product topology and group structure.  Note that $G_1 \times G_2$
is locally compact when $G_1$ and $G_2$ are locally compact.  Suppose
that $d_1(x_1, y_1)$ and $d_2(x_2, y_2)$ are metrics on $G_1$ and
$G_2$, respectively, that determine their given topologies. As in
Section \ref{ultrametrics}, (\ref{D((x_1, x_2), (y_1, y_2)) = d_1(x_1,
  y_1) + d_2(x_2, y_2)}) and (\ref{D'((x_1, x_2), (y_1, y_2)) =
  max(d_1(x_1, y_1), d_2(x_2, y_2))}) define metrics on $G_1 \times
G_2$ corresponding to the product topology.  If $d_1(x_1, y_1)$ and
$d_2(x_2, y_2)$ are both invariant under left or right translations,
then it is easy to see that (\ref{D((x_1, x_2), (y_1, y_2)) = d_1(x_1,
  y_1) + d_2(x_2, y_2)}) and (\ref{D'((x_1, x_2), (y_1, y_2)) =
  max(d_1(x_1, y_1), d_2(x_2, y_2))}) have the same property.

        Now let $I$ be an infinite set, and suppose that for each 
$j \in I$  we have a topological group $G_j$.  It is easy to see that 
the Cartesian product $G = \prod_{j \in I} G_j$ is also a topological
group, with respect to the product topology and group structure
again.  If $G_j$ is compact for each $j \in I$, then $\prod_{j \in I}
G_j$ is compact as well, by Tychonoff's theorem.  If $I$ is countably
infinite and $d_j(x_j, y_j)$ is a metric on $G_j$ that determines the
given topology for each $j \in I$, then we can get a metric on $G$
corresponding to the product topology as in (\ref{d(x, y) = max_{j ge
    1} d_j'(x_j, y_j)}) in Section \ref{ultrametrics}.  This metric
is invariant under left translations on $G$ when $d_j(x_j, y_j)$
is invariant under left translations on $G_j$ for each $j \in I$,
and similarly for right translations.

\section{Uniform continuity}
\label{uniform continuity}

        Let $G$ be a topological group, and let $f$ be a 
real or complex-valued function on $G$.  We say that $f$ is
\emph{left uniformly continuous}\index{uniform continuity} along
a set $A \subseteq G$ if for each $\epsilon > 0$ there is an
open set $U \subseteq G$ such that $e \in U$ and
\begin{equation}
\label{|f(u x) - f(x)| < epsilon}
        |f(u \, x) - f(x)| < \epsilon
\end{equation}
for every $x \in A$ and $u \in U$.  Similarly, we say that $f$ is
\emph{right uniformly continuous}\index{uniform continuity} along $A$
if for every $\epsilon > 0$ there is an open set $U \subseteq G$ such
that $e \in U$ and
\begin{equation}
\label{|f(x u) - f(x)| < epsilon}
        |f(x \, u) - f(x)| < \epsilon
\end{equation}
for every $x \in A$ and $u \in U$.  Of course, these two conditions
are equivalent when $G$ is abelian, and they are analogous to uniform
continuity conditions for functions on metric spaces.  If the topology
on $G$ is determined by a metric $d(x, y)$ that is invariant under
right translations, then left uniform continuity can be reformulated
in terms of a uniform continuity condition with respect to $d(x, y)$.
In the same way, if the topology on $G$ is determined by a metric
$d(x, y)$ that is invariant under left translations, the right uniform
continuity can be reformulated in terms of a uniform continuity condition
with respect to $d(x, y)$.

        If $f$ is continuous on $G$ and $A$ is compact, then $f$
is automatically both left and right uniformly continuous along $A$.
This is analogous to uniform continuity properties of continuous
functions on compact subsets of metric spaces.  To see this, let
$\epsilon > 0$ be given, and for each $p \in A$ let $U(p) \subseteq G$
be an open set such that $e \in U(p)$ and
\begin{equation}
\label{|f(u p) - f(p)| < frac{epsilon}{2}}
        |f(u \, p) - f(p)| < \frac{\epsilon}{2}
\end{equation}
for every $u \in U(p)$.  Using the continuity of multiplication at
$e$, we get for each $p \in A$ an open set $U_1(p) \subseteq G$
such that $e \in U_1(p)$ and
\begin{equation}
\label{U_1(p) U_1(p) subseteq U(p)}
        U_1(p) \, U_1(p) \subseteq U(p).
\end{equation}
Thus $U_1(p) \, p$ is an open set that contains $p$ for each $p \in A$,
and compactness of $A$ implies that there are finitely many points
$p_1, \ldots, p_n$ in $A$ such that
\begin{equation}
\label{A subseteq bigcup_{j = 1}^n U_1(p_j) p_j}
        A \subseteq \bigcup_{j = 1}^n U_1(p_j) \, p_j.
\end{equation}
Put $U = \bigcap_{j = 1}^\infty U_1(p_j)$, which is an open set in $G$
that contains $e$.  If $x \in A$ and $u \in U$, then $x \in U_1(p_j)
\, p_j$ for some $j$, $1 \le j \le n$, and hence
\begin{equation}
\label{u x in U U_1(p_j) p subseteq U_1(p_j U_1(p_j) p_j subseteq U(p_j) p_j}
        u \, x \in U \, U_1(p_j) \, p \subseteq U_1(p_j) \, U_1(p_j) \, p_j
                                                  \subseteq U(p_j) \, p_j.
\end{equation}
It follows that
\begin{equation}
\label{|f(u x) - f(x)| le |f(u x) - f(p_j)| + |f(x) - f(p_j)| < ... = epsilon}
        |f(u \, x) - f(x)| \le |f(u \, x) - f(p_j)| + |f(x) - f(p_j)|
               < \frac{\epsilon}{2} + \frac{\epsilon}{2} = \epsilon,
\end{equation}
using (\ref{|f(u p) - f(p)| < frac{epsilon}{2}}) with $p = p_j$ twice
in the second step.  This implies that $f$ is left uniformly
continuous along $A$, and one can show that $f$ is right uniformly
continuous along $A$ in essentially the same way.  Alternatively, one
can derive the right uniform continuity of $f$ along $A$ from the left
uniform continuity of $f(x^{-1})$ along $A^{-1}$.

        If $f$ is a continuous function on $G$ with compact support,
then $f$ is left and right uniformly continuous along $\supp f$, as
before.  It is not difficult to check that $f$ is actually left and
right uniformly continuous on $G$ under these conditions.  One way to
do this is to take the open sets $U$ in the definition of uniform
continuity to be symmetric about $e$, in the sense that $U^{-1} = U$,
by replacing $U$ with $U \cap U^{-1}$.  Similarly, if $f$ is a continuous
function on $G$ that vanishes at infinity, then $f$ is left and right
uniformly continuous on $G$.  If $\{f_j\}_{j = 1}^\infty$ is a
sequence of left or right uniformly continuous functions on $G$ that
converges uniformly to a function $f$ on $G$, then it is easy to see
that $f$ is also left or right uniformly continuous on $G$, as
appropriate, by standard arguments.  In particular, the spaces of left
or right uniformly continuous functions on $G$ that are bounded on $G$
are closed in $C_b(G)$ with respect to the supremum norm.

        Let $f$ be a continuous function on $G$, and for each $a \in G$
let $L_a(f)$ be the function defined on $G$ by
\begin{equation}
\label{(L_a(f))(x) = f(a^{-1} x)}
        (L_a(f))(x) = f(a^{-1} \, x)
\end{equation}
for each $x \in G$.  If $b$ is another element of $G$, then
\begin{eqnarray}
\label{(L_a(L_b(f)))(x) = ... = f((a b)^{-1} x) = (L_{a b}(f))(x)}
 (L_a(L_b(f)))(x) & = & (L_b(f))(a^{-1} \, x) = f(b^{-1} \, a^{-1} \, x) \\
                  & = & f((a \, b)^{-1} \, x) = (L_{a \, b}(f))(x) \nonumber
\end{eqnarray}
for every $x \in G$.  Similarly, let $R_a(f)$ be the function defined
on $G$ by
\begin{equation}
\label{(R_a(f))(x) = f(x a)}
        (R_a(f))(x) = f(x \, a)
\end{equation}
for every $x \in G$, and observe that
\begin{equation}
\label{(R_a(R_b(f)))(x) = (R_b(f))(x a) = f(x a b) = (R_{a b}(f))(x)}
 (R_a(R_b(f)))(x) = (R_b(f))(x \, a) = f(x \, a \, b) = (R_{a \, b}(f))(x)
\end{equation}
for every $a, b, x \in G$.  Of course, $L_a(f)$ and $R_a(f)$ are
continuous functions on $G$ for every $a \in G$ when $f$ is continuous
on $G$, because of continuity of translations.  The condition that
$f$ be left uniformly continuous on $G$ can be reformulated as saying
that $L_a(f) \to f$ uniformly on $G$ as $a \to e$, and the condition
that $f$ be right uniformly continuous on $G$ is equivalent to asking
that $R_a(f) \to f$ uniformly on $G$ as $a \to e$.

\section{Haar measure}
\label{haar measure}

        Let $G$ be a locally compact topological group.  It is well
known that there is a nonnegative Borel measure $H_L$ on $G$
with suitable regularities properties such that $H_L(K) < \infty$
for every compact set $K$ in $G$, $H_L(U) > 0$ for every nonempty
open set $U$ in $G$, and
\begin{equation}
\label{H_L(a E) = H_L(E)}
        H_L(a \, E) = H_L(E)
\end{equation}
for every Borel set $E \subseteq G$ and $a \in G$.  This measure $H_L$
is said to be a left-invariant \emph{Haar measure}\index{Haar measure}
on $G$, and it is unique in the sense that any other Borel measure on
$G$ with the same properties is equal to a positive real number times
$H_L$.  Although we shall not discuss the proof of existence of Haar
measure here, let us mention some basic examples where it is easy to
identify.  If $G$ is equipped with the discrete topology, then we can
simply take $H_L$ to be counting measure on $G$.  Lebesgue measure on
the real line satisfies the requirements of Haar measure with respect
to addition, and arc length measure on the unit circle satisfies the
requirements of Haar measure with respect to multiplication.  If $G$
is a Lie group, then Haar measure on $G$ can be obtained from a
left-invariant volume form.

        Similarly, a nonnegative Borel measure $H_R$ on a locally compact
group $G$ with suitable regularity properties is a right-invariant Haar
measure on $G$ if $H_R(K) < \infty$ for every compact set $K$ in $G$,
$H_R(U) > 0$ for every nonempty open set $U$ in $G$, and
\begin{equation}
\label{H_R(E b) = H_R(E)}
        H_R(E \, b) = H_R(E)
\end{equation}
for every Borel set $E$ in $G$ and $b \in G$.  Observe that $H_R$ is a
right-invariant Haar measure on $G$ if and only if $H_R(E^{-1})$ is a
left-invariant Haar measure on $G$.  In particular, right-invariant
Haar measure is unique up to multiplication by a positive real number.
Of course, left and right-invariant Haar measures on $G$ are the same
when $G$ is commutative.

        Let ${\bf 1}_E(x)$\index{1_E(x)@${\bf 1}_E(x)$} be the indicator 
function\index{indicator functions} associated to a Borel set $E \subseteq G$, 
which is equal to $1$ when $x \in E$ and to $0$ otherwise.  Observe that
\begin{equation}
\label{{bf 1}_{a E}(x) = {bf 1}_E(a^{-1} x)}
        {\bf 1}_{a \, E}(x) = {\bf 1}_E(a^{-1} \, x)
\end{equation}
for every $a, x \in G$.  If $f$ is a nonnegative Borel measurable
function on $G$ and $H_L$ is a left-invariant Haar measure, then
$f(a^{-1} \, x)$ is Borel measurable for every $a \in G$, and
\begin{equation}
\label{int_G f(a^{-1} x) dH_L(x) = int_G f(x) dH_L(x)}
        \int_G f(a^{-1} \, x) \, dH_L(x) = \int_G f(x) \, dH_L(x),
\end{equation}
as one can see by approximating $f$ by simple functions.  If $f$ is a
real or complex-valued function on $G$ which is integrable with
respect to $H_L$, then one can apply this to $|f(x)|$ to get that
$f(a^{-1} \, x)$ is also integrable with respect to $H_L$ for each $a
\in G$, and that (\ref{int_G f(a^{-1} x) dH_L(x) = int_G f(x)
  dH_L(x)}) still holds.  Similarly,
\begin{equation}
\label{{bf 1}_{E b}(x) = {bf 1}_E(x b^{-1})}
        {\bf 1}_{E \, b}(x) = {\bf 1}_E(x \, b^{-1})
\end{equation}
for every $b, x \in G$.  If $f$ is a nonnegative Borel measurable
function on $G$ and $H_R$ is a right-invariant Haar measure, then $f(x
\, b^{-1})$ is Borel measurable for every $b \in G$, and
\begin{equation}
\label{int_G f(x b^{-1}) dH_R(x) = int_G f(x) dH_R(x)}
        \int_G f(x \, b^{-1}) \, dH_R(x) = \int_G f(x) \, dH_R(x).
\end{equation}
If $f$ is a real or complex-valued function on $G$ which is integrable
with respect to $H_R$, then $f(x \, b^{-1})$ is integrable with
respect to $H_R$ for every $b \in G$, and (\ref{int_G f(x b^{-1})
  dH_R(x) = int_G f(x) dH_R(x)}) still holds.

        Let $H_L$ be a left-invariant Haar measure on $G$, and put
\begin{equation}
\label{I_L(f) = int_G f dH_L}
        I_L(f) = \int_G f \, dH_L
\end{equation}
for each $f \in C_{com}(G)$.  This defines a nonnegative linear
functional on $C_{com}(G)$, and $I_L(f) > 0$ when $f$ is a nonnegative
real-valued continuous function with compact support on $G$ such that
$f(x) > 0$ for some $x \in G$.  As in (\ref{int_G f(a^{-1} x) dH_L(x)
  = int_G f(x) dH_L(x)}),
\begin{equation}
\label{I_L(L_a(f)) = I_L(f)}
        I_L(L_a(f)) = I_L(f)
\end{equation}
for every $f \in C_{com}(G)$ and $a \in G$, where $L_a(f)$ is as in
(\ref{(L_a(f))(x) = f(a^{-1} x)}).  Similarly, if $H_R$ is a
right-invariant Haar measure on $G$, then
\begin{equation}
\label{I_R(f) = int_G f dH_R}
        I_R(f) = \int_G f \, dH_R
\end{equation}
defines a nonnegative linear functional on $C_{com}(G)$ such that
$I_R(f) > 0$ when $f$ is a nonnegative real-valued function with
compact support on $G$ such that $f(x) > 0$ for some $x \in G$, and
\begin{equation}
\label{I_R(R_a(f)) = I_R(f)}
        I_R(R_a(f)) = I_R(f)
\end{equation}
for every $f \in C_{com}(G)$ and $a \in G$, where $R_a(f)$ is as in
(\ref{(R_a(f))(x) = f(x a)}).

        Conversely, suppose that $I_L$ is a nonnegative linear functional 
on $C_{com}(G)$ such that $I_L(f) > 0$ when $f$ is a nonnegative
real-valued continuous function with compact support on $G$, and which
is invariant under left translations in the sense that
(\ref{I_L(L_a(f)) = I_L(f)}) for every $f \in C_{com}(G)$.
Under these conditions, the Riesz representation theorem implies that
there is a unique nonnegative Borel measure $H_L$ with suitable
regularity properties such that (\ref{I_L(f) = int_G f dH_L}) holds
for each $f \in C_{com}(G)$.  It is easy to see that $H_L$ is a
left-invariant Haar measure on $G$ under these conditions, and indeed
the existence of Haar measure is often derived from the existence of
an invariant linear functional on $C_{com}(G)$ of this type.
Similarly, if $I_R$ is a nonnegative linear functional on $C_{com}(G)$
such that $I_R(f) > 0$ when $f$ is a nonnegative real-valued continuous
function on $G$ with compact support such that $f(x) > 0$ for some
$x \in G$, and if $I_R$ is invariant under right translations in the sense
that (\ref{I_R(R_a(f)) = I_R(f)}) holds for every $f \in C_{com}(G)$,
then the Riesz representation theorem leads to a right-invariant Haar
measure on $G$.

        If $G_1$ and $G_2$ are locally compact topological groups,
then $G_1 \times G_2$ is also a locally compact topological group,
with respect to the product group structure and topology.  In this
situation, left or right-invariant Haar measure on $G_1 \times G_2$
basically corresponds to the product of the left and right-invariant
Haar measures on $G_1$ and $G_2$, respectively.  As in Section
\ref{double integrals}, there can be some technical issues related to
this, which are easy to handle when there are countable bases for the
topologies of $G_1$ and $G_2$.  Otherwise, one can get left or
right-invariant Haar integrals as nonnegative linear functionals on
$C_{com}(G_1 \times G_2)$, using the corresponding Haar integrals on
$G_1$ and $G_2$.

        Let $I$ be an infinite set, and let $G_j$ be a compact
topological group for each $j \in I$, so that $G = \prod_{j \in I}
G_j$ is a compact topological group with respect to the product group
structure and topology.  Of course, the Haar measure of $G_j$ is
finite for each $j \in I$, because $G_j$ is compact, and we can
normalize it to be equal to $1$.  Again the Haar measure on $G$
basically corresponds to the product of the Haar measures on the
$G_j$'s, which can also be described by the associated Haar integrals,
as in Section \ref{double integrals}.

\section{Left and right translations}
\label{left, right translations}

        Let $G$ be a compact topological group, and let $H_L$ and
$H_R$ be left and right-invariant Haar measures on $G$, respectively.
Also let $f$ be a continuous function on $G$, so that $f(x \, y)$
is a continuous functions of $(x, y) \in G \times G$, and
\begin{equation}
\label{int_G (int_G f(x y)dH_R(x))dH_L(y) = int_G (int_G f(x y)dH_L(y))dH_R(x)}
        \int_G \Big(\int_G f(x \, y) \, dH_R(x)\Big) \, dH_L(y)
          = \int_G \Big(\int_G f(x \, y) \, dH_L(y)\Big) \, dH_R(x).
\end{equation}
This is a version of Fubini's theorem, which can be verified by
approximating $f(x \, y)$ by finite sums of products of continuous
functions of $x$ and $y$.  At any rate, the invariance of $H_R$ under
right translations implies that
\begin{equation}
\label{int_G f(x y) dH_R(x) = int_G f(x) dH_R(x)}
        \int_G f(x \, y) \, dH_R(x) = \int_G f(x) \, dH_R(x)
\end{equation}
for every $y \in G$.  Thus the left side of (\ref{int_G (int_G f(x
  y)dH_R(x))dH_L(y) = int_G (int_G f(x y)dH_L(y))dH_R(x)}) is equal to
\begin{equation}
\label{H_L(G) int_G f(x) dH_R(x)}
        H_L(G) \, \int_G f(x) \, dH_R(x).
\end{equation}
Similarly, the invariance of $H_L$ under left translations implies that
\begin{equation}
\label{int_G f(x y) dH_L(y) = int_G f(y) dH_L(y)}
        \int_G f(x \, y) \, dH_L(y) = \int_G f(y) \, dH_L(y)
\end{equation}
for each $x \in G$, so that the right side of (\ref{int_G (int_G f(x
  y)dH_R(x))dH_L(y) = int_G (int_G f(x y)dH_L(y))dH_R(x)}) is equal to
\begin{equation}
\label{H_R(G) int_G f(y) dH_L(y)}
        H_R(G) \, \int_G f(y) \, dH_L(y).
\end{equation}
The equality of (\ref{H_L(G) int_G f(x) dH_R(x)}) and (\ref{H_R(G)
  int_G f(y) dH_L(y)}) implies that $H_L$ and $H_R$ are positive
constant multiples of each other, and hence that they are both
invariant under both left and right translations.

        Now let $H_L$ be a left-invariant Haar measure on a locally
compact topological group $G$, and let $a$ be an element of $G$.  It
is easy to see that $H_L(E \, a)$ also satisfies the requirements of
left-invariant Haar measure on $G$, so that there is a positive real
number $\phi_L(a)$ such that
\begin{equation}
\label{H_L(E a) = phi_L(a) H_L(E)}
        H_L(E \, a) = \phi_L(a) \, H_L(E)
\end{equation}
for every Borel set $E \subseteq G$.  If $b$ is another element of $G$, then
\begin{equation}
\label{phi_L(a b) H_L(E) = H_L(E a b) = ... = phi_L(a) phi_L(b) H_L(E)}
        \phi_L(a \, b) H_L(E) = H_L(E \, a \, b) = \phi_L(b) \, H_L(E \, a)
                         = \phi_L(a) \, \phi_L(b) \, H_L(E)
\end{equation}
for every Borel set $E \subseteq G$.  Thus
\begin{equation}
\label{phi_L(a b) = phi_L(a) phi_L(b)}
        \phi_L(a \, b) = \phi_L(a) \, \phi_L(b)
\end{equation}
for every $a, b \in G$, so that $\phi_L$ is a homomorphism from $G$
into the multiplicative group ${\bf R}_+$\index{R_+@${\bf R}_+$} of
positive real numbers.

        If $f$ is a nonnegative Borel measurable function on $G$, then
\begin{equation}
\label{int_G f(x a^{-1}) dH_L(x) = phi_L(a) int_G f(x) dH_L(x)}
 \int_G f(x \, a^{-1}) \, dH_L(x) = \phi_L(a) \, \int_G f(x) \, dH_L(x),
\end{equation}
since one can approximate $f$ by simple functions and use (\ref{{bf
    1}_{E b}(x) = {bf 1}_E(x b^{-1})}) and (\ref{H_L(E a) = phi_L(a)
  H_L(E)}).  If $f$ is a nonnegative real-valued continuous function
with compact support on $G$ such that $f(x) > 0$ for some $x \in G$,
then
\begin{equation}
\label{int_G f(x a^{-1}) dH_L(x) > 0}
        \int_G f(x \, a^{-1}) \, dH_L(x) > 0
\end{equation}
for every $a \in G$, and one can check that
\begin{equation}
\label{lim_{a to e} int_G f(x a^{-1}) dH_L(x) = int_G f(x) dH_L(x)}
 \lim_{a \to e} \int_G f(x \, a^{-1}) \, dH_L(x) = \int_G f(x) \, dH_L(x),
\end{equation}
using the fact that $f$ is right uniformly continuous on $G$, as in
Section \ref{uniform continuity}.  This and (\ref{int_G f(x a^{-1})
  dH_L(x) = phi_L(a) int_G f(x) dH_L(x)}) imply that $\phi_L$ is
continuous at $e$, and hence that $\phi_L$ is continuous on $G$,
because $\phi_L$ is a homomorphism.  

        If $f$ is a nonnegative Borel measurable function on $G$ again, then
\begin{eqnarray}
\label{int_G f(x a^{-1}) phi_L(x)^{-1} dH_L(x) = ...}
\lefteqn{\int_G f(x \, a^{-1}) \, \phi_L(x)^{-1} \, dH_L(x)} \\
& = & \phi_L(a)^{-1} \int_G f(x \, a^{-1}) \, \phi_L(x \, a^{-1})^{-1} \, dH_L(x) 
                                                             \nonumber \\
 & = & \int_G f(x) \, \phi_L(x)^{-1} \, dH_L(x)               \nonumber
\end{eqnarray}
for every $a \in G$.  This uses (\ref{int_G f(x a^{-1}) dH_L(x) =
  phi_L(a) int_G f(x) dH_L(x)}) applied to $f(x) \, \phi_L(x)^{-1}$,
and it would also work when $f$ is an integrable function on $G$ with
respect to $H_L$ with compact support, for instance.  Thus
\begin{equation}
\label{E mapsto int_E phi_L(x)^{-1} dH_L(x)}
        E \mapsto \int_E \phi_L(x)^{-1} \, dH_L(x)
\end{equation}
satisfies the requirements of right-invariant Haar measure on $G$.  Of
course, there are analogous statements for the behavior of
right-invariant Haar measure on a locally compact group under left
translations, which can also be derived from the statements for
left-invariant Haar measure using the mapping $x \mapsto x^{-1}$.

\section{Compact subgroups}
\label{compact subgroups}

        Let $G$ be a locally compact topological group, let $H_L$
be a left-invariant Haar measure on $G$, and let $\phi_L$ be defined
on $G$ as in the previous section.  If $K$ is a compact subgroup of $G$,
then $\phi_L(K)$ is a compact subgroup of ${\bf R}_+$.  It is easy to see
that the only compact subgroup of ${\bf R}_+$ is the trivial subgroup
$\{1\}$, so that
\begin{equation}
        \phi_L(x) = 1
\end{equation}
for every $x \in K$.  This shows that left-invariant Haar measure on
$G$ is invariant under right translations by elements of $K$, and
similarly right-invariant Haar measure on $G$ is invariant under left
translations by elements of $K$.

        Now let $G$ be a topological group, and suppose that $d(x, y)$
is a metric on $G$ that determines the same topology.  Let $K$ be a
compact subgroup of $G$ again, and consider
\begin{equation}
\label{d'(x, y) = sup_{a in K} d(a x, a y)}
        d'(x, y) = \sup_{a \in K} d(a \, x, a \, y)
\end{equation}
for each $x, y \in G$.  Note that this is finite for every $x, y \in
G$, because $K \, x$ and $K \, y$ are compact subsets of $G$, and
hence are bounded with respect to the metric.  One can check that the
supremum is actually attained under these conditions, by standard
arguments using continuity and compactness.  It is easy to see that
(\ref{d'(x, y) = sup_{a in K} d(a x, a y)}) defines a metric on $G$,
which is invariant under left translations by elements of $K$ by
construction.  Of course,
\begin{equation}
        d(x, y) \le d'(x, y)
\end{equation}
for every $x, y \in G$, since we can take $a = e$ in (\ref{d'(x, y) =
  sup_{a in K} d(a x, a y)}).  This implies that every open set in $G$
is an open set with respect to $d'(x, y)$, and we would like to show
that $d'(x, y)$ defines the same topology on $G$.

        It suffices to show that for every $x \in G$ and $\epsilon > 0$
there is an open set $U \subseteq G$ such that $e \in U$ and
\begin{equation}
\label{d'(x, x u) le epsilon}
        d'(x, x \, u) \le \epsilon
\end{equation}
for each $u \in U$.  Equivalently, this means that
\begin{equation}
\label{d(a x, a x u) le epsilon}
        d(a \, x, a \, x \, u) \le \epsilon
\end{equation}
for every $a \in K$ and $u \in U$, which is basically a uniform
continuity condition along $K \, x$, as in Section \ref{uniform
  continuity}.  More precisely, this is the same as saying that the
identity mapping on $G$ is right uniformly continuous along $K \, x$
as a mapping from $G$ as a topological group into $G$ as a metric
space with the metric $d(\cdot, \cdot)$.  This can be verified using
the same type of argument as before, because $K \, x$ is compact, and
the identity mapping on $G$ is continuous as a mapping from $G$ as a
topological group into $G$ as a metric space with the metric $d(\cdot,
\cdot)$ by hypothesis.

        Similarly,
\begin{equation}
\label{d''(x, y) = sup_{b in K} d(x b, y b)}
        d''(x, y) = \sup_{b \in K} d(x \, b, y \, b)
\end{equation}
is a metric on $G$ that is invariant under right translations by
elements of $K$ and determines the same topology on $K$.  If we apply
this to $d'(x, y)$ instead of $d''(x, y)$, then we get a metric
\begin{equation}
\label{d'''(x, y) = sup_{a, b in K} d(a x b, a y b)}
        d'''(x, y) = \sup_{a, b \in K} d(a \, x \, b, a \, y \, b)
\end{equation}
that is invariant under both left and right translations by elements
of $K$, and determines the same topology on $G$.  If $G$ is compact
and metrizable, then it follows that there is a metric on $G$ that
is invariant under left and right translations and determines the same
topology on $G$.  Alternatively, if
\begin{equation}
\label{d(x b, y b) = d(x, y)}
        d(x \, b, y \, b) = d(x, y)
\end{equation}
for some $b \in G$ and every $x, y \in G$, then
\begin{equation}
\label{d'(x b, y b) = d'(x, y)}
        d'(x \, b, y \, b) = d'(x, y)
\end{equation}
for every $x, y \in G$ as well.  Thus one can start with a metric
$d(x, y)$ on $G$ that is invariant under right translations by
elements of $G$, and get a metric $d'(x, y)$ that is invariant under
bith right translations by elements of $G$ and left translations by
elements of $K$.  In the same way, one could start with a metric $d(x,
y)$ that is invariant under left translations by elements of $G$, and
get a metric $d''(x, y)$ that is invariant under left translations by
elements of $G$ and right translations by elements of $K$.  Remember
that for any topological group $G$ with a countable local base for its
topology at the identity element $e$, there is a metric on $G$ that
determines the same topology and is invariant under either left or
right translations, as in Section \ref{metrizability}.

\section{Additional properties}
\label{additional properties}

        Let $G$ be a locally compact topological group, and let $H_L$
be a left-invariant Haar measure on $G$.  Suppose that $H_L(G) < \infty$,
and let us show that $G$ is compact.  Let $U$ be an open set in $G$
such that $e \in U$ and $\overline{U}$ is compact.  Thus $H_L(U) > 0$,
and if $x_1, \ldots, x_n$ are elements of $G$ such that
\begin{equation}
\label{(x_j U) cap (x_l U) = emptyset}
        (x_j \, U) \cap (x_l \, U) = \emptyset
\end{equation}
when $j \ne l$, then
\begin{equation}
\label{n H_L(U) = sum_{j - 1}^n H_L(x_j U) = ... le H_L(G)}
        n \, H_L(U) = \sum_{j - 1}^n H_L(x_j \, U) 
                    = H_L\Big(\bigcup_{j = 1}^n x_j \, U\Big) \le H_L(G).
\end{equation}
This implies that $n$ is bounded by $H_L(G) / H_L(U)$, and we suppose
now that $n$ is the largest possible positive integer for which there
exist $x_1, \ldots, x_n \in G$ such that (\ref{(x_j U) cap (x_l U) =
  emptyset}) holds.

        If $y$ is any element of $G$, then the maximality of $n$ implies
that
\begin{equation}
\label{(y U) cap (x_j U) ne emptyset}
        (y \, U) \cap (x_j \, U) \ne \emptyset
\end{equation}
for some $j$, $1 \le j \le n$.  Equivalently, $y \in x_j \, U \, U^{-1}$,
so that
\begin{equation}
\label{G subseteq bigcup_{j = 1}^n x_j U U^{-1} ...}
        G \subseteq \bigcup_{j = 1}^n x_j \, U \, U^{-1} 
           \subseteq \bigcup_{j = 1}^n x_j \, \overline{U} \, \overline{U}^{-1}.
\end{equation}
By construction, $\overline{U}$ is compact, which implies that
$\overline{U} \, \overline{U}^{-1}$ is compact, and hence that the
right side of (\ref{G subseteq bigcup_{j = 1}^n x_j U U^{-1} ...}) is
compact.  Of course, the right side of (\ref{G subseteq bigcup_{j =
    1}^n x_j U U^{-1} ...}) is contained in $G$, so that they are the
same, and thus $G$ is compact, as desired.

        Let $G$ be a topological group, and let $H$ be a subgroup of $G$.
It is well known that the complement of $H$ in $G$ can be expressed as
a union of cosets of $H$, which are translates of $H$.  In particular,
if $H$ is an open subset of $G$, then every translate of $H$ is an
open set, and hence the complement of $H$ is an open set.  This shows
that open subgroups of $G$ are automatically closed sets.  If $G$ is
connected as a topological space, then it follows that $G$ is the only
open subgroup of itself.

        Suppose that $V$ is an open set in $G$ that contains $e$ and is
symmetric in the sense that $V^{-1} = V$.  Let $V^n$ be $V \, V \cdots V$,
with $n$ $V$'s, or equivalently $V^1 = V$ and $V^{n + 1} = V^n \, V$.
It is easy to see that
\begin{equation}
\label{H = bigcup_{n = 1}^infty V^n}
        H = \bigcup_{n = 1}^\infty V^n
\end{equation}
is a subgroup of $G$, which is also an open set, because $V^n$ is an
open set for each $n$.  Note that
\begin{equation}
        V^n \subseteq \overline{(V^n)} \subseteq V^n \, V = V^{n + 1}
\end{equation}
for each $n$, using (\ref{overline{E} = bigcap E V = bigcap V E}) in
Section \ref{definitions, basic properties} in the second step, so that
\begin{equation}
\label{H = bigcup_{n = 1}^infty overline{(V^n)}}
        H = \bigcup_{n = 1}^\infty \overline{(V^n)}.
\end{equation}
Observe also that $(\overline{V})^n \subseteq \overline{(V^n)}$ for
each $n$, by continuity of multiplication.

        If $G$ is locally compact, then we can choose $V$ so that 
$\overline{V}$ is compact, and hence $(\overline{V})^n$ is compact for 
each $n$.  In particular, $(\overline{V})^n$ is a closed set for each
$n$.  which implies that $\overline{(V^n)} \subseteq
(\overline{V})^n$, because $V^n \subseteq (\overline{V})^n$.  Thus
$\overline{(V^n)} = (\overline{V})^n$ for each $n$ when $\overline{V}$
is compact, and $H = \bigcup_{n = 1}^\infty (\overline{V})^n$
is $\sigma$-compact.

        If a topological space $X$ is $\sigma$-compact, then every
closed set in $X$ is $\sigma$-compact too, because the intersection
of a closed set and a compact set is compact.  If the topology
on $X$ is determined by a metric, then it is well known that every open
set in $X$ can be expressed as a countable union of closed sets.
This was mentioned in Section \ref{metrizability}, where an analogous
argument was given for a topological group with a countable local base
for its topology at the identity element.  It follows that open sets
are also $\sigma$-compact under these conditions.  Note that a
$\sigma$-compact metric space is separable, because compact metric
spaces are separable, and hence has a countable base for its topology.
If a locally compact topological space $X$ has a countable base for
its topology, then $X$ is $\sigma$-compact.  This is because $X$
is covered by open sets contained in compact sets, and the existence
of a countable base for the topology of $X$ implies that this open
covering can be reduced to a subcovering with only finitely or countable
many elements.

        If a locally compact topological group $G$ is $\sigma$-compact,
then left and right-invariant Haar measure on $G$ are both $\sigma$-finite.
Conversely, if left or right-invariant Haar measure on $G$ is
$\sigma$-finite, then $G$ is $\sigma$-compact.  As before, there is an
open subgroup $H$ of $G$ which is $\sigma$-compact, and so it suffices
to show that there are only finitely or countably many left or right 
cosets of $H$ in $G$.  Of course, the left cosets of $H$ are
pairwise-disjoint in $G$, as are the right cosets of $H$.  The main
point is that if a measurable set $E \subseteq G$ has finite left or
right-invariant Haar measure, then the intersection of $E$ with left
or right cosets of $H$ can have positive Haar measure for only
finitely or countably many such cosets.  More precisely, for each
$\epsilon > 0$, there can only be finitely many left or right cosets
of $H$ whose intersection with $E$ has measure at least $\epsilon$,
because $E$ has finite measure.  Applying this to $\epsilon = 1/n$
for each positive integer $n$, it follows that there can only be finitely
or countably many left or right cosets of $H$ whose intersection with
$E$ has positive measure.  If Haar measure on $G$ is $\sigma$-finite,
then there is a sequence $E_1, E_2, E_3, \ldots$ of measurable sets
with finite measure whose union is the whole group.  Because $H$
is an open subgroup, its cosets are nonempty open sets as well,
which have positive Haar measure.  This implies that every coset of $H$
should intersect some $E_j$ in a set of positive measure, and hence
that there are only finitely or countably many cosets of $H$, as desired.

\section{Quotient spaces}
\label{quotient spaces}

        Let $G$ be a group, let $H$ be a subgroup of $G$, and let
$G / H$ be the corresponding quotient space\index{quotient spaces}
of left cosets of $H$ in $G$.  Also let $q$ be the canonical quotient
mapping from $G$ onto $G / H$, which sends each $a \in G$ to the
corresponding left coset $a \, H$.  If $g \in G$, then the left
translation mapping $a \mapsto g \, a$ leads to a natural mapping from
$G / H$ onto itself, which sends a left coset $a \, H$ to $g \, a \,
H$.  Of course, if $H$ is a normal subgroup of $G$, then the quotient
$G / H$ is a group in a natural way, and the quotient mapping $q$ is a
homomorphism from $G$ onto $G / H$.

        Suppose now that $G$ is a topological group, and consider
the corresponding quotient topology on $G / H$.  By definition,
this means that a set $W \subseteq G / H$ is an open set if and only
if $q^{-1}(W)$ is an open set in $G$.  Equivalently, $E \subseteq G /
H$ is a closed set if and only if $q^{-1}(E)$ is a closed set in $G$.
In particular, the quotient mapping $q : G \to G / H$ is automatically
continuous with respect to the quotient topology on $G / H$.  Observe
that
\begin{equation}
\label{q^{-1}(q(A)) = A H}
        q^{-1}(q(A)) = A \, H
\end{equation}
for every $A \subseteq G$.  If $A$ is an open set in $G$, then $A \,
H$ is also an open set in $G$, so that $q(A)$ is an open set in $G /
H$.  This shows that $q$ is an open mapping from $G$ onto $G / H$.

        Suppose from now on in this section that $H$ is a closed
subgroup of $G$.  Of course, $H$ is a coset of itself, and hence an
element of $G / H$.  Because of the way that the quotient topology on
$G / H$ is defined, $H$ is a closed subgroup of $G$ if and only if the
subset of $G / H$ consisting of the one coset $H$ is a closed set.  It
is easy to see that the mappings on $G / H$ corresponding to left
translations on $G$ are homeomorphisms with respect to the quotient
topology, since left translations are homeomorphisms on $G$.  This
implies that every subset of $G / H$ with exactly one element is
closed with respect to the quotient topology.

        Let $a$, $b$ be elements of $G$ such that $a \, H \ne b \, H$,
so that $b^{-1} \, a \not\in H$.  Because $H$ is a closed subgroup of $G$
and hence the complement of $H$ is an open set, the continuity of the
group operations implies that there are open subsets $U$, $V$ of $G$
such that $e \in U, V$ and $b^{-1} \, v^{-1} \, u \, a \not\in H$ for
every $u \in U$ and $v \in V$.  Equivalently, $u \, a \, H \ne v \, b
\, H$ for every $u \in U$ and $v \in V$, which means that
\begin{equation}
\label{q(U a) cap q(V b) = emptyset}
        q(U \, a) \cap q(V \, b) = \emptyset.
\end{equation}
Thus $q(U \, a)$ and $q(V \, b)$ are disjoint open subsets of $G / H$
containing $a \, H$ and $b \, H$, respectively, so that $G / H$ is Hausdorff
with respect to the quotient topology.

        Similarly, let $E$ be a closed set in $G / H$, so that $q^{-1}(E)$
is a closed set in $G$.  If $a \, H \not\in E$, then $a \not\in
q^{-1}(E)$, and there are open subsets $U$, $V$ of $G$ such that $e
\in U, V$ and $v^{-1} \, u \, a \not\in q^{-1}(E)$ for every $u \in U$
and $v \in V$, by the continuity of the group operations.  If $b \, H
\in E$, then $b \, H \subseteq q^{-1}(E)$, and hence $v^{-1} \, u \, a
\not\in b \, H$ for every $u \in U$, and $v \in V$, so that $b^{-1} \,
v^{-1} \, u \, a \not\in H$.  Thus $u \, a \, H \ne v \, b \, H$ for
every $u \in U$ and $v \in V$ when $b \, H \in E$, which is the same as
saying that $q(b) \in E$.  It follows that
\begin{equation}
        q(U \, a) \cap q(V \, q^{-1}(E)) = \emptyset,
\end{equation}
so that $q(U \, a)$ and $q(V \, q^{-1}(E))$ are disjoint open subsets
of $G / H$ that contain $a \, H$ and $E$, respectively.  This implies
that $G / H$ is regular with respect to the quotient topology.

        As in the previous section, if $H$ is an open subgroup of $G$,
then $H$ is also a closed subgroup of $G$.  In this case, the quotient
topology on $G / H$ is the same as the discrete topology.  If $G$ is
locally compact and $H$ is any closed subgroup of $G$, then it is easy
to see that $G / H$ is locally compact also, because $q$ is both
continuous and open.  If $G$ is any topological group with a countable
local base for its topology at $e$ and $H$ is any closed subgroup of
$G$, then one can check that there is also a countable local base for
the quotient topology on $G / H$ at every point.  Of course, it suffices
to have a countable local base for the topology of $G / H$ at the point
corresponding to the coset $H$, since the mappings on $G / H$ corresponding
to left translations on $G$ are homeomorphisms with respect to the
quotient topology.

        Suppose that $d(x, y)$ is a metric on $G$ that determines
the same topology and which is invariant under right translations 
by elements of $H$, so that
\begin{equation}
\label{d(x h, y h) = d(x, y)}
        d(x \, h, y \, h) = d(x, y)
\end{equation}
for every $x, y \in G$ and $h \in H$.  The corresponding quotient
metric\index{quotient metrics} on $G / H$ is defined by
\begin{equation}
\label{d'(a H, b H) = inf {d(a h_1, b h_2) : h_1, h_2 in H}}
        d'(a \, H, b \, H) = \inf \{d(a \, h_1, b \, h_2) : h_1, h_2 \in H\}
\end{equation}
for every $a \, H, b \, H \in G / H$.  Equivalently,
\begin{equation}
\label{d'(a H, b H) = inf {d(a h, b) : h in H} = inf {d(a, b h) : h in H}}
        d'(a \, H, b \, H) = \inf \{d(a \, h, b) : h \in H\} 
                           = \inf \{d(a, b \, h) : h \in H\}
\end{equation}
for every $a \, H, b \, H \in G / H$, because of (\ref{d(x h, y h) = d(x, y)}).
Clearly
\begin{equation}
\label{d'(a H, b H) = d'(b H, a H) ge 0}
        d'(a \, H, b \, H) = d'(b \, H, a \, H) \ge 0
\end{equation}
for every $a \, h, b \, H \in G / H$.  If $a \, H \ne b \, H$, then $b
\not\in a \, H$, and one can check that $d'(a \, H, b \, H) > 0$ using
(\ref{d'(a H, b H) = inf {d(a h, b) : h in H} = inf {d(a, b h) : h in
    H}}) and the fact that $a \, H$ is a closed set, because $H$ is a
closed subgroup.

        Let us check that 
\begin{equation}
\label{d'(a H, c H) le d'(a H, b H) + d'(b H, c H)}
        d'(a \, H, c \, H) \le d'(a \, H, b \, H) + d'(b \, H, c \, H)
\end{equation}
for every $a \, H, b \, H, c \, H \in G / H$, which is to say that
$d'(a \, H, b \, H)$ satisfies the triangle inequality.  By construction,
\begin{equation}
\label{d'(a H, c H) le d(a h_1, c h_2)}
        d'(a \, H, c \, H) \le d(a \, h_1, c \, h_2)
\end{equation}
for every $h_1, h_2 \in H$, and hence
\begin{equation}
\label{d'(a H, c H) le d(a h_1, b) + d(b, c h_2)}
        d'(a \, H, c \, H) \le d(a \, h_1, b) + d(b, c \, h_2),
\end{equation}
by the triangle inequality.  Using (\ref{d'(a H, b H) = inf {d(a h, b)
    : h in H} = inf {d(a, b h) : h in H}}), we can take the infimum of
(\ref{d'(a H, c H) le d(a h_1, b) + d(b, c h_2)}) over $h_1, h_2 \in
H$ to get (\ref{d'(a H, c H) le d'(a H, b H) + d'(b H, c H)}), as
desired.  Thus $d'(a \, H, b \, H)$ defines a metric on $G / H$, and
it is easy to see that the topology on $G / H$ corresponding to $d'(a
\, H, b \, H)$ is the same as the quotient topology.  More precisely,
$q$ maps the open ball in $G$ centered at a point $a$ with radius $r >
0$ with respect to $d(x, y)$ onto the open ball in $G / H$ centered
at $a \, H$ with radius $r$ with respect to $d'(a \, H, b \, H)$.

        If $d(x, y)$ is also invariant under left translations on $G$, 
then it is easy to see that $d'(a \, H, b \, H)$ is invariant under the
induced action of left translations on $G / H$, so that
\begin{equation}
\label{d'(g a H, g b H) = d'(a H, b H)}
        d'(g \, a \, H, g \, b \, H) = d'(a \, H, b \, H)
\end{equation}
for every $a, b, g \in G$.  In particular, if there is a countable
local base for the topology of $G$ at $e$, and if $H$ is compact, then
there is a metric $d(x, y)$ that is invariant under left translations
by elements of $G$, as well as invariant under right translations by
elements of $H$, and which determines the same topology on $G$, as in
Sections \ref{metrizability} and \ref{compact subgroups}.

        The action of $G$ on $G / H$ by left translations corresponds
to a mapping from $(g, a \, H) \in G \times (G / H)$ to $g \, a \, H
\in G / H$.  It is easy to see that this mapping is continuous,
using the quotient topology on $G / H$, and the associated product
topology on $G \times (G / H)$.  If $H$ is a normal subgroup of $G$,
then $G / H$ is a group in a natural way, and the quotient mapping $q$
is a homomorphism from $G$ onto $G / H$.  One can check that
$G / H$ is a topological group with respect to the quotient topology
under these conditions.

\section{Invariant measures}
\label{invariant measures}

        Let $G$ be a locally compact topological group, and let $H$ be 
a closed subgroup of $G$, as in the previous section.  Also let $f(x)$
be a continuous real or complex-valued function on $G$ with compact
support.  Thus $f_x(h) = f(x \, h)$ may be considered as a continuous
function of $h \in H$ with compact support for each $x \in G$.
More precisely,
\begin{equation}
\label{supp f_x subseteq (x^{-1} supp f) cap H}
        \supp f_x \subseteq (x^{-1} \supp f) \cap H
\end{equation}
for each $x \in G$.  If $A$ is a compact subset of $G$, then $A^{-1}
\, \supp f$ is also a compact set, and
\begin{equation}
\label{supp f_x subseteq (A^{-1} supp f) cap H}
        \supp f_x \subseteq (A^{-1} \, \supp f) \cap H
\end{equation}
for each $x \in A$.

        Note that $H$ is a locally compact topological group with respect
to the topology induced by the one on $G$.  Let $f_H(x)$ be the integral
of $f_x(h) = f(x \, h)$ as a function of $h \in H$ with respect to a 
left-invariant Haar measure on $H$.  It is not difficult to check that
$f_H(x)$ is continuous in $x$, using the uniform continuity of $f$
along compact sets, as in Section \ref{uniform continuity}.  This also
uses the fact that if $U \subseteq G$ is an open set such that $x \in U$
and $\overline{U}$ is compact, then $\overline{U}^{-1} \, \supp f$
is a compact set and
\begin{equation}
\label{supp f_y subseteq (overline{U}^{-1} supp f) cap H}
        \supp f_y \subseteq (\overline{U}^{-1} \, \supp f) \cap H
\end{equation}
for every $y \in \overline{U}$, as in (\ref{supp f_x subseteq (A^{-1}
  supp f) cap H}).

        If $a \in H$, then the integral of $f_{x \, a}(h) = f(x \, a \, h)
= f_x(a \, h)$ as a function of $h \in H$ with respect to left-invariant 
Haar measure on $H$ is equal to the integral of $f_x(h) = f(x \, h)$ as
a function of $h$.  This implies that
\begin{equation}
\label{f_H(x a) = f_H(x)}
        f_H(x \, a) = f_H(x)
\end{equation}
for every $x \in G$ and $a \in H$, so that $f_H(x)$ is constant on
left cosets of $H$ in $G$.  Equivalently, there is a function
$\widetilde{f}_H$ on $G / H$ such that
\begin{equation}
\label{f_H = widetilde{f}_H circ q}
        f_H = \widetilde{f}_H \circ q,
\end{equation}
where $q$ is the quotient mapping from $G$ onto $G / H$, as in the
previous section.  It is easy to see that $\widetilde{f}_H$ is
continuous with respect to the usual quotient topology on $G / H$,
because $f_H$ is continuous on $G$.  In addition, $\widetilde{f}_H$
has compact support contained in $q(\supp f)$, because $f$ has compact
support.

        If $H$ is a normal subgroup in $G$, then $G / H$ is a locally
compact topological group as well, as in the previous section.
This permits us to integrate $\widetilde{f}_H$ with respect to a
left-invariant Haar measure on $G / H$.  One can check that this
defines a nonnegative linear functional on $C_{com}(G)$ which is invariant
under left translations, since left translations of $f$ on $G$
correspond to left translations of $\widetilde{f}_H$ on $G / H$.
Thus a left-invariant Haar integral on $G$ may be obtained from
left-invariant Haar integrals on $H$ and $G / H$.

        Now let $H$ be a compact subgroup of $G$ that is not necessarily
normal.  Thus $q^{-1}(q(A)) = A \, H$ is compact for every compact set
$A \subseteq G$ in this case.  If $E \subseteq G / H$ is compact, then
$E$ can be covered by finitely many open sets of the form $q(U)$, where
$U$ is an open set in $G$ such that $\overline{U}$ is compact.  This
implies that $q^{-1}(E)$ is compact in $G$, because it is a closed
set which is contained in the union of finitely many compact sets,
by the preceding remark.  If $f$ is a continuous real or complex-valued
function on $G$ with compact support, then it follows that $f \circ q$
is a continuous function on $G$ with compact support.

        Let $H_L$ be a left-invariant Haar measure on $G$.  
It is easy to see that
\begin{equation}
\label{int_G f(q(x)) dH_L(x)}
        \int_G f(q(x)) \, dH_L(x)
\end{equation}
defines a nonnegative linear functional on $C_{com}(G / H)$ which is
invariant under the action of $G$ on $G / H$ by left translations.
Equivalently,
\begin{equation}
\label{H_L(q^{-1}(E))}
        H_L(q^{-1}(E))
\end{equation}
defines a nonnegative Borel measure on $G / H$ which is invariant under
the action of left translations and has other nice properties.

        As another type of situation, suppose that $H$ is a discrete
subgroup of $G$, in the sense that the topology on $H$ induced by the
one on $G$ is the discrete topology.  This means that $\{e\}$ is a 
relatively open set in $H$, and hence that there is an open set 
$U \subseteq G$ such that $U \cap H = \{e\}$.  Let $U_1 \subseteq G$
be another open set such that $e \in U_1$, $U_1^{-1} = U_1$, and
$U_1 \, U_1 \subseteq U$.  If $x, y \in H$ and $x \ne y$, then
it is easy to see that
\begin{equation}
\label{(x U_1) cap (y U_1) = emptyset and (U_1 x) cap (U_1 y) = emptyset}
        (x \, U_1) \cap (y \, U_1) = \emptyset \quad\hbox{and}\quad
        (U_1 \, x) \cap (U_1 \, y) = \emptyset.
\end{equation}
If $a$ is any element of $G$, then it follows that $a \, U_1$ and $U_1
\, a$ can each contain at most one element of $H$.  In particular,
this implies that $H$ is a closed subgroup of $G$ under these
conditions, because it has no limit points in $G$.  Using the second
part of (\ref{(x U_1) cap (y U_1) = emptyset and (U_1 x) cap (U_1 y) =
  emptyset}), we get that
\begin{equation}
\label{(a U_1 x) cap (a U_1 y) = emptyset}
        (a \, U_1 \, x) \cap (a \, U_1 \, y) = \emptyset
\end{equation}
for every $a \in G$ and $x, y \in H$ with $x \ne y$, which implies that
the restriction of the quotient mapping $q$ to $a \, U_1$ is one-to-one.

        In this case, one normally starts with a Borel measure $\mu$ on $G$
that is invariant under right translations by elements of $H$, and
tries to use it to get a measure on $G / H$ that corresponds to $\mu$
locally under the quotient mapping $q$.  Under suitable compactness or
countability conditions, only finitely or countably many local patches
are needed.  Alternatively, if $G$ and hence $G / H$ are locally
compact, then one can define a nonnegative linear functional on
$C_{com}(G / H)$ using partitions of unity to reduce to the case of
functions supported in a local patch.

\section{Semimetrics}
\label{semimetrics}

        Let $X$ be a set.  A \emph{semimetric}\index{semimetrics}
or \emph{pseudometric}\index{pseudometrics} on $X$ is a nonnegative
real-valued function $d(x, y)$ defined for $x, y \in X$ such that
$d(x, x) = 0$ for every $x \in X$,
\begin{equation}
\label{d(x, y) = d(y, x), 2}
        d(x, y) = d(y, x)
\end{equation}
for every $x, y \in X$, and
\begin{equation}
\label{d(x, z) le d(x, y) + d(y, z), 2}
        d(x, z) \le d(x, y) + d(y, z)
\end{equation}
for every $x, y, z \in X$.  Thus a semimetric is the same as a metric,
but without the requirement that $d(x, y) = 0$ only when $x = y$.

        Similarly, a \emph{seminorm}\index{seminorms} or 
\emph{pseudonorm}\index{pseudonorms} on a real or complex vector space
$V$ is a nonnegative real-valued function $N(v)$ on $V$ such that
\begin{equation}
\label{N(t v) = |t| N(v), 2}
        N(t \, v) = |t| \, N(v)
\end{equation}
for every $v \in V$ and $t \in {\bf R}$ or ${\bf C}$, as appropriate, and
\begin{equation}
\label{N(v + w) le N(v) + N(w), 2}
        N(v + w) \le N(v) + N(w)
\end{equation}
for every $v, w \in V$.  Thus $N(0) = 0$, by (\ref{N(t v) = |t| N(v),
  2}) with $t = 0$, and a seminorm $N(v)$ is a norm when $N(v) > 0$
for each $v \in V$ with $v \ne 0$.  If $N(v)$ is a seminorm on $V$, then
\begin{equation}
\label{d(v, w) = N(v - w), 2}
        d(v, w) = N(v - w)
\end{equation}
defines a semimetric on $V$.

        Let $I$ be a nonempty set, and suppose that for each $j \in I$,
$d_j(x, y)$ is a semimetric on a set $X$.  Let us say that
$\{d_j(x, y)\}_{j \in I}$ is a \emph{nice family of semimetrics}\index{nice 
families of semimetrics} on $X$ if for each $x, y \in X$ with $x \ne y$
there is an $j \in I$ such that $d_j(x, y) > 0$.  In this case, put
\begin{equation}
\label{B_j(x, r) = {x in X : d_j(x, y) < r}}
        B_j(x, r) = \{x \in X : d_j(x, y) < r\}
\end{equation}
for each $x \in X$ and $r > 0$, which is the open ball\index{open
  balls} in $X$ centered at $x$ with radius $r$ associated to
$d_i(\cdot, \cdot)$.  A set $U \subseteq X$ is said to be an open set
with respect to the family of semimetrics $\{d_j(x, y)\}_{j \in I}$ if
for each $x \in U$ there are finitely many indices $j_1, \ldots, j_n
\in I$ and positive real numbers $r_1, \ldots, r_n$ such that
\begin{equation}
\label{bigcap_{l = 1}^n B_{j_l}(x, r_l) subseteq U}
        \bigcap_{l = 1}^n B_{j_l}(x, r_l) \subseteq U.
\end{equation}
It is easy to see that this defines a topology on $X$, which reduces to
the usual topology determined by a metric when $I$ has only one element.

        One can also check that $B_j(x, r)$ is an open set with respect to 
this topology for every $x \in X$, $r > 0$, and $j \in I$, using the
triangle inequality.  If $x, y \in X$ and $x \ne y$, then $d_j(x, y) > 0$ 
for some $j \in I$, and
\begin{equation}
\label{B_j(x, d_j(x, y)/2) cap B_j(y, d_j(x, y)/2) = emptyset}
        B_i(x, d_j(x, y)/2) \cap B_j(y, d_i(x, y)/2) = \emptyset,
\end{equation}
using the triangle inequality again.  Thus the topology on $X$
associated to a nice family of semimetrics is Hausdorff, and one can
show that it is regular as well, in much the same way as for metric
spaces.  Similarly, it is easy to see that
\begin{equation}
\label{f_{p, j}(x) = d_j(x, p)}
        f_{p, j}(x) = d_j(x, p)
\end{equation}
is a continuous real-valued function on $X$ with respect to the
topology associated to this family of semimetrics for each $p \in X$
and $j \in I$, which implies that $X$ is completely regular.

        Conversely, if $X$ is any topological space, and if $f$
is a continuous real or complex-valued function on $X$, then
\begin{equation}
\label{rho_f(x, y) = |f(x) - f(y)|}
        \rho_f(x, y) = |f(x) - f(y)|
\end{equation}
defines a semimetric on $X$.  If a collection of continuous real or
complex-valued functions on $X$ separate points, then the
corresponding collection of semimetrics is a nice family of
semimetrics on $X$.  If $X$ is completely regular, then there is a
nice family of semimetrics on $X$ corresponding to continuous real or
complex-valued functions on $X$ that determines the same topology on
$X$.

        Let $\{d_j(x, y)\}_{j \in I}$ be a nice family of semimetrics
on a set $X$.  If $I$ has only finitely many elements, then it is easy
to see that
\begin{equation}
\label{d(x, y) = max_{j in I} d_j(x, y)}
        d(x, y) = \max_{j \in I} d_j(x, y)
\end{equation}
is a metric on $X$ that determines the same topology.  Suppose now that
$I$ is countably infinite, which we can take to be the set 
${\bf Z}_+$\index{Z_+@${\bf Z}_+$} of positive integers.  Put
\begin{equation}
\label{d_j'(x, y) = min (d_j(x, y), 1/j)}
        d_j'(x, y) = \min (d_j(x, y), 1/j)
\end{equation}
for each $j \in {\bf Z}_+$ and $x, y \in X$, which is also a semimetric
on $X$.  Under these conditions, one can check that
\begin{equation}
\label{d(x, y) = max_{j in {bf Z}_+} d_j'(x, y)}
        d(x, y) = \max_{j \in {\bf Z}_+} d_j'(x, y)
\end{equation}
is a metric on $X$ that determines the same topology.

        If $G$ is a topological group, then a well-known theorem
implies that there is a nice family of left-invariant semimetrics
on $G$ that determine the same topology.  Of course, one could also
get a nice family of right-invariant semimetrics that determine
the same topology.  If $G$ is compact, then one can get a nice family
of right and left-invariant semimetrics on $G$ that determines the
same topology, as in Section \ref{compact subgroups}.  If $H$ is a
closed subgroup of $G$, then one can get semimetrics on $G / H$ from
semimetrics on $G$ that are invariant under right translations
by elements of $H$, as in Section \ref{quotient spaces}.

        Let $V$ be a vector space over the real or complex numbers.
If $V$ is equipped with a topology for which the vector space
operations are continuous and $\{0\}$ is a closed set, then $V$ is
said to be a \emph{topological vector space}.\index{topological vector
  spaces} In particular, a topological vector space is a commutative
topological group with respect to addition.  Let $I$ be a nonempty
set, and suppose that $N_j(v)$ is a seminorm on $V$ for each $j \in
I$.  Let us say that $\{N_j(v)\}_{j \in I}$ is a \emph{nice family of
  seminorms}\index{nice families of seminorms} on $V$ if for each $v
\in V$ with $v \ne 0$ there is a $j \in I$ such that $N_j(v) > 0$.
Equivalently, this means that the corresponding collection of
semimetrics
\begin{equation}
\label{d_j(v, w) = N_j(v - w)}
        d_j(v, w) = N_j(v - w)
\end{equation}
is a nice family of semimetrics on $V$.  One can check that $V$ is a
topological vector space with respect to the topology corresponding to
this family of semimetrics, and more precisely that $V$ is
\emph{locally convex}, in the sense that there is a base for the
topology of $V$ consisting of convex open sets.  Conversely, it is
well known that for every locally convex topological vector space $V$
there is a nice family of seminorms on $V$ that determines the same
topology.

\section{Connectedness}
\label{connectedness}

        If $X$ is any topological space and $x, y \in X$, then put
$x \sim y$ when there is a connected set $E \subseteq X$ such that 
$x, y \in E$.  It is well known and not difficult to check that this
defines an equivalence relation on $X$.  More precisely, this relation
is obviously reflexive and symmetric, and it is transitive because the
union of two connected sets $E_1, E_2 \subseteq X$ is also connected
when $E_1 \cap E_2 \ne \emptyset$.  One can also check that the
equivalence classes in $X$ associated to this equivalence relation are
connected subsets of $X$, known as the \emph{connected
  components}\index{connected components} of $X$.  These are the
maximal connected subsets of $X$, and they are automatically closed
subsets of $X$, because the closure of a connected set is connected.

        Let $G$ be a topological group, and let $x \sim y$ be the
equivalence relation defined on $G$ as in the previous paragraph.
Using continuity of translations, it is easy to see that $x \sim y$
implies that $a \, x \sim a \, y$ and $x \, b \sim y \, b$ for every
$a, b \in G$.  Similarly, $x \sim y$ implies $x^{-1} \sim y^{-1}$,
because of continuity of the mapping $x \mapsto x^{-1}$.  It follows
from these properties that the connected component $H$ of $G$
containing $e$ is a normal subgroup of $G$.  Note that $H$ is closed,
as in the previous paragraph.

        Suppose now that $H$ is any closed connected subgroup of a
topological group $G$.  Let $E$ be a connected subset of $G / H$, and
let us check that $q^{-1}(E)$ is a connected subset of $G$, where $q$
is the usual quotient mapping from $G$ onto $G / H$.  Otherwise, there
are nonempty sets $A, B \subseteq G$ which are separated in the sense that
\begin{equation}
\label{overline{A} cap B = A cap overline{B} = emptyset}
        \overline{A} \cap B = A \cap \overline{B} = \emptyset,
\end{equation}
and for which $A \cup B = q^{-1}(E)$.  If $x \, H \in G / H$, then
$q^{-1}(x \, H) = x \, H$ as a subset of $G$, which is connected
because $H$ is connected.  If $x \, H \in E$, then $x \, H \subseteq
q^{-1}(E) = A \cup B$, so that
\begin{equation}
\label{x H = ((x H) cap A) cup ((x H) cap B)}
        x \, H = ((x \, H) \cap A) \cup ((x \, H) \cap B),
\end{equation}
where $(x \, H) \cap A$ and $(x \, H) \cap B$ are separated, because
$A$ and $B$ are separated.  It follows that $(x \, H) \cap A =
\emptyset$ or $(x \, H) \cap B = \emptyset$, since $x \, H$ is
connected, and hence $x \, H$ is contained in either $A$ or $B$.  Thus
$A = q^{-1}(q(A))$ and $B = q^{-1}(q(B))$, and $q(A)$ and $q(B)$ are
nonempty disjoint subsets of $G / H$ such that $q(A) \cup q(B) =
\emptyset$.  It is not difficult to check that the closure
$\overline{q(A)}$ of $q(A)$ in $G / H$ is equal to $q(\overline{A})$
under these conditions, using the fact that $A = q^{-1}(q(A))$.
Similarly, $\overline{q(B)} = q(\overline{B})$, so that $q(A)$ and
$q(B)$ are separated in $G/H$, since $A$ and $B$ are separated in $G$.
This shows that $E$ is not connected in $G / H$ when $q^{-1}(E)$ is
not connected in $G$, which is the same as saying that $q^{-1}(E)$ is
connected when $E$ is connected, as desired.

        Let $H$ be the connected component of $G$ containing $e$,
which is a closed normal subgroup of $G$, as before.  Note that $x \,
H$ is the connected component of $G$ containing $x$ for every $x \in
G$.  If $E$ is a connected subset of $G / H$, then $q^{-1}(E)$ is a
connected subset of $G$, as in the previous paragraph.  If $E$
contains at least two elements, then $q^{-1}(E)$ contains the union of
two distinct cosets of $H$, contradicting the fact that the cosets of
$H$ are the connected components of $G$.  This shows that $G / H$ is
totally disconnected with respect to the quotient topology, which
means that it does not contain any connected subsets with more than
one element.

        Let $G$ be a topological group again, and suppose that
$K \subseteq G$ is compact, $U \subseteq G$ is an open set, and that
$K \subseteq U$.  Equivalently, the complement $E$ of $U$ is a closed 
set in $G$ which is disjoint from $K$.  As in Section \ref{definitions, 
basic properties}, there is an open set $V \subseteq G$ such that
$e \in V$ and $(K \, V) \cap R = \emptyset$, which means that
\begin{equation}
\label{K V subseteq U}
        K \, V \subseteq U.
\end{equation}
More precisely, this corresponds to (\ref{(K U_1) cap (E V^{-1}) =
  emptyset}), with somewhat different notation, and which actually
gives a stronger conclusion.  We may as well ask also that $V^{-1} =
V$, since otherwise we can replace $V$ with $V \cap V^{-1}$.

        Suppose now that $U$ is a compact open subset of $G$ that contains
$e$ as an element.  The remarks in the preceding paragraph imply that
that there is an open set $V \subseteq G$ such that $e \in V$, $V^{-1} = V$,
and
\begin{equation}
\label{U V subseteq U}
        U \, V \subseteq U.
\end{equation}
Note that $U \, V = U$ under these conditions, because $e \in V$ and
hence $U \subseteq U \, V$.  Similarly, (\ref{U V subseteq U}) implies
that $V \subseteq U$, since $e \in U$.  Plugging this back into
(\ref{U V subseteq U}), we get that $V^2 \subseteq U$.  If $n$ is any
positive integer such that $V^n \subseteq U$, then (\ref{U V subseteq
  U}) implies that $V^{n + 1} \subseteq U$, so that $V^n \subseteq U$
for every $n \in {\bf Z}_+$.  As in Section \ref{additional
  properties}, $H = \bigcup_{n = 1}^\infty V^n$ is an open subgroup of
$G$, which is contained in $U$ in this case.  It follows that $H$ is
also compact under these conditions, because open subgroups are closed
sets, and closed subsets of compact sets are compact.

        A topological space $X$ is said to have topological dimension $0$
if for each point $x \in X$ and open set $W \subseteq X$ with $x \in
W$, there is an open set $U \subseteq X$ such that $x \in U$, $U
\subseteq W$, and $U$ is also a closed set in $X$.  Of course, a
topological group $G$ has topological dimension $0$ if it satisfies
this condition at $x = e$, since translations on $G$ are
homeomorphisms.  If a locally compact topological group $G$ has
topological dimension $0$, then there is a local base for the topology
of $G$ at $e$ consisting of compact open subsets of $G$, and hence a
local base for the topology of $G$ at $e$ consisting of compact open
subgroups of $G$, as in the previous paragraph.  If a topological
space $X$ has topological dimension $0$ and satisfies the first
separation condition, then $X$ is totally disconnected.  Conversely,
it is well known that a locally compact Hausdorff topological space
$X$ has topological dimension $0$ when $X$ is totally disconnected.

        Let $G$ be a topological group, and let $H$ be the connected
component of $G$ containing $e$.  Thus $H$ is a closed normal subgroup 
of $G$, and $G / H$ is totally disconnected, as before.  If $G$ is locally
compact, then $G / H$ is locally compact as well, and hence there is a
local base for the topology of $G / H$ at the identity element consisting
of open subgroups of $G / H$, as in the preceding paragraph.  Of course,
the inverse image of any open subgroup of $G / H$ under the natural
quotient mapping from $G$ onto $G / H$ is an open subgroup of $G$.
If $G$ is locally compact and not connected, then it follows that there
is a proper open subgroup of $G$, because $H$ is a proper subgroup of $G$,
and hence $G / H$ is nontrivial.

\section{Homeomorphism groups}
\label{homeomorphism groups}

        Let $X$ be a topological space, and let $C(X, X)$ be the space
of continuous mappings from $X$ into itself.  This is actually a
semigroup with respect to composition of mappings, with the identity
mapping $\id_X$ on $X$ as the identity element in the semigroup.  The
group $\mathcal{H}(X)$ of all homeomorphisms on $X$ is the same as the
group of invertible elements in $C(X, X)$ as a semigroup.  Let us suppose
from now on in this section that $X$ is compact, and that the topology
on $X$ is determined by a metric $d(x, y)$.  

        In this case, the \emph{supremum metric} on $C(X, X)$ may be 
defined as usual by
\begin{equation}
\label{rho(f, g) = sup_{x in X} d(f(x), g(x))}
        \rho(f, g) = \sup_{x \in X} d(f(x), g(x))
\end{equation}
for each $f, g \in C(X, X)$.  If $h$ is any continuous mapping from
$X$ into itself, then
\begin{equation}
\label{rho(f circ h, g circ h) = sup_{x in X} d(f(h(x)), g(h(x))) le rho(f, g)}
 \rho(f \circ h, g \circ h) = \sup_{x \in X} d(f(h(x)), g(h(x))) \le \rho(f, g)
\end{equation}
for every $f, g \in C(X, X)$.  If we also have that $h(X) = X$, then
we get that
\begin{equation}
\label{rho(f circ h, g circ h) = rho(f, g)}
        \rho(f \circ h, g \circ h) = \rho(f, g)
\end{equation}
for every $f, g \in C(X, X)$.  In particular, the supremum metric is
invariant under right translations on $\mathcal{H}(X)$.

        Let us check that
\begin{equation}
\label{(f, g) mapsto f circ g}
        (f, g) \mapsto f \circ g
\end{equation}
is continuous as a mapping from $C(X, X) \times C(X, X)$ into $C(X,
X)$, using the topology on $C(X, X)$ determined by the supremum
metric, and the associated product topology on $C(X, X) \times C(X,
X)$.  Roughly speaking, to say that this mapping is continuous at a
particular point $(f_0, g_0)$ in $C(X, X) \times C(X, X)$ means that
if $f \in C(X, X)$ is close to $f_0$ and $g \in C(X, X)$ is close to
$g_0$ with respect to the supremum metric, then $f \circ g$ is close
to $f_0 \circ g_0$.  Of course,
\begin{eqnarray}
\label{rho(f circ g, f_0 circ g_0) le ...}
 \rho(f \circ g, f_0 \circ g_0) & \le & \rho(f \circ g, f_0 \circ g)
                                        + \rho(f_0 \circ g, f_0 \circ g_0) \\
        & \le & \rho(f, f_0) + \rho(f_0 \circ g, f_0 \circ g_0), \nonumber
\end{eqnarray}
using (\ref{rho(f circ h, g circ h) = sup_{x in X} d(f(h(x)), g(h(x)))
  le rho(f, g)}) in the second step.  Note that $f_0$ is uniformly
continuous on $X$, because $X$ is compact.  This implies that $f_0
\circ g$ is close to $f_0 \circ g_0$ with respect to the supremum
metric when $g$ is sufficiently close to $g_0$, and hence that $f
\circ g$ is close to $f_0 \circ g_0$ when $f$ is also sufficiently
close to $f_0$, as desired.

        This shows that $C(X, X)$ is a topological semigroup with respect
to the topology determined by the supremum metric, and we would like
to verify that $\mathcal{H}(X)$ is a topological group.  Let us begin
by observing that
\begin{equation}
\label{rho(f^{-1}, {id}_X) = rho(f, {id}_X)}
        \rho(f^{-1}, {\id}_X) = \rho(f, {\id}_X)
\end{equation}
for every $f \in \mathcal{H}(X)$.  This follows from (\ref{rho(f circ
  h, g circ h) = rho(f, g)}), by taking $g = \id_X$ and $h = f^{-1}$.
Hence the mapping $f \mapsto f^{-1}$ is continuous at $\id_X$ on
$\mathcal{H}(X)$, and it remains to check that this mapping is
continuous at any $f_0 \in \mathcal{H}(X)$.  This could be derived
from continuity of $f \mapsto f^{-1}$ at $\id_X$ and continuity of
translations on $\mathcal{H}(X)$, but it is instructive to give a
more explicit argument.  Applying (\ref{rho(f^{-1}, {id}_X) = rho(f,
  {id}_X)}) to $f_0^{-1} \circ f$, we get that
\begin{equation}
\label{rho(f^{-1} circ f_0, {id}_X) = rho(f_0^{-1} circ f, {id}_X)}
        \rho(f^{-1} \circ f_0, {\id}_X) = \rho(f_0^{-1} \circ f, {\id}_X)
\end{equation}
for every $f \in \mathcal{H}(X)$.  We also have that
\begin{equation}
\label{rho(f^{-1}, f_0^{-1}) = ... = rho(f^{-1} circ f_0, {id}_X)}
        \rho(f^{-1}, f_0^{-1}) = \rho(f^{-1} \circ f_0, f_0^{-1} \circ f_0) 
                               = \rho(f^{-1} \circ f_0, {\id}_X),
\end{equation}
as in (\ref{rho(f circ h, g circ h) = rho(f, g)}), so that
\begin{equation}
\label{rho(f^{-1}, f_0^{-1}) = rho(f_0^{-1} circ f, {id}_X)}
        \rho(f^{-1}, f_0^{-1}) = \rho(f_0^{-1} \circ f, {\id}_X)
\end{equation}
for every $f \in \mathcal{H}(X)$.  If $f$ is sufficiently close to
$f_0$ with respect to the supremum metric, then $f_0^{-1} \circ f$ is
close to $f_0^{-1} \circ f_0 = \id_X$, because $f_0$ is uniformly
continuous on $X$.  This implies that $f^{-1}$ is close to $f_0^{-1}$,
by (\ref{rho(f^{-1}, f_0^{-1}) = rho(f_0^{-1} circ f, {id}_X)}),
as desired.

        Suppose that $Y$ is another metric space which is separable,
so that there is a dense set in $Y$ with only finitely or countably
many elements.  It is well known that $C(X, Y)$ is also separable with
respect to the corresponding supremum metric under these conditions.
This uses the fact that continuous mappings from $X$ into $Y$ are
uniformly continuous, because $X$ is compact.  The compactness of $X$
also implies that $X$ can be covered by finitely many balls of
arbitrarily small radius, so that two continuous mappings from $X$
into itself are uniformly close to each other when they are close on
suitable finite subsets of $X$, by uniform continuity.  One can then
use the separability of $Y$ to get that only countably many mappings
are needed to approximate a given mapping from $X$ into $Y$ on a
finite set.  In particular, $C(X, X)$ is separable with respect to
the supremum metric, because compact metric spaces are separable.
Any subset of a separable metric space is also separable
with respect to the induced metric, which implies that
$\mathcal{H}(X)$ is separable as well under these conditions.

        Note that $X$ is complete as a metric space, in the sense that
every Cauchy sequence in $X$ converges to an element of $X$, because
$X$ is compact.  It is well known that $C(X, X)$ is complete with
respect to the supremum metric too.  However, $\mathcal{H}(X)$ is
not normally a closed set in $C(X, X)$, and thus $\mathcal{H}(X)$
is not normally complete with respect to the supremum metric.
The problem is that a sequence of homeomorphisms on $X$ may converge
uniformly to a continuous mapping that is not injective.

        One way to deal with this is to simply use
\begin{equation}
\label{rho(f, g) + rho(f^{-1}, g^{-1})}
        \rho(f, g) + \rho(f^{-1}, g^{-1})
\end{equation}
as the metric on $\mathcal{H}(X)$, which would determine the same
topology on $\mathcal{H}(X)$.  Of course, a sequence $\{f_j\}_{j =
  1}^\infty$ of homeomorphisms on $X$ is a Cauchy sequence with
respect to (\ref{rho(f, g) + rho(f^{-1}, g^{-1})}) if and only if
$\{f_j\}_{j = 1}^\infty$ and $\{f_j^{-1}\}_{j = 1}^\infty$ are Cauchy
sequences with respect to the supremum metric.  In this case, the
completeness of $C(X, X)$ implies that $\{f_j\}_{j = 1}^\infty$
converges uniformly to a continuous mapping $f : X \to X$, and
that $\{f_j^{-1}\}_{j = 1}^\infty$ converges to a continuous mapping
$\widetilde{f} : X \to X$.  Because
\begin{equation}
\label{f_j circ f_j^{-1} = f_j^{-1} circ f_j = {id}_X}
        f_j \circ f_j^{-1} = f_j^{-1} \circ f_j = {\id}_X
\end{equation}
for each $j$, we can use continuity of composition of mappings on
$C(X, X)$ to pass to the limit as $j \to \infty$, to get that
\begin{equation}
\label{f circ widetilde{f} = widetilde{f} circ f = {id}_X}
        f \circ \widetilde{f} = \widetilde{f} \circ f = {\id}_X.
\end{equation}
This implies that $f$ is a homeomorphism of $X$ onto itself, with
$f^{-1} = \widetilde{f}$, and hence that $\{f_j\}_{j = 1}^\infty$
converges to $f$ with respect to (\ref{rho(f, g) + rho(f^{-1},
  g^{-1})}) under these conditions, as desired.

        Alternatively, one can embed $\mathcal{H}(X)$ into
$C(X, X) \times C(X, X)$ using the mapping $f \mapsto (f, f^{-1})$.
The argument in the previous paragraph shows that the image of
$\mathcal{H}(X)$ under this mapping is a closed set in $C(X, X) \times
C(X, X)$ with respect to the product topology corresponding to the
supremum metric on $C(X, X)$.  The completeness of $C(X, X)$ implies
that $C(X, X) \times C(X, X)$ is also complete with respect to an
appropriate product metric, so that $\mathcal{H}(X)$ is complete with
respect to the restriction of such a product metric to the image of
$\mathcal{H}(X)$ in $C(X, X) \times C(X, X)$.

\section{Non-compact spaces}
\label{non-compact spaces}

        Let $(X, d(x, y))$ be a metric space which is not compact.
If $K$ is a nonempty compact subset of $X$ and $f$, $g$ are continuous
mappings from $X$ into itself, then
\begin{equation}
\label{rho_K(f, g) = sup_{x in K} d(f(x), g(x))}
        \rho_K(f, g) = \sup_{x \in K} d(f(x), g(x))
\end{equation}
is finite, because $f(K)$ and $g(K)$ are compact subsets of $X$.  It
is easy to see that this defines a semimetric on $C(X, X)$, and that
the collection of these semimetrics determines a topology on $C(X, X)$
which is Hausdorff, as in Section \ref{semimetrics}.  If $h$ is any
continuous mapping from $X$ into itself, then
\begin{equation}
\label{rho_K(f circ h, g circ h) = rho_{h(K)} (f, g)}
        \rho_K(f \circ h, g \circ h) = \rho_{h(K)} (f, g)
\end{equation}
for every nonempty compact set $K \subseteq X$ and $f, g \in C(X, X)$.
Of course, $h(K)$ is also compact under these conditions, and it
follows that $f \mapsto f \circ h$ is a continuous mapping from
$C(X, X)$ into itself with respect to the topology just defined.

        Suppose from now on in this section that $X$ is locally compact.
We would like to check that $(f, g) \mapsto f \circ g$ is continuous
as a mapping from $C(X, X) \times C(X, X)$ into $C(X, X)$, using the
product topology on $C(X, X) \times C(X, X)$ associated to the
topology on $C(X, X)$ described in the previous paragraph.  Let
$f_0$, $g_0$ in $C(X, X)$ be given, and let us show that $f \circ g$
is close to $f_0 \circ g_0$ in $C(X, X)$ when $f, g \in C(X, X)$ are
sufficiently close to $f_0$, $g_0$, respectively.  More precisely,
if $K$ is any nonempty compact subset of $X$, then we want to show that
$f \circ g$ is uniformly close to $f_0 \circ g_0$ on $K$ when $f$, $g$
are sufficiently close to $f_0$, $g_0$ in $C(X, X)$.  As before,
\begin{eqnarray}
\label{rho_K(f circ g, f_0 circ g_0) le ...}
 \rho_K(f \circ g, f_0 \circ g_0) & \le & \rho_K(f \circ g, f_0 \circ g)
                                   + \rho_K(f_0 \circ g, f_0 \circ g_0) \\
 & = & \rho_{g(K)}(f, f_0) + \rho_K(f_0 \circ g, f_0 \circ g_0), \nonumber
\end{eqnarray}
using (\ref{rho_K(f circ h, g circ h) = rho_{h(K)} (f, g)}) in the
second step.  Of course, $g_0(K)$ is compact in $X$, because $K$ is
compact and $g_0$ is continuous.  Using the local compactness of $X$,
one can get a compact set $L \subseteq X$ that contains $g_0(K)$ in
its interior.  Let us now restrict our attention to $g \in C(X, X)$
that are sufficiently close to $g_0$ on $K$ so that
\begin{equation}
\label{g(K) subseteq L}
        g(K) \subseteq L,
\end{equation}
and hence
\begin{equation}
\label{rho_K(f circ g, f_0 circ g_0) le ..., 2}
        \rho_K(f \circ g, f_0 \circ g_0)
                    \le \rho_L(f, f_0) + \rho_K (f_0 \circ g, f_0 \circ g_0),
\end{equation}
by (\ref{rho_K(f circ g, f_0 circ g_0) le ...}).  Note that $f_0$ is
uniformly continuous on $L$, because $f_0$ is continuous and $L$
is compact.  In order for $f \circ g$ to be close to $f_0 \circ g_0$
on $K$, it suffices that $f$ be sufficiently close to $f_0$ on $L$,
and that $g$ be sufficiently close to $g_0$ on $K$, where the latter
should include (\ref{g(K) subseteq L}) in particular.

        Let us restrict our attention to homeomorphisms on $X$,
and consider the continuity properties of $f \mapsto f^{-1}$ on
$\mathcal{H}(X)$.  If $K$ is a nonempty compact set in $X$ and
$f$ is a homeomorphism on $X$, then
\begin{equation}
\label{rho_K(f^{-1}, {id}_X) = ... = rho_{f^{-1}(K)}({id}_X, f)}
 \rho_K(f^{-1}, {\id}_X) = \rho_{f^{-1}(K)}(f^{-1} \circ f, {\id}_X \circ f) 
                        = \rho_{f^{-1}(K)}({\id}_X, f).
\end{equation}
If $f^{-1}(K)$ is contained in a fixed compact set $L \subseteq X$,
then it follows that $f^{-1}$ is as close as we want to $\id_X$ on $K$
when $f$ is sufficiently close to $\id_X$ on $L$.  Otherwise, in order
to get the continuity of $f \mapsto f^{-1}$ on $\mathcal{H}(X)$, one
can simply use the topology on $\mathcal{H}(X)$ determined by the
collection of supremum seminorms $\rho_K(f, g)$ and $\rho_K(f^{-1},
g^{-1})$ associated to all nonempty compact subsets $K$ of $X$.  
It is easy to see that $(f, g) \mapsto f \circ g$ would still be
continuous as a mapping from $\mathcal{H}(X) \times \mathcal{H}(X)$
into $\mathcal{H}(X)$ with respect to this topology, because of the
discussion in the previous paragraph, and the fact that $(f \circ
g)^{-1} = g^{-1} \circ f^{-1}$.

        If $X = {\bf R}^n$ for some positive integer $n$, for instance,
with the standard metric and topology, then one can avoid this problem.
To see this, let $B_r$ be the closed ball in ${\bf R}^n$ centered at $0$
and with radius $r$ for each nonnegative real number $r$.  If $f$
is a homeomorphism on ${\bf R}^n$ which is sufficiently close to the
identity mapping on $B_{r + 1}$ for some $r > 0$, then well-known results
in topology imply that
\begin{equation}
\label{B_r subseteq f(B_{r + 1})}
        B_r \subseteq f(B_{r + 1}),
\end{equation}
and hence
\begin{equation}
\label{f^{-1}(B_r) subseteq B_{r + 1}}
        f^{-1}(B_r) \subseteq B_{r + 1}.
\end{equation}
This permits one to avoid the problem indicated in the previous
paragraph, to get that $f \mapsto f^{-1}$ is continuous as a mapping
from $\mathcal{H}({\bf R}^n)$ into itself at $\id_{{\bf R}^n}$, with
respect to the topology on $\mathcal{H}({\bf R}^n)$ determined by the
supremum semimetrics $\rho_K(f, g)$.  One can then argue as before
that $f \mapsto f^{-1}$ is continuous at every point in
$\mathcal{H}({\bf R}^n)$, with respect to the same topology.  

        Note that the sequence of semimetrics $\rho_{B_r}(f, g)$ with
$r \in {\bf Z}_+$ determines the same topology on $C({\bf R}^n, {\bf R}^n)$
as the collection of semimetrics $\rho_K(f, g)$ associated to arbitrary
nonempty compact subsets $K$ of ${\bf R}^n$ , because every compact set
$K \subseteq {\bf R}^n$ is contained in $B_r$ for sufficiently large $r$.
This implies that this topology on $C({\bf R}^n, {\bf R}^n)$ can be described
by a single metric, as in Section \ref{semimetrics}.

        Similarly, suppose that the locally compact metric space $X$
is also $\sigma$-compact, which means that $X = \bigcup_{l = 1}^\infty K_l$
for some sequence $K_1, K_2, K_3, \ldots$ of nonempty compact sets.
We may as well ask that $K_l \subseteq K_{l + 1}$ for each $l$,
since we can replace $K_l$ with $\bigcup_{j = 1}^l K_j$ if necessary.
Because $X$ is locally compact, one can also enlarge the $K_l$'s a bit,
so that $K_l$ is contained in the interior of $K_{l + 1}$ for each $l$.
In particular, this implies that the union of the interiors of the
$K_l$'s is equal to $X$.  If $K$ is any compact set in $X$, then it follows
that $K$ is contained in the union of the interiors of $K_l$ for finitely
many $l$, and hence that $K \subseteq K_l$ when $l$ is sufficiently large.

        As before, this implies that the sequence of semimetrics 
$\rho_{K_l}(f, g)$  with $l \in {\bf Z}_+$ determines the same topology 
on $C(X, X)$ as the collection of all semimetrics $\rho_K(f, g)$
associated to arbitary nonempty compact subsets $K$ of $X$ under these
conditions.  In the same way, the countable family of semimetrics of the form
$\rho_{K_l}(f, g)$ and $\rho_{K_l}(f^{-1}, g^{-1})$ with $l \in {\bf Z}_+$
determines the same topology on $\mathcal{H}(X)$ as the collection of all
semimetrics of the form $\rho_K(f, g)$ and $\rho_K(f^{-1}, g^{-1})$ for
some nonempty compact set $K \subseteq X$.  Hence these topologies
on $C(X, X)$ and $\mathcal{H}(X)$ can each be described by a single
metric, as in Section \ref{semimetrics}.

        It is natural to embed $\mathcal{H}(X)$ into $C(X, X) \times C(X, X)$
using the mapping $f \mapsto (f, f^{-1})$, as in the preceding section,
whether or not $X$ is $\sigma$-compact.  As usual, $C(X, X) \times
C(X, X)$ may be equipped with the product topology corresponding to
the topology on $C(X, X)$ determined by the collection of semimetrics
$\rho_K(f, g)$ associated to nonempty compact subsets $K$ of $X$.  By
construction, the topology on $\mathcal{H}(X)$ determined by the
collection of semimetrics $\rho_K(f, g)$ and $\rho_K(f^{-1}, g^{-1})$
associated to nonempty compact subsets $K$ of $X$ corresponds exactly
to the one induced on the image of $\mathcal{H}(X)$ in $C(X, X) \times
C(X, X)$ by the product topology just mentioned.  Note that the image
of $\mathcal{H}(X)$ in $C(X, X) \times C(X, X)$ is a closed set with
respect to the product topology, for basically the same reasons as
before.

        If $X$ is $\sigma$-compact again, then it is easy to see that
$X$ is separable, because compact metric spaces are separable.
This implies that the space $C(K, X)$ of continuous mappings from a
compact metric space $K$ into $X$ is separable with respect to the
supremum metric, as in the previous section.  In particular, this can
be applied to compact sets $K \subseteq X$, using the restriction of
the metric on $X$ to $K$.  If $X$ is locally compact and
$\sigma$-compact, then one can use this to show that $C(X, X)$ is
separable with respect to the usual topology.  This implies that $C(X,
X) \times C(X, X)$ is separable with respect to the product topology,
and hence that $\mathcal{H}(X)$ is separable with respect to the
topology determined by the usual collection of semimetrics of the form
$\rho_K(f, g)$ and $\rho_K(f^{-1}, g^{-1})$.

\section{Groups of isometries}
\label{groups of isometries}

        Let $(X, d(x, y))$ be a metric space again, and let $\mathcal{I}(X, X)$
be the set of mappings $f : X \to X$ which are isometries, so that
\begin{equation}
        d(f(x), f(y)) = d(x, y)
\end{equation}
for every $x, y \in X$.  This is a sub-semigroup of the semigroup
$C(X, X)$ of all continuous mappings from $X$ into itself, which
includes the identity mapping $\id_X$ on $X$.  Similarly, the
collection $\mathcal{I}(X)$ of isometric mappings from $X$ onto itself
is the same as the group of invertible elements in $\mathcal{I}(X,
X)$, and is a subgroup of the group $\mathcal{H}(X)$ of all
homeomorphisms from $X$ onto itself.

        Suppose for the moment that $X$ is compact, and let $\rho(f, g)$
be the supremum metric on $C(X, X)$.  If $h \in \mathcal{I}(X, X)$, then
\begin{equation}
\label{rho(h circ f, h circ g) = rho(f, g)}
        \rho(h \circ f, h \circ g) = \rho(f, g)
\end{equation}
for every $f, g \in C(X, X)$, and in particular the restriction of the
supremum metric to $\mathcal{I}(X)$ is invariant under left as well as
right translations.  Thus
\begin{eqnarray}
\label{rho(f' circ g', f circ g) le ... le rho(f', f) + rho(g, g')}
 \rho(f' \circ g', f \circ g) & \le & \rho(f' \circ g', f \circ g') 
                                        + \rho(f \circ g', f \circ g) \\
                               & \le & \rho(f', f) + \rho(g, g') \nonumber
\end{eqnarray}
for every $f, g, f', g' \in \mathcal{I}(X, X)$, which implies the
continuity of $(f, g) \mapsto f \circ g$ as a mapping from
$\mathcal{I}(X, X) \times \mathcal{I}(X, X)$ into $\mathcal{I}(X, X)$
in a simpler way than before.  Similarly, if $f, g \in \mathcal{I}(X)$, then
\begin{eqnarray}
\label{rho(f^{-1}, g^{-1}) = ... = rho(g, f)}
 \rho(f^{-1}, g^{-1}) & = & \rho(f \circ f^{-1}, f \circ g^{-1}) 
                       = \rho({\id}_X, f \circ g^{-1}) \\
 & = & \rho({\id}_X \circ g, (f \circ g^{-1}) \circ g) = \rho(g, f), \nonumber
\end{eqnarray}
which implies that $f \mapsto f^{-1}$ is an isometry on $\mathcal{I}(X)$
with respect to $\rho(f, g)$.

        Note the set $C_s(X, X)$ of continuous mappings from $X$ onto 
itself is a closed set in $C(X, X)$ with respect to the supremum
metric when $X$ is compact.  One way to see this is to observe first
that $C_s(X, X)$ is the same as the set of continuous mappings $f : X
\to X$ such that $f(X)$ is dense in $X$, because $f(X)$ is compact and
hence closed in $X$.  It is easy to check that the set of continuous
mappings $f : X \to X$ such that $f(X)$ is dense in $X$ is closed in
$C(X, X)$, directly from the definitions.  Alternatively, suppose that
$\{f_j\}_{j = 1}^\infty$ is a sequence of continuous mappings from $X$
onto itself that converges uniformly to a mapping $f$ from $X$ into
itself.  Let $y \in X$ be given, and for each positive integer $j$,
let $x_j$ be an element of $X$ such that $f_j(x_j) = y$.  Because $X$
is compact and hence sequentially compact, there is a subsequence
$\{x_{j_l}\}_{l = 1}^\infty$ of $\{x_j\}_{j = 1}^\infty$ that
converges to an element $x$ of $X$.  This implies that
$\{f_{j_l}(x_{j_l})\}_{l = 1}^\infty$ converges to $f(x)$ in $X$,
since $\{f_{j_l}\}_{l = 1}^\infty$ converges to $f$ uniformly on $X$.
It follows that $f(x) = y$, as desired, because $f_{j_l}(x_{j_l}) = y$
for each $l$ by construction.

        It is easy to see that $\mathcal{I}(X, X)$ is also a closed set
in $C(X, X)$ with respect to the supremum metric.  The usual
Arzela--Ascoli arguments imply that $\mathcal{I}(X, X)$ is actually a
compact set in $C(X, X)$ with respect to the supremum metric when $X$
is compact, since isometric mappings on $X$ are obviously
equicontinuous.  It follows that $\mathcal{I}(X) = \mathcal{I}(X, X)
\cap C_s(X, X)$ is a compact set in $C(X, X)$ too, because $C_s(X, X)$
is a closed set in $C(X, X)$, as in the preceding paragraph.
Alternatively, one can show that $\mathcal{I}(X) = \mathcal{I}(X, X)$
when $X$ is compact.

        Now let $(X, d(x, y))$ be a non-compact metric space, and
let $C(X, X)$ be equipped with the topology determined by the
collection of semimetrics $\rho_K(f, g)$ associated to nonempty
compact subsets $K$ of $X$, as in the previous section.  It is easy to
see that $\mathcal{I}(X, X)$ is a closed set in $C(X, X)$ with respect
to this topology, or even with respect to the topology of pointwise
convergence of mappings from $X$ into itself.  However,
$\mathcal{I}(X)$ may be a proper subset of $\mathcal{I}(X, X)$ when
$X$ is not compact, and the set $C_s(X, X)$ of continuous mappings
from $X$ onto $X$ may not be a closed set in $C(X, X)$.  If $X$ is an
infinite set equipped with the discrete metric, for instance, then
$\mathcal{I}(X)$ consists of all of the one-to-one mappings from $X$
onto itself, while $\mathcal{I}(X, X)$ consists of all one-to-one
mappings from $X$ into itself.  In this case, one can check that
$\mathcal{I}(X, X)$ is the closure of $\mathcal{I}(X)$ in $C(X, X)$,
and in particular that $\mathcal{I}(X)$ is not a closed set in $C(X,
X)$.  Similarly, the set of mappings from $X$ onto itself is dense in
$C(X, X)$ in this situation, and in particular is not closed.  If $X$
is a complete metric space and $f$ is an isometric mapping of $X$ into
itself, then observe that $f(X)$ is a closed set in $X$, because
$f(X)$ is also complete with respect to the restriction of $d(x, y)$
to $x, y \in f(X)$.

        If $f \in \mathcal{I}(X, X)$ and $K$ is a nonempty bounded subset 
of $X$, then $f(K)$ is also a bounded set in $X$, and hence the
semimetric $\rho_K(f, g)$ can be defined on $\mathcal{I}(X, X)$ as in
(\ref{rho_K(f, g) = sup_{x in K} d(f(x), g(x))}).  The collection of
semimetrics $\rho_K(f, g)$ associated to nonempty bounded subsets $K$
of $X$ also defines a topology on $\mathcal{I}(X, X)$, as in Section
\ref{semimetrics}.  One can get the same topology on $\mathcal{I}(X, X)$ 
using the semimetrics associated to the sequence of balls in $X$ centered 
at some fixed point and with radius equal to a positive integer, which
implies that this topology is metrizable, as in Section \ref{semimetrics}.
Of course, compact subsets of $X$ are automatically bounded, so that this
topology is at least as strong as the one determined by the semimetrics
$\rho_K(f, g)$ associated to nonempty compact sets $K \subseteq X$.

        As in (\ref{rho_K(f circ h, g circ h) = rho_{h(K)} (f, g)}),
\begin{equation}
\label{rho_K(f circ h, g circ h) = rho_{h(K)} (f, g), 2}
        \rho_K(f \circ h, g \circ h) = \rho_{h(K)} (f, g)
\end{equation}
for every $f, g, h \in \mathcal{I}(X, X)$ and nonempty bounded set
$K \subseteq X$.  In this case, we also have that
\begin{equation}
\label{rho_K(h circ f, h circ g) = rho_K(f, g)}
        \rho_K(h \circ f, h \circ g) = \rho_K(f, g),
\end{equation}
so that the semimetrics $\rho_K(f, g)$ are invariant under left
translations on $\mathcal{I}(X, X)$.  If $f_0, g_0, f, g \in
\mathcal{I}(X, X)$ and $K$ is a nonempty bounded subset of $X$, then
\begin{eqnarray}
\label{rho_K(f circ g, f_0 circ g_0) le rho_K(g, g_0) + rho_{g_0(K)}(f, f_0)}
 \rho_K(f \circ g, f_0 \circ g_0) & \le & \rho_K(f \circ g, f \circ g_0)
                                       + \rho_K(f \circ g_0, f_0 \circ g_0) \\
                      & \le & \rho_K(g, g_0) + \rho_{g_0(K)}(f, f_0). \nonumber
\end{eqnarray}
As before, this can be used to show that $(f, g) \mapsto f \circ g$ is
continuous as a mapping from $\mathcal{I}(X, X) \times \mathcal{I}(X,
X)$ into $\mathcal{I}(X, X)$, with respect to the topology on
$\mathcal{I}(X, X)$ determined by the semimetrics $\rho_K(f, g)$ on
$\mathcal{I}(X, X)$ associated to nonempty bounded sets $K \subseteq X$.

        Similarly, if $f, f_0 \in \mathcal{I}(X)$ and $K$ is a nonempty
bounded subset of $X$, then
\begin{eqnarray}
\label{rho_K(f^{-1}, g^{-1}) = ... = rho_{g^{-1}(K)}(g, f)}
  \quad \rho_K(f^{-1}, g^{-1}) & = & \rho_K(f \circ f^{-1}, f \circ g^{-1})
                                = \rho_K({\id}_K, f \circ g^{-1})  \\
        & = & \rho_{g^{-1}(K)}({\id}_X \circ g, (f \circ g^{-1}) \circ g)
                                  = \rho_{g^{-1}(K)}(g, f). \nonumber
\end{eqnarray}
Using this, it is easy to see that $f \mapsto f^{-1}$ is continuous as
a mapping from $\mathcal{I}(X)$ to itself, with respect to the
topology determined by the semimetrics $\rho_K(\cdot, \cdot)$
associated to nonempty bounded sets $K \subseteq X$.  More precisely,
this implies that $f \mapsto f^{-1}$ is continuous at every point $g
\in \mathcal{I}(X)$, and hence is continuous as a mapping from
$\mathcal{I}(X)$ into itself.

        It is not difficult to check that the collection of 
$f \in \mathcal{I}(X, X)$ such that $f(X)$ is dense in $X$ is a
closed set in $\mathcal{I}(X, X)$ with respect to the topology on
$\mathcal{I}(X, X)$ determined by the semimetrics $\rho_K(\cdot, \cdot)$
associated to nonempty bounded subsets $K$ of $X$.  Equivalently,
if $f \in \mathcal{I}(X, X)$ can be approximated in this topology by
$g \in \mathcal{I}(X, X)$ such that $g(X)$ is dense in $X$, then
$f(X)$ is dense in $X$ too.  Indeed, let $y \in X$ be given, and
suppose that we are looking for $x \in X$ such that $f(x)$ is close to
$y$.  Let $p$ be any element of $X$, so that we should look for $x \in
X$ such that $d(p, x)$ is close to $d(f(p), y)$.  Suppose that $g \in
\mathcal{I}(X, X)$ has the properties that $g(X)$ is dense in $X$, and
that $g$ is uniformly close to $f$ on the ball $B$ in $X$ centered at
$p$ with radius $d(f(p), y) + 1$, say.  If $x \in X$ and $g(x)$ is
close to $y$, then $d(p, x) = d(g(p), g(x))$ is close to $d(g(p), y)$,
which is less than or equal to $d(f(p), y) + d(f(p), g(p))$.
In particular, if $d(f(p), g(p)) < 1$, then there are points $x \in B$
such that $g(x)$ is as close to $y$ as we want.  If $g$ is also
uniformly close to $f$ on $B$, then we get points $x \in B$ such that
$f(x)$ is arbitrarily close to $y$, as desired.  If $X$ is complete,
then $f(X)$ is a closed set in $X$ for every $f \in \mathcal{I}(X, X)$,
as before.  This implies that $\mathcal{I}(X)$ is a closed set in
$\mathcal{I}(X, X)$ with respect to this topology when $X$ is complete.

        Suppose now that all closed and bounded subsets of $X$ are 
compact.  This is equivalent to asking that closed balls in $X$ be
compact, which implies that $X$ is locally compact, $\sigma$-compact,
and complete.  In this case, the topology on $\mathcal{I}(X, X)$
determined by the collection of semimetrics $\rho_K(f, g)$ associated
to nonempty compact subsets $K$ of $X$ is the same as the topology
determined by the collection of semimetrics $\rho_K(f, g)$ associated
to nonempty bounded sets $K \subseteq X$.  We have already seen that
$\mathcal{I}(X, X)$ is a closed set in $C(X, X)$ with respect to the
topology determined by the semimetrics $\rho_K(f, g)$ associated to
nonempty compact sets $K \subseteq X$, and it follows that
$\mathcal{I}(X)$ is also a closed set in $C(X, X)$ with respect to
this topology, by the remarks in the previous paragraph.

        If $p \in X$ and $A$ is a closed set in $X$, then
\begin{equation}
\label{E(p, A) = {f in mathcal{I}(X, X) : f(p) in A}}
        E(p, A) = \{f \in \mathcal{I}(X, X) : f(p) \in A\}
\end{equation}
is a closed set in $C(X, X)$.  If $A$ is compact, then one can check
that $E(p, A)$ is a compact set in $C(X, X)$, by the usual
Arzela--Ascoli arguments.  More precisely, let $\{f_j\}_{j =
  1}^\infty$ be a sequence of elements of $E(p, A)$, and let $K$ be a
nonempty compact set in $X$.  It is easy to see that there is a
compact set $K_1 \subseteq X$ depending on $K$, $p$, and $A$ such that
$f(K) \subseteq K_1$ for every $f \in E(p, A)$, using the hypothesis
that closed and bounded subsets of $X$ are compact.  Because isometric
mappings are obviously equicontinuous, it follows that there is a
subsequence $\{f_{j_l}\}_{l = 1}^\infty$ of $\{f_j\}_{j = 1}^\infty$
that converges uniformly on $K$ to a mapping $f : K \to X$.  Applying
this to the sequence of closed balls in $X$ centered at $p$ with
positive integer radius, and then using a diagonalization argument,
one can get a subsequence of $\{f_j\}_{j = 1}^\infty$ that converges
uniformly on each such ball, and hence on every compact set in $X$.
This implies that $E(p, A)$ is sequentially compact in $C(X, X)$, and
thus that $E(p, A)$ is compact in $C(X, X)$, because $C(X, X)$ is
metrizable under these conditions.  Alternatively, for each nonempty
compact set $K \subseteq X$, one can show that the set of restrictions
of $f \in E(p, A)$ to $K$ is totally bounded in $C(K, X)$ with respect
to the supremum metric.  This implies that $E(p, A)$ is totally
bounded with respect to a metric on $C(X, X)$ that determines the same
topology, and for which $C(X, X)$ is complete as a metric space.

        It follows that $\mathcal{I}(X, X)$ is locally compact
with respect to this topology under these conditions, which implies
that $\mathcal{I}(X)$ is locally compact too, since it is a
closed set in $\mathcal{I}(X, X)$.

\chapter{Commutative groups}
\label{commutative groups}

\section{Homomorphisms into ${\bf R}$}
\label{homomorphisms into R}

        Let $A$ and $B$ be commutative topological groups, in which 
the group operations are expressed additively.  As usual, a mapping
$\phi : A \to B$ is said to be a \emph{continuous
  homomorphism}\index{continuous homomorphism} if $\phi$ is both a
homomorphism from $A$ into $B$ as commutative groups and a continuous
mapping from $A$ into $B$ as topological spaces.  The collection of
all continuous homomorphisms from $A$ into $B$ will be denoted
$\hom(A, B)$,\index{Hom(A, B)@$\hom(A, B)$} which is easily seen to be
a commutative group with respect to pointwise addition.  Note that the
kernel of a continuous homomorphism is a closed subgroup of the domain.

        Let us now restrict our attention to the case where $B$ is
the additive group of real numbers, with the standard topology.
Observe that the only compact subgroup of ${\bf R}$ is the trivial
subgroup $\{0\}$.  If $\phi$ is a continuous homomorphism from a
compact topological group $A$ into ${\bf R}$, then $\phi(A)$ is a
compact subgroup of ${\bf R}$, and hence $\phi(a) = 0$ for every $a
\in A$.

        Suppose that $\phi$ is a continuous homomorphism from ${\bf R}$
into itself.  It is easy to see that
\begin{equation}
\label{phi(r x) = r phi(x)}
        \phi(r \, x) = r \, \phi(x)
\end{equation}
for every $x \in {\bf R}$ and positive integer $r$, because $r \, x$
is the same as the sum of $r$ $x$'s.  Of course, $\phi(0) = 0$ and
$\phi(-y) = - \phi(y)$ for every $y \in {\bf R}$, so that (\ref{phi(r
  x) = r phi(x)}) holds for every $x \in {\bf R}$ and every integer
$r$.  If $z \in {\bf R}$, $n$ is a positive integer, and $n \, z = 0$,
then $z = 0$, which implies that each $y \in {\bf R}$ is uniquely
determined by $n \, y$ for any positive integer $n$.  Using this, one
can check that (\ref{phi(r x) = r phi(x)}) holds for every $x \in {\bf
  R}$ and rational number $r$.  It follows that (\ref{phi(r x) = r
  phi(x)}) holds for all $x, r \in {\bf R}$, because $\phi$ is
continuous, and the set ${\bf Q}$\index{Q@${\bf Q}$} of rational
numbers is dense in ${\bf R}$.  Hence
\begin{equation}
\label{phi(r) = r phi(1)}
        \phi(r) = r \, \phi(1)
\end{equation}
for every $r \in {\bf R}$.

        If $A_1, A_2, \ldots, A_n$ are finitely many commutative topological
groups, then their Cartesian product $A = A_1 \times A_2 \times \cdots 
\times A_n$ is also a commutative topological group, with respect to the
product group structure and topology.  If $B$ is another commutative
topological group and $\phi_j : A_j \to B$ is a continuous homomorphism
for each $j = 1, 2, \ldots, n$, then
\begin{equation}
\label{phi(a) = phi_1(a_1) + phi_2(a_2) + cdots + phi_n(a_n)}
        \phi(a) = \phi_1(a_1) + \phi_2(a_2) + \cdots + \phi_n(a_n)
\end{equation}
is a continuous homomorphism from $A$ into $B$.  Conversely, it is
easy to see that every continuous homomorphism from $A$ into $B$ is of
this form.  In particular, we can apply this to $A_j = {\bf R}$ for
each $j$, so that $A = {\bf R}^n$ as a topological group with respect
to addition and the standard topology.  If we also take $B = {\bf
  R}$,then we get that the continuous homomorphisms from ${\bf R}^n$
into ${\bf R}$ are exactly the usual linear functionals on ${\bf R}^n$.

        Let $I$ be an infinite set, and suppose that $A_j$ is a commutative
topological group for each $j \in I$.  As before, it is easy to see that
the Cartesian product $A = \prod_{j \in I} A_j$ is also a commutative topological
group, with respect to the product group structure and topology.  Let us 
denote elements of $A$ by $x = \{x_j\}_{j \in I}$, so that $x_j \in A_j$
for each $j \in I$.  If $\phi_l$ is a continuous homomorphism from $A_l$
into a commutative topological group $B$ for some $l \in I$, then
\begin{equation}
\label{widetilde{phi}_l(x) = phi_l(x_l)}
        \widetilde{\phi}_l(x) = \phi_l(x_l)
\end{equation}
is a continuous homomorphism from $A$ into $B$, which is the same as
the canonical coordinate mapping from $A$ onto $A_l$ composed with
$\phi_l$.  Of course, the sum of finitely many homomorphisms of this type
is a continuous homomorphism from $A$ into $B$ as well.

        Suppose that there is an open set $W \subseteq B$ such that
$0 \in W$ and $\{0\}$ is the only subgroup of $B$ contained in $W$.
If $B = {\bf R}$ with the standard topology, for instance, then we
can take $W$ to be any open interval in ${\bf R}$ containing $0$.
Let $\phi$ be a continuous homomorphism from $A$ into $B$, so that
$\phi^{-1}(W)$ is an open set in $A$ that contains $0$.  Because of the
way that the product topology is defined on $A$, there are open sets
$U_j \subseteq A_j$ for each $j \in I$ such that $0 \in A_j$ for each $j$,
\begin{equation}
\label{prod_{j in I} U_j subseteq phi^{-1}(W)}
        \prod_{j \in I} U_j \subseteq \phi^{-1}(W),
\end{equation}
and $U_j = A_j$ for all but finitely many $j$.  If we put $A_j' = A_j$
when $U_j = A_j$ and $A_j' = \{0\}$ otherwise, then $A' = \prod_{j \in
  I} A_j'$ is a subgroup of $A$ contained in $\phi^{-1}(W)$, so that
$\phi(A)$ is a subgroup of $B$ contained in $W$.  This implies that
$\phi(x) = 0$ for every $x \in A'$, and hence that $\phi(x)$ depends
only on $x_j$ for the finitely many $j \in I$ such that $U_j \ne A_j$.
It follows that $\phi$ can be expressed as a finite sum of homomorphisms
of the form (\ref{widetilde{phi}_l(x) = phi_l(x_l)}) under these conditions.

        Let $V$ be a topological vector space over the real numbers.
As in Section \ref{semimetrics}, this means that $V$ is a vector space
over the real numbers equipped with a topology for which the vector
space operations are continuous and $\{0\}$ is a closed set.  In particular,
$V$ is a commmutative topological group with respect to addition.
Suppose that $\lambda$ is a continuous homomorphism from $V$ as a
topological group into ${\bf R}$.  This implies that
\begin{equation}
\label{phi_v(t) = lambda(t v)}
        \phi_v(t) = \lambda(t \, v)
\end{equation}
is a continuous homomorphism from ${\bf R}$ into itself for each $v
\in V$.  As in (\ref{phi(r) = r phi(1)}),
\begin{equation}
\label{phi_v(t) = t phi_v(1)}
        \phi_v(t) = t \, \phi_v(1)
\end{equation}
for each $t \in {\bf R}$, which means that
\begin{equation}
        \lambda(t \, v) = t \, \lambda(v)
\end{equation}
for every $v \in V$ and $t \in {\bf R}$.  Thus continuous homomorphisms
from $V$ as a topological group into ${\bf R}$ are the same as continuous
linear functionals on $V$ as a topological vector space.

\section{Homomorphisms into ${\bf T}$}
\label{homomorphisms into T}

        Let $A$ be a commutative topological group, in which the
group operations are expressed additively.  Remember that the
unit circle ${\bf T}$ in the complex plane is a compact commutative
topological group with respect to multiplication of complex numbers
and the topology induced by the standard topology on ${\bf C}$.
A continuous homomorphism from $A$ into ${\bf T}$ is said to be
a \emph{character}\index{characters} on $A$, and the collection
$\hom(A, {\bf T})$ of characters on $A$ is known as the \emph{dual 
group}\index{dual groups} associated to $A$.  This is a commutative
group with respect to pointwise multiplication of characters, as in 
the previous section, which may be denoted 
$\widehat{A}$.\index{A@$\widehat{A}$}

        The set ${\bf C} \backslash \{0\}$ of nonzero complex
numbers is a commutative topological group with respect to multiplication
and the topology induced by the standard topology on ${\bf C}$.
There is a natural isomorphism between ${\bf C} \backslash \{0\}$
and the product of the unit circle ${\bf T}$ and the multiplicative
group ${\bf R}_+$ of positive real numbers, which sends a nonzero
complex number $z$ to
\begin{equation}
\label{(frac{z}{|z|}, |z|) in {bf T} times {bf R}_+}
        \Big(\frac{z}{|z|}, |z|\Big) \in {\bf T} \times {\bf R}_+.
\end{equation}
More precisely, this is an isomorphism between ${\bf C} \backslash
\{0\}$ and ${\bf R} \times {\bf R}_+$ as topological groups, which
means that it is both a group isomorphism and a homeomorphism.  It
follows that a continuous homomorphism from a commutative topological
group $A$ into ${\bf C} \backslash \{0\}$ is the same as product of
continuous homomorphisms from $A$ into ${\bf T}$ and ${\bf R}_+$.  
The exponential mapping defines an isomorphism from the additive group
${\bf R}$ of real numbers onto the multiplicative group ${\bf R}_+$ of
positive real numbers as topological groups, so that a continuous
homomorphism from $A$ into ${\bf R}_+$ is the same as the exponential
of a continuous homomorphism from $A$ into ${\bf R}$.  If $A$ is
compact, then every continuous homomorphism from $A$ into ${\bf R}$
or equivalently ${\bf R}_+$ is trivial, and hence every continuous homomorphism
from $A$ into ${\bf C} \backslash \{0\}$ actually maps $A$ into ${\bf T}$.
This is basically the same as saying that every compact subgroup of
${\bf C} \backslash \{0\}$ is contained in ${\bf T}$.

        Note that $e_t(x) = \exp (i \, t \, x)$ defines a continuous
homomorphism from ${\bf R}$ into ${\bf T}$ for each $t \in {\bf R}$,
which is surjective when $t \ne 0$.  If $\phi$ is a continuous homomorphism
from a commutative topological group $A$ into ${\bf R}$, then it follows
that $\phi_t(a) = \exp (i \, t \, \phi(a))$ is a continuous
homomorphism from $A$ into ${\bf T}$ for each $t \in {\bf R}$.  These
homomorphisms are all trivial when $A$ is compact.

        If $\psi$ is any continuous mapping from ${\bf R}$ into ${\bf T}$
such that $\psi(0) = 1$, then it is well known that there is a unique
continuous mapping $\phi$ from ${\bf R}$ into itself such that
$\phi(0) = 1$ and
\begin{equation}
\label{psi(x) = exp (i phi(x))}
        \psi(x) = \exp (i \, \phi(x))
\end{equation}
for every $x \in {\bf R}$.  This uses the fact that $\exp (i \, r)$
is a local homeomorphism from ${\bf R}$ onto ${\bf T}$.  Suppose now that
$\psi$ is a continuous homomorphism from ${\bf R}$ into ${\bf T}$,
and let us show that $\phi$ is a homomorphism from ${\bf R}$ into
itself.  Put
\begin{equation}
\label{phi_t(x) = phi(x + y) - phi(y)}
        \phi_y(x) = \phi(x + y) - \phi(y)
\end{equation}
for every $x, y \in {\bf R}$, so that $\phi_y$ is a continuous mapping
from ${\bf R}$ into itself such that $\phi_y(0) = 0$ and
\begin{eqnarray}
\label{exp (i phi_y(x)) = ... = psi(x + y) psi(y)^{-1} = psi(x)}
        \exp (i \, \phi_y(x)) & = & \exp (i \, \phi(x + y) - i \, \phi(y)) \\
   & = & \exp (i \, \phi(x + y)) \, (\exp (i \, \phi(y)))^{-1} \nonumber \\
   & = & \psi(x + y) \, \psi(y)^{-1} = \psi(x) \nonumber
\end{eqnarray}
for every $x, y \in {\bf R}$.  Thus $\phi_y(x)$ satisfies the same
conditions as $\phi(x)$, which implies that $\phi_y(x) = \phi(x)$ for
every $x, y \in {\bf R}$, by uniqueness.  Similarly, $-\phi(-x)$ is a
continuous mapping from ${\bf R}$ into itself which is equal to $0$
at $0$ and satisfies
\begin{equation}
\label{exp (- i phi(-x)) = (exp (i phi(-x))^{-1} = psi(-x)^{-1} = psi(x)}
 \exp (- i \, \phi(-x)) = (\exp (i \, \phi(-x))^{-1} = \psi(-x)^{-1} = \psi(x)
\end{equation}
for every $x \in {\bf R}$, so that $- \phi(-x) = \phi(x)$ for every
$x \in {\bf R}$ by uniqueness of $\phi$ again.  It follows that
$\phi$ is a continuous homomorphism from ${\bf R}$ into itself as
a topological group, and hence there is a $t \in {\bf R}$ such that
\begin{equation}
\label{phi(x) = t x}
        \phi(x) = t \, x
\end{equation}
for every $x \in {\bf R}$, as in the previous section.  This shows that
every continuous homomorphism from ${\bf R}$ into ${\bf T}$ is of the
form $e_t(x) = \exp (i \, t \, x)$ for some $t \in {\bf R}$.

        Now let $\phi$ be a continuous homomorphism from ${\bf T}$
into itself.  Observe that
\begin{equation}
\label{psi(x) = phi(exp (i x))}
        \psi(x) = \phi(\exp (i \, x))
\end{equation}
is a continuous homomorphism from ${\bf R}$ into ${\bf T}$, because
$\exp (i \, x)$ is a continuous homomorphism from ${\bf R}$ into ${\bf T}$.
As in the preceding paragraph, there is a $t \in {\bf R}$ such that
\begin{equation}
\label{psi(x) = exp (i t x)}
        \psi(x) = \exp (i \, t \, x)
\end{equation}
for every $x \in {\bf R}$.  If we take $x = 2 \pi$, then we get that
$\psi(x) = 1$, and hence that $t$ is an integer.  Of course, $\phi(z)
= z^n$ is a continuous homomorphism from ${\bf T}$ into itself for
each integer $n$, and this argument shows that every continuous
homomorphism from ${\bf T}$ into itself is of this form.

        Let $I$ be a nonempty set, and suppose that $A_j$ is a
commutative topological group for each $j \in I$, so that $A =
\prod_{j \in I} A_j$ is a commutative topological group with respect
to the product group structure and topology.  If $l_1, \ldots, l_n$
are finitely many elements of $I$ and $\phi_{l_k}$ is a continuous
homomorphism from $A_{l_k}$ into ${\bf T}$ for each $k = 1, \ldots,
n$, then
\begin{equation}
        \phi(x) = \phi_{l_1}(x_{l_1}) \cdots \phi_{l_n}(x_{l_n})
\end{equation}
is a continuous homomorphism from $A$ into ${\bf T}$.  Conversely,
every continuous homomorphism from $A$ into ${\bf T}$ is of this form,
as in the previous section.  More precisely, note that
\begin{equation}
        W = \{z \in {\bf T} : \re z > 0\}
\end{equation}
is an open set in ${\bf T}$, $1 \in W$, and $\{1\}$ is the only
subgroup of ${\bf T}$ contained in $W$.  This implies that every
continuous homomorphism from $A$ into ${\bf T}$ depends on only
finitely many coordinates, as before.

        Let $V$ be a topological vector space over the real numbers,
and let $\phi$ be a continuous homomorphism from $V$ as a
topological group with respect to addition into ${\bf T}$.
This implies that
\begin{equation}
\label{phi_v(r) = phi(r v)}
        \phi_v(r) = \phi(r \, v)
\end{equation}
is a continuous homomorphism as a function of $r$ from ${\bf R}$
into ${\bf T}$ for each $v \in V$.  Thus there is a real number
$\lambda(v)$ such that
\begin{equation}
\label{phi_v(r) = exp (i r lambda(v))}
        \phi_v(r) = \exp (i \, r \, \lambda(v))
\end{equation}
for each $r \in {\bf R}$, and which is uniquely determined by $v$.  It
is easy to see from the uniqueness of $\lambda(v)$ that $\lambda(t \,
v) = t \, \lambda(v)$ for every $v \in V$ and $t \in {\bf R}$, and that
\begin{equation}
\label{lambda(v + w) = lambda(v) + lambda(w)}
        \lambda(v + w) = \lambda(v) + \lambda(w)
\end{equation}
for every $v, w \in V$, because $\phi$ is a homomorphism with respect to
addition.  This shows that $\lambda$ is a linear functional on $V$, and that
\begin{equation}
\label{phi(v) = exp (i lambda(v))}
        \phi(v) = \exp (i \, \lambda(v))
\end{equation}
for every $v \in V$, by taking $r = 1$ in (\ref{phi_v(r) = exp (i r
  lambda(v))}).

        We would also like to show that $\lambda$ is continuous under
these circumstances.  Of course, it suffices to check that $\lambda$
is continuous at $0$, because $\lambda$ is linear.  It is easy to see
that $\lambda(v)$ is close to $0$ modulo $2 \pi \, {\bf Z}$ when $v$
is close to $0$ in $V$, because $\phi(v)$ is continuous and $\exp (i
\, x)$ is a local homeomorphism from ${\bf R}$ onto ${\bf T}$.  In
order to get that $\lambda(v)$ is close to $0$ when $v$ is close to
$0$, it suffices to use neighborhoods of $0$ in $V$ that are star-like
about $0$.  Conversely, if $\lambda$ is a continuous linear functional
on $V$, then (\ref{phi(v) = exp (i lambda(v))}) defines a continuous
homomorphism from $V$ as a topological group with respect to addition
into ${\bf T}$.

\section{Finite abelian groups}
\label{finite abelian groups}

        Let $n$ be a positive integer, and consider the group
${\bf Z} / n \, {\bf Z}$, which is the quotient of the group ${\bf Z}$
of integers with respect to addition by the subgroup $n \, {\bf Z}$
consisting of integer multiples of $n$.  Thus ${\bf Z} / n \, {\bf Z}$
is the cyclic group of order $n$, which may be considered as a compact
commutative topological group with respect to the discrete topology.
A homomorphism from ${\bf Z} / n \, {\bf Z}$ into ${\bf T}$ is essentially
the same as a homomorphism $\phi$ from ${\bf Z}$ into ${\bf T}$ such that
$\phi(n) = 1$, so that the kernel of $\phi$ contains $n \, {\bf Z}$.
Any homomorphism $\phi$ from ${\bf Z}$ into ${\bf T}$ satisfies
$\phi(j) = \phi(1)^j$ for each $j \in {\bf Z}$, and thus the condition
$\phi(n) = 1$ reduces to $\phi(1)^n = 1$.  It follows that the
dual group associated to ${\bf Z} / n \, {\bf Z}$ is isomorphic to itself.

        Now let $A$ be any finite abelian group, which may be considered 
as a compact commutative topological group with respect to the
discrete topology again.  It is well known that $A$ is isomorphic to a
product of cyclic groups, but let us look at $A$ in terms of
translation operators on functions on $A$.  More precisely, let $V$ be
the vector space of complex-valued functions on $A$, which is a
finite-dimensional vector space with respect to pointwise addition and
scalar multiplication, with dimension equal to the number of elements
of $A$.  Put
\begin{equation}
\label{langle f, g rangle = sum_{x in A} f(x) overline{g(x)}}
        \langle f, g \rangle = \sum_{x \in A} f(x) \, \overline{g(x)}
\end{equation}
for each $f, g \in V$, which defines an inner product on $V$, with the
corresponding norm on $V$ given by
\begin{equation}
\label{||f|| = (sum_{x in A} |f(x)|^2)^{1/2}}
        \|f\| = \Big(\sum_{x \in A} |f(x)|^2\Big)^{1/2}.
\end{equation}
If $a \in A$ and $f \in V$, then let $T_a(f)$ be the function on $A$
defined by
\begin{equation}
\label{(T_a(f))(x) = f(x + a)}
        (T_a(f))(x) = f(x + a)
\end{equation}
for each $x \in V$.  Thus $T_a$ is a one-to-one linear mapping
from $V$ onto itself for each $a \in A$.  It is easy to see that
\begin{equation}
\label{langle T_a(f), T_a(g) rangle = langle f, g rangle}
        \langle T_a(f), T_a(g) \rangle = \langle f, g \rangle
\end{equation}
for every $a \in A$ and $f, g \in V$, so that $T_a$ is a unitary
mapping on $V$ with respect to this inner product.

        If $\phi$ is a homomorphism from $A$ into ${\bf T}$, then
\begin{equation}
\label{(T_a(phi))(x) = phi(x + a) = phi(x) phi(a)}
        (T_a(\phi))(x) = \phi(x + a) = \phi(x) \, \phi(a)
\end{equation}
for every $a, x \in A$, so that $\phi$ is an eigenvector for $T_a$ for
every $a \in A$.  Conversely, suppose that $\psi \in V$ is a nonzero
eigenvector of $T_a$ for every $a \in A$.  This means that $\|\psi\| >
0$, and for each $a \in A$ there is a complex number $\lambda(a)$ such
that
\begin{equation}
\label{T_a(psi) = lambda(a) psi}
        T_a(\psi) = \lambda(a) \, \psi.
\end{equation}
Because $T_a$ is unitary for each $a \in A$, we get that
\begin{equation}
\label{||psi|| = ||T_a(psi)|| = ||lambda(a) psi|| = |lambda(a)| ||psi||}
 \|\psi\| = \|T_a(\psi)\| = \|\lambda(a) \, \psi\| = |\lambda(a)| \, \|\psi\|
\end{equation}
for every $a \in A$, and hence that $|\lambda(a)| = 1$ for each $a \in
A$.  This implies that $|\psi(x)|$ is constant on $A$, and in particular
that $\psi(0) \ne 0$.  Thus
\begin{equation}
\label{phi(x) = psi(0)^{-1} psi(x)}
        \phi(x) = \psi(0)^{-1} \, \psi(x)
\end{equation}
is also an eigenvector of $T_a$ with eigenvalue $\lambda(a)$ for each
$a \in A$, with the additional normalization $\phi(0) = 1$.
The eigenvalue condition implies that
\begin{equation}
\label{phi(x + a) = (T_a(phi))(x) = lambda(a) phi(x)}
        \phi(x + a) = (T_a(\phi))(x) = \lambda(a) \, \phi(x)
\end{equation}
for every $a, x \in A$, and hence that $\lambda(a) = \phi(a)$ for
every $a \in A$, by taking $x = 0$.  It follows that $\phi$ is a
homomorphism from $A$ into ${\bf T}$ under these conditions.

        Suppose that $\phi$ and $\phi'$ are distinct homomorphisms
from $A$ into ${\bf T}$, so that $\phi(a) \ne \phi'(a)$ for some $a \in A$.
Observe that
\begin{equation}
\label{langle phi, phi' rangle = langle T_a(phi), T_a(phi') rangle = ...}
 \langle \phi, \phi' \rangle = \langle T_a(\phi), T_a(\phi') \rangle
              = \phi(a) \, \overline{\phi'(a)} \, \langle \phi, \phi' \rangle,
\end{equation}
because $T_a$ is unitary.  If $\langle \phi, \phi' \rangle \ne 0$,
then it follows that $\phi(a) \, \overline{\phi'(a)} = 1$, which would
contradict the hypothesis that $\phi(a) \ne \phi'(a)$, since
$|\phi'(a)| = 1$.  Thus $\phi$ is orthogonal to $\phi'$, which implies
that the number of distinct homomorphisms from $A$ into ${\bf T}$ is
less than or equal to the dimension of $V$, which is the number of
elements of $A$.  In fact, the number of distinct homomorphisms from
$A$ into ${\bf T}$ is equal to the number of elements of $A$, and they
form an orthogonal basis for $V$.  To see this, remember that a
unitary transformation on a finite-dimensional complex inner product
space can be diagonalized in an orthogonal basis.  This implies that
for each $a \in A$, $V$ can be expressed as an orthogonal direct sum
of eigenspaces for $T_a$.  We also have that $T_a \circ T_b = T_b
\circ T_a$ for every $a, b \in A$, because $A$ is abelian, which
implies that $T_b$ maps the eigenspaces of $T_a$ into themselves for
each $a, b \in A$.  This permits one to get an orthogonal basis for
$V$ consisting of eigenvectors for $T_a$ for every $a \in A$,
which can be normalized to get homomorphisms from $A$ into ${\bf T}$
as in the preceding paragraph.

\section{Discrete groups}
\label{discrete groups}

        If $\phi$ is a homomorphism from the group ${\bf Z}$ of
integers with respect to addition into any group $G$, then
\begin{equation}
\label{phi(j) = phi(1)^j}
        \phi(j) = \phi(1)^j
\end{equation}
for each $j \in {\bf Z}$.  Thus $\phi$ is uniquely determined by
$\phi(1)$, and for each $g \in G$ there is a homomorphism $\phi$ from
${\bf Z}$ into $G$ such that $g = \phi(1)$.  If $G$ is an abelian
group, then the collection of homomorphisms from ${\bf Z}$ into $G$ is
a group, and it follows that this group is isomorphic to $G$.
In particular, the group of homomorphisms from ${\bf Z}$ into ${\bf T}$
is isomorphic to ${\bf T}$.

        Let $A$ be a nonempty set, and let ${\bf T}^A$ be the collection 
of all functions on $A$ with values in ${\bf T}$, which is a
commutative group with respect to pointwise multiplication.  This is
the same as the Cartesian product of copies of ${\bf T}$ indexed by
$A$, which is a compact Hausdorff space with respect to the product
topology associated to the standard topology on ${\bf T}$.  As usual,
${\bf T}^A$ is a topological group.

        Now let $A$ be any abelian group, equipped with the discrete 
topology.  It is easy to see that the group $\hom(A, {\bf T})$ of 
homomorphisms from $A$ into  ${\bf T}$ is a closed subgroup of ${\bf T}^A$.  
Thus $\hom(A, {\bf T})$ is compact with respect to the topology induced by 
the product topology on ${\bf T}^A$.

        Let $E$ be a nonempty subset of $A$, and let ${\bf T}^E$ be the 
group of functions on $E$ with values in ${\bf T}$, as before.  There
is a natural mapping $\pi_E$ from ${\bf T}^A$ onto ${\bf T}^E$, which
sends a ${\bf T}$-valued function on $A$ to its resriction to $E$, and
which is clearly a continuous homomorphism.  Thus $\pi_E$ maps 
$\hom(A, {\bf T})$ onto a compact subgroup of ${\bf T}^E$.  If $E$ is a 
set of generators of $A$, so that every element of $A$ can be expressed 
as a finite sum of elements of $E$ and their inverses, then the restriction 
of $\pi_E$ to $\hom(A, {\bf T})$ is a homeomorphism onto its image in 
${\bf T}^E$.

        In particular, if $A = {\bf Z}$, then we can take $E = \{1\}$,
and identify ${\bf T}^E$ with ${\bf T}$.  In this case, $\pi_E$ maps
$\hom({\bf Z}, {\bf T})$ onto ${\bf T}$, and the restriction of
$\pi_E$ to $\hom({\bf Z}, {\bf T})$ corresponds exactly to the isomorphism
between $\hom({\bf Z}, {\bf T})$ and ${\bf T}$ mentioned earlier.
The topology on $\hom({\bf Z}, {\bf T})$ induced by the product topology
on ${\bf T}^{\bf Z}$ corresponds to the standard topology on ${\bf T}$.
Similarly, if $A = {\bf Z}^n$ for some positive integer $n$, then we can
take $E$ to be the set of $n$ standard generators of ${\bf Z}^n$.
This leads to an isomorphism between $\hom({\bf Z}^n, {\bf T})$ and
${\bf T}^n$, for which the standard topology on ${\bf T}^n$ corresponds
to the topology on $\hom({\bf Z}^n, {\bf T})$ induced by the product
topology on ${\bf T}^{({\bf Z}^n)}$.

        Let $A$ be any abelian group with the discrete topology again.
Suppose that $B$ is a subgroup of $A$, $\phi$ is a homomorphism from
$B$ into ${\bf T}$, and that $x \in A \backslash B$.  We would like to
extend $\phi$ to a homomorphism from the subgroup $B(x)$ of $A$
generated by $B$ and $x$ into ${\bf T}$.  As usual, $j \cdot x$
can be defined for every $j \in {\bf Z}$, as the sum of $j$ $x$'s when 
$j \in {\bf Z}_+$, and so on.  If $j \cdot x \not\in B$ for any 
$j \in {\bf Z}_+$, then every element of $B(x)$ has a unique representation
as $b + j \cdot x$ for some $b \in B$ and $j \in {\bf Z}$, so that
$B(x)$ is isomorphic to $B \times {\bf Z}$.  In this case, one can extend
$\phi$ to $B(x)$ by putting $\phi(b + j \cdot x) = \phi(b)$ for every
$b \in B$ and $j \in {\bf Z}$.  Otherwise, if $j \cdot x \in B$
for some $j \in {\bf Z}_+$, then we let $n$ be the smallest positive
integer with this property.  Thus $n \cdot x \in B$, so that $\phi(n \cdot x)$
is already defined, and we choose $\zeta \in {\bf T}$ such that $\zeta^n =
\phi(n \cdot x)$.  In this case, we put
\begin{equation}
\label{phi(b + j cdot x) = phi(b) zeta^j}
        \phi(b + j \cdot x) = \phi(b) \, \zeta^j
\end{equation}
for every $b \in B$ and $j \in {\bf Z}$, which one can check is a
well-defined homomorphism from $B(x)$ into ${\bf T}$.  This uses the
fact that $j \cdot x \in B$ for some $j \in {\bf Z}$ if and only if
$j$ is an integer multiple of $n$.

        If $A$ is finitely generated, or generated by $B$ and finitely
many additional elements, then one can repeat the process to get an
extension of $\phi$ to a homomorphism from $A$ into ${\bf T}$.
Similarly, if $A$ can be generated by $B$ and countably many additional
elements, then an increasing sequence of these extensions leads to an
extension of $\phi$ to a homomorphism from $A$ into ${\bf T}$.
Otherwise, one can use Zorn's lemma or the Hausdorff maximality principle
to argue that there is a maximal extension of $\phi$ to a homomorphism
from subgroup of $A$ into ${\bf T}$, and then use the preceding construction
to show that a maximal extension is defined on all of $A$.

        In particular, if $a$ is any nonzero element of $A$, then one
can take $B$ to be the subgroup of $A$ generated by $a$, and it is easy
to see that there is a homomorphism $\phi$ from $B$ into ${\bf T}$ such
that $\phi(a) \ne 1$.  The previous discussion then leads to an extension
of $\phi$ to a homomorphism from $A$ into ${\bf T}$.

\section{Compact groups}
\label{compact groups}

        Let $A$ be a compact abelian topological group.  Remember
that $\{1\}$ is the only subgroup of ${\bf T}$ contained in the set
$W$ of $z \in {\bf T}$ with $\re z > 0$, for instance.  If $\phi$ is a
continuous homomorphism from $A$ into ${\bf T}$ such that $\re \phi(a)
> 0$ for every $a \in A$, then $\phi(a)$ is a subgroup of ${\bf T}$
contained in $W$, and hence $\phi(a) = 1$ for every $a \in A$.  In
particular, this holds when
\begin{equation}
\label{|phi(a) - 1| < 1}
        |\phi(a) - 1| < 1
\end{equation}
for every $a \in A$.  Suppose that $\phi_1$, $\phi_2$ are continuous
homomorphisms from $A$ into ${\bf T}$ such that
\begin{equation}
\label{|phi_1(a) - phi_2(a)| < 1}
        |\phi_1(a) - \phi_2(a)| < 1
\end{equation}
for every $a \in A$.  Applying the previous argument to $\phi(a) =
\phi_1(a) \, \phi_2(a)^{-1}$, we get that $\phi_1(a) = \phi_2(a)$ for
every $a \in A$.  Because of this, it is natural to consider the group
$\hom(A, {\bf T})$ of continuous homomorphisms from $A$ into ${\bf T}$
as being equipped with the discrete topology in this case.

        Let $H$ be Haar measure on $A$, which we normalize so that
$H(A) = 1$, and let $\phi$ be a continuous homomorphism from $A$ into
${\bf T}$.  Observe that
\begin{equation}
\label{int_A phi(x) dH(x) = int_A phi(x + a) dH(x) = phi(a) int_A phi(x) dH(x)}
        \int_A \phi(x) \, dH(x) = \int_A \phi(x + a) \, dH(x) 
                                = \phi(a) \, \int_A \phi(x) \, dH(x)
\end{equation}
for every $a \in A$.  If $\phi(a) \ne 1$ for some $a \in A$, then it
follows that
\begin{equation}
\label{int_A phi(x) dH(x) = 0}
        \int_A \phi(x) \, dH(x) = 0.
\end{equation}
Suppose that $\phi_1$, $\phi_2$ are distinct continuous homomorphisms
from $A$ into ${\bf T}$, so that $\phi_1(a) \ne \phi_2(a)$ for some $a
\in A$.  Applying the previous argument to $\phi(x) = \phi_1(x) \,
\phi_2(x)^{-1} = \phi_1(x) \, \overline{\phi_2(x)}$, we get that
\begin{equation}
\label{int_A phi_1(x) overline{phi_2(x)} dH(x) = 0}
        \int_A \phi_1(x) \, \overline{\phi_2(x)} \, dH(x) = 0.
\end{equation}
This shows that distinct characters on $A$ are orthogonal with respect
to the standard $L^2$ integral inner product
\begin{equation}
\label{langle f, g rangle = int_A f(x) overline{g(x)} dH(x)}
        \langle f, g \rangle = \int_A f(x) \, \overline{g(x)} \, dH(x).
\end{equation}
The normalization $H(A) = 1$ implies that characters on $A$ have
$L^2$ norm equal to $1$, so that they are actually orthonormal in $L^2(A)$.

        Let $C(A)$ be the space of continuous complex-valued functions
on $A$, and let $\mathcal{E}$ be the linear subspace of $C(A)$ consisting
of finite linear combinations of characters on $A$.  Note that this contains 
constant functions on $A$, since the constant function equal to $1$ on $A$
is a character.  If $f \in \mathcal{E}$, then $\overline{f} \in \mathcal{E}$,
because the complex conjugate of a character on $A$ is also a character.
Similarly, if $f, g \in \mathcal{E}$, then $f \, g \in \mathcal{E}$,
because the product of two characters on $A$ is a character as well.
If $\mathcal{E}$ separates points in $A$, then the Stone--Weierstrass
theorem implies that $\mathcal{E}$ is dense in $C(A)$ with respect to
the supremum norm.  Of course, $\mathcal{E}$ separates points in $A$
if and only if the set of characters on $A$ separates points.  In this
case, it follows that the characters on $A$ form an orthonormal basis
for $L^2(A)$, since $C(A)$ is dense in $L^2(A)$.

        If $A$ is a finite abelian group, then $C(A)$ is spanned by
the characters on $A$, as in Section \ref{finite abelian groups}.
If $A = {\bf T}$ with the standard topology, then the characters
on $A$ are given by $z \mapsto z^n$ for $n \in {\bf Z}$, which
obviously separate points on ${\bf T}$.  Let $I$ be a nonempty set,
let $A_i$ be a compact abelian topological group for each $i \in I$,
and suppose that characters on $A_i$ separate points for each $i \in I$.
Under these conditions, $A = \prod_{i \in I} A_i$ is a compact abelian
group with respect to the product topology and group structure, and it
is easy to see that characters on $A$ separate points too.

        Let $A$ be any compact abelian topological group again, and
let $E_0$ be a subgroup of the group of characters on $A$.  Also let
$\mathcal{E}_0$ be the linear subspace of $C(A)$ spanned by $E_0$,
which is a subalgebra of $C(A)$ that contains the constant functions
and is invariant under complex conjugation, as before.  If $E_0$
separates points in $A$, then $\mathcal{E}_0$ separates points in $A$,
and hence $\mathcal{E}_0$ is dense in $C(A)$ with respect to the
supremum norm, by the Stone--Weierstrass theorem.  If $\phi$ is a
character on $A$ not in $E_0$, then $\phi$ is orthogonal to the
elements of $\mathcal{E}_0$ in $L^2(A)$, as in (\ref{int_A phi_1(x)
  overline{phi_2(x)} dH(x) = 0}).  This shows that every character
on $A$ is an element of $E_0$ when $E_0$ separates points in $A$.

        Suppose that $A$ is a compact subgroup of a commutative topological
group $B$, with the induced topology, and that characters on $B$ separate
points in $B$.  If $E_0$ is the subgroup of $\hom(A, {\bf T})$ consisting of
restrictions of characters on $B$ to $A$, then $E_0$ separates points in $A$,
and hence every character on $A$ is in $E_0$, by the previous argument.

        Suppose now that $A$ is an abelian group with the discrete
topology, and let $\widehat{A}$ be the group of characters on $A$.
As in the preceding section, $\widehat{A}$ may be considered as a
closed subgroup of ${\bf T}^A$ with respect to the product topology,
and thus as a compact abelian topological group.  If $a \in A$
and $\phi \in \widehat{A}$, then put
\begin{equation}
\label{Psi_a(phi) = phi(a)}
        \Psi_a(\phi) = \phi(a),
\end{equation}
so that $\Psi_a$ defines a continuous homomorphism from $\widehat{A}$
into ${\bf T}$ for each $a \in A$.  If $\phi_1$ and $\phi_2$ are
distinct elements of $\widehat{A}$, then $\phi_1(a) \ne \phi_2(a)$ for
some $a \in A$, which means that the collection of $\Psi_a$'s with $a
\in A$ separates points on $\widehat{A}$.  It is easy to see that the
collection of $\Psi_a$'s with $a \in A$ forms a subgroup of the group
of characters on $\widehat{A}$, and in fact $a \mapsto \Psi_a$ defines
a homomorphism from $A$ into the group of characters on $\widehat{A}$.
The discussion in the previous paragraphs implies that every continuous
homomorphism from $\widehat{A}$ into ${\bf T}$ is of the form $\Psi_a$
for some $a \in A$.  If $a \in A$ and $a \ne 0$, then we saw in the
preceding section that $\phi(a) \ne 0$ for some $\phi \in \widehat{A}$.
This shows that $a \mapsto \Psi_a$ is an isomorphism from $A$ onto the
group of characters on $\widehat{A}$ in this case.

\section{The dual topology}
\label{dual topology}

        Let $X$ be a topological space, and let $C(X)$ be the space
of continuous complex-valued functions on $X$.  If $f \in C(X)$ and $K$ 
is a nonempty compact subset of $X$, then put
\begin{equation}
\label{||f||_K = sup_{x in K} |f(x)|}
        \|f\|_K = \sup_{x \in K} |f(x)|.
\end{equation}
It is easy to see that this defines a seminorm on $C(X)$, known as the
\emph{supremum seminorm}\index{supremum seminorms} associated to $K$.
The collection of these seminorms determines a topology on $C(X)$
which makes $C(X)$ into a topological vector space, as in Section
\ref{semimetrics}.  Of course, if $X$ is compact, then this is the
same as the topology on $C(X)$ determined by the supremum norm.

        In addition to being a vector space with respect to pointwise
addition and scalar multiplication, $C(X)$ is a commutative algebra
with respect to pointwise multiplication.  Observe that
\begin{equation}
\label{||f g||_K le ||f||_K ||g||_K}
        \|f \, g\|_K \le \|f\|_K \, \|g\|_K
\end{equation}
for every $f, g \in C(X)$ and nonempty compact set $K \subseteq X$.
Using this, one can check that pointwise multiplication of functions
defines a continuous mapping from $C(X) \times C(X)$ into $C(X)$,
so that $C(X)$ is actually a topological algebra.

        Let $C(X, {\bf T})$ be the space of continuous mappings from $X$
into ${\bf T}$.  This is a commutative group with respect to pointwise
multiplication, and a topological group with respect to the topology
induced by the one on $C(X)$ just defined.  Note that $C(X, {\bf T})$
is a closed set in $C(X)$ with respect to this topology.  If $f, g \in
C(X)$, $h \in C(X, {\bf T})$, and $K \subseteq X$ is a nonempty
compact set, then
\begin{equation}
\label{||f h - g h||_K = ||f - g||_K}
        \|f \, h - g \, h\|_K = \|f - g\|_K.
\end{equation}
In particular, the corresponding semimetric
\begin{equation}
\label{d_K(f, g) = ||f - g||_K}
        d_K(f, g) = \|f - g\|_K
\end{equation}
is invariant under translations on $C(X, {\bf T})$ as a group with
respect to pointwise multiplication.

        Now let $A$ be a commutative topological group, and let 
$\widehat{A} = \hom(A, {\bf T})$ be the dual group of continuous 
homomorphisms from $A$ into ${\bf T}$.  Thus $\widehat{A} \subseteq
C(A)$, and we consider the topology on $\widehat{A}$ induced by the
topology on $C(A)$ determined by the supremum seminorms associated to
nonempty compact subsets of $A$, as in the previous paragraphs.  More
precisely, $\widehat{A}$ is a subgroup of $C(A, {\bf T})$, and hence a
topological group with respect to the induced topology.  It is easy to
see that $\widehat{A}$ is a closed set in $C(A)$ with respect to this
topology.

        If $A$ is discrete, then the compact subsets of $A$ are the same 
as the finite subsets of $A$, and this topology on $\widehat{A}$ is
the same as the one induced by the product topology on ${\bf T}^A$ as
in Section \ref{discrete groups}.  If $A$ is compact, then this
topology on $C(A)$ is the same as the one determined by the supremum
norm, and the induced topology on $\widehat{A}$ is discrete, as in the
preceding section.  If $A = {\bf R}$ as a commutative group with
respect to addition and equipped with the standard topology, then we
have seen in Section \ref{homomorphisms into T} that every continuous
homomorphism from ${\bf R}$ into ${\bf T}$ is of the form $e_t(x) =
\exp (i \, t \, x)$ for some $t \in {\bf R}$.  Thus $t \mapsto e_t$
defines a group isomorphism from ${\bf R}$ onto $\widehat{\bf R}$, and
one can check that this is also a homeomorphism with respect to the
standard topology on ${\bf R}$ and the topology induced on
$\widehat{\bf R}$ by the one on $C({\bf R})$ as before.

        Suppose that $A$ and $B$ are commutative topological groups,
and consider $A \times B$ as a topological group with respect to the
product topology and group structure.  This is also a commutative
topological group, and we have seen that the corresponding dual group 
is isomorphic in a natural way to $\widehat{A} \times \widehat{B}$.
It is easy to see that the topology on $\widehat{A \times B}$ induced
by the usual topology on $C(A \times B)$ is the same as the product
topology on $\widehat{A} \times \widehat{B}$ associated to the
topologies induced on $\widehat{A}$ and $\widehat{B}$ by those on
$C(A)$ and $C(B)$, respectively.  This uses the fact that every
compact subset of $A \times B$ is contained in $H \times K$ for some
compact sets $H \subseteq A$ and $K \subseteq B$.  In particular, the
dual of ${\bf R}^n$ is isomorphic to itself as a topological group
with the standard topology for each positive integer $n$.

          Let $X$ be a topological space again, and suppose that $X$
is $\sigma$-compact.  This means that there is a sequence $K_1, K_2,
K_3, \ldots$ of compact subsets of $X$ whose union is equal to $X$.
We may as well ask that $K_l \ne \emptyset$ and $K_l \subseteq K_{l +
  1}$ for each $l$, by replacing $K_l$ with the union of $K_1, \ldots,
K_l$ if necessary.  If $X$ is also locally compact, then every compact
set in $X$ is contained in the interior of another compact set.  Using
this, one can modify the $K_l$'s again to get that $K_l$ is contained
in the interior of $K_{l + 1}$ for each $l$.  In particular, this
implies that $X$ is equal to the union of the interiors of the
$K_l$'s.  If $K \subseteq X$ is any compact set, then it follows that
$K$ is contained in the union of the interiors of finitely many
$K_l$'s, and hence that $K \subseteq K_l$ for some $l$.  This implies
that the topology on $C(X)$ determined by the supremum seminorms
associated to nonempty compact subsets of $X$ is the same as the
topology determined by the supremum seminorms associated to the
$K_l$'s under these conditions.  It follows that this topology on
$C(X)$ is metrizable, as in Section \ref{semimetrics}.

\section{Equicontinuity}
\label{equicontinuity}

        Let $A$ be a locally compact abelian topological group, and let
$\widehat{A} = \hom(A, {\bf T})$ be the corresponding dual group, with the
topology induced by the usual one on $C(A)$, as in the previous
section.  A set $E \subseteq \widehat{A}$ is said to be
\emph{equicontinuous}\index{equicontinuity} if for each $\epsilon > 0$
there is an open set $U \subseteq A$ such that $0 \in U$ and
\begin{equation}
\label{|phi(x) - 1| < epsilon}
        |\phi(x) - 1| < \epsilon
\end{equation}
for every $\phi \in E$ and $x \in U$.  Although this type of condition
might normally be described as ``equicontinuity of $E$ at $0$'', it
implies equicontinuity of $E$ at every point in $A$, because $E
\subseteq \widehat{A}$.  Similarly, the remarks in this section may be
considered as special cases of standard arguments about collections of
continuous functions on locally compact Hausdorff topological spaces,
with simplifications resulting from the group structure and the
restriction to characters on $A$.

        Suppose that $E \subseteq \widehat{A}$ is compact, and let us
check that $E$ is equicontinuous.  Let $\epsilon > 0$ be given, and for
each $\phi \in E$, let $U(\phi)$ be an open set in $A$ such that 
$0 \in U(\phi)$ and
\begin{equation}
\label{sup {|phi(x) - 1| : x in overline{U(phi)}} < epsilon}
        \sup \{|\phi(x) - 1| : x \in \overline{U(\phi)}\} < \epsilon.
\end{equation}
Such an open set $U(\phi)$ exists because $\phi$ is continuous at
$0$, and we may also ask $\overline{U(\phi)}$ to be compact, since $A$
is locally compact.  The latter condition implies the set $W(\phi)$ of
$\psi \in \widehat{A}$ that satisfy
\begin{equation}
\label{sup {|psi(x) - 1| : x in overline{U(phi)}} < epsilon}
        \sup \{|\psi(x) - 1| : x \in \overline{U(\phi)}\} < \epsilon
\end{equation}
is an open set in $\widehat{A}$.  Of course, $\phi \in W(\phi)$ by
construction, so that the collection of $W(\phi)$ with $\phi \in E$
is an open covering of $E$ in $\widehat{A}$.  Thus there are finitely
many elements $\phi_1, \ldots, \phi_n$ of $E$ such that $E \subseteq
\bigcup_{j = 1}^n W(\phi_j)$, because $E$ is compact in $\widehat{A}$
by hypothesis.  If we put $U = \bigcap_{j = 1}^n U(\phi_j)$, then it
follows that $U$ is an open set in $A$ such that $0 \in U$ and
\begin{equation}
\label{sup {|psi(x) - 1| : x in overline{U}} < epsilon}
        \sup \{|\psi(x) - 1| : x \in \overline{U}\} < \epsilon
\end{equation}
for every $\psi \in E$, as desired.  Alternatively, the following is
perhaps a slightly more conventional version of the same type of
argument.  If $V$ is an open set in $A$ such that $0 \in V$ and $K =
\overline{V}$ is compact, then $E$ is totally bounded with respect to
$\|\cdot \|_K$, in the sense that that every element of $E$ can be
uniformly approximated on $K$ by elements of a finite subset of $E$.
This permits the equicontinuity of $E$ to be reduced to the
equicontinuity of a finite set of continuous functions, in essentially
the same way as before.

        Conversely, if $E \subseteq \widehat{A}$ is closed and
equicontinuous, then $E$ is compact.  To see this, let ${\bf T}^A$ be
the set of all functions on $A$ with values in ${\bf T}$, as in
Section \ref{discrete groups}.  As before, ${\bf T}^A$ is a compact
commutative topological group with respect to the product topology,
which contains the group of all homomorphisms from $A$ into ${\bf T}$
as a closed subgroup.  In particular, there is a natural inclusion
mapping of $\widehat{A}$ into ${\bf T}^A$ which is a continuous
homomorphism, because finite subsets of $A$ are compact.  Suppose that
$\psi$ is a homomorphism from $A$ into ${\bf T}$ which is in the
closure of $E$ in ${\bf T}^A$ with respect to the product topology.
It is easy to see that $\psi$ is continous at $0$, because $E$ is
equicontinuous.  This implies that $\psi$ is continuous at every point
in $A$, since $\psi$ is a homomorphism, and hence that $\psi \in
\widehat{A}$.  One can also use the equicontinuity of $E$ to show that
$\psi$ is in the closure of $E$ with respect to the topology on
$\widehat{A} \subseteq C(A)$.  More precisely, this means that $\psi$
can be approximated by elements of $E$ uniformly on compact subsets of
$A$, using the equicontinuity of $E$ to reduce to approximations
of $\psi$ by elements of $E$ on finite subsets of $A$.  This implies
that $E$ is a closed subset of ${\bf T}^A$ with respect to the
product topology, and hence that $E$ is compact with respect to the
product topology.  In order to show that $E$ is compact as a subset
of $\widehat{A}$, it suffices to check that the topology induced on
$E$ by the one on $\widehat{A} \subseteq C(A)$ is the same as the one
induced on $E$ by the product topology on ${\bf T}^A$.  As before,
this uses the equicontinuity of $E$, to show that $\phi, \phi' \in E$
are uniformly close on a compact set $K \subseteq A$ when $\phi$ and $\phi'$
are sufficiently close on a suitable finite subset of $K$.

        If $A$ is $\sigma$-compact, then the topology on $C(A)$ 
that we are using can be defined by a sequence of seminorms, as in the
previous section.  This implies that the same topology on $C(A)$ can
be defined by a translation-invariant metric, as in Section
\ref{semimetrics}.  One can also check that $C(A)$ is complete with
respect to such a metric, by standard arguments.  It is well known
that a subset of a complete metric space is compact if and only if it
is closed and totally bounded.  If $E \subseteq \widehat{A}$ is
equicontinuous, then one can show that $E$ is totally bounded in
$C(A)$, which basically means that $E$ is totally bounded with respect
to $\|f\|_K$ for each nonempty compact set $K \subseteq A$.  This uses
the fact that $\phi, \phi' \in E$ are uniformly close on $K$ when
$\phi$ and $\phi'$ are sufficiently close on a suitable finite subset
of $K$, as in the preceding paragraph.  Of course, if $E$ is
relatively closed in $\widehat{A}$, then $E$ is a closed set in
$C(A)$, because $\widehat{A}$ is a closed set in $C(A)$.  This leads
to a somwhat simpler proof of the compactness of $E \subseteq
\widehat{A}$ when $E$ is closed and equicontinuous and $A$ is
$\sigma$-compact.

\section{Local compactness of $\widehat{A}$}
\label{local compactness of widehat{A}}

        Let $A$ be a locally compact abelian topological group again.
Remember that the dual group $\widehat{A} = \hom(A, {\bf T})$ is also a
topological group with respect to the topology induced by the usual 
one on $C(A)$, as before.  Let $U$ be an open set in $A$ such that $0
\in U$ and $\overline{U}$ is compact, and put
\begin{equation}
\label{B_{overline{U}} = {phi in widehat{A} : sup |phi(x) - 1| < 1}}
        B_{\overline{U}} = \bigg\{\phi \in \widehat{A} : 
                            \sup_{x \in \overline{U}} |\phi(x) - 1| < 1\bigg\}.
\end{equation}
Thus $B_{\overline{U}}$ is an open set in $\widehat{A}$ that contains
the identity element of $\widehat{A}$, which is the character equal to
$1$ at every point in $A$.  In order to show that $\widehat{A}$ is
locally compact as well, it suffices to check that
\begin{equation}
\label{overline{B}_{overline{U}} = {phi in widehat{A} : sup |phi(x) - 1| le 1}}
        \overline{B}_{\overline{U}} = \bigg\{\phi \in \widehat{A} :
                         \sup_{x \in \overline{U}} |\phi(x) - 1| \le 1\bigg\}
\end{equation}
is compact in $\widehat{A}$.

        Of course, $\overline{B}_{\overline{U}}$ is a closed set in
$\widehat{A}$, by construction.  It remains to show that
$\overline{B}_{\overline{U}}$ is equicontinuous, because of the
discussion in the previous section.  Put $U_0 = U$ for convenience,
and let $U_1$ be an open set in $A$ such that $0 \in U_1$ and $U_1 +
U_1 \subseteq U_0$.  Continuing in this way, for each positive integer
$l$ there is an open set $U_l$ in $A$ such that $0 \in U_l$ and
\begin{equation}
\label{U_l + U_l subseteq U_{l - 1}}
        U_l + U_l \subseteq U_{l - 1}.
\end{equation}
If $x \in U_l$, then $2^j \, x \in U_{l - j}$ for $j = 1, \ldots, l$.
Of course, $\phi(n \, x) = \phi(x)^n$ for every $x \in A$, $\phi \in
\widehat{A}$, and $n \in {\bf Z}$, so that
\begin{equation}
\label{|phi(x)^{2^j} - 1| le 1}
        |\phi(x)^{2^j} - 1| \le 1
\end{equation}
for every $x \in U_l$, $\phi \in \overline{B}_{\overline{U}}$, and
$j = 0, 1, \ldots, l$.  It is easy to see that there is a sequence
$\{r_l\}_{l = 0}^\infty$ of positive real numbers converging to $0$
such that
\begin{equation}
\label{|phi(x) - 1| le r_l}
        |\phi(x) - 1| \le r_l
\end{equation}
for every $x \in U_l$, $\phi \in \overline{B}_{\overline{U}}$, and $l
\ge 0$ under these conditions.  This implies that
$\overline{B}_{\overline{U}}$ is equicontinuous, as desired.

        If $A$ is compact, then we can take $U = A$ in the previous 
discussion.  In this case, $B_{\overline{U}}$ and
$\overline{B}_{\overline{U}}$ contain only the trivial character on
$A$, which corresponds to the fact that $\widehat{A}$ is discrete.
Similarly, if $A$ is discrete, then we can take $U = \{0\}$,
so that $B_{\overline{U}} = \overline{B}_{\overline{U}} = \widehat{A}$,
which we have already seen is compact in this situation.

        Note that $\overline{B}_{\overline{U}}$ is a compact set in
$\widehat{A}$ for any open set $U \subseteq A$ with $0 \in U$,
by the same arguments as before.  The additional hypothesis that
$\overline{U}$ be compact is only needed to get $B_{\overline{U}}$
to be an open set in $\widehat{A}$.  If $U$ is an open subgroup in $A$,
then $U$ is also a closed set in $A$, and
\begin{equation}
\label{overline{B}_{overline{U}} = {phi in widehat{A} : phi(x) = 1, x in U}}
        \overline{B}_{\overline{U}} = \bigg\{\phi \in \widehat{A} :
                           \phi(x) = 1 \hbox{ for every } x \in U\bigg\},
\end{equation}
for the usual reasons.  This shows that $\overline{B}_{\overline{U}}$
is a compact subgroup of $\widehat{A}$ under these conditions.
Similiarly, if $K$ is any nonempty compact subset of $A$, then
\begin{equation}
\label{B_K = {phi in widehat{A} : sup_{x in K} |phi(x) - 1| < 1}}
 B_K = \bigg\{\phi \in \widehat{A} : \sup_{x \in K} |\phi(x) - 1| < 1\bigg\}
\end{equation}
is an open subset of $\widehat{A}$ that contains the identity element
in $\widehat{A}$.  If $K$ is a compact subgroup of $A$, then it
follows that
\begin{equation}
\label{B_K = {phi in widehat{A} : phi(x) = 1 for every x in K}}
 B_K = \{\phi \in \widehat{A} : \phi(x) = 1 \hbox{ for every } x \in K\},
\end{equation}
is an open subgroup of $\widehat{A}$.  In particular, if $K = U$ is a 
compact open subgroup of $A$, then $B_K = \overline{B}_{\overline{U}}$
is a compact open subgroup of $\widehat{A}$.

\section{Some additional properties}
\label{some additional properties}

        Let $X$ be a metric space, and let $C(X)$ be the space of
continuous complex-valued functions on $X$, with the usual topology
determined by the collection of supremum seminorms associated to
nonempty compact subsets of $X$.  If $X$ is compact, then one can
simply use the supremum norm on $C(X)$, and it is well known that
$C(X)$ is separable.  This uses the fact that continuous functions on
compact metric spaces are uniformly continuous.  Similarly, if $X$ is
locally compact and $\sigma$-compact, then we have seen that the usual
topology on $C(X)$ can be described by a metric, and one can also show
that $C(X)$ is separable.  This can be obtained from the case of
compact metric spaces, applied to an increasing sequence of compact
subsets of $X$, the union of whose interiors is $X$.

        Now let $A$ be a commutative topological group.  If $A$ is
compact, then we have seen that the topology on $\widehat{A} = \hom(A,
{\bf T})$ determined by the supremum norm on $C(A)$ is discrete.  If
$A$ is compact and metrizable, then $C(A)$ is separable, as in the
previous paragraph, which implies that $\widehat{A}$ can have only
finitely or countable many elements.  If $A$ is locally compact,
$\sigma$-compact, and metrizable, then $C(A)$ is still metrizable and
separable, which implies that $\widehat{A}$ is metrizable and
separable with respect to the induced topology.

        Alternatively, suppose that $A$ has a countable local base for its
topology at $0$.  This means that there is a sequence $U_1, U_2, U_3,
\ldots$ of open subsets of $A$ such that $0 \in U_j$ for each $j$, and
for each open set $V$ in $A$ with $0 \in V$ there is a positive
integer $j$ such that $U_j \subseteq V$.  If $B_j =
\overline{B}_{\overline{U_j}}$ is as in the previous section, then
$B_j$ is a compact set in $\widehat{A}$ for each $j$.  It is easy to
see that each $\phi \in \widehat{A}$ is contained in $B_j$ for some
$j$, so that $\widehat{A} = \bigcup_{j = 1}^\infty B_j$.  This shows
that $\widehat{A}$ is $\sigma$-compact when $A$ is metrizable.  In
particular, this gives another way to see that $\widehat{A}$ has only
finitely many elements when $A$ is compact and metrizable, so that
$\widehat{A}$ is discrete.  If $A$ is locally compact,
$\sigma$-compact, and metrizable, then it follows that $\widehat{A}$
is metrizable and $\sigma$-compact, and hence separable.

        More precisely, it is easy to see that compact metric spaces are 
totally bounded and hence separable, which implies that
$\sigma$-compact metric spaces are separable too.  In the other
direction, if $Y$ is a separable metric space, then it is well known
that there is a countable base for the topology of $Y$.  If $Y$ is
also locally compact, then it follows that $Y$ is $\sigma$-compact,
because every open covering of $Y$ can be reduced to a subcovering
with only finitely or countable many elements, by Lindel\"of's theorem.
If $A$ is a locally compact commutative topological group which is
$\sigma$-compact and metrizable, then $\widehat{A}$ is metrizable and
separable, by the argument at the beginning of the section.  In this
case, $\widehat{A}$ is also locally compact, as in the previous
section, and one can use separability and metrizability of
$\widehat{A}$ to get that $\widehat{A}$ is $\sigma$-compact.

        Similarly, let $G$ be a topological group, and let $H$ be an open 
subgroup of $G$.  If $G$ is separable, then it is easy to see that
there are only finitely or countably many left or right cosets of $H$
in $G$.  In particular, if $G$ is locally compact, then there is an
open subgroup $H$ of $G$ which is $\sigma$-compact, as in Section
\ref{additional properties}.  If $G$ is locally compact and separable,
then it follows that $G$ is $\sigma$-compact as well, because $G$ is
the union of finitely or countably many translates of $H$.  If $A$ is
a locally compact commutative topological group which is
$\sigma$-compact and metrizable, then $\widehat{A}$ is a locally
compact commutative topological group which is also separable, as
before, and one can apply the preceding argument with $G =
\widehat{A}$ to get that $\widehat{A}$ is $\sigma$-compact.

        Let $A$ be a commutative topological group again, and let $H$
and $K$ be nonempty compact subsets of $A$.  Note that
\begin{equation}
\label{H + K = {x + y : x in H, y in K}}
        H + K = \{x + y : x \in H, \, y \in K\}
\end{equation}
is also a compact set in $A$, because of continuity of addition on
$A$, and the fact that $H \times K$ is a compact subset of $A \times A$.  
If $\phi, \psi \in \widehat{A}$, $x \in H$, and $y \in K$, then
\begin{eqnarray}
\label{|phi(x + y) - psi(x + y)| le |phi(x) - psi(x)| + |phi(y) - psi(y)|}
 \lefteqn{|\phi(x + y) - \psi(x + y)| = 
|\phi(x) \, \phi(y) - \psi(x) \, \psi(y)|} \\
          & & \le |\phi(x) \, \phi(y) - \psi(x) \, \phi(y)| 
                  + |\psi(x) \, \phi(y) - \psi(x) \, \psi(y)| \nonumber \\
           & & = |\phi(x) - \psi(x)| + |\phi(y) - \psi(y)|. \nonumber
\end{eqnarray}
This implies that
\begin{equation}
\label{||phi - psi||_{H + K} le ||phi - psi||_H + ||phi - psi||_K}
        \|\phi - \psi\|_{H + K} \le \|\phi - \psi\|_H + \|\phi - \psi\|_K
\end{equation}
for every $\phi, \psi \in \widehat{A}$, where $\|f\|_E$ denotes the
supremum seminorm of $f \in C(A)$ associated to a nonempty compact set
$E \subseteq A$, with $E = H$, $K$, or $H + K$.

        If $E$ is any subset of $A$, then put
\begin{equation}
\label{-E = {-x : x in E}}
        -E = \{-x : x \in E\}.
\end{equation}
This was denoted $E^{-1}$ in Chapter \ref{topological groups}, where the
group operations were expressed multiplicatively.  Also put
\begin{equation}
\label{E(n) = {sum_{j = 1}^n x_j : x_j in E for each j = 1, ldots, n}}
        E(n) = \bigg\{\sum_{j = 1}^n x_j : x_j \in E \hbox{ for each } 
                                               j = 1, \ldots, n\bigg\},
\end{equation}
for each positive integer $n$, which was denoted $E^n$ in Chapter
\ref{topological groups}.  If $E$ is a nonempty compact set in $A$,
then $E(n)$ is also compact for each $n$, because of continuity of
addition on $A$, as before.  In this case,
\begin{equation}
\label{||phi - psi||_{E(n)} le n ||phi - psi||_E}
        \|\phi - \psi\|_{E(n)} \le n \, \|\phi - \psi\|_E
\end{equation}
for every $\phi, \psi \in \widehat{A}$, as one can see by applying
(\ref{||phi - psi||_{H + K} le ||phi - psi||_H + ||phi - psi||_K})
repeatedly.

        Suppose now that $A$ is locally compact, and let $V$ be an 
open set in $A$ such that $0 \in V$, $-V = V$, and $\overline{V}$ is
compact.  Thus $\bigcup_{n = 1}^\infty V(n)$ is an open subgroup of
$A$, as in Section \ref{additional properties}.  Let us ask also that
$A = \bigcup_{n = 1}^\infty V(n)$, which holds automatically when $A$
is connected.  If $K \subseteq A$ is nonempty and compact, then 
$K \subseteq V(n) \subseteq \overline{V}(n)$ for some $n$, and hence
\begin{equation}
\label{||phi - psi||_K le ... le n ||phi - psi||_{overline{V}}}
        \|\phi - \psi\|_K \le \|\phi - \psi\|_{\overline{V}(n)} 
                           \le n \, \|\phi - \psi\|_{\overline{V}}
\end{equation}
for every $\phi, \psi \in \widehat{A}$.  It follows that $\|\phi -
\psi\|_{\overline{V}}$ is a metric on $\widehat{A}$ that determines
the same topology on $\widehat{A}$ as the one induced by the usual
topology on $C(A)$ under these conditions.

\chapter{The Fourier transform}
\label{fourier transform}

\section{Integrable functions}
\label{integrable functions, 2}

         Throughout this chapter, we let $A$ be a locally compact abelian 
topological group, and we let $H$ be a translation-invariant Haar measure
on $A$.  If $A$ is discrete, then it is customary to take $H$ to be
counting measure on $A$, instead of some other multiple of counting
measure, and if $A$ is compact, then we can normalize $H$ so that
$H(A) = 1$.  The \emph{Fourier transform}\index{fourier transforms} of
a complex-valued integrable function $f$ on $A$ is the function
$\widehat{f}$\index{f(phi)@$\widehat{f}(\phi)$} defined on the dual
group $\widehat{A}$ by
\begin{equation}
\label{widehat{f}(phi) = int_A f(x) overline{phi(x)} dH(x)}
        \widehat{f}(\phi) = \int_A f(x) \, \overline{\phi(x)} \, dH(x)
\end{equation}
for each $\phi \in \widehat{A}$.  If $A$ is the unit circle ${\bf T}$,
for instance, then we have seen that $\widehat{A}$ is isomorphic to
the group ${\bf Z}$ with respect to addition, and
(\ref{widehat{f}(phi) = int_A f(x) overline{phi(x)} dH(x)}) reduces to
the usual definition (\ref{widehat{f}(n) = frac{1}{2 pi} int_{bf T}
  f(z) overline{z}^n |dz|}) of the Fourier coefficients of $f$.  As in
(\ref{|widehat{f}(n)| le frac{1}{2 pi} int_{bf T} |f(z)| |dz|}),
we have that
\begin{equation}
\label{|widehat{f}(phi)| le int_A |f(x)| dH(x)}
        |\widehat{f}(\phi)| \le \int_A |f(x)| \, dH(x)
\end{equation}
for every $\phi \in \widehat{A}$.

        Let $T_a(f)$ be the integrable function on $A$ defined for
each $a \in A$ by
\begin{equation}
\label{(T_a(f))(x) = f(x + a), 2}
        (T_a(f))(x) = f(x + a).
\end{equation}
The Fourier transform of $T_a(f)$ is given by
\begin{eqnarray}
\label{widehat{(T_a(f))}(phi)  = ... = phi(a) widehat{f}(phi)}
\widehat{(T_a(f))}(\phi) & = & \int_A f(x + a) \, \overline{\phi(x)} \, dH(x) \\
     & = & \int_A f(x) \, \overline{\phi(x - a)} \, dH(x) \nonumber \\
     & = & \phi(a) \, \int_A f(x) \, \overline{\phi(x)} \, dH(x)
              = \phi(a) \, \widehat{f}(\phi),                \nonumber
\end{eqnarray}
using the fact that $\phi$ is a homomorphism from $A$ into ${\bf T}$
in the third step.

        If $E$ is a Borel set in $A$, then $-E = \{-x : x \in E\}$
is also a Borel set in $A$, because $x \mapsto -x$ is a homeomorphism
from $A$ onto itself.  Let us check that
\begin{equation}
\label{H(-E) = H(E)}
        H(-E) = H(E).
\end{equation}
Because $H(-E)$ is a translation-invariant measure on $A$ that
satisfies the same conditions as Haar measure, it is equal to a
constant $c > 0$ times $H(E)$, and we would like to show that $c = 1$.
If $U$ is an open set in $A$ such that $0 \in U$ and $\overline{U}$ is
compact, then $V = U \cap (-U)$ has the same properties, and also
satisfies $-V = V$.  Thus $H(V)$ is positive, finite, and $H(-V) =
H(V)$, which implies that $c = 1$, as desired.

        Let $f$ be an integrable function on $A$ again, and consider
\begin{equation}
\label{g(x) = overline{f(-x)}}
        g(x) = \overline{f(-x)}.
\end{equation}
The Fourier transform of $g$ is given by
\begin{eqnarray}
\label{widehat{g}(phi) = ... = overline{widehat{f}(phi)}}
 \widehat{g}(\phi) & = & \int_A \overline{f(-x)} \, \overline{\phi(x)} \, dH(x)
                                                                            \\
   & = & \overline{\int_A f(-x) \, \phi(x) \, dH(x)}   \nonumber \\
   & = & \overline{\int_A f(x) \, \phi(-x) \, dH(x)}   \nonumber \\
   & = & \overline{\int_A f(x) \, \overline{\phi(x)} \, dH(x)}
           = \overline{\widehat{f}(\phi)}              \nonumber
\end{eqnarray}
for every $\phi \in \widehat{A}$.

\section{Complex Borel measures}
\label{complex borel measures}

        Let $\mu$ be a regular complex Borel measure on $A$.
The \emph{Fourier transform}\index{Fourier transforms} of $\mu$ is the
function $\widehat{\mu}$ defined on the dual group $\widehat{A}$ by
\begin{equation}
\label{widehat{mu}(phi) = int_A overline{phi(x)} d mu(x)}
        \widehat{\mu}(\phi) = \int_A \overline{\phi(x)} \, d\mu(x)
\end{equation}
for every $\phi \in \widehat{A}$.  If
\begin{equation}
\label{mu(E) = int_E f(x) dH(x)}
        \mu(E) = \int_E f(x) \, dH(x)
\end{equation}
for some integrable function $f$ on $A$ and every Borel set $E \subseteq A$,
then $\widehat{\mu} = \widehat{f}$.  As before, this definition reduces to
the previous one for the Fourier coefficients of a Borel measure on the unit
circle when $A = {\bf T}$.  We also have that
\begin{equation}
\label{|widehat{mu}(phi)| le |mu|(A)}
        |\widehat{\mu}(\phi)| \le |\mu|(A)
\end{equation}
for each $\phi \in \widehat{A}$, where $|\mu|$ is the total variation
measure on $A$ associated to $\mu$.

        Suppose that $K \subseteq A$ is nonempty and compact, and let $\mu_K$
be the Borel measure on $A$ defined by
\begin{equation}
\label{mu_K(E) = mu(E cap K)}
        \mu_K(E) = \mu(E \cap K)
\end{equation}
for each Borel set $E \subseteq A$.  Thus
\begin{equation}
\label{widehat{mu_K}(phi) = int_K overline{phi(x)} d mu(x)}
        \widehat{\mu_K}(\phi) = \int_K \overline{\phi(x)} \, d\mu(x)
\end{equation}
for each $\phi \in \widehat{\phi}$.  If $\phi, \psi \in \widehat{A}$, then
\begin{eqnarray}
\label{|widehat{mu_K}(phi) - widehat{mu_K}(psi)| le ...}
        |\widehat{\mu_K}(\phi) - \widehat{\mu_K}(\psi)|
            & = & \biggl|\int_K \overline{\phi(x)} \, d\mu(x) - 
                     \int_K \overline{\psi(x)} \, d\mu(x) \biggr| \\
 & \le & \Big(\sup_{x \in K} |\phi(x) - \psi(x)|\Big) \, |\mu|(K). \nonumber
\end{eqnarray}
This implies that $\widehat{\mu_K}(\phi)$ is continuous with respect
to the topology induced on $\widehat{A}$ by the usual one on $C(A)$,
defined by the supremum seminorms associated to nonempty compact
subsets of $A$.  More precisely, $\widehat{\mu_K}(\phi)$ is uniformly
continuous on $\widehat{A}$ as a topological group, because
(\ref{|widehat{mu_K}(phi) - widehat{mu_K}(psi)| le ...}) implies that
\begin{equation}
\label{|widehat{mu_K}(phi) - widehat{mu_K}(psi)| le ..., 2}
        |\widehat{\mu_K}(\phi) - \widehat{\mu_K}(\psi)|
 \le \Big(\sup_{x \in K} |\phi(x) \, \psi(x)^{-1} - 1|\Big) \, |\mu|(K)
\end{equation}
for every $\phi, \psi \in \widehat{A}$.

        To say that $\mu$ is a regular complex Borel measure on $A$
means that $|\mu|$ is a regular Borel measure on $A$.  This implies
that for every $\epsilon > 0$ there is a compact set $K \subseteq A$ 
such that
\begin{equation}
\label{|mu|(A backslash K) < epsilon}
        |\mu|(A \backslash K) < \epsilon.
\end{equation}
It follows that
\begin{equation}
\label{|widehat{mu}(phi) - widehat{mu_K}(phi)| = ... < epsilon}
        |\widehat{\mu}(\phi) - \widehat{\mu_K}(\phi)|
           = \biggl|\int_{A \backslash K} \overline{\phi(x)} \, d\mu(x)\biggr|
               \le |\mu|(A \backslash K) < \epsilon
\end{equation}
for every $\phi \in \widehat{A}$, so that $\widehat{\mu}$ can be
uniformly approximated on $\widehat{A}$ by functions of the form
$\widehat{\mu_K}$, where $K \subseteq A$ is compact.  This shows that
$\widehat{\mu}$ is uniformly continuous on $\widehat{A}$, since it can
be approximated uniformly on $\widehat{A}$ by uniformly continuous
functions.  Of course, if $A$ is compact, then $\widehat{A}$ is
discrete, and every function on $\widehat{A}$ is uniformly continuous
trivially.

\section{Vanishing at infinity}
\label{vanishing at infinity}

        Let $f$ be an integrable function on $A$, and let $T_a(f)$
be defined for $a \in A$ as in (\ref{(T_a(f))(x) = f(x + a), 2}) 
in Section \ref{integrable functions, 2}.  Thus the Fourier transform
of $T_a(f)$ is equal to $\phi(a) \, \widehat{f}(\phi)$, as in 
(\ref{widehat{(T_a(f))}(phi)  = ... = phi(a) widehat{f}(phi)}),
so that $(\phi(a) - 1) \, \widehat{f}(\phi)$ is the Fourier transform
of $T_a(f) - f$.  In particular,
\begin{equation}
\label{|phi(a) - 1| |widehat{f}(phi)| le int_A |f(x + a) - f(x)| dH(x)}
 |\phi(a) - 1| \, |\widehat{f}(\phi)| \le \int_A |f(x + a) - f(x)| \, dH(x)
\end{equation}
for every $a \in A$ and $\phi \in \widehat{A}$.  If $f$ is a
continuous function on $A$ with compact support, then we have seen
that $f$ is uniformly continuous on $A$, and hence
\begin{equation}
\label{int_A |f(x + a) - f(x)| dH(x) to 0}
        \int_A |f(x + a) - f(x)| \, dH(x) \to 0
\end{equation}
as $a \to 0$ in $A$.  This also holds for any integrable function $f$
on $A$, since we can approximate $f$ by continuous functions with
compact support with respect to the $L^1$ norm.

        Let $\epsilon > 0$ be given, and let $U$ be an open set in $A$
such that $0 \in U$ and
\begin{equation}
\label{int_A |f(x + a) - f(x)| dH(x) le epsilon}
        \int_A |f(x + a) - f(x)| \, dH(x) \le \epsilon
\end{equation}
for every $a \in \overline{U}$.  Combining this with (\ref{|phi(a) -
  1| |widehat{f}(phi)| le int_A |f(x + a) - f(x)| dH(x)}), we get that
\begin{equation}
\label{|phi(a) - 1| |widehat{f}(phi)| le epsilon}
        |\phi(a) - 1| \, |\widehat{f}(\phi)| \le \epsilon
\end{equation}
for every $a \in \overline{U}$ and $\phi \in \widehat{A}$.  Let
$\overline{B}_{\overline{U}}$ be as in (\ref{overline{B}_{overline{U}}
  = {phi in widehat{A} : sup |phi(x) - 1| le 1}}) in Section
\ref{local compactness of widehat{A}}, which is a compact subset of
$\widehat{A}$.  If $\phi \in \widehat{A} \backslash
\overline{B}_{\overline{U}}$, then
\begin{equation}
\label{|phi(a) - 1| > 1}
        |\phi(a) - 1| > 1
\end{equation}
for some $a \in \overline{U}$, and (\ref{|phi(a) - 1|
  |widehat{f}(phi)| le epsilon}) implies that
\begin{equation}
\label{|widehat{f}(phi)| < epsilon}
        |\widehat{f}(\phi)| < \epsilon.
\end{equation}
This shows that $\widehat{f}$ vanishes at infinity on $\widehat{A}$.
If $A$ is discrete, then $\widehat{A}$ is compact, and every function
on $A$ has this property trivially.  Of course, if $A$ is discrete,
then one can take $U = \{0\}$, and (\ref{int_A |f(x + a) - f(x)| dH(x)
  le epsilon}) is trivial too.

        Let us give another proof of this when $A$ is compact, which
is like the one for $A = {\bf T}$ in (\ref{lim_{|n| to infty}
  |widehat{f}(n)| = 0}) in Section \ref{basic notions}.  
If $f \in L^2(A)$, then
\begin{equation}
\label{widehat{f}(phi) = langle f, phi rangle}
        \widehat{f}(\phi) = \langle f, \phi \rangle
\end{equation}
for every $\phi \in \widehat{A}$, using the standard integral inner
product on $L^2(A)$.  Because characters on $A$ are orthonormal in
$L^2(A)$, we get that $\widehat{f}(\phi)$ is square-summable on
$\widehat{A}$, with
\begin{equation}
\label{sum_{phi in widehat{A}} |widehat{f}(phi)|^2 le int_A |f(x)|^2 dH(x)}
 \sum_{\phi \in \widehat{A}} |\widehat{f}(\phi)|^2 \le \int_A |f(x)|^2 \, dH(x).
\end{equation}
In particular, $\widehat{f}(\phi)$ vanishes at infinity on
$\widehat{A}$ when $f \in L^2(A)$.  If $f$ is an integrable function
on $A$, then we can get the same conclusion by approximating $f$ by
square-integrable functions on $A$ with respect to the $L^1$ norm,
using also (\ref{|widehat{f}(phi)| le int_A |f(x)| dH(x)}).

\section{Convolution of integrable functions}
\label{convolution of integrable functions}

        Let $f$, $g$ be nonnegative Borel measurable functions on $A$,
and put
\begin{equation}
\label{(f * g)(x) = int_A f(x - y) g(y) dH(y)}
        (f * g)(x) = \int_A f(x - y) \, g(y) \, dH(y).
\end{equation}
Because of Fubini's theorem, we get that
\begin{eqnarray}
\label{int_A (f * g)(x) dH(x) = int_A int_A f(x - y) g(y) dH(y) dH(x) = ...}
 \int_A (f * g)(x) \, dH(x) & = & \int_A \int_A f(x - y) \, g(y) dH(y) \, dH(x) 
                                                                        \\
 & = & \int_A \int_A f(x - y) \, g(y) \, dH(x) \, dH(y) \nonumber \\
 & = & \int_A \int_A f(x) \, g(y) \, dH(x) \, dH(y) \nonumber \\
 & = & \Big(\int_A f(x) \, dH(x)\Big) \, \Big(\int_A g(y) \, dH(y)\Big),
                                                              \nonumber
\end{eqnarray}
using also translation-invariance of Haar measure in the third step.
In particular, if $f$ and $g$ are integrable on $A$, then $f * g$ is
also integrable on $A$, and hence finite almost everywhere.

        More precisely, $f(x - y)$ is a Borel measurable function on
$A \times A$ when $f$ is Borel measurable on $A$, because $(x, y) \mapsto
x - y$ is a continuous mapping from $A \times A$ into $A$.  This implies
that $f(x - y) \, g(y)$ is also Borel measurable on $A \times A$.  As
in Section \ref{double integrals}, there are some additional
technicalities related to the way that the product measure is defined
on $A \times A$.  If there is a countable base for the topology of
$A$, then every open covering of $A$ can be reduced to a subcovering
with only finitely or countable many elements.  In this case,
the local compactness of $A$ implies that $A$ is $\sigma$-compact,
and hence that Haar measure on $A$ is $\sigma$-finite.  This permits one
to use the standard construction of the product measure on $A \times A$.
One also has a countable base for the topology of $A \times A$, consisting
of products of basic open subsets of $A$, so that open subsets of
$A \times A$ can be expressed as countable unions of products of
open subsets of $A$.  This implies that open subsets of $A \times A$
are measurable with respect to the standard product measure construction,
and hence that Borel subsets of $A \times A$ are measurable with respect
to the product measure as well.  Alternatively, the product measure
can be defined as a Borel measure with suitable regularity properties.

        If $f$ and $g$ are integrable complex-valued functions on $A$,
then one would like to define their convolution\index{convolutions} $f * g$
in the same way.  One can first apply the previous discussion to
$|f|$ and $|g|$, to get that
\begin{equation}
\label{int_A |f(x - y)| |g(y)| dH(y) < infty}
        \int_A |f(x - y)| \, |g(y)| \, dH(y) < \infty
\end{equation}
for almost every $x \in A$ with respect to $H$.  This means that
$f(x - y) \, g(y)$ is an integrable function of $y$ for almost
every $x \in A$, so that (\ref{(f * g)(x) = int_A f(x - y) g(y) dH(y)})
is defined for almost every $x \in A$.  Of course,
\begin{equation}
        |(f * g)(x)| \le \int_A |f(x - y)| \, |g(y)| \, dH(y)
\end{equation}
when $(f * g)(x)$ is defined, and we can integrate this in $x$ and
interchange the order of integration as in (\ref{int_A (f * g)(x)
  dH(x) = int_A int_A f(x - y) g(y) dH(y) dH(x) = ...})  to get that
\begin{equation}
\label{int_A |(f * g)(x)| dH(x) le (int_A |f(x)| dH(x)) (int_A |g(y)| dH(y))}
 \int_A |(f * g)(x)| \, dH(x) \le \Big(\int_A |f(x)| \, dH(x)\Big) \,
                                   \Big(\int_A |g(y)| \, dH(y)\Big).
\end{equation}
This shows that $f * g$ is integrable on $A$ when $f$ and $g$ are
integrable.

        It is easy to see that
\begin{equation}
\label{(f * g)(x) = (g * f)(x)}
        (f * g)(x) = (g * f)(x)
\end{equation}
when (\ref{int_A |f(x - y)| |g(y)| dH(y) < infty}) holds, using the
change of variables $y \mapsto x - y$.  Similarly, one can check that
\begin{equation}
        (f * g) * h = f * (g * h)
\end{equation}
for every $f, g, h \in L^1(A)$.  Let us show that
\begin{equation}
\label{widehat{(f * g)}(phi) = widehat{f}(phi) widehat{g}(phi)}
        \widehat{(f * g)}(\phi) = \widehat{f}(\phi) \, \widehat{g}(\phi)
\end{equation}
for every $\phi \in \widehat{A}$ when $f$ and $g$ are integrable functions
on $A$.  By the definition of the Fourier transform,
\begin{eqnarray}
\label{widehat{(f * g)}(phi) = ...}
 \widehat{(f * g)}(\phi) & = & \int_A (f * g)(x) \, \overline{\phi(x)} \, dH(x)
                                                                           \\
 & = & \int_A \int_A f(x - y) \, g(y) \, \overline{\phi(x)} \, dH(y) \, dH(x)
                                                               \nonumber \\
 & = & \int_A \int_A f(x - y) \, \overline{\phi(x - y)} \, g(y) \, 
                                \overline{\phi(y)} \, dH(y) \, dH(x), \nonumber
\end{eqnarray}
using the fact that $\phi$ is a character on $A$ in the last step.
Interchanging the order of integration as in (\ref{int_A (f * g)(x)
  dH(x) = int_A int_A f(x - y) g(y) dH(y) dH(x) = ...}), and using
translation-invariance of the resulting integral in $x$, we get that
this is equal to
\begin{equation}
\label{... = widehat{f}(phi) widehat{g}(phi)}
        \Big(\int_A f(x) \, \overline{\phi(x)} \, dH(x)\Big) \,
         \Big(\int_A g(y) \, \overline{\phi(y)} \, dH(y)\Big)
            = \widehat{f}(\phi) \, \widehat{g}(\phi),
\end{equation}
as desired.

\section{Convolution of other functions}
\label{convolution of other functions}

        If $f$ and $g$ are continuous complex-valued functions
with compact support on $A$, then $(f * g)(x)$ is defined for every $x
\in A$, and indeed it can be defined in terms of the Haar integral as
a nonnegative linear functional on $C_{com}(A)$.\index{convolutions}
One can also check that $f * g$ is a continuous function on $A$ with
compact support, using the fact that continuous functions on $A$ are
uniformly continuous on compact sets.  As in Section \ref{double
  integrals}, interchanging the order of integration of a continuous
function on $A \times A$ with compact support can be seen in a more
elementary way, which can be used to derive the same properties of $f
* g$ as in the previous section.  Of course, integrable functions on
$A$ can be approximated by continuous functions on $A$ with compact
support with respect to the $L^1$ norm, which gives another way to
look at the convolution of integrable functions on $A$.

        Suppose now that $f \in L^p(A)$ and $g \in L^q(A)$, where 
$1 < p, q < \infty$ are conjugate exponents, so that $1/p + 1/q = 1$.
H\"older's inequality implies that
\begin{eqnarray}
\label{int_A |f(x - y)| |g(y)| dH(y) le ...}
\lefteqn{\int_A |f(x - y)| \, |g(y)| \, dH(y)}     \\
         & \le & \Big(\int_A |f(x - y)|^p \, dH(y)\Big)^{1/p} \, 
                               \Big(\int_A |g(y)|^q\Big)^{1/q} \nonumber \\
           & = & \Big(\int_A |f(y)|^p \, dH(y)\Big)^{1/p} \, 
                       \Big(\int_A |g(y)|^q \, dH(y)\Big)^{1/q}, \nonumber
\end{eqnarray}
for every $x \in A$, using the change of variables $y \mapsto x - y$
in the second step.  Thus $(f * g)(x)$\index{convolutions} can be
defined as in (\ref{(f * g)(x) = int_A f(x - y) g(y) dH(y)}) for every
$x \in A$, and satisfies
\begin{equation}
\label{|(f * g)(x)| le (int_A |f(y)|^p dH)^{1/p} (int_A |g(z)|^q dH)^{1/q}}
        |(f * g)(x)| \le \Big(\int_A |f(y)|^p \, dH(y)\Big)^{1/p} \,
                          \Big(\int_A |g(z)|^q \, dH(z)\Big)^{1/q}.
\end{equation}
One can also check that $f * g$ is a continuous function on $A$ that
vanishes at infinity under these conditions, by approximating $f$ and
$g$ by continuous functions with compact support.  Note that (\ref{(f
  * g)(x) = (g * f)(x)}) holds in this case as well, for the same
reasons as before.

        Similarly, if $f$ is an integrable function on $A$ and $g$
is a bounded Borel measurable function on $A$, then
\begin{eqnarray}
\label{int_A |f(x - y)| |g(y)| dH(y) le ..., 2}
 \int_A |f(x - y)| \, |g(y)| \, dH(y) & \le &
  \Big(\int_A |f(x - y)| \, dH(y)\Big) \, \Big(\sup_{z \in A} |g(z)|\Big) \\
 & = & \Big(\int_A |f(y)| \, dH(y)\Big) \, \Big(\sup_{z \in A} |g(z)|\Big)
                                                                  \nonumber
\end{eqnarray}
for every $x \in A$.  This implies that $(f *
g)(x)$\index{convolutions} can be defined as in (\ref{(f * g)(x) =
  int_A f(x - y) g(y) dH(y)}) for every $x \in A$ again, and satisfies
\begin{equation}
\label{|(f * g)(x)| le (int_A |f(y)| dH(y)) (sup_{z in A} |g(z)|)}
        |(f * g)(x)| \le \Big(\int_A |f(y)| \, dH(y)\Big) \, 
                            \Big(\sup_{z \in A} |g(z)|\Big).
\end{equation}
If $f$ is a continuous function on $A$ with compact support, then $f$
is uniformly continuous on $A$, and it is easy to see that $f * g$ is
uniformly continuous on $A$ too.  If $f$ is an integrable function on
$A$, then one can approximate $f$ by continuous functions with compact
support, to get that $f * g$ is uniformly continuous on $A$.  Note
that $f * g$ is constant when $g$ is constant, so that $f * g$ may not
vanish at infinity when $A$ is not compact.

        As in Section \ref{integrable functions, 2}, the uniqueness
of Haar measure can be used to show that $H(-E) = H(E)$ for every
Borel set $E \subseteq A$.  Alternatively, we can choose Haar measure
to have this property by replacing $H(E)$ with $(H(E) + H(-E))/2$, if
necessary.  Let $f$ and $g$ be continuous functions on $A$ with
compact support again, and suppose that $H'$ is another Haar measure
on $A$.  If $f * g$ is the convolution of $f$ and $g$ with respect to
$H$, then
\begin{eqnarray}
\label{int_A (f * g)(x) dH'(x) = ...}
 \int_A (f * g)(x) \, dH'(x) & = & 
                    \int_A \int_A f(x - y) \, g(y) \, dH(y) \, dH'(x) \\
 & = & \int_A \int_A f(x - y) \, g(y) \, dH'(x) \, dH(y) \nonumber \\
 & = & \int_A \int_A f(x) \, g(y) \, dH'(x) \, dH(y) \nonumber \\
 & = & \Big(\int_A f(x) \, dH'(x)\Big) \, \Big(\int_A g(y) \, dH(y)\Big)
                                                      \nonumber
\end{eqnarray}
as in (\ref{int_A (f * g)(x) dH(x) = int_A int_A f(x - y) g(y) dH(y)
  dH(x) = ...}).  Because $H$ is invariant under translations and the
mapping $x \mapsto -x$ on $A$, we still have that $f * g = g * f$.
This implies that
\begin{eqnarray}
\label{int_A (f * g)(x) dH'(x) = int_A (g * f)(x) dH'(x) = ...}
 \int_A (f * g)(x) \, dH'(x) & = & \int_A (g * f)(x) \, dH'(x) \\
 & = & \Big(\int_A g(x) \, dH'(x)\Big) \, \Big(\int_A f(y) \, dH(y)\Big),
                                                             \nonumber
\end{eqnarray}
as in (\ref{int_A (f * g)(x) dH'(x) = ...}), and hence
\begin{eqnarray}
\label{(int_A f(x) dH'(x)) (int_A g(y) dH(y)) ...}
\lefteqn{\Big(\int_A f(x) \, dH'(x)\Big) \, \Big(\int_A g(y) \, dH(y)\Big)} \\
& = & \Big(\int_A g(x) \, dH'(x)\Big) \, \Big(\int_A f(y) \, dH(y)\Big)
                                                               \nonumber
\end{eqnarray}
for every $f, g \in C_{com}(A)$.  If $g$ is a nonnegative real-valued
continuous function on $A$ with compact support such that $g(x) > 0$
for some $x \in A$, then the integrals of $g$ with respect to $H$ and
$H'$ are both positive and finite.  This shows the integrals of $f$
with respect to $H$ and $H'$ are the same up to multiplication by a
positive constant, so that $H'$ is a constant multiple of $H$.

\section{Convolution of measures}
\label{convolution of measures}

        Let $\mu$ and $\nu$ be regular complex Borel measures on $A$,
and let $\mu \times \nu$ be the corresponding product regular Borel measure
on $A \times A$.  As usual, this is all a bit simpler when there is a
countable base for the topology of $A$, so that the standard construction
of product measures can be used.  Otherwise, one can get $\mu \times \nu$
from the corresponding bounded linear functional on $C_0(A \times A)$,
as in Section \ref{double integrals}.  If $E$ is a Borel set in $A$, then
\begin{equation}
\label{E' = {(x, y) in A times A : x + y in E}}
        E' = \{(x, y) \in A \times A : x + y \in E\}
\end{equation}
is a Borel set in $A \times A$, because $(x, y) \mapsto x + y$
is a continuous mapping from $A \times A$ into $A$.  It is easy to see that
\begin{equation}
\label{(mu * nu)(E) = (mu times nu)(E')}
        (\mu * \nu)(E) = (\mu \times \nu)(E')
\end{equation}
defines a complex Borel measure on $A$, known as the
convolution\index{convolutions} of $\mu$ and $\nu$.

        Similarly, if $\phi$ is a complex-valued Borel measurable function
on $A$, then $\phi(x + y)$ is Borel measurable on $A \times A$, again because
$(x, y) \mapsto x + y$ is continuous.  If $\phi$ is also bounded on $A$, 
then it follows that
\begin{equation}
\label{int_A phi d(mu * nu) = int_A int_A phi(x + y) dmu(x) dnu(y)}
 \int_A \phi \, d(\mu * \nu) = \int_A \int_A \phi(x + y) \, d\mu(x) \, d\nu(y).
\end{equation}
More precisely, this is equivalent to (\ref{(mu * nu)(E) = (mu times
  nu)(E')}) when $\phi$ is the indicator function associated to the
Borel set $E$, and otherwise one can reduce to that case by
approximating $\phi$ by simple functions.

        Alternatively, if $A$ is compact, then one can use this to define
the bounded linear functional on $C(A)$ that corresponds to $\mu *
\nu$ in terms of the bounded linear functionals on $C(A)$
corresponding to $\mu$ and $\nu$.  If $A$ is not compact, then one can
still try to define the bounded linear functional on $C_0(A)$
corresponding to $\mu * \nu$ in this way, but $\phi(x + y)$ does not
vanish at infinity on $A \times A$ when $\phi(x) \ne 0$ for some $x
\in A$, and so it is better to be able to integrate at least bounded
continuous functions on $A \times A$.  If $\mu$ or $\nu$ has compact
support, then one can extend the corresponding linear functional on
$C_0(A)$ to $C(A)$, and avoid this problem.  In particular, if $\mu$
and $\nu$ both have compact support, then one can get a linear
functional on $C(A \times A)$ from the linear functionals on $C(A)$
associated to $\mu$ and $\nu$.  If $\mu$ and $\nu$ do not have compact
support, then there are still relatively simple approximation
arguments for extending these linear functionals to bounded continuous
functions.

        At any rate, an advantage to defining the convolution in terms 
of bounded linear functionals on $C_0(A)$ is that one automatically
gets the regularity of the corresponding Borel measures.  Otherwise,
one can show directly that $\mu * \nu$ is regular when $\mu$ and $\nu$
are regular, using the regularity of $\mu \times \nu$ on $A \times A$.
More precisely, one can begin by reducing to the case where $\mu$ and
$\nu$ are real-valued and nonnegative, using the Jordan decomposition.
There is another trick, which is to first check that $\mu * \nu$ is
inner regular, and then get outer regularity by taking complements.
Let $E$ be a Borel set in $A$, and suppose that $H \subseteq A \times A$
is a compact set such that $H \subseteq E'$.  Note that
\begin{equation}
\label{K = {x + y : (x, y) in H}}
        K = \{x + y : (x, y) \in H\}
\end{equation}
is a compact set in $A$, because $(x, y) \mapsto x + y$ is continuous.
By construction, $K \subseteq E$ and $H \subseteq K' \subseteq E'$.
If $H$ approximates $E'$ well with respect to $\mu \times \nu$, then
$K$ approximates $E$ well with respect to $\mu * \nu$, as desired.

        Suppose that $f$ and $g$ are integrable functions on $A$, and that
\begin{equation}
\label{mu(E) = int_E f(x) dH(x), nu(E) = int_E g(x) dH(x)}
        \mu(E) = \int_E f(x) \, dH(x), \quad \nu(E) = \int_E g(x) \, dH(x)
\end{equation}
for all Borel sets $E \subseteq A$.  Let us check that $\mu * \nu$
corresponds to $f * g$ in the same way, so that
\begin{equation}
\label{(mu * nu)(E) = int_E (f * g)(x) dH(x)}
        (\mu * \nu)(E) = \int_E (f * g)(x) \, dH(x)
\end{equation}
for all Borel sets $E \subseteq A$.  Equivalently,
\begin{equation}
\label{int_A phi d(mu * nu) = int_A phi(x) (f * g)(x) dH(x)}
        \int_A \phi \, d(\mu * \nu) = \int_A \phi(x) \, (f * g)(x) \, dH(x)
\end{equation}
for every bounded complex-valued Borel measurable function $\phi$ on
$A$.  To see this, observe that
\begin{eqnarray}
\label{int_A phi(x) (f * g)(x) dH(x) = ...}
        \quad  \int_A \phi(x) \, (f * g)(x) \, dH(x) & = &
           \int_A \int_A \phi(x) \, f(x - y) \, g(y) \, dH(y) \, dH(x) \\
 & = & \int_A \int_A \phi(x) \, f(x - y) \, g(y) \, dH(x) \, dH(y) \nonumber \\
 & = & \int_A \int_A \phi(x + y) \, f(x) \, g(y) \, dH(x) \, dH(y), \nonumber
\end{eqnarray}
using the definition of $f * g$ in the first step, Fubini's theorem in
the second step, and translation-invariance of Haar measure in the
third step.  This is the same as the right side of (\ref{int_A phi
  d(mu * nu) = int_A int_A phi(x + y) dmu(x) dnu(y)}) in this case,
as desired.

        It is easy to see that $\mu * \nu = \nu * \mu$ for any two
regular complex Borel measures $\mu$, $\nu$ on $A$, using
commutativity of addition.  One can also show that convolution is
associative, using associativity of addition and Fubini's theorem.
More precisely, if $E$ is a Borel set in $A$, then
\begin{equation}
\label{E'' = {(x, y, z) in A times A times A : x + y + z in E}}
        E'' = \{(x, y, z) \in A \times A \times A : x + y + z \in E\}
\end{equation}
is a Borel set in $A \times A \times A$, and the convolution of three
measures can be described in terms of the measure of $E''$ with
respect to the correspondng product measure on $A \times A \times A$.
Similarly, if $\phi$ is a bounded Borel measurable function on $A$,
then the integral of $\phi$ with respect to the convolution of three
measures on $A$ can be expressed in terms of a triple integral of
$\phi(x + y + z)$ on $A \times A \times A$.  Using Fubini's theorem,
one can check that the two different ways of convolving three measures
give the same results.

        Remember that $|\mu|$ denotes the total variation measure
corresponding to a complex Borel measure $\mu$ on $A$, and that
$\|\mu\| = |\mu|(A)$ defines a norm on the vector space of complex
Borel measures on $A$.  If $\mu$, $\nu$ are regular complex Borel measures
on $A$, then one can check that
\begin{equation}
\label{|(mu * nu)(E)| le (|mu| * |nu|)(E)}
        |(\mu * \nu)(E)| \le (|\mu| * |\nu|)(E)
\end{equation}
for every Borel set $E \subseteq A$.  This implies that
\begin{equation}
\label{|mu * nu|(E) le (|mu| * |nu|)(E)}
        |\mu * \nu|(E) \le (|\mu| * |\nu|)(E),
\end{equation}
because of the way that the total variation measure is defined.  In
particular,
\begin{equation}
\label{|mu * nu|(A) le (|mu| * |nu|)(A) = |mu|(A) cdot |nu|(A)}
        |\mu * \nu|(A) \le (|\mu| * |\nu|)(A) = |\mu|(A) \cdot |\nu|(A),
\end{equation}
using the fact that $A' = A \times A$ in the lest step.  Thus
\begin{equation}
\label{||mu * nu|| le ||mu|| ||nu||}
        \|\mu * \nu\| \le \|\mu\| \, \|\nu\|,
\end{equation}
which extends the fact that the $L^1$ norm of the convolution of two
integrable functions is less than or equal to the product of the $L^1$
norms of the two functions, as in (\ref{int_A |(f * g)(x)| dH(x) le
  (int_A |f(x)| dH(x)) (int_A |g(y)| dH(y))}) in Section
\ref{convolution of integrable functions}.

        Let $\delta_a$\index{Dirac mass@Dirac mass $\delta_a$} be the 
Dirac mass at the point $a \in A$, which is the measure defined on $A$
by putting $\delta_a(E)$ equal to $1$ when $a \in E$ and equal to $0$
otherwise.  Note that $\|\delta_a\| = 1$ for each $a \in A$, and that
\begin{equation}
\label{delta_0 * mu = mu * delta_0 = mu}
        \delta_0 * \mu = \mu * \delta_0 = \mu
\end{equation}
for every complex Borel measure $\mu$ on $A$.  By construction,
$\widehat{\delta_a}(\phi) = \overline{\phi(a)}$ for every $a \in A$
and $\phi \in \widehat{A}$.  If $\mu$ and $\nu$ are complex Borel
measures on $A$, then
\begin{eqnarray}
\label{widehat{(mu * nu)}(phi) = ... = widehat{mu}(phi) widehat{nu}(phi)}
 \widehat{(\mu * \nu)}(\phi) & = & \int_A \overline{\phi} \, d(\mu * \nu)
          = \int_A \int_A \overline{\phi(x + y)} \, d\mu(x) \, d\nu(y)  \\
 & = & \int_A \int_A \overline{\phi(x)} \, \overline{\phi(y)} \, d\mu(x) \, 
                                                       d\nu(y) \nonumber \\
 & = & \Big(\int_A \overline{\phi(x)} \, d\mu(x)\Big) \,
                              \Big(\int_A \overline{\phi(y)} \, d\nu(y)\Big)
                  = \widehat{\mu}(\phi) \, \widehat{\nu}(\phi) \nonumber
\end{eqnarray}
for every $\phi \in \widehat{A}$, extending the analogous fact
(\ref{widehat{(f * g)}(phi) = widehat{f}(phi) widehat{g}(phi)}) for
integrable functions.  In particular, this is consistent with
(\ref{delta_0 * mu = mu * delta_0 = mu}), because
$\widehat{\delta_0}(\phi) = 1$ for every $\phi \in \widehat{A}$.

\section{Functions and measures}
\label{functions, measures}

        If $f$ is an integrable complex-valued function on $A$
and $\nu$ is a regular complex Borel measure on $A$, then we would
like to define their convolution as a function on $A$ by
\begin{equation}
\label{(f * nu)(x) = int_A f(x - y) d nu(y)}
        (f * \nu)(x) = \int_A f(x - y) \, d\nu(y).
\end{equation}
Let us start with the case where $f$ are $\nu$ are real-valued and
nonnegative, so that (\ref{(f * nu)(x) = int_A f(x - y) d nu(y)})
makes sense as a nonnegative extended real number.  Using Fubini's
theorem and translation-invariance of Haar measure, we get that
\begin{eqnarray}
\label{int_A (f * nu)(x) dH(x) = ... = (int_A f(x) dH(x)) nu(A)}
\int_A (f * \nu)(x) \, dH(x) & = & \int_A \int_A f(x - y) \, d\nu(y) \, dH(x) \\
         & = & \int_A \int_A f(x - y) \, dH(x) \, d\nu(y) \nonumber \\
         & = & \int_A \int_A f(x) \, dH(x) \, d\nu(y)     \nonumber \\
         & = & \Big(\int_A f(x) \, dH(x)\Big) \, \nu(A).  \nonumber
\end{eqnarray}
This shows that $f * \nu$ is integrable on $A$ with respect to Haar
measure when $f$ is a nonnegative integrable function on $A$ and $\nu$
is a finite nonnegative regular Borel measure on $A$.

        If $f$ and $\nu$ are complex-valued, then we can apply the previous
argument to $|f|$ and $|\nu|$, which implies in particular that
\begin{equation}
\label{int_A |f(x - y)| d |nu|(y) < infty}
        \int_A |f(x - y)| \, d|\nu|(y) < \infty
\end{equation}
for almost every $x \in A$ with respect to Haar measure.  Thus $(f *
\nu)(x)$ may be defined as in (\ref{(f * nu)(x) = int_A f(x - y) d
  nu(y)}) for almost every $x \in A$ with respect to $H$, and satisfies
\begin{equation}
\label{|(f * nu)(x)| le int_A |f(x - y)| d|nu|(y)}
        |(f * \nu)(x)| \le \int_A |f(x - y)| \, d|\nu|(y).
\end{equation}
Integrating this as before, we get that
\begin{eqnarray}
\label{int_A |(f * nu)(x)| dH(x) le ... = (int_A |f(x)| dH(x)) |nu|(A)}
 \int_A |(f * \nu)(x)| \, dH(x) & \le & \int_A \int_A |f(x - y)| \, d|\nu|(y)
                                                                  \, dH(x) \\
         & = & \Big(\int_A |f(x)| \, dH(x)\Big) \, |\nu|(A),     \nonumber
\end{eqnarray}
so that $f * \nu$ is also integrable with respect to $H$ on $A$.  If
$f(x) = 0$ for almost every $x \in A$ with respect to $H$, then it
follows that $(f * \nu)(x) = 0$ almost everywhere on $A$ with respect
to $H$ too.

        Let $g$ be a complex-valued integrable function on $A$,
and consider
\begin{equation}
\label{nu(E) = int_E g(y) dH(y)}
        \nu(E) = \int_E g(y) \, dH(y)
\end{equation}
as a complex Borel measure on $A$.  In this case, $f * \nu$ is the
same as $f * g$.  Now let $\nu$ be any complex regular Borel measure
on $A$ again, and consider
\begin{equation}
\label{mu(E) = int_E f(x) dH(x), 2}
        \mu(E) = \int_E f(x) \, dH(x)
\end{equation}
as another Borel measure on $A$.  We would like to check that
\begin{equation}
\label{(mu * nu)(E) = int_E (f * nu)(x) dH(x)}
        (\mu * \nu)(E) = \int_E (f * \nu)(x) \, dH(x)
\end{equation}
for every Borel set $E \subseteq A$, or equivalently that
\begin{equation}
\label{int_A phi d(mu * nu) = int_A phi(x) (f * nu)(x) dH(x)}
        \int_A \phi \, d(\mu * \nu) = \int_A \phi(x) \, (f * \nu)(x) \, dH(x)
\end{equation}
for every bounded complex-valued Borel measurable function $\phi$ on $A$.
As in (\ref{int_A phi(x) (f * g)(x) dH(x) = ...}), we have that
\begin{eqnarray}
\label{int_A phi(x) (f * nu)(x) dH(x) = ...}
        \int_A \phi(x) \, (f * \nu)(x) \, dH(x) & = &
         \int_A \int_A \phi(x) \, f(x - y) \, d\nu(y) \, dH(x) \\
 & = & \int_A \int_A \phi(x) \, f(x - y) \, dH(x) \, d\nu(y) \nonumber \\
 & = & \int_A \int_A \phi(x + y) \, f(x) \, dH(x) \, d\nu(y), \nonumber
\end{eqnarray}
which is the same as the right side of (\ref{int_A phi d(mu * nu) =
  int_A int_A phi(x + y) dmu(x) dnu(y)}) in this case, as desired.

        If $f \in C_{com}(A)$, then $(f * \nu)(x)$ is defined for every
$x \in A$, and it is easy to see that $f * \nu$ is continuous on $A$ too, 
using the fact that $f$ is uniformly continuous on $A$.  If $\nu$ has
compact support on $A$, so that there is a compact set $K \subseteq A$
such that $|\nu|(A \backslash K) = 0$, then $f * \nu$ has compact
support as well.  Otherwise, the regularity of $\nu$ implies that for
each $\epsilon > 0$ there is a compact set $K(\epsilon) \subseteq A$
such that $|\nu|(A \backslash K(\epsilon)) < \epsilon$, and one can
use this to check that $f * \nu \in C_0(A)$ when $f \in C_{com}(A)$.

        Let $f$ be a nonnegative real-valued Borel measurable function
on $A$ again, and let $\nu$ be a nonnegative real-valued regular Borel
measure on $A$, so that $(f * \nu)(x)$ is defined as a nonnegative
extended real number for each $x \in A$ as before.  If $1 < p < \infty$, 
then
\begin{equation}
\label{((f * nu)(x))^p le nu(A)^{p - 1} int_A f(x - y)^p d nu(y)}
        ((f * \nu)(x))^p \le \nu(A)^{p - 1} \, \int_A f(x - y)^p \, d\nu(y)
\end{equation}
for every $x \in A$, by Jensen's or H\"older's inequality.  Thus
\begin{eqnarray}
\label{int_A ((f * nu)(x))^p dH(x) le ... = nu(A)^p (int_A f(x)^p dH(x))}
        \int_A ((f * \nu)(x))^p \, dH(x) & \le &
         \nu(A)^{p - 1} \, \int_A \int_A f(x - y)^p \, d\nu(y) \, dH(x) \\
& = & \nu(A)^{p - 1} \, \int_A \int_A f(x - y)^p \, dH(x) \, d\nu(y) \nonumber \\
& = & \nu(A)^p \, \Big(\int_A f(x)^p \, dH(x)\Big), \nonumber
\end{eqnarray}
by Fubini's theorem and translation-invariance of Haar measure.  If $f
\in L^p(A)$, then it follows that $f * \nu \in L^p(A)$ too, and in
particular that $f * \nu(x) < \infty$ for almost every $x \in A$ with
respect to Haar measure.

        If $f$ is a complex-valued function in $L^p(A)$, $1 < p < \infty$,
and if $\nu$ is a complex regular Borel measure on $A$, then we can
apply the previous argument to $|f|$ and $|\nu|$, to get that 
$(|f| * |\nu|)(x) < \infty$ for almost every $x \in A$ with respect to $H$.
Thus $(f * \nu)(x)$ can be defined for almost every $x \in A$ as before,
and satisfies (\ref{|(f * nu)(x)| le int_A |f(x - y)| d|nu|(y)}).
The analogue of (\ref{int_A ((f * nu)(x))^p dH(x) le ... = nu(A)^p
  (int_A f(x)^p dH(x))}) for $|f|$ and $|\nu|$ implies that
$f * \nu \in L^p(A)$, with
\begin{equation}
\label{(int_A |(f * nu)(x)|^p dH(x))^{1/p} le ...}
        \Big(\int_A |(f * \nu)(x)|^p \, dH(x)\Big)^{1/p} 
         \le |\nu|(A) \, \Big(\int_A |f(x)|^p \, dH(x)\Big)^{1/p}.
\end{equation}
If $f(x) = 0$ for almost every $x \in A$ with respect to $H$, then it
follows that $(f * \nu)(x) = 0$ almost everywhere on $A$ too, as in
the $p = 1$ case.

        Let $f$ be a bounded complex-valued Borel measurable function
on $A$, and let $\nu$ be a complex regular Borel measure on $A$ again.
Under these conditions, $(f * \nu)(x)$ is defined for every $x \in A$,
and satisfies
\begin{equation}
\label{|(f * nu)(x)| le ... le (sup_{y in A} |f(y)|) |nu|(A)}
        |(f * \nu)(x)| \le \int_A |f(x - y)| \, d|\nu|(y)
                        \le \Big(\sup_{y \in A} |f(y)|\Big) \, |\nu|(A)
\end{equation}
for every $x \in A$.  If $f$ is bounded and continuous on $A$, and if
$\nu$ has compact support in $A$, then it is easy to see that $f *
\nu$ is also continuous on $A$, using the uniform continuity of $f$ on
compact subsets of $A$.  Actually, if $\nu$ has compact support in
$A$, then $f * \nu$ is defined and continuous on $A$ for every
continuous function $f$ on $A$, by the same argument.  If $f$ is
bounded and continuous on $A$ and $\nu$ does not have compact support,
then one can still check that $f * \nu$ is continuous on $A$, using
the regularity of $\nu$ to approximate the relevant integrals by ones
over compact sets.  Similarly, if $f \in C_0(A)$, then $f * \nu \in
C_0(A)$ as well.  If $f$ is bounded and uniformly continuous on $A$,
then $f * \nu$ is uniformly continuous on $A$ too, by a more direct
version of the same type of argument.

        If $\phi \in \widehat{A}$, then $\phi$ is bounded and continuous
on $A$ in particular, so that $\phi * \nu$ is defined as in the
previous paragraph.  In this case, we have that
\begin{eqnarray}
\label{(phi * nu)(x) = ... = widehat{nu}(phi) phi(x)}
        (\phi * \nu)(x) & = & \int_A \phi(x - y) \, d\nu(y)
                          = \int_A \phi(x) \, \phi(-y) \, d\nu(y) \\
                        & = & \phi(x) \int_A \overline{\phi(y)} \, d\nu(y)
                          = \widehat{\nu}(\phi) \, \phi(x) \nonumber
\end{eqnarray}
for every $x \in A$.

\section{Density in $C_0(\widehat{A})$}
\label{density in C_0(widehat{A})}

        Let $a \in A$ and $\phi \in \widehat{A}$ be given, and let
$f$ be a nonnegative real-valued integrable function on $A$ supported
in a small neighborhood $U$ of $-a$ in $A$ such that
\begin{equation}
\label{int_A f(x) dH(x) = 1}
        \int_A f(x) \, dH(x) = 1.
\end{equation}
If $U$ is a sufficiently small neighborhood of $-a$ in $A$, then it is
easy to see that $\widehat{f}(\phi)$ is as close to $\phi(a)$ as we
want, because $\phi$ is continuous at $a$.  Now let $\mathcal{E}$ be
the collection of functions on $\widehat{A}$ of the form $\widehat{f}$
for some complex-valued integrable function $f$ on $A$.  This is a
subalgebra of the algebra $C_0(\widehat{A})$ of all continuous
complex-valued functions that vanish at infinity on $\widehat{A}$,
since the Fourier transform of the convolution of two integrable
functions is the same as the product of the corresponding Fourier
transforms.  We have also seen that the complex conjugate of every
element of $\mathcal{E}$ is an element of $\mathcal{E}$ as well.  If
$\phi \in \widehat{A}$, then there are $f \in L^1(A)$ such that
$\widehat{f}(\phi)$ approximates $\phi(0) = 1$, as before, and hence
$\widehat{f}(\phi) \ne 0$.  Similarly, if $\phi, \psi \in \widehat{A}$
and $\phi \ne \psi$, then $\phi(a) \ne \psi(a)$ for some $a \in A$,
and there are $f \in L^1(A)$ such that $\widehat{f}(\phi)$,
$\widehat{f}(\psi)$ approximate $\phi(a)$, $\psi(a)$, respectively, as
before.  In particular,
\begin{equation}
\label{widehat{f}(phi) ne widehat{f}(psi)}
        \widehat{f}(\phi) \ne \widehat{f}(\psi)
\end{equation}
for some $f \in L^1(A)$, which implies that $\mathcal{E}$ separates
points in $\widehat{A}$.

        A version of the Stone--Weierstrass theorem implies that
$\mathcal{E}$ is dense in $C_0(\widehat{A})$ with respect to the
supremum norm.  More precisely, if $\widehat{A}$ is compact, then
$C_0(\widehat{A}) = C(A)$, and one can use the standard version of the
Stone--Weierstrass theorem.  Otherwise, let $\widehat{A}_1$ be the
one-point compactification of $\widehat{A}$, so that elements of
$C_0(\widehat{A})$ can be identified with continuous complex-valued
functions on $\widehat{A}_1$ that vanish at the point at infinity.
Also let $\mathcal{E}_1$ be the collection of functions on
$\widehat{A}_1$ that can be expreseed as the sum of a constant
function on $\widehat{A}_1$ and a function that corresponds to an
element of $\mathcal{E}$.  One can check that this is a subalgebra of
the algebra $C(\widehat{A}_1)$ of all continuous complex-valued
functions on $\widehat{A}_1$ that separates points on $\widehat{A}_1$
and is invariant under complex conjugation.  Thus the standard version
of the Stone--Weierstrass theorem implies that $\mathcal{E}_1$ is
dense in $C(\widehat{A}_1)$ with respect to the supremum norm.  In
particular, continuous functions on $\widehat{A}_1$ that vanish at the
point at infinity can be approximated by elements of $\mathcal{E}_1$
uniformly on $\widehat{A}_1$, in which case the approximations should
also be small at the point at infinity.  Using this, it is easy to see
that $\mathcal{E}$ is dense in $C_0(\widehat{A})$ with respect to the
supremum norm, as desired.

        Of course, if $A$ is discrete, then $\widehat{A}$ is compact, 
and $C_0(\widehat{A}) = C(\widehat{A})$.  In this case, $\mathcal{E}$
contains the constant functions on $\widehat{A}$, which are the
Fourier transforms of functions on $A$ supported at $0$.  Similarly,
for each $a \in A$,
\begin{equation}
\label{Psi_a(phi) = phi(a), 2}
        \Psi_a(\phi) = \phi(a)
\end{equation}
is the Fourier transform of the function on $A$ equal to $1$ at $-a$
and to $0$ at every other point in $A$, and hence is in $\mathcal{E}$.
Let $\mathcal{E}_0$ be the collection of functions on $\widehat{A}$
which can be expressed as finite linear combinations of functions of
the form (\ref{Psi_a(phi) = phi(a), 2}) for some $a \in A$.  This is
also a subalgebra of $C(\widehat{A})$ that contains the constant
functions, separates points in $\widehat{A}$, and which is invariant
under complex conjugation.  The density of $\mathcal{E}_0$ in
$C(\widehat{A})$ implicitly came up before, in Section \ref{compact
  groups}.  Equivalently, $\mathcal{E}_0$ consists of the Fourier
transforms of functions on $A$ with finite support, which are dense in
$L^1(A)$ when $A$ is discrete.

        If $A$ is compact, then $\widehat{A}$ is discrete, and every
element of $\widehat{A}$ may be considered as an integrable function
on $A$.  The Fourier transform of $\phi \in \widehat{A}$ satisfies
$\widehat{\phi}(\psi) = 0$ for every $\psi \in \widehat{A}$ such that
$\phi \ne \psi$, because of the usual orthogonality properties of
characters on $A$.  If Haar measure $H$ on $A$ is normalized so that
$H(A) = 1$, then we also have that $\widehat{\phi}(\phi) = 1$.  Thus
the Fourier transform maps finite linear combinations of characters on
$A$ to functions with finite support on $\widehat{A}$, and every
function on $\widehat{A}$ with finite support is of this form.  In
this situation, $C_0(\widehat{A})$ is the same as the space
$c_0(\widehat{A})$ of all complex-valued functions on $\widehat{A}$
that vanish at infinity, and functions on $\widehat{A}$ with finite
support are clearly dense in $c_0(\widehat{A})$ with respect to the
supremum norm.

\chapter{Banach algebras}
\label{banach algebras}

\section{Definitions and examples}
\label{definitions, examples}

        Let $\mathcal{A}$ be a vector space over the real or complex numbers, 
and suppose that there is a bilinear mapping that sends $(x, y) \in
\mathcal{A} \times \mathcal{A}$ to an element $x \, y$ of
$\mathcal{A}$.  This means that $x \mapsto x \, y$ is a linear mapping
from $\mathcal{A}$ into itself for each $y \in \mathcal{A}$, and that
$y \mapsto x \, y$ is a linear mapping from $\mathcal{A}$ into itself
for each $x \in \mathcal{A}$.  If this bilinear mapping satisfies the
associative law
\begin{equation}
\label{(x y) z = x (y z)}
        (x \, y) \, z = x \, (y \, z)
\end{equation}
for every $x, y, z \in \mathcal{A}$, then $\mathcal{A}$ is said to be an
\emph{algebra}\index{algebras} over ${\bf R}$ or ${\bf C}$, as appropriate.
If we also have that
\begin{equation}
\label{x y = y x}
        x \, y = y \, x
\end{equation}
for every $x, y \in \mathcal{A}$, then $\mathcal{A}$ is said to be a
\emph{commutative algebra}.\index{commutative algebras} Suppose that
$\mathcal{A}$ is an algebra over ${\bf R}$ or ${\bf C}$ which is
equipped with a norm $\|x\|$ such that
\begin{equation}
\label{||x y|| le ||x|| ||y||}
        \|x \, y\| \le \|x\| \, \|y\|
\end{equation}
for every $x, y \in \mathcal{A}$.  It is easy to see that this implies
that multiplication on $\mathcal{A}$ is continuous as a mapping from
$\mathcal{A} \times \mathcal{A}$ into $\mathcal{A}$.  If $\mathcal{A}$
is complete with respect to the metric associated to the norm, then
$\mathcal{A}$ is said to be a \emph{Banach algebra}.\index{Banach
  algebras} In many situations, there may be a nonzero element $e$ of
$\mathcal{A}$ which is the multiplicative identity element, which
means that
\begin{equation}
\label{e x = x e = x}
        e \, x = x \, e = x
\end{equation}
for every $x \in \mathcal{A}$.  Note that $e$ is unique when it
exists, and that (\ref{||x y|| le ||x|| ||y||}) and (\ref{e x = x e =
  x}) imply that $\|e\| \ge 1$.  It is customary to ask that $\|e\| =
1$ under these conditions.

        If $X$ is a nonempty topological space, then the space $C_b(X)$ 
of bounded continuous real or complex-valued functions on $X$ is a
commutative Banach algebra with respect to pointwise addition and
multiplication of functions, and using the supremum norm.  In
particular, the constant function ${\bf 1}_X$ equal to $1$ at every
point in $X$ is the multiplicative identity element in $C_b(X)$.  If
$X$ is a locally compact Hausdorff topological space which is not
compact, then the space $C_0(X)$ of continuous real or complex-valued
functions on $X$ that vanish at infinity is a Banach algebra without a
multiplicative identity element.  Of course, if $X$ is equipped with
the discrete topology, then $C_b(X)$ is the same as $\ell^\infty(X)$,
and $C_0(X)$ is the same as $c_0(X)$.

        Let $V$ be a vector space over the real or complex numbers, 
and equipped with a norm $\|v\|_V$.  The space $\mathcal{BL}(V) =
\mathcal{BL}(V, V)$\index{BL(V)@$\mathcal{BL}(V)$} of bounded linear
mappings from $V$ into itself is an algebra with composition of linear
mappings as multiplication.  As in (\ref{||T_2 circ T_1||_{op, 13} le
  ||T_1||_{op, 12} ||T_2||_{op, 23}}) in Section \ref{bounded linear
  mappings}, the operator norm $\|T\|_{op}$ on $\mathcal{BL}(V)$
corresponding to the norm $\|v\|_V$ on $V$ satisfies (\ref{||x y|| le
  ||x|| ||y||}).  If $V$ is complete with respect to $\|v\|_V$, then
$\mathcal{BL}(V)$ is complete with respect to the operator norm, as
mentioned in Section \ref{bounded linear mappings}, and hence
$\mathcal{BL}(V)$ is a Banach algebra.  The identity mapping $I = I_V$
on $V$ is the multiplicative identity element in $\mathcal{BL}(V)$,
which is nonzero as long as $V \ne \{0\}$, in which case $\|I\|_{op} =
1$.

        Let $A$ be a locally compact commutative topological group, with
a suitable choice of Haar measure.  The space $L^1(A)$ of integrable
complex-valued functions on $A$ is a commutative Banach algebra, with
convolution as multiplication.  Similarly, the space of regular
complex Borel measures on $A$ is a commutative Banach algebra with
respect to convolution.  The Fourier transform defines a homomorphism
from $L^1(A)$ with convolution as multiplication into
$C_0(\widehat{A})$, where $\widehat{A}$ is the dual group of
continuous homomorphisms from $A$ into ${\bf T}$, as usual.  The
Fourier transform also defines a homomorphism from the algebra of
complex regular Borel measures on $A$ with respect to convolution into
$C_b(\widehat{A})$.

\section{Invertibility}
\label{invertibility}

        Let $\mathcal{A}$ be a real or complex algebra with a nonzero
multiplicative identity element $e$.  An element $x$ of $\mathcal{A}$
is said to be \emph{invertible}\index{invertible elements of an
  algebra} if there is an element $x^{-1}$ of $\mathcal{A}$ such that
\begin{equation}
\label{x^{-1} x = x x^{-1} = e.}
        x^{-1} \, x = x \, x^{-1} = e.
\end{equation}
It is easy to see that $x^{-1}$ is unique when it exists, in which
case $x^{-1}$ is also invertible, with $(x^{-1})^{-1} = x$.  If $x$
and $y$ are invertible elements of $\mathcal{A}$, then their product
$x \, y$ is invertible in $\mathcal{A}$ too, and the inverse is given by
\begin{equation}
\label{(x y)^{-1} = y^{-1} x^{-1}}
        (x \, y)^{-1} = y^{-1} \, x^{-1}.
\end{equation}
Thus the invertible elements of $\mathcal{A}$ form a group with
respect to multiplication.

        If $w \in \mathcal{A}$ is invertible and $w$ commutes 
with $z \in \mathcal{A}$, so that $w \, z = z \, w$, then $z$ commutes
with $w^{-1}$ as well.  In particular, if $x$ and $y$ are any two
commuting elements of $\mathcal{A}$ whose product $x \, y$ is
invertible in $\mathcal{A}$, then $x \, y$ commutes with both $x$ and
$y$, and hence $(x \, y)^{-1}$ commutes with both $x$ and $y$.  Under
these conditions, it follows that $x$ and $y$ are both invertible in
$\mathcal{A}$, with
\begin{equation}
\label{x^{-1} = (x y)^{-1} y, y^{-1} = (x y)^{-1} x}
        x^{-1} = (x \, y)^{-1} \, y, \quad y^{-1} = (x \, y)^{-1} \, x.
\end{equation}
Note that this does not necessarily work when $x$ and $y$ do not commute.

        Suppose now that $\mathcal{A}$ is equipped with a norm $\|\cdot \|$ 
that satisfies (\ref{||x y|| le ||x|| ||y||}) for every $x, y \in
\mathcal{A}$.  If $x, y \in \mathcal{A}$ are invertible, then
\begin{equation}
\label{x^{-1} - y^{-1} = ... = x^{-1} (y - x) y^{-1}}
        x^{-1} - y^{-1} = x^{-1} \, y \, y^{-1} - x^{-1} \, x \, y^{-1} 
                        = x^{-1} \, (y - x) \, y^{-1},
\end{equation}
and hence
\begin{equation}
\label{||x^{-1} - y^{-1}|| le ||x^{-1}|| ||x - y|| ||y^{-1}||}
        \|x^{-1} - y^{-1}\| \le \|x^{-1}\| \, \|x - y\| \, \|y^{-1}\|.
\end{equation}
In particular,
\begin{equation}
\label{||y^{-1}|| le ||x^{-1}|| ||x - y|| ||y^{-1}|| + ||x^{-1}||}
        \|y^{-1}\| \le \|x^{-1}\| \, \|x - y\| \, \|y^{-1}\| + \|x^{-1}\|.
\end{equation}
If $\|x^{-1}\| \, \|x - y\| < 1$, then
\begin{equation}
\label{(1 - ||x^{-1}|| ||x - y||) ||y^{-1}|| le ||x^{-1}||}
        (1 - \|x^{-1}\| \, \|x - y\|) \, \|y^{-1}\| \le \|x^{-1}\|
\end{equation}
implies that
\begin{equation}
\label{||y^{-1}|| le frac{||x^{-1}||}{1 - ||x^{-1}|| ||x - y||}}
        \|y^{-1}\| \le \frac{\|x^{-1}\|}{1 - \|x^{-1}\| \, \|x - y\|}.
\end{equation}
Combining this with (\ref{||x^{-1} - y^{-1}|| le ||x^{-1}|| ||x - y||
  ||y^{-1}||}), we get that
\begin{equation}
\label{||x^{-1} - y^{-1}|| le ...}
        \|x^{-1} - y^{-1}\| \le
           \frac{\|x^{-1}\|^2 \, \|x - y\|}{1 - \|x^{-1}\| \, \|x - y\|}
\end{equation}
when $\|x^{-1}\| \, \|x - y\| < 1$.  If $\|x^{-1}\| \, \|x - y\| \le
1/2$, for instance, then it follows that
\begin{equation}
\label{||x^{-1} - y^{-1}|| le 2 ||x^{-1}||^2 ||x - y||}
        \|x^{-1} - y^{-1}\| \le 2 \, \|x^{-1}\|^2 \, \|x - y\|.
\end{equation}
This shows that $x \mapsto x^{-1}$ is a continuous mapping on the set
of invertible elements of $\mathcal{A}$ with respect to the metric
associated to the norm, so that the group of invertible elements of
$\mathcal{A}$ is actually a topological group with respect to the
topology induced by this metric.

        If $a$ is any element of $\mathcal{A}$ and $j$ is a positive integer,
then we let $a^j$ be the product of $j$ $a$'s, as usual, so that $a^1
= a$ and $a^j = a \, a^{j - 1}$ when $j \ge 2$.  It is customary to
put $a^0 = e$ for every $a \in \mathcal{A}$ when there is a nonzero
multiplicative identity element $e$ in $\mathcal{A}$.  Because of
(\ref{||x y|| le ||x|| ||y||}), we have that
\begin{equation}
\label{||a^j|| le ||a||^j}
        \|a^j\| \le \|a\|^j
\end{equation}
for each $j \ge 1$, and we ask that $\|e\| = 1$, as in the previous
section.  This implies that $\sum_{j = 0}^\infty a^j$ converges
absolutely in $\mathcal{A}$ when $\|a\| < 1$, and hence that $\sum_{j
  = 0}^\infty a^j$ converges in $\mathcal{A}$ when $\mathcal{A}$ is a
Banach algebra.  Note that
\begin{equation}
\label{||sum_{j = 0}^infty a^j|| le ... = frac{1}{1 - ||a||}}
        \biggl\|\sum_{j = 0}^\infty a^j\biggr\| \le \sum_{j = 0}^\infty \|a^j\| 
                           \le \sum_{j = 0}^\infty \|a\|^j = \frac{1}{1 - \|a\|}
\end{equation}
under these conditions.  A standard computation shows that
\begin{equation}
\label{(e - a) (sum_{j = 0}^n a^j) = ... = e - a^{n + 1}}
 (e - a) \, \Big(\sum_{j = 0}^n a^j\Big) = \Big(\sum_{j = 0}^n a^j\Big) \, (e - a)
                                                = e - a^{n + 1}
\end{equation}
for each nonnegative integer $n$, and of course $a^{n + 1} \to 0$ in
$\mathcal{A}$ as $n \to \infty$ when $\|a\| < 1$, by (\ref{||a^j|| le
  ||a||^j}).  Thus
\begin{equation}
\label{(e - a) (sum_{j = 0}^infty a^j) = (sum_{j = 0}^infty a^j) (e - a) = e}
        (e - a) \, \Big(\sum_{j = 0}^\infty a^j\Big)
                  = \Big(\sum_{j = 0}^\infty a^j\Big) \, (e - a) = e,
\end{equation}
so that $e - a$ is invertible in $\mathcal{A}$ when $\|a\| < 1$ and
$\mathcal{A}$ is a Banach algebra, with inverse equal to $\sum_{j =
  0}^\infty a^j$.  If $x$ is an invertible element of $\mathcal{A}$,
and $y \in \mathcal{A}$ satisfies $\|x^{-1}\| \, \|x - y\| < 1$, then
we can apply the previous argument to $a = x^{-1} \, (x - y)$, to get
that $e - x^{-1} \, (x - y)$ is invertible in $\mathcal{A}$ when
$\mathcal{A}$ is a Banach algebra.  This implies that
\begin{equation}
\label{y = x - (x - y) = x (e - x^{-1} (x - y))}
        y = x - (x - y) = x \, (e - x^{-1} \, (x - y))
\end{equation}
is invertible in $\mathcal{A}$ under these conditions, so that the set
of invertible elements of $\mathcal{A}$ is an open subset of
$\mathcal{A}$ when $\mathcal{A}$ is a Banach algebra.

\section{Spectrum and spectral radius}
\label{spectrum, spectral radius}

        Let $\mathcal{A}$ be a Banach algebra over the real or complex
numbers with a nonzero multiplicative identity element $e$.  The
\emph{spectrum}\index{spectrum} of an element $x$ of $\mathcal{A}$ is
defined to be the set $\sigma(x)$\index{sigma(x)@$\sigma(x)$} of
$\lambda \in {\bf R}$ or ${\bf C}$, as appropriate, such that $x -
\lambda \, e$ is not invertible in $\mathcal{A}$.  If $|\lambda| >
\|x\|$, so that $\|\lambda^{-1} \, x\| = |\lambda|^{-1} \, \|x\| < 1$,
then $e - \lambda^{-1} \, x$ is invertible in $\mathcal{A}$, as in the
previous section, and hence $\lambda \not\in \sigma(x)$.
Equivalently,
\begin{equation}
\label{|lambda| le ||x||}
        |\lambda| \le \|x\|
\end{equation}
for every $\lambda \in \sigma(x)$.  We also know from the previous
section that that the set of invertible elements of $\mathcal{A}$ is
an open set, which implies that the set of $\lambda \in {\bf R}$ or
${\bf C}$ such that $\lambda \not\in \sigma(x)$ is an open set, so
that $\sigma(x)$ is a closed set in ${\bf R}$ or ${\bf C}$, as
appropriate.

        If $\mathcal{A}$ is a complex Banach algebra, then a famous
theorem states that $\sigma(x)$ is nonempty for every $x \in
\mathcal{A}$.  To see this, suppose for the sake of a contradiction
that $x - \lambda \, e$ is invertible for every $\lambda \in {\bf C}$.
The main idea is that $(x - \lambda \, e)^{-1}$ should be holomorphic
as a function of $\lambda$ on the complex plane with values in
$\mathcal{A}$.  In particular, one can develop the theory of
holomorphic functions with values in a complex Banach space, as well
as other complex topological vector spaces.  To avoid technicalities
about holomorphic vector-valued functions, one can use bounded linear
functionals on $\mathcal{A}$ to reduce to the case of ordinary
complex-valued holomorphic functions.  More precisely, one can show
that $\phi((x - \lambda \, e)^{-1})$ is a holomorphic complex-valued
function of $\lambda$ on the complex plane for each bounded linear
functional $\phi$ on $\mathcal{A}$.  Because
\begin{equation}
\label{(x - lambda e)^{-1} = - lambda^{-1} (e - lambda^{-1} x)^{-1} to 0}
 (x - \lambda \, e)^{-1} = - \lambda^{-1} \, (e - \lambda^{-1} \, x)^{-1} \to 0
\end{equation}
as $|\lambda| \to \infty$, we get that $\phi((x - \lambda \, e)^{-1})
\to 0$ as $|\lambda| \to 0$ for every bounded linear functional $\phi$
on $\mathcal{A}$.  This implies that $\phi((x - \lambda \, e)^{-1}) =
0$ for every $\lambda \in {\bf C}$ and bounded linear functional
$\phi$ on $\mathcal{A}$, by well-known theorems in complex analysis.
Using the Hahn-Banach theorem, it follows that $(x - \lambda \,
e)^{-1} = 0$ for every $\lambda \in {\bf C}$, which is a
contradiction, since the inverse of any invertible element of
$\mathcal{A}$ is nonzero.

        Let $\mathcal{A}$ be a real or complex Banach algebra again.
If $a \in \mathcal{A}$ satisfies $\|a^l\| < 1$ for some positive
integer $l$, then $e - a$ is invertible in $\mathcal{A}$.  One way to
see this is to observe that $\sum_{j = 0}^\infty a^j$ also converges
in $\mathcal{A}$ in this situation, and that the sum is the inverse of
$e - a$.  Alternatively, one can use the previous argument to get that
$e - a^l$ is invertible in $\mathcal{A}$, and then apply (\ref{(e - a)
  (sum_{j = 0}^n a^j) = ... = e - a^{n + 1}}) with $n = l - 1$ to obtain
\begin{equation}
\label{(e - a) (sum_{j = 0}^{l - 1} a^j) = e - a^l}
        (e - a) \, \Big(\sum_{j = 0}^{l - 1} a^j\Big) = e - a^l.
\end{equation}
This implies that $e - a$ is invertible, because $e - a^l$ is
invertible and $e - a$ commutes with $\sum_{j = 0}^{l - 1} a^j$, as in
the preceding section.

        Let $x$ be any element of $\mathcal{A}$, and put
\begin{equation}
\label{r(x) = inf_{l ge 1} ||x^l||^{1/l}}
        r(x) = \inf_{l \ge 1} \|x^l\|^{1/l},
\end{equation}
where more precisely the infimum is taken over all positive integers
$l$.  If $\lambda$ is a real or complex number, as appropriate, such
that $|\lambda| > r(x)$, then $|\lambda| > \|x^l\|^{1/l}$ for some
positive integer $l$, and hence $|\lambda^l| > \|x^l\|$.
Equivalently, $\|(\lambda^{-1} \, x)^l\| < 1$, which implies that $e -
\lambda^{-1} \, x$ is invertible in $\mathcal{A}$, as in the previous
paragraph.  Thus $x - \lambda \, e$ is invertible when $|\lambda| >
r(x)$, which means that
\begin{equation}
\label{|lambda| le r(x)}
        |\lambda| \le r(x)
\end{equation}
for every $\lambda \in \sigma(x)$.

        Another famous theorem states that
\begin{equation}
\label{r(x) = sup {|lambda| : lambda in sigma(x)}}
        r(x) = \sup \{|\lambda| : \lambda \in \sigma(x)\}
\end{equation}
for every $x \in \mathcal{A}$ when $\mathcal{A}$ is a complex Banach
algebra.  To see this, let $\rho > 0$ be defined by
\begin{equation}
\label{sup {|lambda| : lambda in sigma(x)} = 1/rho}
        \sup \{|\lambda| : \lambda \in \sigma(x)\} = 1/\rho,
\end{equation}
so that $\rho = +\infty$ when $\sigma(x) = \{0\}$.  Thus $x - \lambda
\, e$ is invertible in $\mathcal{A}$ when $|\lambda| > 1/\rho$, which
is the same as saying that $e - \lambda^{-1} \, x$ is invertible when
$|\lambda^{-1}| < \rho$.  If we put $\zeta = \lambda^{-1}$, then we
get that $e - \zeta \, x$ is invertible when $\zeta \in {\bf C}$
satisfies $0 < |\zeta| < \rho$, and of course this holds trivially
when $\zeta = 0$.

        As before, the main idea is to look at $(e - \zeta \, x)^{-1}$
as a holomorphic $\mathcal{A}$-valued function on the open disk where
$|\zeta| < \rho$.  In particular, we know from the previous section
that $(e - \zeta \, x)^{-1}$ is a continuous function on this open
disk, which implies that it is bounded on the compact sub-disk where
$|\zeta| \le \rho_1$ for any $\rho_1 < \rho$.  We also know that
\begin{equation}
\label{(e - zeta x)^{-1} = sum_{j = 0}^infty zeta^j x^j}
        (e - \zeta \, x)^{-1} = \sum_{j = 0}^\infty \zeta^j \, x^j
\end{equation}
when $|\zeta| \, \|x\| < 1$, which is the power series expansion for
$(e - \zeta \, x)^{-1}$ as a function of $\zeta$ at $0$.  If this were
an ordinary complex-valued holomorphic function, then we could
represent the coefficients of this power series expansion in terms of
suitable integrals of the function on the circle where $|\zeta| =
\rho_1$.  This type of result can be extended to vector-valued
functions, or one can consider $\phi((e - \zeta \, x)^{-1})$ as a
complex-valued holomorphic function of $\zeta$ on the open disk where
$|\zeta| < \rho$ for each bounded linear functional $\phi$ on
$\mathcal{A}$.  Either way, one can use integral expressions for the
power series coefficients to show that for each $\rho_1 < \rho$ there
is a positive real number $C(\rho_1)$ such that
\begin{equation}
\label{rho_1^j ||x^j|| le C(rho_1)}
        \rho_1^j \, \|x^j\| \le C(\rho_1)
\end{equation}
for every positive integer $j$.  This also uses the Hahn--Banach
theorem, if one considers $\phi((e - \zeta \, x)^{-1})$ as a
holomorphic function of $\zeta$ for each bounded linear functional
$\phi$ on $\mathcal{A}$.  Equivalently, we get that
\begin{equation}
\label{||x^j||^{1/j} le rho_1^{-1} C(rho_1)^{1/j}}
        \|x^j\|^{1/j} \le \rho_1^{-1} \, C(\rho_1)^{1/j}
\end{equation}
for each $j \ge 1$ and $\rho_1 \in (0, \rho)$, which implies that
\begin{equation}
\label{r(x) le rho_1^{-1}}
        r(x) \le \rho_1^{-1}
\end{equation}
for each $\rho_1 \in (0, \rho)$, because $C^{1/j} \to 1$ as $j \to
\infty$ for every $C > 0$.  Taking $\rho_1 \to \rho$, we get that
\begin{equation}
\label{r(x) le rho^{-1}}
        r(x) \le \rho^{-1},
\end{equation}
which is exactly what we wanted, since the opposite inequality follows
from (\ref{|lambda| le r(x)}).

        Note that $\|x^n\|^{1/n} \to r(x)$ as $n \to \infty$ for
every $x \in \mathcal{A}$.  This can be derived from the previous
argument in the complex case, but we can also verify it more directly
from the definitions in both the real and complex cases.  If $l$ is
any positive integer, then any other positive integer $n$ can be
expressed as $j \, l + k$, where $j$ and $k$ are nonnegative integers
and $k < l$.  Thus we get that
\begin{equation}
\label{||x^n|| = ||x^{j l + k}|| le ||x^l||^j ||x||^k}
        \|x^n\| = \|x^{j \, l + k}\| \le \|x^l\|^j \, \|x\|^k,
\end{equation}
which implies that
\begin{equation}
\label{||x^n||^{1/n} le ... = (||x^l||^{1/l})^{1 - (k/n)} ||x||^{k/n}}
        \|x^n\|^{1/n} \le (\|x^l\|^{1/l})^{jl/n} \, \|x\|^{k/n} 
                         = (\|x^l\|^{1/l})^{1 - (k/n)} \, \|x\|^{k/n}.
\end{equation}
Because $0 \le k < l$, we get that
\begin{equation}
\label{limsup_{n to infty} ||x^n||^{1/n} le ||x^l||^{1/l}}
        \limsup_{n \to \infty} \|x^n\|^{1/n} \le \|x^l\|^{1/l}
\end{equation}
for each positive integer $l$.  Now we can take the infimum over $l$,
to obtain that
\begin{equation}
\label{limsup_{n to infty} ||x^n||^{1/n} le r(x)}
        \limsup_{n \to \infty} \|x^n\|^{1/n} \le r(x).
\end{equation}
This implies that $\|x^n\|^{1/n} \to r(x)$ as $n \to \infty$, as
desired, since $\|x^n\|^{1/n} \ge r(x)$ for each $n \ge 1$ by
definition of $r(x)$.

\section{Maximal ideals and homomorphisms}
\label{maximal ideals, homomorphisms}

        Let $\mathcal{A}$ be a commutative ideal over the real or complex
numbers.  As usual, a linear subspace $\mathcal{I}$ of $\mathcal{A}$
is said to be an \emph{ideal}\index{ideals} if for every $x \in
\mathcal{I}$ and $a \in \mathcal{A}$ we have that $a \, x \in
\mathcal{I}$.  Thus $\mathcal{A}$ is automatically an ideal in itself,
and an ideal $\mathcal{I}$ in $\mathcal{A}$ is said to be
\emph{proper} if $\mathcal{I} \ne \mathcal{A}$.  A proper ideal
$\mathcal{I}$ in $\mathcal{A}$ is said to be \emph{maximal}\index{maximal
ideals} if $\mathcal{A}$ and $\mathcal{I}$ are the only ideals in
$\mathcal{A}$ that contain $\mathcal{I}$.

        Suppose for the moment that $\mathcal{A}$ has a nonzero 
multiplicative identity element $e$.  Note that an ideal $\mathcal{I}$
in $\mathcal{A}$ is proper if and only if $e \not\in \mathcal{I}$.  It
is well known that every proper ideal in $\mathcal{A}$ is contained in
a maximal ideal, as one can show using Zorn's lemma or the Hausdorff
maximality principle.  The main point is that the union of a chain of
ideals in $\mathcal{A}$ is also an ideal in $\mathcal{A}$, and in fact
the union of a chain of proper ideals in $\mathcal{A}$ is a proper
ideal in $\mathcal{A}$, because it does not contain $e$ as an element.

        If $\mathcal{A}$ is a Banach algebra, then the closure of every
ideal in $\mathcal{A}$ is an ideal in $\mathcal{A}$ too.  If
$\mathcal{I}$ is a proper ideal in $\mathcal{A}$ and $\mathcal{A}$ has
a nonzero multiplicative identity element $e$, then $\mathcal{I}$ does
not contain any invertible elements of $\mathcal{A}$.  This implies
that $e - x \not\in \mathcal{I}$ for every $x \in \mathcal{A}$ with
$\|x\| < 1$, since $e - x$ is invertible in $\mathcal{A}$ when $\|x\| < 1$,
as in Section \ref{invertibility}.  It follows that $e \not\in
\overline{I}$, so that the closure $\overline{I}$ of $\mathcal{I}$ is
also a proper ideal in $\mathcal{A}$.  In particular, maximal ideals
are automatically closed in $\mathcal{A}$ under these conditions.

        Let $\phi$ be a homomorphism from $\mathcal{A}$ into the real
or complex numbers, as appropriate.  This means that $\phi$ is a linear
functional on $\mathcal{A}$ such that
\begin{equation}
\label{phi(x y) = phi(x) phi(y)}
        \phi(x \, y) = \phi(x) \, \phi(y)
\end{equation}
for every $x, y \in \mathcal{A}$.  Of course, the kernel of $\phi$ is
an ideal in $\mathcal{A}$.  If $\phi(x) \ne 0$ for some $x \in
\mathcal{A}$, then $\phi$ maps $\mathcal{A}$ onto ${\bf R}$ or ${\bf
  C}$, as appropriate, and the kernel of $\phi$ is a maximal ideal in
$\mathcal{A}$.  If $\mathcal{A}$ has a nonzero multiplicative identity
element $e$ and $\phi(x) \ne 0$ for some $x \in \mathcal{A}$, then
$\phi(e) = 1$, and $\phi(y) \ne 0$ for every invertible element $y$
of $\mathcal{A}$.

        Suppose that $\mathcal{A}$ is a Banach algebra with a nonzero
multiplicative identity element $e$, and that $\phi$ is a homomorphism
from $\mathcal{A}$ into ${\bf R}$ or ${\bf C}$, as appropriate, which
is not identically $0$.  If $x \in \mathcal{A}$ and $\|x\| < 1$, then
$e - x$ is invertible in $\mathcal{A}$, as in Section \ref{invertibility},
and hence $\phi(e - x) \ne 0$.  Thus $\phi(x) \ne 1$ when $\|x\| < 1$,
and one can apply this argument to $t \, x$ for each $t \in {\bf R}$
or ${\bf C}$, as appropriate, such that $|t| \le 1$, to get that
$|\phi(x)| < 1$ when $\|x\| < 1$.  This implies that
\begin{equation}
\label{|phi(x)| le ||x||}
        |\phi(x)| \le \|x\|
\end{equation}
for every $x \in \mathcal{A}$.

        There is a version of this argument that does not require a
multiplicative identity element.  If $x \in \mathcal{A}$ and $\|x\| < 1$,
then the series $\sum_{j = 1}^\infty x^j$ converges in $\mathcal{A}$
when $\mathcal{A}$ is a Banach algebra, for the same reasons as in
Section \ref{invertibility}.  Put $y = \sum_{j = 1}^\infty x^j$,
and observe that
\begin{equation}
\label{x y = sum_{j = 1}^infty x^{j + 1} = sum_{j = 2}^infty x^j = y - x}
 x \, y = \sum_{j = 1}^\infty x^{j + 1} = \sum_{j = 2}^\infty x^j = y - x.
\end{equation}
Let $\phi$ be a homomorphism from $\mathcal{A}$ into ${\bf R}$ or ${\bf
  C}$, as appropriate, so that
\begin{equation}
\label{phi(x) phi(y) = phi(x y) = phi(y) - phi(x)}
        \phi(x) \, \phi(y) = \phi(x \, y) = \phi(y) - \phi(x).
\end{equation}
If $\phi(x) = 1$, then this implies that $\phi(y) = \phi(y) - 1$,
which is impossible.  Thus we get that $\phi(x) \ne 1$ when $\|x\| <
1$.  As before, this implies that (\ref{|phi(x)| le ||x||}) holds for
every $x \in \mathcal{A}$.

        Let $\mathcal{B}$ be a complex Banach algebra with a nonzero
multiplicative identity element $e$.  If every nonzero element of
$\mathcal{B}$ is invertible, then a famous theorem states that $\mathcal{B}$
is isomorphic to the field of complex numbers.  To see this, let
$x \in \mathcal{A}$ be given, and let $\lambda \in {\bf C}$ be an
element of the spectrum of $x$, whose existence was discussed in the
previous section.  Thus $x - \lambda \, e$ is not invertible in
$\mathcal{A}$, which implies that $x = \lambda \, e$ in this situation.
This shows that every element of $\mathcal{B}$ can be expressed as a
complex multiple of $e$, as desired.

        Suppose now that $\mathcal{A}$ is a complex commutative Banach
algebra with nonzero multiplicative identity element $e$, and let
$\mathcal{I}$ be a proper closed ideal in $\mathcal{A}$.  The quotient
$\mathcal{A} / \mathcal{I}$ is a complex commutative algebra with a
nonzero multiplicative identity element in a natural way, and it can
be shown that $\mathcal{A} / \mathcal{I}$ is also a Banach algebra
with respect to the corresponding quotient norm.  If $\mathcal{I}$
is a maximal ideal in $\mathcal{A}$, then $\mathcal{A} / \mathcal{I}$
is a field, and hence is isomorphic to the field of complex numbers,
as in the previous paragraph.  It follows that every maximal ideal in
$\mathcal{A}$ is the kernel of a homomorphism from $\mathcal{A}$
onto ${\bf C}$ under these conditions.

\section{Homomorphisms on $L^1(A)$}
\label{homomorphisms on L^1(A)}

        Let $A$ be a locally compact commutative topological group,
and let $H$ be a Haar measure on $A$.  As before, the space $L^1(A)$
of complex-valued integrable functions on $A$ is a commutative Banach
algebra with respect to convolution and the $L^1$ norm.  If $\phi \in
\widehat{A}$, then the mapping from $f \in L^1(A)$ to
$\widehat{f}(\phi)$ is a homomorphism from $L^1(A)$ as a Banach
algebra into the field of complex numbers, as in (\ref{widehat{(f *
    g)}(phi) = widehat{f}(phi) widehat{g}(phi)}) in Section
\ref{convolution of integrable functions}.  Note that
$\widehat{f}(\phi) \ne 0$ when $f$ is a nonnegative real-valued
function supported on a sufficiently small neighborhood of $0$ in $A$
with integral equal to $1$, as in Section \ref{density in
  C_0(widehat{A})}, because $\phi(0) = 1$ and $\phi$ is continuous on
$A$.  Similarly, if $\phi$ and $\phi'$ are distinct elements of
$\widehat{A}$, then $\phi(a) \ne \phi'(a)$ for some $a \in A$, and it
is easy to see that $\widehat{f}(\phi) \ne \widehat{f}(\phi')$ when
$f$ is a nonnegative real-valued function on $A$ supported in a
sufficiently small neighborhood of $a$ with integral equal to $1$.
Thus we get a natural one-to-one mapping from $\widehat{A}$ into the
set of nonzero homomorphisms on $L^1(A)$.  We would like to show that
every nonzero homomorphism on $L^1(A)$ is of the form $f \mapsto
\widehat{f}(\phi)$ for some $\phi \in \widehat{A}$, so that this
mapping is a surjection.

        Suppose for the moment that $A$ is discrete, so that we can take
$H$ to be counting measure on $A$, and $L^1(A)$ is the same as $\ell^1(A)$.
For each $a \in A$, let $\delta_a(x)$ be the function on $A$ equal to
$1$ when $x = a$ and to $0$ otherwise.  In this situation, $L^1(A)$ has a
multiplicative identity element given by $\delta_0$, and
\begin{equation}
\label{delta_a * delta_b = delta_{a + b}}
        \delta_a * \delta_b = \delta_{a + b}
\end{equation}
for every $a, b \in A$.  Let $\Phi$ be a nonzero homomorphism on $L^1(A)$,
and put
\begin{equation}
\label{phi(a) = overline{Phi(a)}}
        \phi(a) = \overline{\Phi(a)}
\end{equation}
for each $a \in A$.  Thus $\phi(0) = 1$, because $\Phi(\delta_0) = 1$,
and
\begin{equation}
\label{phi(a + b) = phi(a) phi(b)}
        \phi(a + b) = \phi(a) \, \phi(b)
\end{equation}
for every $a, b \in A$.  This implies that $\phi(a) \ne 0$ for each $a
\in A$, by taking $b = -a$, and hence that $\phi$ is a homomorphism
from $A$ into the multiplicative group of nonzero complex numbers.  We
also have that
\begin{equation}
\label{|phi(a)| = |Phi(delta_a)| le 1}
        |\phi(a)| = |\Phi(\delta_a)| \le 1
\end{equation}
for every $a \in A$, by (\ref{|phi(x)| le ||x||}) in the previous
section, and because $\delta_a$ has $L^1$ norm equal to $1$ for each
$a \in A$.  Applying this to $-a$, we get that $|\phi(a)| = 1$ for
every $a \in A$, so that $\phi$ is a homomorphism from $A$ into the
unit circle ${\bf T}$.  By construction, $\Phi(f) = \widehat{f}(\phi)$
when $f = \delta_a$ for any $a \in A$, and this also holds when $f$
has finite support on $A$, by linearity.  It follows that $\Phi(f) =
\widehat{f}(\phi)$ for every $f \in L^1(A)$, because functions with
finite support on $A$ are dense in $L^1(A)$, and both $\Phi$ and $f
\mapsto \widehat{f}(\phi)$ are bounded linear functionals on $L^1(A)$.

        Suppose now that $A$ is $\sigma$-compact, and let $\Phi$ be a
nonzero homomorphism from $L^1(A)$ into ${\bf C}$.  Thus $\Phi$ is a
bounded linear functional on $L^1(A)$ with dual norm less than or
equal to $1$, by (\ref{|phi(x)| le ||x||}) in the previous section.
Because $A$ is $\sigma$-compact and hence is $\sigma$-finite with
respect to Haar measure, the Riesz representation theorem implies that
there is a bounded complex-valued Borel measurable function $\phi$
on $A$ with essential supremum norm less than or equal to $1$ such that
\begin{equation}
\label{Phi(f) = int_A f(x) overline{phi(x)} dH(x)}
        \Phi(f) = \int_A f(x) \, \overline{\phi(x)} \, dH(x)
\end{equation}
for every $f \in L^1(A)$.  If $f, g \in L^1(A)$, then
\begin{equation}
\label{Phi(f * g) = Phi(f) Phi(g) = Phi(f) int_A g(y) overline{phi(y)} dH(y)}
        \Phi(f * g) = \Phi(f) \, \Phi(g)
                    = \Phi(f) \, \int_A g(y) \, \overline{\phi(y)} \, dH(y)
\end{equation}
and
\begin{eqnarray}
\label{Phi(f * g) = int_A (f * g)(x) overline{phi(x)} dH(x) = ...}
 \Phi(f * g) & = & \int_A (f * g)(x) \, \overline{\phi(x)} \, dH(x) \\
 & = & \int_A \int_A f(x - y) \, g(y) \, \overline{\phi(x)} \, dH(y) \, dH(x).
                                                             \nonumber
\end{eqnarray}
Interchanging the order of integration and comparing the result with
(\ref{Phi(f * g) = Phi(f) Phi(g) = Phi(f) int_A g(y) overline{phi(y)} dH(y)}),
we get that
\begin{equation}
\label{Phi(f) overline{phi(y)} = ... = Phi(f_y)}
 \Phi(f) \, \overline{\phi(y)} = \int_A f(x - y) \, \overline{\phi(x)} \, dH(x)
                           = \Phi(f_y)
\end{equation}
for almost every $y \in A$, where $f_y(x) = f(x - y)$.

        By standard arguments, $y \mapsto f_y$ is a continuous mapping
from $A$ into $L^1(A)$ for each $f \in L^1(A)$.  More precisely, this
follows from uniform continuity when $f$ is a continuous function on $A$
with compact support, and otherwise one can get this for an arbitrary
integrable function $f$ on $A$ by approximating $f$ by a continuous
function with compact support with respect to the $L^1$ norm.  This
implies that $\Phi(f_y)$ is a continuous function of $y$ on $A$,
because $\Phi$ is a bounded linear functional on $L^1(A)$.  If we apply
this to $f \in L^1(A)$ such that $\Phi(f) \ne 0$, then we get that
$\phi(y)$ is equal to a continuous function of $y$ on $A$ almost
everywhere.  Thus we may as well suppose that $\phi$ is continuous
on $A$, and hence that (\ref{Phi(f) overline{phi(y)} = ... =
  Phi(f_y)}) holds for every $f \in L^1(A)$ and every $y \in A$.

        In particular, $\phi(0) = 1$, since (\ref{Phi(f) overline{phi(y)} 
= ... = Phi(f_y)}) holds for $y = 0$ and $f \in L^1(A)$ such that
$\Phi(f) \ne 0$.  If $y, z \in A$, then we can apply (\ref{Phi(f) 
overline{phi(y)} = ... = Phi(f_y)}) to $y + z$ instead of $y$ to get that
\begin{equation}
\label{Phi(f) overline{phi(y + z)} = Phi(f_{y + z})}
        \Phi(f) \, \overline{\phi(y + z)} = \Phi(f_{y + z})
\end{equation}
for every $f \in L^1(A)$.  Similarly, we can apply (\ref{Phi(f) 
overline{phi(y)} = ... = Phi(f_y)}) to $f_y$ instead of $f$ and $z$
in place of $y$ to get that
\begin{equation}
\label{Phi(f_{y + z}) = ... = Phi(f) overline{phi(y)} overline{phi(z)}}
        \Phi(f_{y + z}) = \Phi(f_y) \, \overline{\phi(z)}
         = \Phi(f) \, \overline{\phi(y)} \, \overline{\phi(z)},
\end{equation}
using (\ref{Phi(f) overline{phi(y)} = ... = Phi(f_y)}) again in the
second step.  Combining (\ref{Phi(f) overline{phi(y + z)} = Phi(f_{y +
    z})}) and (\ref{Phi(f_{y + z}) = ... = Phi(f) overline{phi(y)}
  overline{phi(z)}}), we get that
\begin{equation}
\label{phi(y + z) = phi(y) phi(z)}
        \phi(y + z) = \phi(y) \, \phi(z)
\end{equation}
for every $y, z \in A$, by taking $f \in L^1(A)$ such that $\Phi(f) \ne 0$.
As before, this implies that $\phi(y) \ne 0$ for every $y \in A$, by
taking $z = -y$.

        Remember that $|\phi(y)| \le 1$ for almost every $y \in A$.
By replacing $\phi$ with a continuous function on $A$, we get that
$|\phi(y)| \le 1$ for every $y \in A$.  This implies that $|\phi(y)| =
1$ for every $y \in A$, because $|\phi(-y)| \le 1$ and $\phi(-y) =
1/\phi(y)$.  Thus $\phi$ is a continuous homomorphism from $A$ into
${\bf T}$ and $\Phi(f) = \widehat{f}(\phi)$ for every $f \in L^1(A)$, 
as desired.

        Essentially the same argument can be used when $A$ is not
$\sigma$-compact, with some adjustments.  The main point is that there
is an open subgroup $B$ of $A$ which is $\sigma$-compact, as in
Section \ref{additional properties}.  Thus $A$ is partitioned into the
cosets of $B$, each of which is $\sigma$-finite with respect to Haar
measure.  If $a + B$ is any coset of $B$ in $A$, then we can apply
the Riesz representation theorem to the restriction of $\Phi$
to the subspace of $L^1(A)$ consisting of integrable functions
supported on $a + B$.  This leads to the same type of representation
for $\Phi$ as before, since only finitely or countable many cosets of
$B$ in $A$ are needed at each step.

\section{The weak$^*$ topology}
\label{weak^* topology}

        Let $V$ be a vector space over the real or complex numbers
with a norm $\|v\|$, and let $V^*$ be the corresponding dual space of
bounded linear functionals on $V$, with the dual norm $\|\lambda\|_*$.
Observe that
\begin{equation}
\label{N_v^*(lambda) = |lambda(v)|}
        N_v^*(\lambda) = |\lambda(v)|
\end{equation}
defines a seminorm on $V^*$ for each $v \in V$.  The collection of
these seminorms $N_v^*$, $v \in V$, is a nice collection of seminorms
on $V^*$ in the sense discussed in Section \ref{semimetrics}, and thus
defines a topology on $V^*$, known as the \emph{weak$^*$
  topology}.\index{weak^* topology@weak$^*$ topology} It is easy
to see that the closed unit ball
\begin{equation}
\label{B_{V^*} = {lambda in V^* : ||lambda||_* le 1}}
        B_{V^*} = \{\lambda \in V^* : \|\lambda\|_* \le 1\}
\end{equation}
in $V^*$ is a closed set in $V^*$ with respect to the weak$^*$ topology.
A famous theorem of Banach and Alaoglu\index{Banach--Alaoglu theorem}
states that $B_{V^*}$ is actually compact with respect to the weak$^*$
topology on $V^*$.  If $V$ is separable, then one can also show that
the topology on $B_{V^*}$ induced by the weak$^*$ topology on $V^*$
is metrizable.  The main point is that the topology on $B_{V^*}$
induced by the weak$^*$ topology on $V^*$ is the same as the
topology determined by the seminorms $N_v^*$ corresponding to a
dense set of $v \in V$, or even a set of $v \in V$ whose linear
span is dense in $V$ with respect to the norm $\|\cdot \|$.

        Now let $\mathcal{A}$ be a Banach algebra over the real or 
complex numbers, and let $\mathcal{H}(\mathcal{A})$ be the collection
of homomorphisms from $\mathcal{A}$ into ${\bf R}$ or ${\bf C}$, as
appropriate.  Thus $\mathcal{H}(\mathcal{A})$ is contained in the
closed unit ball $B_{\mathcal{A}^*}$ in the dual $\mathcal{A}^*$ of
$\mathcal{A}$, as in Section \ref{maximal ideals, homomorphisms}.  It
is not difficult to check that $\mathcal{H}(\mathcal{A})$ is a closed
set in $\mathcal{A}^*$ with respect to the weak$^*$ topology, and
hence is compact with respect to the weak$^*$ topology, by the
Banach--Alaoglu theorem.  If $\mathcal{A}$ has a nonzero
multiplicative identity element $e$, then
\begin{equation}
\label{mathcal{H}_1(mathcal{A}) = {lambda in cal{H}(cal{A}) : lambda(e) = 1}}
        \mathcal{H}_1(\mathcal{A}) = \{\lambda \in \mathcal{H}(\mathcal{A}) : 
                                                              \lambda(e) = 1\}
\end{equation}
is the same as the set of nonzero homomorphisms on $\mathcal{A}$.  In
this case, $\mathcal{H}_1(\mathcal{A})$ is also a closed set in
$\mathcal{A}^*$ with respect to the weak$^*$ topology, and hence is
compact with respect to the weak$^*$ topology, by the Banach--Alaoglu
theorem again.

        Let $A$ be a locally compact commutative topological group,
and let $H$ be a Haar measure on $A$.  As usual, the space $L^1(A)$ of
complex-valued integrable functions on $A$ is a commutative Banach
algebra with respect to convolutions.  If $\phi \in \widehat{A}$, then
$f \mapsto \widehat{f}(\phi)$ is a nonzero homomorphism on $L^1(A)$,
and we saw in the preceding section that every nonzero homomorphism
from $L^1(A)$ into the field of complex numbers is of this form.  
This defines a mapping from $\widehat{A}$ into the set
$\mathcal{H}(L^1(A))$ of complex homomorphisms on $L^1(A)$, and one
can check that this mapping is continuous with respect to the topology
on $\mathcal{H}(L^1(A))$ induced by the weak$^*$ topology on the dual of
$L^1(A)$ and the topology defined on $\widehat{A}$ previously.
This is very easy to do when $A$ is compact, and otherwise the main
point is that integrable functions on $A$ can be approximated by functions
with compact support with respect to the $L^1$ norm.

        Of course, if $A$ is compact, then the topology defined on
$\widehat{A}$ previously is the same as the discrete topology, and so 
any mapping from $\widehat{A}$ into another topological space is continuous.
In this case, one can show that the topology on $\widehat{A}$ correspomding
to the weak$^*$ topology on the dual of $L^1(A)$ is also discrete.
This uses the fact that $\widehat{A} \subseteq L^1(A)$ when $A$ is
compact, and the orthogonality of distinct elements of $\widehat{A}$
with respect to the standard $L^2$ inner product on $A$.  If $A$
is discrete, then $\widehat{A}$ is compact, and $L^1(A)$ has a
multiplicative identity element, as in the preceding section.
It is easy to see that the topology defined previously on $\widehat{A}$
is the same as the one induced by the weak$^*$ topology on $L^1(A)$
in this situation as well, basically because compact subsets of
$A$ are finite when $A$ is discrete.

\section{Comparing topologies on $\widehat{A}$}
\label{comparing topologies on widehat{A}}

        Let $A$ be a locally compact commutative topological group
with a Haar measure $H$, as before.  If $\phi \in \widehat{A}$,
then $f \mapsto \widehat{f}(\phi)$ is a bounded linear functional
on $L^1(A)$, which leads to a natural mapping from $\widehat{A}$ into
the dual of $L^1(A)$.  It is easy to see that this mapping is continuous
with respect to the topology on $\widehat{A}$ defined previously and
the weak$^*$ topology on $L^1(A)^*$, as mentioned in the preceding section.
In fact, this mapping is a homeomorphism onto its image, so that the usual
topology on $\widehat{A}$ is the same as the one induced by the weak$^*$
topology on $L^1(A)^*$.

        To see this, let $\phi \in \widehat{A}$ be given, as well as a
nonempty compact set $K \subseteq A$ and $\epsilon > 0$.  Thus
\begin{equation}
\label{{psi in widehat{A} : sup_{x in K} |psi(x) - phi(x)| < epsilon}}
        \bigg\{\psi \in \widehat{A} : \sup_{x \in K} |\psi(x) - \phi(x)| 
                                                         < \epsilon\bigg\}
\end{equation}
is a basic open set in $\widehat{A}$ containing $\phi$.  In order to show
that this contains a relative neighborhood of $\phi$ in $\widehat{A}$
with respect to the topology induced by the weak$^*$ topology on $L^1(A)^*$,
it suffices to check that there are finitely many integrable functions
$f_1, \ldots, f_l$ on $A$ and finitely many positive real numbers
$r_1, \ldots, r_n$ such that
\begin{equation}
\label{{psi in widehat{A} : ... for j = 1, ldots, l}}
 \{\psi \in \widehat{A} : |\widehat{f_j}(\psi) - \widehat{f_j}(\phi)| < r_j
                                          \hbox{ for } j = 1, \ldots, l\}
\end{equation}
is contained in (\ref{{psi in widehat{A} : sup_{x in K} |psi(x) -
    phi(x)| < epsilon}}).

        Let $f$ be a nonnegative continuous real-valued function on $A$
with compact support and $\int_A f \, dH = 1$.  If the support of $f$
is contained in a sufficiently small neighborhood of $0$ in $A$, then
\begin{equation}
\label{|widehat{f}(phi) - 1| = ... < frac{1}{4}}
        |\widehat{f}(\phi) - 1| = 
 \biggl|\int_A f(x) \, \overline{\phi(x)} \, dH(x) - 1\biggr| < \frac{1}{4},
\end{equation}
because $\phi$ is continuous on $A$ and $\phi(0) = 1$.  In particular,
\begin{equation}
\label{|widehat{f}(phi)| > frac{3}{4}}
        |\widehat{f}(\phi)| > \frac{3}{4}.
\end{equation}
Let us now fix such a function $f$ for the rest of the argument.

        Let $\eta$ be a positive real number less than or equal to $1/4$.
As one of the conditions on $\psi \in \widehat{A}$ as in (\ref{{psi in
    widehat{A} : ... for j = 1, ldots, l}}), we ask that
\begin{equation}
\label{|widehat{f}(psi) - widehat{f}(phi)| < eta le frac{1}{4}}
        |\widehat{f}(\psi) - \widehat{f}(\phi)| < \eta \le \frac{1}{4}.
\end{equation}
In particular, this implies that
\begin{equation}
\label{|widehat{f}(psi)| > frac{1}{2}}
        |\widehat{f}(\psi)| > \frac{1}{2}.
\end{equation}

        As in Section \ref{integrable functions, 2}, put $(T_a(f))(x) =
f(x + a)$ for each $a, x \in A$, so that 
\begin{equation}
\label{(widehat{T_a(f)})(psi) = psi(a) widehat{f}(psi)}
        (\widehat{T_a(f)})(\psi) = \psi(a) \, \widehat{f}(\psi)
\end{equation}
for every $\psi \in \widehat{A}$.  If $\widehat{f}(\psi)$ is
sufficiently close to $\widehat{f}(\phi)$, and
$(\widehat{T_a(f)})(\psi)$ is sufficiently close to
$(\widehat{T_a(f)}(\phi)$, then it follows that $\psi(a)$ is as close
as one wants to $\phi(a)$.  Of course,
\begin{eqnarray}
\label{|(widehat{T_a(f)})(psi) - (widehat{T_b(f)})(psi)| le ...}
        |(\widehat{T_a(f)})(\psi) - (\widehat{T_b(f)})(\psi)|
          & \le & \int_A |(T_a(f))(x) - (T_b(f))(x)| \, dH(x) \\
           & = & \int_A |f(x + a) - f(x + b)| \, dH(x) \nonumber \\
           & = & \int_A |f(x + a - b) - f(x)| \, dH(x) \nonumber
\end{eqnarray}
for every $\psi \in \widehat{A}$ and $a, b \in A$.  Because $f$ is a
continuous function with compact support on $A$, and hence $f$ is
uniformly continuous, the right side of (\ref{|(widehat{T_a(f)})(psi)
  - (widehat{T_b(f)})(psi)| le ...}) tends to $0$ as $a - b \to 0$ in
$A$.  If $(\widehat{T_a(f)})(\psi)$ is sufficiently close to
$(\widehat{T_a(f)})(\phi)$, then it follows that $(\widehat{T_b(f)})(\psi)$
is as close to $(\widehat{T_b(f)})(\phi)$ as one wants when $a - b$
is sufficiently close to $0$ in $A$.

        If $\widehat{f}(\psi)$ is sufficiently close to $\widehat{f}(\phi)$,
and if $(\widehat{T_a(f)})(\psi)$ is sufficiently close to
$(\widehat{T_a(f)})(\phi)$, then we get that $\psi(b)$ is as close as
we want to $\phi(b)$ when $a - b$ is sufficiently close to $0$ in $A$.
Note that this works uniformly over $\psi \in \widehat{A}$, in the sense
that how close $a - b$ should be to $0$ in $A$ does not depend on $\psi$.
If $K \subseteq A$ is compact, then $K$ can be covered by neighborhoods
of finitely many of its elements to which the previous statement applies.
In order to ensure that $\psi$ is uniformly close to $\phi$ on $K$, it is 
enough to know that $\widehat{f}(\psi)$ is sufficiently close to 
$\widehat{f}(\phi)$, and that $(\widehat{T_a(f)})(\psi)$ is sufficiently 
close to $(\widehat{T_a(f)})(\phi)$ for finitely many $a \in K$, as desired.

        Suppose now that $A$ is also $\sigma$-compact, which implies
that there is a sequence $K_1, K_2, K_3, \ldots$ of compact subsets of
$A$ such that $K_j$ is contained in the interior of $K_{j + 1}$ for
each $j$ and $\bigcup_{j = 1}^\infty K_j = A$, and hence that every
compact subset of $A$ is contained in $K_j$ for some $j$.  If $A$ is
metrizable too, then one can check that $L^1(A)$ is separable, using
the fact that continuous functions with compact support on $A$ are
dense in $L^1(A)$.  This implies that the closed unit ball in
$L^1(A)^*$ is metrizable with respect to the weak$^*$ topology, as in
the previous section.  In particular, it follows that
$\mathcal{H}(L^1(A))$ is metrizable with respect to the weak$^*$
topology in this case.

\section{Involutions}
\label{involutions}

        Let $\mathcal{A}$ be an algebra over the real or complex numbers.
A mapping
\begin{equation}
\label{x mapsto x^*}
        x \mapsto x^*
\end{equation}
from $\mathcal{A}$ into itself is said to be an
\emph{involution}\index{involutions} if it satisfies the following
conditions.  In the real case, (\ref{x mapsto x^*}) should be linear,
while in the complex case, (\ref{x mapsto x^*}) should be
conjugate-linear.  This means that
\begin{equation}
\label{(x + y)^* = x^* + y^*}
        (x + y)^* = x^* + y^*
\end{equation}
for every $x, y \in \mathcal{A}$ in both cases, and that
\begin{equation}
\label{(t x)^* = t x^*}
        (t \, x)^* = t \, x^*
\end{equation}
for every $x \in \mathcal{A}$ and $t \in {\bf R}$ in the real case, and
\begin{equation}
\label{(t x)^* = overline{t} x^*}
        (t \, x)^* = \overline{t} \, x^*
\end{equation}
for every $x \in \mathcal{A}$ and $t \in {\bf C}$ in the complex case.
In both cases, (\ref{x mapsto x^*}) should also satisfy
\begin{equation}
\label{(x y)^* = y^* x^*}
        (x \, y)^* = y^* \, x^*
\end{equation}
and
\begin{equation}
\label{(x^*)^* = x}
        (x^*)^* = x
\end{equation}
for every $x, y \in \mathcal{A}$.  Note that (\ref{(x^*)^* = x})
implies that (\ref{x mapsto x^*}) maps $\mathcal{A}$ onto itself.  If
$\mathcal{A}$ has a nonzero multiplicative identity element $e$, then
it is easy to see that $e^* = e$.  In this situation, if $x \in
\mathcal{A}$ is invertible, then $x^*$ is invertible in $\mathcal{A}$
as well, and
\begin{equation}
\label{(x^*)^{-1} = (x^{-1})^*}
        (x^*)^{-1} = (x^{-1})^*.
\end{equation}
If $\mathcal{A}$ is equipped with a norm $\|x\|$, then it would be
nice to have that
\begin{equation}
\label{||x^*|| = ||x||}
        \|x^*\| = \|x\|
\end{equation}
for every $x \in \mathcal{A}$.  In particular, this would imply that
(\ref{x mapsto x^*}) is a continuous mapping from $\mathcal{A}$ onto
itself.

        Suppose that $\lambda$ is a homomorphism from $\mathcal{A}$ into
the real or complex numbers, as appropriate.  If (\ref{x mapsto x^*})
is an involution on $\mathcal{A}$, then it follows that $\lambda(x^*)$
is also a homomorphism on $\mathcal{A}$ in the real case, and that
$\overline{\lambda(x^*)}$ is a homomorphism on $\mathcal{A}$ in the
complex case.  If $\mathcal{A}$ has a nonzero multiplicative identity
element $e$, then it is easy to see that $\sigma(x^*) = \sigma(x)$ for
every $x \in \mathcal{A}$ in the real case, and 
\begin{equation}
\label{sigma(x^*) = {overline{mu} : mu in sigma(x)}}
        \sigma(x^*) = \{\overline{\mu} : \mu \in \sigma(x)\}
\end{equation}
for every $x \in \mathcal{A}$ in the complex case.  If $\mathcal{A}$
is equipped with a norm $\|x\|$ that satisfies (\ref{||x^*|| =
  ||x||}), then $r(x^*) = r(x)$ for each $x \in \mathcal{A}$, where
$r(x)$ is as in (\ref{r(x) = inf_{l ge 1} ||x^l||^{1/l}}) in Section
\ref{spectrum, spectral radius}.

        Consider the algebra $C_b(X)$ of bounded continuous complex-valued 
functions on a nonempty topological space $X$.  The mapping that sends
$f \in C_b(X)$ to its complex conjugate defines an involution on
$C_b(X)$ that preserves the supremum norm.  Similarly, if $X$ is a
locally compact Hausdorff topological space that is not compact, then
complex-conjugation defines an involution on $C_0(X)$ that preserves
the supremum norm.  Now let $A$ be a locally compact commutative topological
group, and let $H$ be a Haar measure on $A$.  Thus the space $L^1(A)$
of complex-valued integrable functions on $A$ is a commutative Banach
algebra with respect to convolution.  Although complex-conjugation
also defines a norm-preserving involution on $L^1(A)$, it is customary
to use instead the mapping that sends $f \in L^1(A)$ to
$\overline{f(-x)}$.  One can check that this defines a norm-preserving
involution on $L^1(A)$ as well, which corresponds to complex
conjugation of the Fourier transform of $f$, as in Section
\ref{integrable functions, 2}.  If $V$ is a real or complex Hilbert
space, then the adjoint of a bounded linear operator on $V$ defines an
involution on the algebra of bounded linear operators on $V$ that
preserves the operator norm.

        Let $\mathcal{A}$ be a complex Banach algebra with a nonzero 
multiplicative identity element $e$, a norm $\|x\|$, and an involution
(\ref{x mapsto x^*}).  If
\begin{equation}
\label{||x x^*|| = ||x||^2}
        \|x \, x^*\| = \|x\|^2
\end{equation}
for each $x \in \mathcal{A}$, then $\mathcal{A}$ is said to be a
\emph{$C^*$-algebra}.  The algebra $C_b(X)$ of bounded continuous
complex-valued functions on a nonempty topological space $X$ is a
$C^*$-algebra, using the supremum norm and complex-conjugation as the
involution.  The algebra of bounded linear operators on a nontrivial
complex Hilbert space $V$ is a $C^*$-algebra too, using the operator
norm and the involution defined by the adjoint.  If $\mathcal{A}$ is a
$C^*$-algebra, then
\begin{equation}
\label{||x||^2 = ||x x^*|| le ||x|| ||x^*||}
        \|x\|^2 = \|x \, x^*\| \le \|x\| \, \|x^*\|
\end{equation}
for each $x \in \mathcal{A}$, and hence $\|x\| \le \|x^*\|$.  Applying
this to $x^*$, we get that
\begin{equation}
        \|x^*\| \le \|(x^*)^*\| = \|x\|,
\end{equation}
so that (\ref{||x^*|| = ||x||}) holds for every $x \in \mathcal{A}$.
Note that
\begin{equation}
\label{(x x^*)^* = (x^*)^* x^* = x x^*}
        (x \, x^*)^* = (x^*)^* \, x^* = x \, x^*
\end{equation}
for every $x \in \mathcal{A}$, and that $\|y^2\| = \|y\|^2$ for every
$y \in \mathcal{A}$ with $y^* = y$.

        Applying the previous remark to $y = x \, x^*$, we get that
\begin{equation}
\label{||x x^* x x^*|| = ||x x^*||^2 = ||x||^4}
        \|x \, x^* \, x \, x^*\| = \|x \, x^*\|^2 = \|x\|^4
\end{equation}
for every $x \in \mathcal{A}$, using also (\ref{||x x^*|| = ||x||^2})
in the last step.  If $x$ commutes with $x^*$, then it follows that
\begin{equation}
\label{||x||^4 = ||x^2 (x^*)^2|| le ||x^2|| ||(x^*)^2||}
        \|x\|^4 = \|x^2 \, (x^*)^2\| \le \|x^2\| \, \|(x^*)^2\|.
\end{equation}
Of course, $(x^*)^2 = (x^2)^*$, so that $\|(x^*)^2\| = \|(x^2)^*\| =
\|x^2\|$.  Combining this with (\ref{||x||^4 = ||x^2 (x^*)^2|| le
  ||x^2|| ||(x^*)^2||}), we get that
\begin{equation}
\label{||x||^4 le ||x^2||^2}
        \|x\|^4 \le \|x^2\|^2,
\end{equation}
which is to say that $\|x\|^2 \le \|x^2\|$.  This implies that
\begin{equation}
\label{||x^2|| = ||x||^2}
        \|x^2\| = \|x\|^2
\end{equation}
for every $x \in \mathcal{A}$ such that $x \, x^* = x^* \, x$, since
$\|x^2\| \le \|x\|^2$ by the definition of a Banach algebra.

        Suppose that $x \in \mathcal{A}$ satisfies
\begin{equation}
\label{||x^n|| = ||x||^n}
        \|x^n\| = \|x\|^n
\end{equation}
for some positive integer $n$, and let us check that $\|x^l\| =
\|x\|^l$ for each $l \in {\bf Z}_+$ with $l < n$.  By hypothesis,
\begin{equation}
\label{||x||^n = ||x^n|| = ||x^l x^{n - l}|| le ||x||^l ||x||^{n - l}}
        \|x\|^n = \|x^n\| = \|x^l \, x^{n - l}\| \le \|x\|^l \, \|x\|^{n - l}
\end{equation}
for each $l < n$, using the definition of a Banach algebra in the last
step.  Thus $\|x\|^l \le \|x^l\|$ when $l < n$, and hence $\|x^l\| =
\|x\|^l$, because $\|x^l\| \le \|x\|^l$ for each $l$.

        If $x \in \mathcal{A}$ and $x \, x^* = x^* \, x$, then $x^k$ has
the same property for each positive integer $k$.  Applying 
(\ref{||x^2|| = ||x||^2}) to $x^k$, we get that $\|x^{2 \, k}\| = \|x^k\|^2$
for each $k \in {\bf Z}_+$, which we can use repeatedly to obtain that
(\ref{||x^n|| = ||x||^n}) holds when $n = 2^r$ for some $r \in {\bf
  Z}_+$.  This implies that (\ref{||x^n|| = ||x||^n}) holds for every
$n \in {\bf Z}_+$, by the remarks in the previous paragraph.  Thus
\begin{equation}
\label{r(x) = ||x||}
        r(x) = \|x\|,
\end{equation}
where $r(x)$ is as in (\ref{r(x) = inf_{l ge 1} ||x^l||^{1/l}}) in
Section \ref{spectrum, spectral radius}.  Because $\mathcal{A}$ is
complex, we can combine this with (\ref{r(x) = sup {|lambda| : lambda
    in sigma(x)}}), to get that
\begin{equation}
\label{||x|| = sup {|mu| : mu in sigma(x)}}
        \|x\| = \sup \{|\mu| : \mu \in \sigma(x)\}.
\end{equation}
Note that the supremum on the right side of (\ref{||x|| = sup {|mu| :
    mu in sigma(x)}}) is attained, since $\sigma(x)$ is compact.  If
$\mathcal{A}$ is a commutative $C^*$-algebra, then it follows that for
each $x \in \mathcal{A}$ there is a homomorphism $\lambda$ from
$\mathcal{A}$ into ${\bf C}$ such that $|\lambda(x)| = \|x\|$.

\chapter{Operators on $L^2$}
\label{operators on L^2}

\section{Convolution operators}
\label{convolution operators}

        Let $A$ be a locally compact commutative topological group
with a Haar measure $H$.  If $\theta$ is a complex-valued integrable
function on $A$ with respect to $H$, then
\begin{equation}
\label{C_theta(f) = theta * f}
        C_\theta(f) = \theta * f
\end{equation}
defines a bounded linear operator on the space $L^2(A)$ of
complex-valued square-integrable functions on $A$ with respect to $H$,
with operator norm less than or equal to the $L^1$ norm of $\theta$.
This follows from the discussion in Chapter \ref{functions, measures},
applied to the measure $\nu(E) = \int_E \theta \, dH$.  Using Fubini's
theorem, it is easy to see that the adjoint of $C_\theta$ with respect
to the standard inner product
\begin{equation}
\label{langle f, g rangle = int_A f(x) overline{g(x)} dH(x), 2}
        \langle f, g \rangle = \int_A f(x) \, \overline{g(x)} \, dH(x)
\end{equation}
on $L^2(A)$ is given by $C_\theta^* = C_{\widetilde{\theta}}$, where
$\widetilde{\theta}(x) = \overline{\theta(-x)}$.

        Suppose now that $A$ is compact, and let us normalize $H$
as usual so that $H(A) = 1$.  If $\theta$ is a continuous function on
$A$, then one can check that $C_\theta$ is a compact linear operator
on $L^2(A)$, in the sense that $C_\theta$ can be approximated in the
operator norm by finite rank operators on $L^2(A)$.  This uses the
fact that $\theta$ is uniformly continuous on $A$, and it also holds
for integrable functions $\theta$ on $A$, since they can be
approximated by continuous functions with respect to the $L^1$ norm.
Note that $C_\theta$ is normal as well, which means that it commutes
with its adjoint, since its adjoint is a convolution operator too, as
in the previous paragraph.

        It follows from well-known results about compact normal
linear operators on Hilbert spaces that for each $\theta \in L^1(A)$,
there is an orthonormal basis for $L^2(A)$ consisting of eigenvectors
for $C_\theta$.  Put
\begin{equation}
\label{E(lambda, theta) = {f in L^2(A) : C_theta(f) = lambda f}}
        E(\lambda, \theta) = \{f \in L^2(A) : C_\theta(f) = \lambda \, f\}
\end{equation}
for each $\lambda \in {\bf C}$ and $\theta \in L^1(A)$, which is the
eigenspace in $L^2(A)$ corresponding to the eigenvalue $\lambda$ and
the linear operator $C_\theta$.  It is also well known that
$E(\lambda, \theta)$ is a finite-dimensional linear subspace of
$L^2(A)$ when $\lambda \ne 0$, because $C_\theta$ is a compact linear
operator on $L^2(A)$.

        Let $T_a$ be the usual translation operator on $L^2(A)$ for
each $a \in A$, so that $(T_a(f))(x) = f(x + a)$ for every $f \in
L^2(A)$.  Thus $T_a$ is a unitary operator on $L^2(A)$ for each $a \in
A$, and it is easy to see that $T_a$ commutes with $C_\theta$ for
every $\theta \in L^1(A)$.  This implies that $T_a$ maps $E(\lambda,
\theta)$ into itself for each $a \in A$, $\lambda \in {\bf C}$, and
$\theta \in L^1(A)$, and in fact that
\begin{equation}
\label{T_a(E(lambda, theta)) = E(lambda, theta)}
        T_a(E(\lambda, \theta)) = E(\lambda, \theta)
\end{equation}
for every $a \in A$, $\lambda \in {\bf C}$, and $\theta \in L^1(A)$,
since the previous statement also applies to $T_{-a} = T_a^{-1}$.
Suppose that $\lambda \ne 0$, so that $E(\lambda, \theta)$ is a
finite-dimensional linear subspace of $L^2(A)$, as before.  The
restriction of $T_a$ to $E(\lambda, \theta)$ is a unitary operator on
$E(\lambda, \theta)$ with respect to the restriction of the standard
inner product on $L^2(A)$ to $E(\lambda, \theta)$ for each $a \in A$,
and so for each $a \in A$ there is an orthonormal basis for
$E(\lambda, \theta)$ consisting of eigenvectors for $T_a$.  Because
$A$ is commutative, and hence $T_a$ commutes with $T_b$ for every $a,
b \in A$, it is well known that there is an orthonormal basis for
$E(\lambda, \theta)$ whose elements are eigenvectors for $T_a$ for
each $a \in A$ simultaneously.  This uses the fact that the
eigenspaces for $T_a$ are invariant under $T_b$ for every $a, b \in
A$, since $T_a$ and $T_b$ commute with each other.

\section{Simultaneous eigenfunctions}
\label{simultaneous eigenfunctions}

        Let $A$ be a compact commutative topological group, as in the
previous section.  Also let $f \in L^2(A)$ be an eigenvector for $T_a$
for each $a \in A$ with $L^2$ norm equal to $1$, so that for each $a
\in A$ there is a complex number $\mu(a)$ such that
\begin{equation}
\label{T_a(f) = mu(a) f}
        T_a(f) = \mu(a) \, f.
\end{equation}
In particular,
\begin{equation}
\label{mu(a) = langle T_a(f), f rangle}
        \mu(a) = \langle T_a(f), f \rangle
\end{equation}
for each $a \in A$, which implies that $\mu(a)$ is a continuous
function of $a \in A$, since $a \mapsto T_a(f)$ is continuous as a
mapping from $A$ into $L^2(A)$.  As usual, the continuity of $a
\mapsto T_a(f)$ follows from the uniform continuity of $f$ when $f$ is
a continuous function on $A$, and otherwise one can approximate any $f
\in L^2(A)$ by continuous functions on $A$ with respect to the $L^2$
norm.

        Note that $\mu(0) = 1$ and $|\mu(a)| = 1$ for each $a \in A$, 
because $T_0$ is the identity operator and $T_a$ is unitary for every
$a \in A$.  Similarly,
\begin{equation}
\label{mu(a + b) = mu(a) mu(b)}
        \mu(a + b) = \mu(a) \, \mu(b)
\end{equation}
for all $a, b \in A$, because $T_a \circ T_b = T_{a + b}$.  This shows
that $\mu$ defines a continuous character on $A$.

        If we put $g(x) = \mu(x)^{-1} \, f(x)$, then it is easy to
see that $T_a(g) = g$ for every $a \in A$.  More precisely, this means
that $T_a(g)$ is equal to $g$ as elements of $L^2(A)$, so that for
each $a \in A$ we have that $g(x + a) = g(x)$ for almost every $x \in
A$ with respect to Haar measure.  If $f$ is continuous on $A$, then
$g$ is continuous on $A$ too, and it follows that $g(x + a) = g(x)$
for every $a, x \in A$, which is to say that $g$ is constant on $A$.
Otherwise, one can use Fubini's theorem to argue that for almost every
$x \in A$, $g(x + a) = g(x)$ for almost every $a \in A$.  As soon as
this holds for any $x \in A$, we get that $g$ is equal to a constant
almost everywhere on $A$.  Alternatively, if for each $a \in A$, $g(x+
a) = g(x)$ for almost every $x \in A$, then the convolution of $g$
with any integrable function on $A$ is constant on $A$.  This also
implies that $g$ is equal to a constant almost everywhere on $A$, by
approximating $g$ with respect to the $L^2$ norm by convolutions of
$g$ with functions supported in small neighborhoods of $0$ in $A$.

        Thus $f$ is equal to a constant multiple of $\mu$ almost everywhere
on $A$.  If $\theta \in L^1(A)$, $\lambda \in {\bf C}$, $\lambda \ne
0$, and $E(\lambda, \theta) \ne \{0\}$, then it follows that there is
an orthonormal basis for $E(\lambda, \theta)$ consisting of characters
on $A$.  Conversely, if $\phi \in \widehat{A}$, then $C_\theta(\phi) =
\theta * \phi = \widehat{\theta}(\phi) \, \phi$, as in (\ref{(phi *
  nu)(x) = ... = widehat{nu}(phi) phi(x)}) in Section \ref{functions,
  measures}.

        If $\theta \in L^1(A)$ satisfies $\widehat{\theta}(\phi) = 0$
for every $\phi \in \widehat{A}$, then the previous discussion implies
that $\lambda = 0$ is the only eigenvalue of $C_\theta$, and hence
that $C_\theta(f) = 0$ for every $f \in L^2(A)$.  This implies that
$\theta = 0$ under these conditions, since there are continuous
functions $f$ on $A$ supported in a small neighborhood of $0$ such
that $\theta * f$ approximates $\theta$ with respect to the $L^1$
norm.

        In particular, if $\theta \in L^2(A)$ is orthogonal to every
$\phi \in \widehat{A}$, then $\theta = 0$.  This implies that finite
linear combinations of characters on $A$ are dense in $L^2(A)$.

        Suppose now that $f$ is a continuous function on $A$, and let us
show that $f$ can be approximated uniformly on $A$ by finite linear
combinations of characters on $A$.  As usual, $f$ can be approximated
uniformly on $A$ by functions of the form $f * h$, where $h$ is a
continuous function supported in a small neighborhood of $0$ in $A$.
If $g \in L^2(A)$ approximates $f$ with respect to the $L^2$ norm,
then $g * h$ approximates $f * h$ uniformly on $A$ for any fixed $h
\in L^2(A)$, because the supremum norm of the convolution of two
functions in $L^2(A)$ is less than or equal to the product of their
$L^2$ norms.  If $g$ is a finite linear combination of characters on
$A$, then $g * h$ is also a finite linear combination of characters on
$A$, since the convolution of an integrable function on $A$ with a
character is equal to a constant multiple of that character, as
before.  Because $f$ can be approximated by finite linear combinations
of characters with respect to the $L^2$ norm, $f * h$ can be
approximated uniformly by finite linear combinations of characters for
each $h \in L^2(A)$, and hence $f$ can be approximated uniformly by
finite linear combinations of characters, as desired.

        If $\mu$ is a complex regular Borel measure on $A$ such that
$\widehat{\mu}(\phi) = 0$ for every $\phi \in \widehat{A}$, then
\begin{equation}
\label{int_A f d mu = 0}
        \int_A f \, d\mu = 0
\end{equation}
for every $f \in \widehat{A}$, and hence for every function $f$ on $A$
which can be expressed as a finite linear combination of characters on
$A$.  This implies that (\ref{int_A f d mu = 0}) holds for every
continuous function $f$ on $A$, because every continuous function on
$A$ can be approximated by finite linear combinations of characters
uniformly on $A$, as in the previous paragraph.  It follows that $\mu
= 0$ on $A$ when $\widehat{\mu}(\phi) = 0$ for every $\phi \in
\widehat{A}$.

\section{Non-compact groups}
\label{non-compact groups}

        Let $A$ be a locally compact commutative topological group again.
Even if $A$ is not compact, $\theta \mapsto C_\theta$ defines a
bounded linear mapping from $L^1(A)$ into the space
$\mathcal{BL}(L^2(A))$ of bounded linear operators on $L^2(A)$.  This
mapping is also a homomorphism from $L^1(A)$ as an algebra with
respect to convolution into $\mathcal{BL}(L^2(A))$ as an alegbra with
respect to composition of operators.  Thus the collection
$\mathcal{A}_0$ of bounded linear operators on $L^2(A)$ of the form
$C_\theta$ for some $\theta \in L^1(A)$ is a commutative subalgebra of
$\mathcal{BL}(L^2(A))$.  More precisely, $\mathcal{A}_0$ is a
$*$-subalgebra of $\mathcal{BL}(L^2(A))$, because $C_\theta^* \in
\mathcal{A}_0$ for every $\theta \in L^1(A)$, by the remarks at the
beginning of the chapter.  

        Let $\mathcal{A}_1$ be the linear span of $\mathcal{A}_0$ and 
the identity operator $I$ on $L^2(A)$ in $\mathcal{BL}(L^2(A))$, which
is also a $*$-subalgebra of $\mathcal{BL}(L^2(A))$.  Note that $I \in
\mathcal{A}_0$ when $A$ is discrete, in which case $\mathcal{A}_1 =
\mathcal{A}_0$.  The closure $\mathcal{A}$ of $\mathcal{A}_1$ in
$\mathcal{BL}(L^2(A))$ is a commutative Banach subalgebra of
$\mathcal{BL}(L^2(A))$ that contains the identity operator and is a
$*$-subalgebra of $\mathcal{BL}(L^2(A))$, which implies that
$\mathcal{A}$ is a $C^*$-algebra.  It follows that for each $T \in
\mathcal{A}$ there is a complex homomorphism $\Phi$ on $\mathcal{A}$
such that $|\Phi(T)|$ is equal to the operator norm $\|T\|_{op}$ of
$T$ on $L^2(A)$, as mentioned at the end of Section \ref{involutions}.

        If $\Phi$ is any complex homomorphism on $\mathcal{A}$, then
$\theta \mapsto \Phi(C_\theta)$ defines a complex homomorphism on $L^1(A)$.
As in Section \ref{homomorphisms on L^1(A)}, if $\Phi(C_\theta) \ne 0$
for some $\theta \in L^1(A)$, then there is a continuous homomorphism
$\phi : A \to {\bf T}$ such that
\begin{equation}
\label{Phi(C_theta) = int_A theta overline{phi} dH = widehat{theta}(phi)}
        \Phi(C_\theta) = \int_A \theta(x) \, \overline{\phi(x)} \, dH(x)
                       = \widehat{\theta}(\phi)
\end{equation}
for every $\theta \in L^1(A)$, where $H$ denotes Haar measure on $A$,
as usual.  If $\theta \in L^1(A)$ is not equal to $0$ almost
everywhere on $A$ with respect to $H$, then $C_\theta(f) \ne 0$ for
some $f \in L^2(A)$, as in the previous section.  This implies that
there is a complex homomorphism $\Phi$ on $\widehat{A}$ such that
$\Phi(C_\theta) \ne 0$, as in the preceding paragraph, and hence that
$\widehat{\theta}(\phi) \ne 0$ for some $\phi \in \widehat{A}$.

        If $\nu$ is a nonzero complex regular Borel measure on $A$, 
then it is easy to see that there is a continuous complex-valued
function $f$ on $A$ with compact support such that $(f * \nu)(0) \ne
0$.  It follows that $(f * \nu)(x) \ne 0$ for every $x$ in a
neighborhood of $0$ in $A$, because $f * \nu$ is continuous on $A$, as
in Section \ref{functions, measures}.  The convolution $f * \nu$ may
also be considered as an integrable function on $A$ with respect to
Haar measure $H$, which is thus not equal to $0$ almost everywhere on
$A$ with respect to $H$.  The discussion in the preceding paragraph
implies that the Fourier transform of $f * \nu$ is not identically
zero on $\widehat{A}$, and in particular that the Fourier transform of
$\nu$ is not identically zero on $\widehat{A}$.

\section{Some connections with duality}
\label{some connections with duality}

        Let $A$ be a locally compact commutative topological group.
If $a, b \in A$ and $a \ne b$, then there is a $\phi \in \widehat{A}$
such that $\phi(a) \ne \phi(b)$.  More precisely, if $A$ is discrete,
then this follows from the discussion at the end Section \ref{discrete
  groups}.  If $A$ is compact, then this can be derived from the fact
that every continuous complex-valued function on $A$ can be
approximated uniformly by a finite linear combination of characters,
as in Section \ref{simultaneous eigenfunctions}.  Of course, this also
uses the fact that continuous functions on $A$ separate points, by
Urysohn's lemma.  If $A$ is any locally compact commutative
topological group, then consider the Borel measure on $A$ given by
\begin{equation}
\label{mu_{a, b} = delta_a - delta_b}
        \mu_{a, b} = \delta_a - \delta_b,
\end{equation}
where $\delta_a$, $\delta_b$ are the Dirac masses at $a$, $b$, respectively.
Thus
\begin{equation}
\label{widehat{mu_{a, b}}(phi) = overline{phi(a)} - overline{phi(b)}}
        \widehat{\mu_{a, b}}(\phi) = \overline{\phi(a)} - \overline{\phi(b)}
\end{equation}
for every $\phi \in \widehat{A}$, and $\mu_{a, b} \not\equiv 0$ when
$a \ne b$.  The discussion in the previous section implies that
$\widehat{\mu_{a, b}}(\phi) \ne 0$ for some $\phi \in \widehat{A}$,
which means exactly that $\phi(a) \ne \phi(b)$, as desired.

        Suppose for the moment that $A$ is a closed subgroup of a compact
abelian topological group $B$, with the induced topology, so that $A$
is also compact.  As in the previous paragraph, characters on $B$
separate points.  This implies that every character on $A$ is the
restriction to $A$ of a character on $B$, by the discussion in Section
\ref{compact groups}.  Equivalently, every character on $A$ can be
extended to a character on $B$.  Similarly, if $B$ is a discrete
abelian group and $A$ is any subgroup of $B$, then every character on
$A$ can be extended to a character on $B$, as in Section \ref{discrete
  groups}.  In fact, if $A$ is a closed subgroup of a locally compact
abelian topological group $B$, then it is well known that every
character on $A$ can be extended to $B$, but we shall not get into
this here.  As another variant of this, suppose that $A$ is an open
subgroup of an abelian topological group $B$, which implies that $A$
is a closed subgroup of $B$.  If $\phi$ is a character on $A$, then
$\phi$ can be extended to a homomorphism from $B$ into ${\bf T}$,
as in the case where $B$ is discrete.  Under these conditions,
any extension of $\phi$ to a homomorphism from $B$ into ${\bf T}$
is continuous, because $\phi$ is continuous on $A$ and $A$ is an
open subgroup of $B$.

        Suppose now that $A$ and $B$ are locally compact commutative
topological groups, and that $h$ is a continuous homomorphism from $A$
into $B$.  This leads to a natural \emph{dual homomorphism}\index{dual
  homomorphisms} $\widehat{h}$ from $\widehat{B}$ into $\widehat{A}$,
which sends a character $\phi$ on $B$ to $\phi \circ h$ as a character
on $A$.  It is easy to see that $\widehat{h}$ is continuous as a
mapping from $\widehat{B}$ into $\widehat{A}$ with respect to their
corresponding dual topologies, basically because $h$ maps compact
subsets of $A$ to compact subsets of $B$.  If $A$ is a subgroup of $B$
with the induced topology, and if $h$ is the obvious inclusion mapping
that sends each element of $A$ to itself as an element of $B$, then
$\widehat{h}$ sends each character $\phi$ on $B$ to its restriction to
$A$, as in the previous paragraph.

        Alternatively, let $A$ be a locally compact commutative topological
group, and let $C$ be a closed subgroup of $A$.  Thus the quotient
$B = A / C$ is also a locally compact abelian topological group with
respect to the quotient topology, and the natural quotient mapping
$h : A \to B$ is a continuous homomorphism.  If $\phi$ is a character
on $B$, then $\widehat{h}(\phi) = \phi \circ h$ is a character on $A$
whose kernel contains $C$ as a subgroup.  Conversely, if $\psi$ is a
homomorphism from $A$ into ${\bf T}$ whose kernel contains $C$ as a
subgroup, then $\psi$ can be expressed as $\phi \circ h$ for a unique
homomorphism $\phi$ from $B$ into ${\bf T}$.  If $\psi$ is also continuous
as a mapping from $A$ into ${\bf T}$, then one check that $\phi$ is 
continuous as a mapping from $B$ into ${\bf T}$, because of the way that
the quotient topology is defined.

        It follows that $\widehat{h}$ defines a one-to-one mapping from
$\widehat{B}$ onto the subgroup of $\widehat{A}$ consisting of characters
on $A$ whose kernels contain $C$.  It is easy to see that this is a
closed subgroup of $\widehat{A}$ with respect to the usual dual
topology on $\widehat{A}$.  We have already seen that $\widehat{h}$ is
continuous, and in this case one can check that $\widehat{h}$ is
actually a homeomorphism from $\widehat{B}$ onto its image in
$\widehat{A}$ in this case, with respect to the topology induced by
the one one $\widehat{A}$.  To see this, the main point is that every
compact set in $B$ is contained in the image of a compact set in $A$
under $h$.  More precisely, every point in $B$ is contained in an open
set that is contained in the image of a compact set in $A$, because
$A$ is locally compact and $h$ is an open mapping.  Hence every
compact set $K \subseteq B$ is contained in the union of finitely many
open sets, each of which is contained in the image of a compact set in
$A$ under $h$, so that $K$ is contained in the image of a compact set
$E$ in $A$ under $h$.  Of course, $h^{-1}(K)$ is a closed set in $A$,
and so $E \cap h^{-1}(K)$ is a compact set in $A$ which is mapped onto
$K$ by $h$.

        Remember that the dual group $\widehat{A}$ associated to a 
locally compact abelian topological group $A$ is also a locally
compact abelian topological group in a natural way, as in Section
\ref{local compactness of widehat{A}}.  Thus the \emph{second
  dual}\index{second duals}
$\widehat{\widehat{A}}$\index{A'@$\widehat{\widehat{A}}$} may be
defined as usual as the group of continuous homomorphisms from
$\widehat{A}$ into the multiplicative group ${\bf T}$ of complex
numbers with modulus $1$.  This is a locally compact commutative
topological group as well, using the topology determined on
$\widehat{\widehat{A}}$ by the one on $\widehat{A}$ in the same way
that the topology on $\widehat{A}$ was defined in terms of the
topology on $A$.  In particular, if $A$ is discrete, then
$\widehat{A}$ is compact, and hence $\widehat{\widehat{A}}$ is
discrete too.  Similarly, if $A$ is compact, then $\widehat{A}$ is
discrete, which implies that $\widehat{\widehat{A}}$ is compact.

        Let $a \in A$ be given, and put
\begin{equation}
\label{Psi_a(phi) = phi(a), 3}
        \Psi_a(\phi) = \phi(a)
\end{equation}
for each $\phi \in \widehat{A}$, as in (\ref{Psi_a(phi) = phi(a)})
in Section \ref{compact groups}.  This defines a homomorphism from
$\widehat{A}$ into ${\bf T}$, since the group structure on
$\widehat{A}$ is defined by pointwise multiplication.  It is easy to
see that $\Psi_a$ is also continuous with respect to the usual topology
on $\widehat{A}$, basically because $\{a\}$ is a nonempty compact
subset of $A$.  Thus $\Psi_a \in \widehat{\widehat{A}}$ for each $a \in A$,
and in fact the mapping $a \mapsto \Psi_a$ is a homomorphism from $A$
into $\widehat{\widehat{A}}$, because
\begin{equation}
\label{Psi_{a + b}(phi) = phi(a + b) = phi(a) phi(b) = Psi_a(phi) Psi_b(phi)}
        \Psi_{a + b}(\phi) = \phi(a + b) = \phi(a) \, \phi(b) 
                           = \Psi_a(\phi) \, \Psi_b(\phi)
\end{equation}
for every $a, b \in A$.  The fact that elements of $\widehat{A}$
separate points on $A$ says exactly that $a \mapsto \Psi_a$ is a
one-to-one mapping from $A$ into $\widehat{\widehat{A}}$.

        If $A$ is discrete, then every element of $\widehat{\widehat{A}}$
is of the form $\Psi_a$ for some $a \in A$, as discussed at the end of
Section \ref{compact groups}.  If $A$ is any locally compact abelian 
topological group, then one can check that $a \mapsto \Psi_a$ is continuous
as a mapping from $A$ into $\widehat{\widehat{A}}$.  Of course, this is
trivial when $A$ is discrete, and otherwise one can use the fact that
compact subsets of $\widehat{A}$ are equicontinuous on $A$, as in Section
\ref{equicontinuity}.  A well-known theorem states that $a \mapsto \Psi_a$
is a homeomorphism from $A$ onto $\widehat{\widehat{A}}$ for any locally
compact abelian topological group $A$.  If $A$ is compact, then this can
be seen as follows.

        Let $\widehat{\widehat{A}}_1$ be the subgroup of 
$\widehat{\widehat{A}}$ consisting of characters on $\widehat{A}$
of the form $\Psi_a$ for some $a \in A$.  If $A$ is compact, then
$\widehat{\widehat{A}}_1$ is a compact subgroup of $\widehat{\widehat{A}}$,
because $a \mapsto \Psi_a$ is a continuous mapping from $A$ into
$\widehat{\widehat{A}}$.  In this case, it follows immediately that
$a \mapsto \Psi_a$ is a homeomorphism from $A$ onto $\widehat{\widehat{A}}_1$
with respect to the topology induced by the one on $\widehat{\widehat{A}}$,
because a one-to-one continuous mapping from a compact topological
space onto a Hausdorff space is always a homeomorphism.  Thus it remains
to show that $\widehat{\widehat{A}}$ is equal to $\widehat{\widehat{A}}_1$
when $A$ is compact.

        Note that $\widehat{\widehat{A}}_1$ is a closed subgroup of
$\widehat{\widehat{A}}$ when $A$ is compact, because $\widehat{\widehat{A}}_1$
is compact, as before.  Thus $\widehat{\widehat{A}} / \widehat{\widehat{A}}_1$
is also a compact abelian topological group under these conditions.
If $\widehat{\widehat{A}}_1 \ne \widehat{\widehat{A}}$, then 
$\widehat{\widehat{A}} / \widehat{\widehat{A}}_1$ is nontrivial,
and hence there is a nontrivial continuous character on 
$\widehat{\widehat{A}} / \widehat{\widehat{A}}_1$, by the remarks
at the beginning of the section.  Equivalently, this means that there
is a nontrivial continuous character on $\widehat{\widehat{A}}$ whose
kernel contains $\widehat{\widehat{A}}_1$.  However, $\widehat{\widehat{A}}$
is the dual of the discrete abelian group $\widehat{A}$, which implies
that every continuous character on $\widehat{\widehat{A}}$ is given by
evaluation at an element of $\widehat{A}$, as before.  If $\phi \in
\widehat{A}$ has the property that $\Psi_a(\phi) = \phi(a) = 1$ for
every $a \in A$, then $\phi$ is the identity element of $\widehat{A}$,
and hence the character on $\widehat{\widehat{A}}$ corresponding to
evaluation at $\phi$ is trivial as well.  This implies that there is
no nontrivial continuous character on $\widehat{\widehat{A}}$ whose 
kernel contains $\widehat{\widehat{A}}_1$, so that $\widehat{\widehat{A}}_1
= \widehat{\widehat{A}}$, as desired.

\chapter{$p$-Adic numbers}
\label{p-adic numbers}

\section{The $p$-adic absolute value}
\label{p-adic absolute value}

        Let $p$ be a prime number, and let $x$ be a rational number.
The \emph{$p$-adic absolute value}\index{p-adic absolute value@$p$-adic
absolute value} $|x|_p$ of $x$ is defined as follows.  If $x = 0$,
then $|x|_p = 0$.  Otherwise, if $x = p^j \, a / b$ for some integers
$a$, $b$, and $j$, where $a, b \ne 0$ and $a$, $b$ are not divisible
by $p$, then we put
\begin{equation}
\label{|x|_p = p^{-j}}
        |x|_p = p^{-j}.
\end{equation}

        It is easy to see that
\begin{equation}
\label{|x + y|_p le max(|x|_p, |y|_p)}
        |x + y|_p \le \max(|x|_p, |y|_p)
\end{equation}
for every $x, y \in {\bf Q}$.  More precisely, if $x = p^j \, a / b$
and $y = p^j \, c / d$, where $a$, $b$, $c$, $d$, and $j$ are
integers, and $b, d \ne 0$ are not divisible by $p$, then
\begin{equation}
\label{x + y = p^j (frac{a d + b c}{b d})}
        x + y = p^j \Big(\frac{a \, d + b \, c}{b \, d}\Big)
\end{equation}
where $b \, d \ne 0$ is not divisible by $p$.  This implies that $|x +
y|_p \le p^{-j}$ under these conditions, and (\ref{|x + y|_p le
  max(|x|_p, |y|_p)}) follows by taking $j$ as large as possible.

        Like the ordinary absolute value on ${\bf R}$, we also have that
\begin{equation}
\label{|x y|_p = |x|_p |y|_p}
        |x \, y|_p = |x|_p \, |y|_p
\end{equation}
for every $x, y \in {\bf Q}$.  The \emph{$p$-adic metric}\index{p-adic
metric@$p$-adic metric} on ${\bf Q}$ is defined by
\begin{equation}
\label{d_p(x, y) = |x - y|_p}
        d_p(x, y) = |x - y|_p.
\end{equation}
This is actually an ultrametric on ${\bf Q}$, because
\begin{equation}
\label{d_p(x, z) le max(d_p(x, y), d_p(y, z))}
        d_p(x, z) \le \max(d_p(x, y), d_p(y, z))
\end{equation}
for every $x, y, z \in {\bf Q}$, by (\ref{|x + y|_p le max(|x|_p,
  |y|_p)}).

        Note that $|x|_p \le 1$ for every integer $x$.  More precisely,
a rational number $x$ satisfies $|x|_p \le 1$ if and only if $x$ can
be expressed as $a / b$, where $a$ and $b$ are integers, $b \ne 0$,
and $b$ is not divisible by $p$.  Because $p$ is prime, the ring ${\bf
  Z} / p \, {\bf Z}$ of integers modulo $p$ is a field, and hence there
is an integer $c$ such that $b \, c \equiv 1$ modulo $p$.  Thus
$x = (a \, c) / (b \, c) = (a \, c) / (1 - p \, l)$ for some integer $l$.
As usual,
\begin{equation}
\label{sum_{j = 0}^n (p l)^j = frac{1 - (p l)^{n + 1}}{1 - p l}}
        \sum_{j = 0}^n (p \, l)^j = \frac{1 - (p \, l)^{n + 1}}{1 - p \, l}
\end{equation}
for each integer $n \ge 0$.  In this context, $(p \, l)^{n + 1} \to 0$
as $n \to \infty$ with respect to the $p$-adic metric, so that
\begin{equation}
\label{a c sum_{j = 0}^n (p l)^n to frac{a c}{1 - p l} = x}
        a \, c \, \sum_{j = 0}^n (p \, l)^n \to \frac{a \, c}{1 - p \, l} = x
\end{equation}
as $n \to \infty$ with respect to the $p$-adic metric.  This shows that
every $x \in {\bf Q}$ with $|x|_p \le 1$ can be approximated by integers
with respect to the $p$-adic metric.

\section{Completion}
\label{completion}

        In the same way that the real numbers ${\bf R}$ can be obtained
by completing the rational numbers ${\bf Q}$ as a metric space with
respect to the standard metric, the \emph{$p$-adic
  numbers}\index{p-adic numbers@$p$-adic numbers} ${\bf
  Q}_p$\index{Q_p@${\bf Q}_p$} are obtained by completing ${\bf Q}$
with respect to the $p$-adic metric.  Thus ${\bf Q}_p$ can be defined
initially as a complete metric space, with an isometric embedding of
${\bf Q}$ with the $p$-adic metric onto a dense subset of ${\bf Q}_p$.
Let us identify ${\bf Q}$ with this dense subset of ${\bf Q}_p$, and
let $d_p(x, y)$ be the metric on ${\bf Q}_p$, which extends the
$p$-adic metric on ${\bf Q}$.  As usual, ${\bf Q}_p$ can be described
in terms of Cauchy sequences in ${\bf Q}$ with respect to the $p$-adic
metric, and ${\bf Q}_p$ is uniquely determined up to isometric
equivalence.

        Using the description of ${\bf Q}_p$ in terms of Cauchy sequences
in ${\bf Q}$ with respect to the $p$-adic metric or the other properties
of ${\bf Q}_p$ mentioned in the previous paragraph, one can check that the
$p$-adic metric is an ultrametric on ${\bf Q}_p$, and that $d_p(x, y)$
is an integer power of $p$ for each $x, y \in {\bf Q}_p$ with $x \ne y$.
Addition and multiplication can also be extended to ${\bf Q}_p$, so that
${\bf Q}_p$ becomes a commutative ring, with $0$ and $1$ as its additive
and multiplicative identity elements.  The $p$-adic absolute value $|x|_p$
extends to $x \in {\bf Q}_p$ as well, and is an integer power of $p$
when $x \ne 0$.  It is easy to see that (\ref{|x + y|_p le max(|x|_p, |y|_p)}),
(\ref{|x y|_p = |x|_p |y|_p}), and (\ref{d_p(x, y) = |x - y|_p}) continue
to hold on ${\bf Q}_p$.  If $x \in {\bf Q}_p$ and $x \ne 0$, then one
can show that $x$ has a multiplicative inverse in ${\bf Q}_p$, so that
${\bf Q}_p$ is a field.

        If $(M, d(x, y))$ is any metric space and $\{x_j\}_{j = 1}^\infty$
is a Cauchy sequence of elements of $M$, then
\begin{equation}
\label{lim_{j to infty} d(x_j, x_{j + 1}) = 0}
        \lim_{j \to \infty} d(x_j, x_{j + 1}) = 0.
\end{equation}
Although this condition is not normally sufficient to imply that
$\{x_j\}_{j = 1}^\infty$ is a Cauchy sequence in $M$, this does work
when $d(x, y)$ is an ultrametric on $M$.  In particular, if
$\{a_j\}_{j = 1}^\infty$ is a sequence of elements of ${\bf Q}_p$
that converges to $0$ with respect to the $p$-adic metric, then the
partial sums $s_n = \sum_{j = 1}^n a_j$ satisfy
\begin{equation}
\label{s_{n + 1} - s_n = a_{n + 1} to 0}
        s_{n + 1} - s_n = a_{n + 1} \to 0
\end{equation}
in ${\bf Q}_p$ as $n \to \infty$, and hence $\{s_n\}_{n = 1}^\infty$
is a Cauchy sequence in ${\bf Q}_p$.  Thus an infinite series $\sum_{j
  = 1}^\infty a_j$ of $p$-adic numbers converges when $\{a_j\}_{j =
  1}^\infty$ converges to $0$ in ${\bf Q}_p$, because ${\bf Q}_p$
is complete.

        It is easy to see that addition and multiplication define
continuous mappings from ${\bf Q}_p \times {\bf Q}_p$ into ${\bf Q}_p$,
using the product topology on ${\bf Q}_p \times {\bf Q}_p$ corresponding
to the topology on ${\bf Q}_p$ determined by the $p$-adic metric.
In particular, ${\bf Q}_p$ is a commutative topological group with
respect to addition, since the additive inverse of $x \in {\bf Q}_p$
is the same as multiplying $x$ by $-1$.

\section{$p$-Adic integers}
\label{p-adic integers}

        The set ${\bf Z}_p$\index{Z_p@${\bf Z}_p$} of \emph{$p$-adic 
integers}\index{p-adic integers@$p$-adic integers} may be defined by
\begin{equation}
\label{{bf Z}_p = {x in {bf Q}_p : |x|_p le 1}}
        {\bf Z}_p = \{x \in {\bf Q}_p : |x|_p \le 1\}.
\end{equation}
If $x, y \in {\bf Z}_p$, then $x + y \in {\bf Z}_p$ and $x \, y \in
{\bf Z}_p$, by (\ref{|x + y|_p le max(|x|_p, |y|_p)}) and (\ref{|x
  y|_p = |x|_p |y|_p}), so that ${\bf Z}_p$ is a sub-ring of ${\bf
  Q}_p$.  Note that ${\bf Z} \subseteq {\bf Z}_p$, and that ${\bf
  Z}_p$ is a closed set in ${\bf Q}_p$ with respect to the $p$-adic
metric.  Equivalently, ${\bf Z}_p$ is the same as the closure of ${\bf
  Z}$ in ${\bf Q}_p$.  To see this, let $x \in {\bf Z}_p$ be given, so
that $|x|_p \le 1$.  Because ${\bf Q}$ is dense in ${\bf Q}_p$, for
each $\epsilon > 0$ there is a $y \in {\bf Q}$ such that $|x - y|_p <
\epsilon$.  If $\epsilon \le 1$, then we also get that $|y|_p \le 1$,
by the ultrametric version of the triangle inequality.  As in Section
\ref{p-adic absolute value}, $y$ can be approximated by integers with
respect to the $p$-adic metric, and hence $x$ can be too, as desired.

        If $j$ is an integer, then put
\begin{equation}
\label{p^j {bf Z} = {p^j x : x in {bf Z}}}
        p^j \, {\bf Z} = \{p^j \, x : x \in {\bf Z}\},
\end{equation}
and similarly
\begin{equation}
\label{p^j {bf Z}_p = {p^j x : x in {bf Z}_p}}
        p^j \, {\bf Z}_p = \{p^j \, x : x \in {\bf Z}_p\}.
\end{equation}
Thus $p^j \, {\bf Z}$ and $p^j \, {\bf Z}_p$ are sub-rings of ${\bf Q}$
and ${\bf Q}_p$, respectively, and they are ideals in ${\bf Z}$ and
${\bf Z}_p$, respectively, when $j \ge 0$.  Equivalently,
\begin{equation}
\label{p^j {bf Z}_p = {y in {bf Q}_p : |y|_p le p^{-j}}}
        p^j \, {\bf Z}_p = \{y \in {\bf Q}_p : |y|_p \le p^{-j}\}
\end{equation}
for each $j$.  If $j \ge 0$, then the quotients ${\bf Z} / p^j \, {\bf
  Z}$ and ${\bf Z}_p / p^j \, {\bf Z}_p$ can be defined as commutative
rings in the usual way.  Because ${\bf Z} \subseteq {\bf Z}_p$ and
$p^j \, {\bf Z} \subseteq {\bf Z}_p$, there is a natural homomorphism
from ${\bf Z} / p^j \, {\bf Z}$ into ${\bf Z}_p / p^j \, {\bf Z}_p$.
It is easy to see that this homomorphism is a surjection, using the
fact that ${\bf Z}$ is dense in ${\bf Z}_p$.  One can also check that
this homomorphism is injective, and hence an isomorphism.  The main
point is that if $x \in {\bf Z} \cap (p^j \, {\bf Z}_p)$, then $|x|_p
\le p^{-j}$, and hence $x \in p^j \, {\bf Z}_p$.

        Let $(M, d(x, y))$ be any metric space for the moment.
Remember that a set $E \subseteq M$ is said to be
\emph{bounded}\index{bounded sets} if $E$ is contained in a ball in
$M$, and that $E$ is said to be 
\emph{totally bounded}\index{totally bounded sets} in $M$ if for 
each $\epsilon > 0$, $E$ is contained in the union of finitely many
balls of radius $\epsilon$ in $M$.  It is easy to see that totally
bounded subsets of $M$ are bounded, and that compact subsets of $M$
are totally bounded.  If $M$ is complete, then it is well known that
$E \subseteq M$ is compact if and only if $E$ is closed and totally
bounded.  Note that bounded subsets of the real line with the standard
metric are totally bounded.

        If $j$ is a positive integer, then the natural isomorphism
between ${\bf Z} / p^j \, {\bf Z}$ and ${\bf Z}_p / p^j \, {\bf Z}_p$
implies that the latter has exactly $p^j$ elements.  This implies that
${\bf Z}_p$ is the union of $p^j$ translates of $p^j \, {\bf Z}_p$,
and hence that ${\bf Z}_p$ is the union of $p^j$ closed balls of
radius $p^{-j}$ in ${\bf Q}_p$ for each positive integer $j$.
Thus ${\bf Z}_p$ is a compact set in ${\bf Q}_p$ with respect to
the $p$-adic metric, because ${\bf Q}_p$ is complete and ${\bf Z}_p$
is closed and totally bounded.  It follows that ${\bf Q}_p$ is
locally compact, since translated of ${\bf Z}_p$ in ${\bf Q}_p$
are compact as well.  Observe that $p^k \, {\bf Z}_p$ is compact
in ${\bf Q}_p$ for each integer $k$ too, which implies that closed
and bounded subsets of ${\bf Q}_p$ are compact.

\section{Haar measure on ${\bf Q}_p$}
\label{haar measure on Q_p}

        Let $H$ be Haar measure on ${\bf Q}_p$, normalized so that
$H({\bf Z}_p) = 1$.  If $j$ is a positive integer, the it follows
from the discusion in the previous section that ${\bf Z}_p$ is the
union of $p^j$ pairwise-disjoint translates of $p^j \, {\bf Z}_p$.
This implies that
\begin{equation}
\label{H(p^j {bf Z}_p) = p^{-j}}
        H(p^j \, {\bf Z}_p) = p^{-j}.
\end{equation}
This equation also holds when $j \le 0$, in which case $p^j \, {\bf Z}_p$
is the union of $p^{-j}$ pairwise-disjoint translates of ${\bf Z}_p$.
Note that $p^j \, {\bf Z}_p$ is an open and closed subgroup of ${\bf Q}_p$
for each integer $j$, and that the restriction of $H$ to $p^j \, {\bf Z}_p$
is a Haar measure on $p^j \, {\bf Z}_p$ as a compact topological group.

        With this normalization, Haar measure on ${\bf Q}_p$ is the same
as $1$-dimensional Hausdorff measure on ${\bf Q}_p$.  This can be used
as a way to construct Haar measure on ${\bf Q}_p$, which is analogous
to the standard construction of Lebesgue measure on the real line.
Of course, any Hausdorff measure on ${\bf Q}_p$ is automatically
invariant under translations, because the $p$-adic metric on ${\bf
  Q}_p$ is invariant under translations by construction.  However,
$1$-dimensional Hausdorff measure on ${\bf Q}_p$ has the additional
property that it is finite for bounded subsets of ${\bf Q}_p$, and it
is positive for nonempty open subsets of ${\bf Q}_p$.

        Alternatively, if $f$ is a continuous real or complex-valued
function on ${\bf Z}_p$, then one can define a Haar integral of $f$, 
as follows.  We have seen that ${\bf Z}_p$ is the union of $p^j$ 
pairwise-disjoint translates of $p^j \, {\bf Z}_p$ for each positive
integer $j$.  Using this, one can define Riemann sums by averaging
$f$ over $p^j$ points in ${\bf Z}_p$, with exactly one point in each
of the translates of $p^j \, {\bf Z}_p$ in ${\bf Z}_p$.  Because $f$
is continuous on ${\bf Z}_p$ and ${\bf Z}_p$ is compact, $f$ is also
uniformly continuous on ${\bf Z}_p$, and one can use this to show that
any sequence of Riemann sums defined as before is a Cauchy sequence
as $j \to \infty$.  Thus the Riemann sums converge in ${\bf R}$ or
${\bf C}$, as appropriate, as $j \to \infty$, and the Haar integral
of $f$ is defined to be the limit of the Riemann sums.  One can also
use the uniform continuity of $f$ on ${\bf Z}_p$ to show that the limit
of the Riemann sums does not depend on which points in the translates of
$p^j \, {\bf Z}_p$ in ${\bf Z}_p$ were used in the definition of the 
Riemann sums.  One could instead consider upper and lower Riemann sums
for continuous real-valued functions on ${\bf Z}_p$, and show that they
give the same result in the limit, using uniform continuity again.
At any rate, this defines a Haar integral as a nonnegative linear
functional on continuous functions on ${\bf Z}_p$, and one can check that
it is invariant under translations on ${\bf Z}_p$.  As usual, one can then
get Haar measure on ${\bf Z}_p$ as a regular Borel measure from the Haar
integral using the Riesz representation theorem.

        Similarly, one can define a Haar integral for continuous functions
on $p^l \, {\bf Z}_p$ for each integer $l$.  If $l < 0$, then $p^l \,
{\bf Z}_p$ is the union of $p^{-l}$ pairwise-disjoint translates of
${\bf Z}_p$, and the Haar integral of a function on $p^l \, {\bf Z}_p$
is the same as the sum of $p^{-l}$ Haar integrals over ${\bf Z}_p$.
If $f$ is a continuous real or complex-valued function on ${\bf Q}_p$
with compact support, then the support of $f$ is contained in $p^l \,
{\bf Z}_p$ for some integer $l$, and the Haar integral of $f$ as a
function on ${\bf Q}_p$ can be defined as the Haar integral of $f$
over $p^l \, {\bf Z}_p$ for such an $l$.  One can check that this
defines a nonnegative linear functional on the space of continuous
functions on ${\bf Q}_p$ with compact support which is invariant under
translations, so that Haar measure on ${\bf Q}_p$ corresponds to this
Haar integral as in the Riesz representation theorem.

\section{${\bf Q}_p / {\bf Z}_p$}
\label{Q_p / Z_p}

        Consider ${\bf Q}_p / {\bf Z}_p$ as a commutative group with
respect to addition, which is the quotient of the commutative group
${\bf Q}_p$ of $p$-adic numbers with respect to addition by the
subgroup ${\bf Z}_p$ of $p$-adic integers.  This is analogous to ${\bf
  Q} / {\bf Z}$ as the quotient of the commutative group ${\bf Q}$ of
rational numbers with respect to addition by the subgroup ${\bf Z}$ of
integers.  Because of the natural inclusions of ${\bf Q}$ in ${\bf Q}_p$
and ${\bf Z}$ in ${\bf Z}_p$, there is a natural homomorphism from
${\bf Q} / {\bf Z}$ into ${\bf Q}_p / {\bf Z}_p$.  It is easy to see
that this homomorphism is a surjection, because ${\bf Q}$ is dense in
${\bf Q}_p$, and ${\bf Z}_p$ is an open subgroup of ${\bf Q}_p$.
However, this homomorphism is not an injection, because ${\bf Z}$
is a proper subgroup of ${\bf Q} \cap {\bf Z}_p$.

        To deal with this, let $B_p$ be the subgroup of ${\bf Q}$
consisting of rational numbers of the form $p^{-j} \, a$, where $a$
and $j$ are integers.  It is easy to see that $B_p$ is a subgroup of
${\bf Q}$ with respect to addition, and that $B_p$ is dense in ${\bf
  Q}_p$, because ${\bf Z}$ is dense in ${\bf Z}_p$.  Of course, ${\bf
  Z} \subseteq B_p$, and so the quotient $B_p / {\bf Z}$ makes sense
and is a subgroup of ${\bf Q} / {\bf Z}$.  As before, the restriction
of the natural homomorphism from ${\bf Q} / {\bf Z}$ onto ${\bf Q}_p /
{\bf Z}_p$ also maps $B_p / {\bf Z}$ onto ${\bf Q}_p / {\bf Z}_p$,
because $B_p$ is dense in ${\bf Q}_p$ and ${\bf Z}_p$ is an open
subgroup of ${\bf Q}_p$.  This actually defines an isomorphism from
$B_p / {\bf Z}$ onto ${\bf Q}_p / {\bf Z}_p$, because $B_p \cap {\bf
  Z}_p = {\bf Z}$.

        Let $B_{p, j} = p^{-j} \, {\bf Z}$ be the subgroup of $B_p$ 
consisting of rational numbers of the form $p^{-j} \, a$ with $a \in
{\bf Z}$ for each nonnegative integer $j$, so that $B_{p, j} \subseteq
B_{p, j + 1}$ for each $j$, and $B_p = \bigcup_{j = 0}^\infty B_{p,
  j}$.  In particular, ${\bf Z} \subseteq B_{p, j}$ for each $j \ge
0$, so that the quotient $B_{p, j} / {\bf Z}$ makes sense and is a
subgroup of $B_p / {\bf Z}$.  Put $C_p = B_p / {\bf Z}$ and $C_{p, j}
= B_{p, j} / {\bf Z}$ for each nonnegative integer $j$, and observe
that $C_{p, j} \subseteq C_{p, j + 1}$ for each $j \ge 0$ and that
$C_p = \bigcup_{j = 0}^\infty C_{p, j}$.  By construction, $C_{p, j}
= p^{-j} \, {\bf Z} / {\bf Z}$ is isomorphic as a group to
${\bf Z} / p^j \, {\bf Z}$ for each $j \ge 0$.  Similarly,
$C_{p, j + 1}$ is isomorphic to ${\bf Z} / p^{j + 1} \, {\bf Z}$,
and $C_{p, j} \subseteq C_{p, j + 1}$ corresponds to the subgroup
$p \, {\bf Z} / p^{j + 1} \, {\bf Z}$ of ${\bf Z} / p^{j + 1} \, {\bf Z}$
under this isomorphism.

        Of course, we can also think of $B_p$ as a subgroup of the
group ${\bf R}$ of real numbers with respect to addition, and we can
think of $C_p = B_p / {\bf Z}$ as a subgroup of ${\bf R} / {\bf Z}$,
which is isomorphic to the multiplicative group ${\bf T}$ of complex
numbers with modulus equal to $1$.  If $j$ is a positive integer,
then $C_{p, j}$ corresponds to the subgroup of ${\bf T}$ consisting
of complex numbers $z$ such that $z$ to the power $p^j$ is equal to $1$.
Thus $C_p$ corresponds to the subgroup of ${\bf T}$ consisting of
complex numbers $z$ such that $z$ to the power $p^j$ is equal to $1$
for some nonnegative integer $j$.  Note that the quotient topology on
${\bf Q}_p / {\bf Z}_p$ is the same as the discrete topology,
because ${\bf Z}_p$ is an open subgroup of ${\bf Q}_p$.
By constrast, the subgroup of ${\bf T}$ corresponding to $C_p$
is dense in ${\bf T}$ with respect to the standard topology.

\section{Coherent sequences}
\label{coherent sequences}

        As in Section \ref{p-adic integers}, there is a natural ring
isomorphism between ${\bf Z}_p / p^j \, {\bf Z}_p$ and ${\bf Z} / p^j
\, {\bf Z}$ for each positive integer $j$.  Because $p^{j + 1} \, {\bf
  Z} \subseteq p^j \, {\bf Z}$ and $p^{j + 1} \, {\bf Z}_p \subseteq
p^j \, {\bf Z}_p$, there are also natural ring homomorphisms
from ${\bf Z} / p^{j + 1} \, {\bf Z}$ onto ${\bf Z} / p^j \, {\bf Z}$
and from ${\bf Z}_p / p^{j + 1} \, {\bf Z}$ onto ${\bf Z}_p / p^j \,
{\bf Z}_p$ for each $j$.  It is easy to see that these homomorphisms
correspond to each other under the isomorphisms mentioned earlier.

        Consider the Cartesian product $X = \prod_{j = 1}^\infty ({\bf Z} /
p^j \, {\bf Z})$, so that the elements of $X$ are sequences $x = 
\{x_j\}_{j = 1}^\infty$ with $x_j \in {\bf Z} / p^j \, {\bf Z}$
for each $j$.  Thus $X$ is a commutative ring with respect to coordinatewise
addition and multiplication, and a compact Hausdorff topological space with
respect to the product topology associated to the discrete topology on the
finite set ${\bf Z} / p^j \, {\bf Z}$ for each $j$.  Of course, the ring
operations on $X$ are continuous with respect to this topology, so that $X$
is a topological ring.

        Let us say that $x \in X$ is a \emph{coherent sequence}\index{coherent
sequences} if $x_j$ is the image of $x_{j + 1}$ under the natural
homomorphism from ${\bf Z} / p^{j + 1} \, {\bf Z}$ onto ${\bf Z} / 
p^j \, {\bf Z}$ for each $j \ge 1$.  Let $Y$ be the set of coherent
sequences in $X$, which is a closed sub-ring of $X$.

        Let $q_j$ be the natural quotient ring homomorphism from ${\bf Z}$
onto ${\bf Z} / p^j \, {\bf Z}$ for each positive integer $j$.  Put
\begin{equation}
\label{q(a) = {q_j(a)}_{j = 1}^infty}
        q(a) = \{q_j(a)\}_{j = 1}^\infty
\end{equation}
for each $a \in {\bf Z}$, so that $q$ defines a ring homomorphism from
${\bf Z}$ into $X$.  It is easy to see that $q$ is one-to-one, and that
$q(a)$ is a coherent sequence in $X$ for each $a \in {\bf Z}$.  Thus
$q$ maps ${\bf Z}$ into $Y$, and one can check that $q({\bf Z})$ is dense
in $Y$, which is to say that the closure of $q({\bf Z})$ in $X$ with
respect to the product topology is equal to $Y$.  More precisely, for each
$y \in Y$ and positive integer $l$, there is an $a \in {\bf Z}$ such that
$q_l(a) = y_l$, which implies that $q_j(a) = y_j$ when $j \le l$,
by coherence.

        Similarly, there is a natural ring homomorphism $q'$
from ${\bf Z}_p$ onto ${\bf Z} / p^j \, {\bf Z}$ for each positive
integer $j$, which is the composition of the quotient homomorphism from 
${\bf Z}_p$ onto ${\bf Z}_p / p^j \, {\bf Z}_p$ with the usual isomorphism
from ${\bf Z}_p / p^j \, {\bf Z}_p$ onto ${\bf Z} / p^j \, {\bf Z}$.
Thus
\begin{equation}
\label{q'(a) = {q'_j(a)}_{j = 1}^infty}
        q'(a) = \{q'_j(a)\}_{j = 1}^\infty
\end{equation}
defines a ring homomorphism from ${\bf Z}_p$ into $X$.  Note that the
restriction of $q'_j$ to ${\bf Z}$ is equal to $q_j$ for
each $j$, and hence the restriction of $q'$ to ${\bf Z}$ is
equal to $q$.  As before, it is easy to see that $q'$ is a
one-to-one mapping from ${\bf Z}_p$ into the sub-ring $Y$ of coherent
sequences in $X$.

        Because $p^j \, {\bf Z}_p$ is an open subset of ${\bf Z}_p$
for each $j$, $q'_j$ is continuous as a mapping from ${\bf
  Z}_p$ into ${\bf Z} / p^j \, {\bf Z}$ equipped with the discrete
topology.  This implies that $q'$ is continuous as a
mapping from ${\bf Z}_p$ into $X$ with the corresponding product
topology.  In particular, $q'$ maps ${\bf Z}_p$ onto a
compact subset of $X$, since ${\bf Z}_p$ is compact.  Thus
$q'({\bf Z}_p)$ is a closed set in $X$, which is contained
in $Y$ and contains $q({\bf Z})$.  It follows that $q'({\bf Z}_p)
= Y$, because $Y$ is the closure of $q({\bf Z})$ in $X$, although this
can also be verified more directly from the definitions.

        It is easy to see that $q'$ is actually a homeomorphism
from ${\bf Z}_p$ onto $Y$ with respect to the topology on $Y$ induced
by the product topology on $X$, because the open sub-rings $p^j \,
{\bf Z}_p$ of ${\bf Z}_p$ form a local base for the topology of ${\bf
  Z}_p$ at $0$.  Thus $q'$ defines an isomorphism from
${\bf Z}_p$ onto $Y$ as topological rings, and as topological groups
with respect to addition in particular.

\section{Characters}
\label{characters}

        If $\phi$ is a continuous group homomorphism from ${\bf Z}_p$
as a topological group with respect to addition into the
multiplicative group ${\bf T}$ of complex numbers with modulus $1$,
then there is a positive integer $j$ such that $\phi$ maps $p^j \,
{\bf Z}_p$ into the set of $z \in {\bf T}$ such that $\re z > 0$.  As
usual, this implies that $\phi$ maps $p^j \, {\bf Z}_p$ onto the
trivial subgroup $\{1\}$ of ${\bf T}$, since this is the only subgroup
of ${\bf T}$ contained in the set of $z \in {\bf T}$ with $\re z > 0$.
It follows that $\phi$ can be expressed as the composition of the
natural quotient homomorphism from ${\bf Z}_p$ onto ${\bf Z}_p / p^j
\, {\bf Z}_p$ with a group homomorphism from ${\bf Z}_p / p^j \, {\bf
  Z}_p$ into ${\bf T}$.  Conversely, the composition of the natural
quotient homomorphism from ${\bf Z}_p$ onto ${\bf Z}_p / p^j \, {\bf
  Z}_p$ with a group homomorphism from ${\bf Z}_p / p^j \, {\bf Z}_p$
into ${\bf T}$ is a continuous group homomorphism from ${\bf Z}_p$
into ${\bf T}$, because $p^j \, {\bf Z}_p$ is an open subgroup of
${\bf Z}_p$ for each positive integer $j$.

        As in Section \ref{p-adic integers}, ${\bf Z}_p / p^j \, {\bf Z}_p$
is isomorphic as a ring and hence as a group with respect to addition
to ${\bf Z} / p^j \, {\bf Z}$ for every positive integer $j$.  The
group of homomorphisms from ${\bf Z} / p^j \, {\bf Z}$ into ${\bf T}$
is isomorphic to ${\bf Z} / p^j \, {\bf Z}$, as discussed at the
beginning of Section \ref{finite abelian groups}.  Similarly, one can
check that the group of continuous group homomorphisms from ${\bf
  Z}_p$ into ${\bf T}$ is isomorphic to the group $C_p$ from Section
\ref{Q_p / Z_p}.  More precisely, for each positive integer $j$, the
subgroup $C_{p, j}$ of $C_p$ described in Section \ref{Q_p / Z_p}
corresponds exactly to the group of group homomorphisms from ${\bf Z}_p$
into ${\bf T}$ that send $p^j \, {\bf Z}_p$ to the trivial subgroup
$\{1\}$ of ${\bf T}$.  Equivalently, $C_{p, j}$ can be identified with
the group of group homomorphisms from ${\bf Z} / p^j \, {\bf Z}$ into
${\bf T}$ for each $j$.

        Alternatively, let $B(x, y)$ be the image of the product $x \, y$
of $x, y \in {\bf Q}_p$ in ${\bf Q}_p / {\bf Z}_p$, so that $\phi_y(x)
= B(x, y)$ defines a group homomorphism from ${\bf Z}_p$ into ${\bf
  Q}_p / {\bf Z}_p$ for each $y \in {\bf Q}_p$.  As in Section
\ref{Q_p / Z_p}, ${\bf Q}_p / {\bf Z}_p$ is isomorphic to $C_p = B_p /
{\bf Z}$ as a commutative group with respect to addition, which may be 
considered as a subgroup of ${\bf R} / {\bf Z}$, which is isomorphic
to the multiplicative group ${\bf T}$ of $z \in {\bf C}$ with $|z| =
1$.  Thus $\phi_y$ leads to a group homomorphism from ${\bf Z}_p$ into
${\bf T}$, which is easily seen to be continuous, because $B(x, y) =
0$ when $x \, y \in {\bf Z}_p$.  In particular, $B(x, y) = 0$ when
both $x$ and $y$ are in ${\bf Z}_p$, which implies that $\phi_y$
really only depends on the image of $y$ in ${\bf Q}_p / {\bf Z}_p$.
This leads to a homomorphism from ${\bf Q}_p / {\bf Z}_p \cong C_p$
into the dual group of ${\bf Z}_p$, and one can check that this
homomorphism is actually an isomporhism.

        Similarly, for each $x \in {\bf Z}_p$, $y \mapsto B(x, y)$
is a group homomorphism from ${\bf Q}_p$ into ${\bf Q}_p / {\bf Z}_p$
with respect to addition, whose kernel contains ${\bf Z}_p$.  Thus we
can identify this homomorphism with a homomorphism from ${\bf Q}_p /
{\bf Z}_p$ into itself, which leads to a homomorphism from ${\bf Q}_p
/ {\bf Z}_p$ into ${\bf T}$ as before.  One can check directly that
every homomorphism from ${\bf Q}_p / {\bf Z}_p$ into ${\bf T}$ can be
represented by an element of ${\bf Z}_p$ in this way, using the
isomorphism between ${\bf Q}_p / {\bf Z}_p$ and $C_p$ discussed in
Section \ref{Q_p / Z_p}.  More precisely, remember that $C_p$ is the
union of the groups $C_{p, j} \cong {\bf Z} / p^j \, {\bf Z}$ for $j =
1, 2, 3, \ldots$, where $C_{p, j} \subseteq C_{p, j + 1}$ for each
$j$.  A homomorphism on $C_p$ is basically the same as a sequence of
homomorphisms on the $C_{p, j}$'s, with the additional condition that
the $j$th homomorphism on $C_{p, j}$ be equal to the restriction of
the $(j + 1)$th homomorphism on $C_{p, j + 1}$ to $C_{p, j}$.  The
dual of $C_{p, j}$ is isomorphic to ${\bf Z} / p^j \, {\bf Z}$ in the
usual way, and one can use this to get an isomorphism between the dual
of ${\bf Q}_p / {\bf Z}_p$ and ${\bf Z}_p$ as a group with respect to
addition.  This also uses the description of ${\bf Z}_p$ in terms of
coherent sequences in the previous section, where the coherence
condition corresponds exactly to the compatibility condition for
homomorphisms on the $C_{p, j}$'s mentioned earlier.  In particular,
one can use this to show that every element of the dual of $C_p \cong
{\bf Q}_p / {\bf Z}_p$ can be represented in terms of $y \mapsto B(x,
y)$ for some $x \in {\bf Z}_p$, as before.  Moreover, one can check
that this isomorphism between the dual of ${\bf Q}_p / {\bf Z}_p$ and
${\bf Z}_p$ is a homeomorphism with respect to the usual topology on
${\bf Z}_p$ and the topology on the dual of ${\bf Q}_p / {\bf Z}_p$ as
a discrete group.

        Of course, $\phi_y(x) = B(x, y)$ also defines a homomorphism
from ${\bf Q}_p$ as a group with respect to addition into ${\bf Q}_p /
{\bf Z}_p$ for each $y \in {\bf Q}_p$.  If $|y|_p = p^j$ for some
integer $j$, then the kernel of $\phi_y$ is equal to $p^j \, {\bf
  Z}_p$, and otherwise $\phi_y \equiv 0$ when $y = 0$.  As before,
${\bf Q}_p / {\bf Z}_p$ is isomorphic to $C_p = B_p / {\bf Z}$, which
may be considered as a subgroup of ${\bf R} / {\bf Z} \cong {\bf T}$.
This permits $\phi_y$ to be interpreted as a homomorphism from ${\bf
  Q}_p$ into ${\bf T}$ for each $y \in {\bf Q}_p$, and it is easy to
see that this homomorphism is continuous, because of the previous
description of its kernel.

        Conversely, suppose that $\phi$ is a continuous homomorphism
from ${\bf Q}_p$ into ${\bf T}$.  As usual, the continuity condition
implies that there is an integer $l$ such that $\re \phi(x) > 0$
for every $x \in p^l \, {\bf Z}_p$, and hence that $\phi(x) = 1$
for every $x \in p^l \, {\bf Z}_p$, because $\{1\}$ is the only
subgroup of ${\bf T}$ contained in the right half-plane.  Thus $\phi$
is the same as the composition of the natural quotient mapping from
${\bf Q}_p$ onto ${\bf Q}_p / p^l \, {\bf Z}_p$ with a homomorphism
from ${\bf Q}_p / p^l \, {\bf Z}_p$ into ${\bf T}$.  If $l = 0$, then
we have seen that $\phi$ can be represented by $\phi_y$ for some $y
\in {\bf Z}_p$, and otherwise one can reduce to that case using
multiplication in ${\bf Q}_p$ to get that $\phi$ can be represented by
$\phi_y$ for some $y \in p^{-l} \, {\bf Z}_p$.

        Remember that $p^k \, {\bf Z}_p$ is a compact subgroup of ${\bf Q}_p$
with respect to addition for each integer $k$, and that every compact
subset of ${\bf Q}_p$ is contained in $p^k \, {\bf Z}_p$ when $-k$ is
sufficiently large.  If $\phi$ and $\psi$ are continuous homomorphisms
from ${\bf Q}_p$ into ${\bf T}$ that are uniformly close to each other
on $p^k \, {\bf Z}_p$ for some integer $k$, then $\phi(x) = \psi(x)$
for every $x \in p^k \, {\bf Z}_p$, by the usual arguments.  This
simplifies the description of the topology on the dual of ${\bf Q}_p$,
so that a pair of continuous homomorphisms $\phi, \psi : {\bf Q}_p \to
{\bf T}$ are close to each other with respect to this topology when
$\phi(x) = \psi(x)$ for every $x \in p^k \, {\bf Z}_p$ for some
negative integer $k$ such that $-k$ is large.  Using this and the
discussion in the previous paragraphs, it follows that ${\bf Q}_p$ is
isomorphic to its dual group, by an isomorphism that is also a
homeomorphism with respect to the appropriate topologies.

\chapter{$r$-Adic integers and solenoids}
\label{r-adic integers, solenoids}

\section{$r$-Adic absolute values}
\label{r-adic absolute values}

        Let $r = \{r_j\}_{j = 1}^\infty$ be a sequence of positive
integers such that $r_j \ge 2$ for each $j$, and put $R_0 = 1$ and
$R_l = \prod_{j =  1}^l r_j$ for each positive integer $l$.  Also
put $l(0) = +\infty$, and for each nonzero integer $a$, let $l(a)$
be the largest nonnegative integer such that $a$ is an integer
multiple of $R_{l(a)}$.  Equivalently, $l(a) + 1$ is the smallest
positive integer such that $a$ is not an integer multiple of
$R_{l(a) + 1}$ when $a \ne 0$.  It is easy to see that $l(-a) = l(a)$,
\begin{equation}
\label{l(a + b) ge min(l(a), l(b))}
        l(a + b) \ge \min(l(a), l(b))
\end{equation}
and
\begin{equation}
\label{l(a b) ge max(l(a), l(b))}
        l(a \, b) \ge \max(l(a), l(b))
\end{equation}
for every $a, b \in {\bf Z}$.

        Let $t = \{t_l\}_{l = 0}^\infty$ be a strictly decreasing
sequence of positive real numbers that converges to $0$, and put
\begin{equation}
\label{|a|_r = t_{l(a)}}
        |a|_r = t_{l(a)}
\end{equation}
for each nonzero integer $a$, and $|0|_r = 0$.  Thus $|-a|_r = |a|_r$,
\begin{equation}
\label{|a + b|_r le max(|a|_r, |b|_r)}
        |a + b|_r \le \max(|a|_r, |b|_r),
\end{equation}
and
\begin{equation}
\label{|a b|_r le min(|a|_r, |b|_r)}
        |a \, b|_r \le \min(|a|_r, |b|_r)
\end{equation}
for every $a, b \in {\bf Z}$, by the corresponding statements for $l(a)$
in the previous paragraph.  If we put
\begin{equation}
\label{d_r(a, b) = |a - b|_r}
        d_r(a, b) = |a - b|_r
\end{equation}
for each $a, b \in {\bf Z}$, then it follows that
\begin{equation}
\label{d_r(a, c) le max(d_r(a, b), d_r(b, c))}
        d_r(a, c) \le \max(d_r(a, b), d_r(b, c))
\end{equation}
for every $a, b, c \in {\bf Z}$, and hence that $d_r(a, b)$ is an
ultrametric on ${\bf Z}$.  Let us call $|a|_r$ the \emph{$r$-adic
  absolute value}\index{r-adic absolute values@$r$-adic absolute
  values} of $a$ associated to $t = \{t_l\}_{l = 0}^\infty$, and
$d_r(a, b)$ the corresponding \emph{$r$-adic metric}\index{r-adic
  metrics@$r$-adic metrics} on ${\bf Z}$.

        Suppose for the moment that $r = \{r_j\}_{j = 1}^\infty$
is a constant sequence, so that $r_j = r_1$ for each $j$, and $R_l =
(r_1)^l$ for every $l \ge 0$.  In this case, (\ref{l(a b) ge max(l(a),
  l(b))}) can be improved to
\begin{equation}
\label{l(a b) ge l(a) + l(b)}
        l(a \, b) \ge l(a) + l(b)
\end{equation}
for every $a, b \in {\bf Z}$.  If $t_{k + l} \le t_k \, t_l$ for
every $k, l \ge 0$, then we also get that
\begin{equation}
\label{|a b|_r le |a|_r |b|_r}
        |a \, b|_r \le |a|_r \, |b|_r
\end{equation}
for every $a, b \in {\bf Z}$.  If $r_1 = p$ is a prime number,
then we have that
\begin{equation}
\label{l(a b) = l(a) + l(b)}
        l(a \, b) = l(a) + l(b)
\end{equation}
for every $a, b \in {\bf Z}$, and $|a|_r$ is the same as the usual
$p$-adic abslolute value of $a$ when $t_l = p^{-l}$ for each $l \ge 0$.

        Perhaps the simplest choice of $t = \{t_l\}_{l = 0}^\infty$
for an arbitrary $r = \{r_j\}_{j = 1}^\infty$ is given by $t_l =
1/R_l$, which reduces to $t_l = p^{-l}$ when $r_j = p$ for each $j$.
However, for many purposes the choice of $t = \{t_l\}_{l = 0}^\infty$
does not really matter, because the corresponding ultrametrics $d_r(a,
b)$ on ${\bf Z}$ will be topologically equivalent and invariant under
translations.  More precisely, suppose that $t = \{t_l\}_{l = 0}^\infty$
and $t' = \{t_l'\}_{l = 0}^\infty$ are strictly decreasing sequences
of positive real numbers that converge to $0$, and let $d_r(a, b)$ and
$d_r'(a, b)$ be the corresponding ultrametrics on ${\bf Z}$,
associated to the same sequence $r = \{r_j\}_{j = 1}^\infty$.
Under these conditions, it is easy to see that the identity mapping on
${\bf Z}$ is uniformly continuous as a mapping from ${\bf Z}$ equipped
with $d_r(a, b)$ onto ${\bf Z}$ equipped with $d_r'(a, b)$, and as a
mapping from ${\bf Z}$ equipped with $d_r'(a, b)$ onto ${\bf Z}$
equipped with $d_r(a, b)$.

\section{$r$-Adic integers}
\label{r-adic integers}

        Let $r = \{r_j\}_{j = 1}^\infty$ and $R_l = \prod_{j = 1}^l r_j$
be as in the previous section, so that $R_l \, {\bf Z}$ is an ideal
in ${\bf Z}$ as a commutative ring, and $R_{l + 1} \, {\bf Z} \subseteq
R_l \, {\bf Z}$ for each $l$.  Thus the quotient ${\bf Z} / R_l \, {\bf Z}$
is also a commutative ring for each $l$, and there is a natural ring
homomorphism from ${\bf Z} / R_{l + 1} \, {\bf Z}$ onto
${\bf Z} / R_l \, {\bf Z}$.

        Consider the Cartesian product $X = \prod_{l = 1}^\infty 
({\bf Z} / R_l \, {\bf Z})$, whose elements are sequences $x = 
\{x_l\}_{l = 1}^\infty$ with $x_l \in {\bf Z} / R_l \, {\bf Z}$ 
for each $l$.  Observe that $X$ is a commutative ring with respect to
coordinatewise addition and multiplication, as well as a compact
Hausdorff topological space with respect to the product topology
associated to the discrete topologies on the finite sets ${\bf Z} /
R_l \, {\bf Z}$ for each $l$.  More precisely, $X$ is a topological
ring, since the ring operations are continuous with respect to the
product topology on $X$.

        Let $\delta_l(x_l, y_l)$ be the discrete metric on ${\bf Z} /
R_l \, {\bf Z}$ for each $l$, so that $\delta_l(x_l, y_l)$ is equal to
$1$ when $x_l \ne y_l$ and to $0$ when $x_l = y_l$.  Let $t = 
\{t_l\}_{l = 0}^\infty$ be as in the previous section, and put
\begin{equation}
\label{delta_l'(x_l, y_l) = ... = t_{l - 1} delta_l(x_l, y_l)}
        \delta_l'(x_l, y_l) = \min(\delta_l(x_l, y_l), t_{l - 1}) 
                            = t_{l - 1} \, \delta_l(x_l, y_l).
\end{equation}
As in Section \ref{ultrametrics},
\begin{equation}
\label{delta(x, y) = max_{l ge 1} delta_l'(x_l, y_l)}
        \delta(x, y) = \max_{l \ge 1} \delta_l'(x_l, y_l)
\end{equation}
defines an ultrametric on $X$ for which the corresponding topology is
the same as the product topology associated to the discrete topology
on ${\bf Z} / R_l \, {\bf Z}$ for each $l$.  Equivalently, if $x \ne
y$, then
\begin{equation}
\label{delta(x, y) = t_{l(x, y)}}
        \delta(x, y) = t_{l(x, y)},
\end{equation}
where $l(x, y)$ is the smallest nonnegative integer $l$ such that
$x_{l + 1} \ne y_{l + 1}$, which is the same as the largest nonnegative
integer such that $x_j = y_j$ for each $j \le l(x, y)$.  If $x = y$,
then we can put $l(x, y) = +\infty$.

        Let $q_l$ be the natural quotient ring homomorphism from ${\bf Z}$
onto ${\bf Z} / R_l \, {\bf Z}$ for each $l$, and put
\begin{equation}
\label{q(a) = {q_l(a)}_{l = 1}^infty, 2}
        q(a) = \{q_l(a)\}_{l = 1}^\infty
\end{equation}
for each $a \in {\bf Z}$.  This defines a ring homomorphism from ${\bf
  Z}$ into $X$, and it is easy to see that the kernel of this
homomorphism is equal to $\{0\}$, because $R_l \to +\infty$ as $l \to
\infty$.  One can also check that
\begin{equation}
\label{delta(q(a), q(b)) = d_r(a, b)}
        \delta(q(a), q(b)) = d_r(a, b)
\end{equation}
for every $a, b \in {\bf Z}$, where $d_r(a, b)$ is as in (\ref{d_r(a,
  b) = |a - b|_r}).  More precisely,
\begin{equation}
\label{l(q(a), q(b)) = l(a - b)}
        l(q(a), q(b)) = l(a - b)
\end{equation}
for every $a, b \in {\bf Z}$, where $l(a - b)$ is as defined in the
preceding section.

        As before, $x = \{x_l\}_{l = 1}^\infty \in X$ is said to be a 
\emph{coherent sequence}\index{coherent sequences} if $x_l$ is the image
of $x_{l + 1}$ under the natural homomorphism from ${\bf Z} / R_{l +
  1} \, {\bf Z}$ onto ${\bf Z} / R_l \, {\bf Z}$ for each $l$.  It is
easy to see that the collection $Y$ of coherent sequences in $X$ is a
closed sub-ring of $X$, and that $q({\bf Z}) \subseteq Y$.  One can
also check that $Y$ is equal to the closure of $q({\bf Z})$ in $X$.

        Of course, $X$ is complete as a metric space with respect to
$\delta(x, y)$, since it is compact.  This can also be verified more 
directly from the definitions, and it follows that $Y$ is complete as 
a metric space with respect to the restriction of $\delta(x, y)$ to
$x, y \in X$ as well, because $Y$ is closed in $X$.  Thus the
completion ${\bf Z}_r$\index{Z_r@${\bf Z}_r$} of ${\bf Z}$ as a metric
space with respect to the $r$-adic metric $d_r(a, b)$ can be
identified with $Y$ with the restriction of $\delta(x, y)$ to $x, y
\in Y$.  The elements of ${\bf Z}_r$ may be referred to as
\emph{$r$-adic integers}.\index{r-adic integers@$r$-adic integers}

        The identification of ${\bf Z}_r$ with $Y$ shows that addition
and multiplication of integers can be extended in a nice way to
$r$-adic integers, so that ${\bf Z}_r$ is a compact topological ring.
The $r$-adic metric $d_r(a, b)$ has a nice extension to an ultrametric
${\bf Z}_r$ too, which corresponds to $\delta(x, y)$ on $Y$.
Similarly, the $r$-adic absolute value $|a|_r$ can be extended to
${\bf Z}_r$, with properties like those discussed in the previous
section.  Although the $r$-adic absolute value and metric depend on
the choice of $t = \{t_l\}_{l = 0}^\infty$, different choices of $t$
leads to equivalent completions ${\bf Z}_r$ of ${\bf Z}$ as a
topological ring.  In particular, the definition of $Y$ in terms of
coherent sequences does not depend on the choice of $t$.

        Let $k$ be a positive integer, and let $Y_k$ be the set of 
$y = \{y_l\}_{l = 1}^\infty \in Y$ such that $y_k = 0$ in ${\bf Z} / 
R_k \, {\bf Z}$.  This implies that $y_l = 0$ in ${\bf Z} / R_l \, 
{\bf Z}$ when $l \le k$, because of the coherence condition.  It is 
easy to see that $Y_k$ is a closed ideal in $Y$, which is also
relatively open in $Y$.  Equivalently, $Y_k$ is the same as the
closure of $q(R_k \, {\bf Z})$ in $X$, and one can check that $Y /
Y_k$ is isomorphic as a ring to ${\bf Z} / R_k \, {\bf Z}$.

        Let $H$ be Haar measure on $Y$, normalized so that $H(Y) = 1$.
Observe that
\begin{equation}
\label{H(Y_k) = 1/R_k}
        H(Y_k) = 1/R_k
\end{equation}
for each positive integer $k$, because $Y$ can be expressed as the
union of $R_k$ pairwise-disjoint translates of $Y_k$.  If $t_l =
1/R_l$ for each $l \ge 0$, then one can verify that $H$ is the same as
the $1$-dimensional Hausdorff measure on $Y$ corresponding to the
ultrametric $\delta(x, y)$ on $Y$.  Alternatively, a
translation-invariant Haar integral can be defined for continuous real
or complex-valued functions on $Y$ by approximation by suitable
Riemann sums, and Haar measure on $Y$ can be derived from this using
the Riesz representation theorem, as usual.

        Suppose that $\phi$ is a continuous group homomorphism from $Y$
as a group with respect to addition into the multiplicative group ${\bf T}$
of complex numbers with modulus equal to $1$.  Continuity of $\phi$ at 
$0 \in Y$ implies that there is a positive integer $k$ such that
$\re \phi(y) > 0$ for every $y \in Y_k$, and hence that $\phi(y) = 1$
for every $y \in Y_k$, because $\{1\}$ is the only subgroup of ${\bf T}$
contained in the right half-plane.  Thus $\phi$ can be expressed as
the composition of the natural quotient mapping from $Y$ onto
$Y / Y_k \cong {\bf Z} / R_k \, {\bf Z}$ with a homomorphism from
$Y / Y_k$ into ${\bf T}$, and conversely any homomorphism from $Y$
into ${\bf T}$ of this type is continuous.

        More precisely, for each nonnegative integer $l$, let $B_{r, l}$
be the subgroup of the group ${\bf Q}$ of rational numbers with
respect to addition consisting of rational numbers of the form $a /
R_l$, where $a$ is an integer.  Thus $B_{r, 0} = {\bf Z}$ and $B_{r,
  l} \subseteq B_{r, l + 1}$ for each $l \ge 0$, amd we put $B_r =
\bigcup_{l = 0}^\infty B_{r, l}$, which is also a subgroup of ${\bf
  Q}$ with respect to addition.  Similarly, put $C_{r, l} = B_{r, l} /
       {\bf Z}$ and $C_r = B_r / {\bf Z}$, so that $C_{r, l} \subseteq
       C_{r, l + 1} \subseteq C_r$ for each $l$ and $C_r = \bigcup_{l
         = 0}^\infty C_{r, l}$.  One can check that the dual of
$Y \cong {\bf Z}_r$ is isomorphic to $C_r$ as a discrete group,
and that the dual of $C_r$ is isomorphic to $Y$.

\section{Some solenoids}
\label{solenoids}

        Let us continue with the notations and hypotheses in the 
previous two sections.  Of course, $R_l \, {\bf Z}$ is a subgroup of
the real line ${\bf R}$ as a group with respect to addition for each
nonnegative integer $l$.  The quotient ${\bf R} / R_l \, {\bf Z}$ can
be defined as a group as well as a topological space and even a
$1$-dimensional smooth manifold, which is equivalent to the unit
circle ${\bf T}$ in the usual way.

        Let $\widetilde{X}$ be the Cartesian product $\prod_{l = 0}^\infty
({\bf R} / R_l \, {\bf Z})$, so that $\widetilde{X}$ consists of the
sequences $x = \{x_l\}_{l = 0}^\infty$ such that $x_l \in {\bf R} /
R_l \, {\bf Z}$ for each $l \ge 0$.  Thus $\widetilde{X}$ is a commutative
group with respect to coordinatewise addition, and a compact Hausdorff
topological space with respect to the product topology associated to the
quotient topology on ${\bf R} / R_l \, {\bf Z}$ for each $l$.  It is easy
to see that the group operations on $\widetilde{X}$ are continuous,
so that $\widetilde{X}$ is a topological group, which is isomorphic to the
product of a sequence of copies of the unit circle ${\bf T}$.

        Because $R_{l + 1} = R_l \, r_{l + 1}$ and hence $R_{l + 1} \, 
{\bf Z} \subseteq R_l \, {\bf Z}$ for each $l \ge 0$, there is a natural
group homomorphism from ${\bf R} / R_{l + 1} \, {\bf Z}$ onto ${\bf R}
/ R_l \, {\bf Z}$.  The kernel of this homomorphism is equal to $R_l
\, {\bf Z} / R_{l + 1} \, {\bf Z} \cong {\bf Z} / r_{l + 1} \, {\bf
  Z}$ as a subgroup of ${\bf R} / R_{l + 1} \, {\bf Z}$.  This
homomorphism is also a local diffeomorphism from ${\bf R} / R_{l + 1}
\, {\bf Z}$ onto ${\bf R} / R_l \, {\bf Z}$ as $1$-dimensional smooth
manifolds.

        Let us say that $x = \{x_l\}_{l = 0}^\infty \in \widetilde{X}$
is a \emph{coherent sequence}\index{coherent sequences} if $x_l$ is
the image of $x_{l + 1}$ under the natural homomorphism from 
${\bf R} / R_{l + 1} \, {\bf Z}$ onto ${\bf R} / R_l \, {\bf Z}$
described in the preceding paragraph for each $l \ge 0$.  It is easy
to see that the collection $\widetilde{Y}$ of coherent sequences in
$\widetilde{X}$ is a closed subgroup of $\widetilde{X}$.

        Let $\widetilde{q}_l$ be the canonical quotient mapping from
${\bf R}$ onto ${\bf R} / R_l \, {\bf Z}$ for each $l \ge 0$,
which is both a group homomorphism and a local diffeomorphism from
${\bf R}$ onto ${\bf R} / R_l \, {\bf Z}$ as $1$-dimensional smooth
manifolds.  Put
\begin{equation}
\label{widetilde{q}(a) = {widetilde{q}_l(a)}_{l = 0}^infty}
        \widetilde{q}(a) = \{\widetilde{q}_l(a)\}_{l = 0}^\infty
\end{equation}
for each $a \in {\bf R}$, so that $\widetilde{q}$ is a group homomorphism
from ${\bf R}$ into $\widetilde{X}$ which is also continuous with respect
to the product topology on $\widetilde{X}$ mentioned earlier.  It is easy 
to see that the kernel of $\widetilde{q}$ is trivial, and hence that
$\widetilde{q}$ is one-to-one, because $R_l \to \infty$ as $l \to \infty$.

        By construction, $\widetilde{q}(a)$ is a coherent sequence for
each $a \in {\bf R}$, so that $\widetilde{q}({\bf R}) \subseteq
\widetilde{Y}$.  If $y \in \widetilde{Y}$ and $k$ is a nonnegative
integer, then there is an $a \in {\bf R}$ such that $\widetilde{q}_k(a)
= y_k$, which implies that $\widetilde{q}_l(a) = y_l$ for each $l \le
k$, by coherence.  This shows that $\widetilde{Y}$ is the closure of
$\widetilde{q}({\bf R})$ in $\widetilde{X}$.  In particular, it follows
that $\widetilde{Y}$ is connected, since the closure of a connected set
is connected.

        Let $k$ be a nonnegative integer, and let $\widetilde{Y}_k$ be 
the set of $y \in \widetilde{Y}$ such that $y_k = 0$ in ${\bf R} / R_k
\, {\bf Z}$.  Note that $\widetilde{Y}_k$ is a closed subgroup of
$\widetilde{Y}$, and that $y_l = 0$ in ${\bf R} / R_l \, {\bf Z}$ when
$y \in \widetilde{Y}_k$ and $l \le k$, by coherence.  If $y \in
\widetilde{Y}_0$, then $y_0 = 0$ in ${\bf R} / {\bf Z}$, and hence $y_l
\in {\bf Z} / R_l \, {\bf Z}$ for each $l \ge 1$, again by coherence.
Thus $y' = \{y_l\}_{l = 1}^\infty$ is an element of the compact
commutative topological group $Y$ defined in the previous section, and
it is easy to see that $y \mapsto y'$ is an isomorphism between
$\widetilde{Y}_0$ and $Y \cong {\bf Z}_r$ as topological groups.
Under this isomorphism, $\widetilde{Y}_k$ corresponds exactly to the
subgroup $Y_k$ of $Y$ discussed in the previous section for each
positive integer $k$.

        Let $\pi_k$ be the $k$th coordinate projection from $\widetilde{X}$
onto ${\bf R} / R_k \, {\bf Z}$ for each $k \ge 0$, so that $\pi_k(x)
= x_k$ for every $x = \{x_l\}_{l = 0}^\infty \in \widetilde{X}$.  Thus
$\pi_k$ is a continuous homomorphism from $\widetilde{X}$ onto ${\bf R}
/ R_k \, {\bf Z}$, and it is easy to see that $\pi_k(\widetilde{Y}) =
{\bf R} / R_k \, {\bf Z}$ too.  Note that the kernel of the restriction
of $\pi_k$ to $\widetilde{Y}$ is equal to $\widetilde{Y}_k$.

        Suppose that $\phi$ is a continuous homomorphism from $\widetilde{Y}$
into the multiplicative group ${\bf T}$ of complex numbers with
modulus equal to $1$.  Because $\phi$ is continuous at $0$, there is a
nonnegative integer $k$ such that $\re \phi(y) > 0$ for every $y \in
\widetilde{Y}_k$, which implies that $\phi(\widetilde{Y}_k) = \{1\}$, as
usual.  Under these conditions, one can check that $\phi$ is the
composition of the restriction of $\pi_k$ to $\widetilde{Y}_k$ with a
continuous homomorphism from ${\bf R} / R_k \, {\bf Z}$ onto ${\bf
  T}$.  Conversely, every homomorphism from $\widetilde{Y}$ into ${\bf
  T}$ of this form is continuous.

        As in the previous section, let $B_{r, k}$ be the subgroup of
${\bf Q}$ with respect to addition consisting of integer multiples of
$1/R_k$, and put $B_r = \bigcup_{k = 0}^\infty B_{r, k}$.  Using the
discussion of continuous homomorphisms from $\widetilde{Y}$ into ${\bf T}$
in the preceding paragraph, one can show that the dual of $\widetilde{Y}$
is isomorphic to $B_r$.  Similarly, one can show directly that the
dual of $B_r$ is isomorphic to $\widetilde{Y}$, by parameterizing the
dual of $B_{r, k}$ by ${\bf R} / R_k \, {\bf Z}$ for each $k$.
Of course, $B_{r, k}$ is isomorphic as a group to ${\bf Z}$ for each $k$,
and ${\bf R} / R_k \, {\bf Z}$ is isomorphic as a topological group to
${\bf R} / {\bf Z}$, but one should be a it careful about how they
fit together to get $B_r$ and $\widetilde{Y}$.

\section{Complexification}
\label{complexification}

        Again we continue with the notation and hypotheses in the previous 
sections.  Let us now consider the complex plane ${\bf C}$ as a
commutative topological group with respect to addition, which contains
the real line ${\bf R}$ as a closed subgroup, and also $R_l \, {\bf
  Z}$ for each nonnegative integer $l$.  The quotient ${\bf C} / R_l
\, {\bf Z}$ may be considered as a topological group as well as a
Riemann surface.  The complex exponential function determines a
homomorphic group isomorphism between ${\bf C} / {\bf Z}$ and the
multiplicative group of nonzero complex numbers, and similarly
${\bf C} / R_l \, {\bf Z}$ is equivalent to ${\bf C} \backslash \{0\}$
for each $l$.

        Let $\widetilde{X}^{\bf C}$ be the Cartesian product 
$\prod_{l = 0}^\infty ({\bf C} / R_l \, {\bf Z})$, which is the set
of all sequences $x = \{x_l\}_{l = 0}^\infty$ with $x_l \in {\bf C} /
R_l \, {\bf Z}$ for each $l \ge 0$.  As usual, $\widetilde{X}^{\bf C}$
is a commutative topological group with respect to coordinatewise
addition and the product topology associated to the quotient topology
on ${\bf C} / R_l \, {\bf Z}$ for each $l$.  However,
$\widetilde{X}^{\bf C}$ is not locally compact, because ${\bf C} / R_l
\, {\bf Z}$ is not compact for any $l$.  Note that $\widetilde{X}$ as
defined in the preceding section may be considered as a closed
subgroup of $\widetilde{X}^{\bf C}$.

        As before, there is a natural homomorphism from ${\bf C} / 
R_{l + 1} \, {\bf Z}$ onto ${\bf C} / R_l \, {\bf Z}$ for each $l$.
because $R_{l + 1} \, {\bf Z} \subseteq R_l \, {\bf Z}$.  This
homomorphism is also a holomorphic local diffeomorphism from ${\bf C}
/ R_{l + 1} \, {\bf Z}$ onto ${\bf C} / R_l \, {\bf Z}$ as Riemann
surfaces for each $l$.  The kernel of this homomorphism is equal to
$R_l \, {\bf Z} / R_{l + 1} \, {\bf Z} \cong {\bf Z} / r_{l + 1} \,
{\bf Z}$ as a subgroup of ${\bf C} / R_{l + 1} \, {\bf Z}$.

        A sequence $x = \{x_l\}_{l = 0}^\infty \in \widetilde{X}^{\bf C}$
is said to be a \emph{coherent sequence}\index{coherent sequences}
if $x_l$ is the image of $x_{l + 1}$ under the natural homomorphism
from ${\bf C} / R_{l + 1} \, {\bf Z}$ onto ${\bf C} / R_l \, {\bf Z}$
for each $l$, and the collection $\widetilde{Y}^{\bf C}$ of coherent
sequences in $\widetilde{X}^{\bf C}$ is a closed subgroup of
$\widetilde{X}^{\bf C}$.  Note that the set $\widetilde{Y}$ of
coherent sequences in $\widetilde{X}$ is the same as the intersection
of $\widetilde{Y}^{\bf C}$ with $\widetilde{X}$.  Of course,
${\bf C} / R_l \, {\bf Z}$ is isomorphic as a topological group to
the product of ${\bf R} / R_l \, {\bf Z}$ and ${\bf R}$ for each $l$,
and similarly $\widetilde{Y}^{\bf C}$ is isomorphic as a topological
group to the product of $\widetilde{Y}$ and ${\bf R}$.  In particular,
$\widetilde{Y}^{\bf C}$ is locally compact, because $\widetilde{Y}$
is compact and ${\bf R}$ is locally compact.

        Let $\widetilde{q}_l^{\bf C}$ be the canonical quotient mapping
from ${\bf C}$ onto ${\bf C} / R_l \, {\bf Z}$ for each $l$, which is
a group homomorphism and a holomorphic local diffeomorphism from
${\bf C}$ onto ${\bf C} / R_l \, {\bf Z}$ as Riemann surfaces for
each $l$.  Put
\begin{equation}
\label{widetilde{q}^{bf C}(a) = {widetilde{q}_l^{bf C}(a)}_{l = 0}^infty}
 \widetilde{q}^{\bf C}(a) = \{\widetilde{q}_l^{\bf C}(a)\}_{l = 0}^\infty
\end{equation}
for each $a \in {\bf C}$, which defines a continuous group
homomorphism from ${\bf C}$ into $\widetilde{X}^{\bf C}$.  More
precisely, $\widetilde{q}^{\bf C}$ is a one-to-one mapping from ${\bf
  C}$ into $\widetilde{Y}^{\bf C}$, and the closure of
$\widetilde{q}^{\bf C}({\bf C})$ in $\widetilde{X}^{\bf C}$ is equal
to $\widetilde{Y}^{\bf C}$.  This is basically the same as for
$\widetilde{Y}$, because $\widetilde{Y}^{\bf C}$ is isomorphic to
$\widetilde{Y} \times {\bf R}$ as a topological group.

        Let $\pi_k^{\bf C}$ be the $k$th coordinate projection from
$\widetilde{X}^{\bf C}$ onto $C / R_k \, {\bf Z}$ for each nonnegative
integer $k$, which sends $x = \{x_l\}_{l = 0}^\infty \in
\widetilde{X}^{\bf C}$ to $x_k$.  Note that $\pi_k^{\bf
  C}(\widetilde{Y}^{\bf C}) = {\bf C} / R_k \, {\bf Z}$ for each $k$,
and that the kernel of the restriction of $\pi_k^{\bf C}$ to
$\widetilde{Y}^{\bf C}$ is equal to the subgroup $\widetilde{Y}_k$ of
$\widetilde{Y} \subseteq \widetilde{Y}^{\bf C}$ defined in the
previous section.

        Suppose that $\phi$ is a continuous homomorphism from
$\widetilde{Y}^{\bf C}$ into the group ${\bf C} \backslash \{0\}$ of 
nonzero complex numbers with respect to multiplication.  As usual,
continuity of $\phi$ at $0$ implies that there is a nonnegative
integer $k$ such that
\begin{equation}
\label{|phi(y) - 1| < 1/2}
        |\phi(y) - 1| < 1/2
\end{equation}
for every $y \in \widetilde{Y}_k$, and hence that
$\phi(\widetilde{Y}_k) = \{1\}$, because $\phi(\widetilde{Y}_k)$ is a
subgroup of ${\bf C} \backslash \{0\}$.  Using this, one can check
that $\phi$ is the composition of the restriction of $\pi_k^{\bf C}$
to $\widetilde{Y}^{\bf C}$ with a continuous homomorphism $\phi_k :
{\bf C} / R_k \, {\bf Z} \to {\bf C} \backslash \{0\}$.  Thus
\begin{equation}
\label{phi circ widetilde{q}^{bf C} = ... = phi_k circ widetilde{q}_k^{bf C}}
        \phi \circ \widetilde{q}^{\bf C} 
               = \phi_k \circ \pi_k^{\bf C} \circ \widetilde{q}^{\bf C}
               = \phi_k \circ \widetilde{q}_k^{\bf C},
\end{equation}
so that $\phi \circ \widetilde{q}^{\bf C}$ is holomorphic as a mapping
from ${\bf C}$ into ${\bf C} \backslash \{0\}$ if and only if $\phi_k$
is holomorphic as a mapping from ${\bf C} / R_k \, {\bf Z}$ as a Riemann
surface into ${\bf C} \backslash \{0\}$.

        Of course, the continuous homomorphisms from ${\bf T}$ into
${\bf C} \backslash \{0\}$ all map ${\bf T}$ into itself, and are of the
form $z \mapsto z^j$ for some integer $j$.  Similarly, $z \mapsto z^j$
defines a holomorphic homomorphism from ${\bf C} \backslash \{0\}$
into itself for each integer $j$, and every holomorphic homomorphism
from ${\bf C} \backslash \{0\}$ into itself is of this form.  Because
${\bf C} / R_k \, {\bf Z}$ is holomorphically isomorphic to ${\bf C}
\backslash \{0\}$ for each $k$, one can use this to characterize the
continuous homomorphisms $\phi : \widetilde{Y}^{\bf C} \to {\bf C}
\backslash \{0\}$ such that $\phi \circ \widetilde{q}^{\bf C} : {\bf
  C} \to {\bf C} \backslash \{0\}$ is holomorphic.  In particular,
every continuous homomorphism from ${\bf T}$ into itself has a unique
extension to a holomorphic homomorphism from ${\bf C} \backslash
\{0\}$ into itself, and every continuous homomorphism from
$\widetilde{Y}$ into ${\bf T}$ has a unique extension to a continuous
homomorphism from $\widetilde{Y}^{\bf C}$ into ${\bf C} \backslash
\{0\}$ whose composition with $\widetilde{q}^{\bf C}$ is holomorphic.

\chapter{Compactifications}
\label{compactifications}

\section{Compactifications and duality}
\label{compactifications, duality}

        Let $A$ and $B$ be commutative topological groups, and let $h$
be a continuous homomorphism from $A$ into $B$.  Also let
$\widehat{A}$, $\widehat{B}$ be the corresponding dual groups of
continuous homomorphisms from $A$, $B$ into the multiplicative group
${\bf T}$ of complex numbers with modulus equal to $1$, respectively.
If $\phi$ is a continuous homomorphism from $B$ into ${\bf T}$, then
$\phi \circ h$ is a continuous homomorphism from $A$ into ${\bf T}$,
and $\widehat{h}(\phi) = \phi \circ h$ defines a homomorphism from
$\widehat{B}$ into $\widehat{A}$, as in Section \ref{some connections
  with duality}.  We shall be especially interested in the case where
$h(A)$ is dense in $B$, which implies that $\widehat{h}$ is
one-to-one.  Let us also restrict our attention from now on to the
case where $B$ is compact, so that $\widehat{B}$ is discrete.

        Conversely, let $C_d$ be a subgroup of $\widehat{A}$, equipped with
the discrete topology.  Thus the dual $\widehat{C_d}$ of $C_d$ consists
of all homomorphisms from $C_d$ into ${\bf T}$, and is compact with
respect to the usual dual topology.  If $a \in A$ and $\phi \in C_d$,
then put
\begin{equation}
\label{Psi_a(phi) = phi(a), 4}
        \Psi_a(\phi) = \phi(a),
\end{equation}
as in Sections \ref{compact groups} and \ref{some connections with duality}.
It is easy to see that $\Psi_a$ defines a homomorphism from $C_d$ into
${\bf T}$ for each $a \in A$, which is automatically continuous, since
$C_d$ is discrete.  As usual, $a \mapsto \Psi_a$ defines a homomorphism
from $A$ into $\widehat{C_d}$, basically because each $\phi \in C_d$
is a homomorphism on $A$.  Note that $a \mapsto \Psi_a$ is continuous
as a mapping from $A$ into $\widehat{C_d}$, with respect to the usual
dual topology on $\widehat{C_d}$.  This uses from the fact that
(\ref{Psi_a(phi) = phi(a), 4}) is continuous as a function of $a \in A$
for each $\phi \in C_d$, since $C_d \subseteq \widehat{A}$.  If $C_d$
separates points in $A$, then $a \mapsto \Psi_a$ is injective as a
mapping from $A$ into $\widehat{C_d}$.

        Let $B_0$ be the subgroup of $B = \widehat{C_d}$ consisting
of homomorphisms from $C_d$ into ${\bf T}$ of the form $\Psi_a$ for
some $a \in A$, and let $B_1$ be the closure of $B_1$ in
$\widehat{C_d}$.  We would like to check that $B_1 = \widehat{C_d}$,
so that $B_0$ is dense in $\widehat{C_d}$.  Of course, $\widehat{C_d}
/ B_1$ is a compact commutative topological group.  If $\widehat{C_d}
/ B_1$ is not the trivial group, then there is a nontrivial
homomorphism from $\widehat{C_d} / B_1$ into ${\bf T}$.  This would
imply that there is a nontrivial homomorphism from $\widehat{C_d}$
into ${\bf T}$ whose kernel contains $B_1$, by composing the previous
homomorphism with the standard quotient mapping from $\widehat{C_d}$
onto $\widehat{C_d} / B_1$.  However, we have seen that every
continuous homomorphism from $\widehat{C_d}$ into ${\bf T}$
corresponds to evaluation at some $\phi \in C_d$.  If evaluation at
$\phi \in C_d$ contains $B_1$ in its kernel, then (\ref{Psi_a(phi) =
  phi(a), 4}) is equal to $1$ for every $a \in A$, and hence $\phi$ is
the trivial character on $A$.  It follows that there is no nontrivial
continuous homomorphism from $\widehat{C_d}$ into ${\bf T}$ that
contains $B_1$ in its kernel, as desired.

        As a basic class of examples, take $A = {\bf Z}$ and $B = {\bf T}$.
Let $z$ be an element of ${\bf T}$, and put $h_z(j) = z^j$ for each $j
\in {\bf Z}$.  This defines a homomorphism from ${\bf Z}$ into ${\bf
  T}$, and every homomorphism from ${\bf Z}$ into ${\bf T}$ is of this
form.  Of course, $h_z({\bf Z})$ is a finite subgroup of ${\bf T}$ when
$z$ corresponds to an angle which is a rational multiple of $\pi$, and
it is well known that $h_z({\bf Z})$ is dense in ${\bf T}$ otherwise.
In this situation, $\widehat{A} \cong {\bf T}$. $\widehat{B} \cong
{\bf Z}$, and the dual homomorphism $\widehat{h_z}$ maps $\widehat{B}$
onto the subgroup of $\widehat{A}$ generated by $h_z$ as an element of
$\widehat{A}$.

\section{Almost periodic functions}
\label{almost periodic functions}

        Let $A$ be a commutative topological group again, and let
$C_b(A)$ be the space of bounded continuous complex-valued functions
$f$ on $A$, equipped with the supremum norm $\|f\|_{sup}$.  Put
\begin{equation}
\label{f_a(x) = f(x + a)}
        f_a(x) = f(x + a)
\end{equation}
for each $a, x \in A$ and $f \in C_b(A)$, and
\begin{equation}
\label{T(f) = {f_a : a in A}}
        T(f) = \{f_a : a \in A\},
\end{equation}
which is a subset of $C_b(A)$.  If $T(f)$ is totally bounded in
$C_b(A)$ with respect to the supremum norm, then $f$ is said to be
\emph{almost periodic}\index{almost periodic functions} on $A$.
Equivalently, $f$ is almost periodic on $A$ if the closure
$\overline{T(f)}$ of $T(f)$ in $C_b(A)$ is compact, since $C_b(A)$ is
complete with respect to the supremum norm.

        If $f$ is constant on $A$, then $f_a = f$ for each $a \in A$,
$T(f) = \{f\}$, and $f$ is obviously almost periodic.  Similarly,
if $f$ is a continuous homomorphism from $A$ into ${\bf T}$, then $f_a
= f(a) \, f$ for each $a \in A$, and it is easy to see that $f$ is
almost periodic on $A$, basically because the unit circle is compact.
If $A$ is compact, then every continuous function $f$ on $A$ is
uniformly continuous, which implies that $a \mapsto f_a$ is continuous
as a mapping from $A$ into $C_b(A)$.  It follows that $T(f)$ is a
compact set in $C_b(A)$ in this case, and hence that $f$ is almost
periodic.

        If $E$, $E'$ are totally bounded subsets of $C_b(A)$,
then it is easy to see that
\begin{equation}
\label{E + E' = {g + g' : g in E, g' in E'}}
        E + E' = \{g + g' : g \in E, \, g' \in E'\}
\end{equation}
and
\begin{equation}
\label{E cdot E' = {g g' : g in E, g' in E'}}
        E \cdot E' = \{g \, g' : g \in E, \, g' \in E'\}
\end{equation}
are totally bounded in $C_b(A)$ as well.  This uses the continuity of
addition and multiplication on $C_b(A)$, and the fact that totally
bounded sets in $C_b(A)$ are bounded in the second case.  It follows
that sums and products of almost periodic functions on $A$ are almost
periodic, since
\begin{equation}
\label{T(f + f') subseteq T(f) + T(f')}
        T(f + f') \subseteq T(f) + T(f')
\end{equation}
and
\begin{equation}
\label{T(f f') subseteq T(f) cdot T(f')}
        T(f \, f') \subseteq T(f) \cdot T(f')
\end{equation}
for every $f, f' \in C_b(A)$.  Thus the space $\mathcal{AP}(A)$ of
almost periodic functions on $A$ is a subalgebra of $C_b(A)$, and one
can also check that $\mathcal{AP}(A)$ is a closed set in $C_b(A)$ with
respect to the supremum metric.  Note that the complex conjugate of an
almost periodic function on $A$ is almost periodic too.

        If $E \subseteq C_b(A)$ is any totally bounded set with respect to 
the supremum norm, then $E$ is equicontinuous at $0$, in the sense that
for each $\epsilon > 0$ there is an open set $U \subseteq A$ such that
$0 \in U$ and
\begin{equation}
\label{|g(x) - g(0)| < epsilon}
        |g(x) - g(0)| < \epsilon
\end{equation}
for every $x \in U$ and $g \in E$.  More precisely, if $E$ is totally
bounded in $C_b(A)$, then the elements of $E$ can be approximated
uniformly on $A$ by finitely many continuous functions.  This permits
the equicontinuity of $E$ at $0$ to be derived from the continuity of
finitely many functions at $0$, by standard arguments.  In particular,
if $f \in \mathcal{AP}(A)$, then this can be applied to $E = T(f)$, to
get that $f$ is uniformly continuous on $A$.

        Let $B$ be another commutative topological group, and suppose that
$h$ is a continuous homomorphism from $A$ into $B$.  Observe that the mapping
from $g \in C_b(B)$ to $g \circ h \in C_b(A)$ is a bounded linear mapping
with respect to the supremum norms on $C_b(A)$ and $C_b(B)$, with operator
norm equal to $1$, since constant functions on $B$ are sent to constant
functions on $A$ with the same constant value.  If $h(A)$ is dense in $B$,
then $g \mapsto g \circ h$ is an isometric linear embedding of $C_b(B)$
into $C_b(A)$.  At any rate, if $E \subseteq C_b(B)$ is totally bounded
with respect to the supremum norm, then it follows that the set
\begin{equation}
\label{E_h = {g circ h : g in E}}
        E_h = \{g \circ h : g \in E\}
\end{equation}
is totally bounded in $C_b(A)$ with respect to the supremum norm.

        If $g \in C_b(B)$ is almost periodic on $B$, then it is easy to see
that $g \circ h$ is almost periodic on $A$.  More precisely, let $E$ be
the set of translates of $g$ on $B$, which is totally bounded in $C_b(B)$
by hypothesis.  Thus the set $E_h$ of compositions of elements of $E$
with $h$ is totally bounded in $C_b(A)$, as in the previous paragraph.
Because $h : A \to B$ is a homomorphism, every translate of $g \circ h$
on $A$ can be expressed as the composition of a translate of $g$ on $B$
with $h$.  This implies that the collection of translates of $g \circ h$ 
on $A$ is a subset of $E_h$, and hence that the collection of translates
of $g \circ h$ on $A$ is totally bounded in $C_b(A)$, as desired.

        If $B$ is compact, then every continuous function on $B$ is
almost periodic, as before.  In this case, we have also seen that
every continuous function on $B$ can be uniformly approximated by
finite linear combinations of characters on $B$.  It follows that if
$g$ is a continuous function on $B$, then $g \circ h$ can be uniformly
approximated on $A$ by finite linear combinations of characters on
$A$.

        Conversely, let $f$ be any almost periodic function on $A$,
and let $K$ be the closure of $T(f)$ in $C_b(A)$.  Thus $K$ is compact
with respect to the supremum metric on $C_b(A)$.  If $r \in A$, then
consider the mapping that sends $g \in C_b(A)$ to its translate
$g_r(x) = g(x + r)$.  This is an isometric linear mapping from $C_b(A)$
onto itself for each $r \in A$, which sends $f_a$ to $f_{a + r}$ for each
$a \in A$.  It follows that $g \mapsto g_r$ sends $T(f)$ onto itself for
each $r \in A$, which implies that it also sends $K$ onto itself for each 
$r \in A$.

        We have seen that $f$ is uniformly continuous on $A$, so that
for each $\epsilon > 0$ there is an open set $U(\epsilon) \subseteq A$
such that $0 \in U(\epsilon)$ and
\begin{equation}
\label{|f(x + y) - f(x)| < epsilon}
        |f(x + y) - f(x)| < \epsilon
\end{equation}
for every $x \in A$ and $y \in U(\epsilon)$.  This implies that each
translate $f_a$ of $f$ satisfies the same condition, so that
\begin{equation}
\label{|f_a(x + y) - f_a(x)| < epsilon}
        |f_a(x + y) - f_a(x)| < \epsilon
\end{equation}
for every $a, x \in A$ and $y \in U(\epsilon)$.  Similarly, if $g \in K$,
then
\begin{equation}
\label{|g(x + y) - g(x)| le epsilon}
        |g(x + y) - g(x)| \le \epsilon
\end{equation}
for every $x \in A$ and $y \in U(\epsilon)$.  This shows that the
functions in $K$ are uniformly equicontinuous on $A$.

        As before, $g \mapsto g_r$ defines a linear isometry from $C_b(A)$
onto itself for each $r \in A$.  More precisely, the collection of
linear isometries from $C_b(A)$ onto itself is a group with respect
to composition, and the mapping from $r \in A$ to $g \mapsto g_r$
is a homomorphism from $A$ into this group.  We have also seen that
$g \mapsto g_r$ maps $K$ onto itself for each $r \in A$, so that we get
a homomorphism from $A$ into the group $\mathcal{I}(K)$ of isometries
from $K$ onto itself, with respect to the restriction of the supremum
metric on $C_b(A)$ to $K$.

        As in Section \ref{groups of isometries}, $\mathcal{I}(K)$ is a
compact topological group with respect to an appropriate metric, which
is the supremum metric on the space of continuous mappings from $K$
into itself.  It is easy to see that the homomorphism from $A$ into
$\mathcal{I}(K)$ described in the preceding paragraph is continuous
with respect to this metric, because of the uniform equicontinuity
of the elements of $K$ as functions on $A$.  Let $B$ be the closure of
the image of $A$ in $\mathcal{I}(K)$ under this homomorphism, which
is a compact commutative group with respect to the induced topology.

        Put $\lambda(g) = g(0)$ for each $g \in C_b(A)$, which is a 
bounded linear functional on $C_b(A)$.  In particular, the restriction of
$\lambda$ to $K$ is a continuous function on $K$ with respect to the
supremum metric.  Remember that $f \in T(f) \subseteq K$, so that if
$R \in \mathcal{I}(K)$, then $R(f) \in K$ too.  The mapping from $R \in 
\mathcal{I}(K)$ to $R(f) \in K$ is continuous with respect to the
corresponding supremum metrics on $K$ and $\mathcal{I}(K)$, and hence
\begin{equation}
\label{Lambda(R) = lambda(R(f)) = (R(f))(0)}
        \Lambda(R) = \lambda(R(f)) = (R(f))(0)
\end{equation}
is a continuous function on $\mathcal{I}(K)$.  If we apply this to the
$R \in \mathcal{I}(K)$ that corresponds to $g \mapsto g_r$ for some $r
\in A$, then we get that
\begin{equation}
\label{Lambda(R) = (R(f))(0) = f_r(0) = f(r)}
        \Lambda(R) = (R(f))(0) = f_r(0) = f(r).
\end{equation}
This shows that $f$ can be expressed as the composition of the
continuous function $\Lambda$ on $B$ with the natural continuous
homomorphism from $A$ into $B$.  It follows from this and the earlier
discussion that $f$ can be uniformly approximated by finite sums of
characters on $A$, because $B$ is compact.

\section{Spaces of continuous functions}
\label{spaces of continuous functions}

        Let $X$, $Y$ be (nonempty) topological spaces, and suppose 
that $h$ is a continuous mapping from $X$ into $Y$.  Remember that
the spaces $C_b(X)$, $C_b(Y)$ of bounded continuous complex-valued 
functions on $X$, $Y$, respectively, are commutative Banach
algebras with respect to pointwise addition and multiplication
and the corresponding supremum norms.  If $g \in C_b(Y)$, then $g \circ h
\in C_b(X)$, and in fact $g \mapsto g \circ h$ defines a bounded linear
mapping from $C_b(Y)$ into $C_b(X)$ which is also an algebra homomorphism.
More precisely, the operator norm of this mapping with respect to the
supremum norms on $C_b(X)$, $C_b(Y)$ is equal to $1$, as one can see
by considering constant functions $g$ on $Y$.

        If $h(X)$ is dense in $Y$, then $g \mapsto g \circ h$ is an
isometric linear mapping from $C_b(Y)$ into $C_b(X)$.  In this case,
\begin{equation}
\label{{g circ h : g in C_b(Y)}}
        \{g \circ h : g \in C_b(Y)\}
\end{equation}
is actually a closed subalgebra of $C_b(X)$ with respect to the
supremum norm.  More precisely, (\ref{{g circ h : g in C_b(Y)}}) is
complete with respect to the supremum norm on $X$ under these
conditions, because $C_b(Y)$ is complete with respect to the supremum
norm on $Y$.  This implies that (\ref{{g circ h : g in C_b(Y)}}) is a closed
set in $C_b(X)$ with respect to the supremum norm, by standard arguments.

        Now let $\mathcal{A}$ be any closed subalgebra of $C_b(X)$ that
contains the constant functions.  Thus $\mathcal{A}$ is also a
commutative Banach algebra, and the constant function ${\bf 1}_X$ on
$X$ equal to $1$ at each point is the multiplicative identity element
in $\mathcal{A}$.  If $\lambda$ is a homomorphism from $\mathcal{A}$
into the complex numbers such that $\lambda(f) \ne 0$ for some $f \in
\mathcal{A}$, then $\lambda({\bf 1}_X) = 1$, and $\lambda$ is a bounded
linear functional on $\mathcal{A}$ with respect to the supremum norm,
with dual norm equal to $1$.  Let $\mathcal{H}_1(\mathcal{A})$ be the 
collection of nonzero complex homomorphisms on $\mathcal{A}$, which is 
a subset of the unit ball of the dual $\mathcal{A}^*$ of $\mathcal{A}$
as a Banach space.  As before, $\mathcal{H}_1(\mathcal{A})$ is a closed
set in $\mathcal{A}^*$ with respect to the weak$^*$ topology, and hence
is compact with respect to the weak$^*$ topology, by the Banach--Alaoglu
theorem.

        If $p \in X$, then
\begin{equation}
\label{lambda_p(f) = f(p)}
        \lambda_p(f) = f(p)
\end{equation}
defines a nonzero complex homomorphism on $C_b(X)$, and the restriction
of $\lambda_p$ to $\mathcal{A}$ is an element of $\mathcal{H}_1(\mathcal{A})$.
Consider the mapping from $p \in X$ to the restriction of $\lambda_p$
to $\mathcal{A}$, as a mapping from $X$ into
$\mathcal{H}_1(\mathcal{A})$.  It is easy to see that this mapping is
continuous with respect to the weak$^*$ topology on
$\mathcal{H}_1(\mathcal{A}) \subseteq \mathcal{A}^*$, because
(\ref{lambda_p(f) = f(p)}) is a continuous function of $p$ for each $f
\in C_b(X)$.  Of course, this mapping from $X$ into
$\mathcal{H}_1(\mathcal{A})$ is injective if and only if $\mathcal{A}$
separates points on $X$.

        Note that a bounded continuous function $f$ on $X$ is invertible
in $C_b(X)$ if and only if $f(x) \ne 0$ for each $x \in X$ and
$1/f(x)$ is bounded on $X$.  Equivalently, $f$ is invertible in
$C_b(X)$ when there is a $\delta > 0$ such that $|f(x)| \ge \delta$ for
every $x \in X$.  Suppose that $\mathcal{A}$ has the additional
property that for each $f \in \mathcal{A}$, the complex conjugate
$\overline{f}$ of $f$ is an element of $\mathcal{A}$ as well.  In this
case, if $f \in \mathcal{A}$ is invertible in $C_b(X)$, then $1/f \in
\mathcal{A}$ too.  Indeed, if $f \in \mathcal{A}$, then $\overline{f}
\in \mathcal{A}$, and hence $|f|^2 \in \mathcal{A}$.  If $f$ is
invertible in $C_b(X)$, then $|f|^2$ is also invertible in $C_b(X)$,
and $1/f = \overline{f}/|f|^2$.  Thus it suffices to show that
$1/|f|^2 \in \mathcal{A}$, which basically means that one can reduce
to the case of positive real-valued functions on $X$.  If $|f(x)| \le
1$ for each $x \in X$, then the supremum norm of $1 - |f(x)|^2$ is
strictly less than $1$ when $f$ is invertible on $X$, which implies that
\begin{equation}
\label{|f|^2 = 1 - (1 - |f|^2)}
        |f|^2 = 1 - (1 - |f|^2)
\end{equation}
is invertible in $\mathcal{A}$, as desired.  Of course, it is easy to
reduce to the case where the supremum norm of $f$ is equal to $1$.

        Let $\mathcal{H}_1(\mathcal{A}, X)$ be the subset of 
$\mathcal{H}_1(\mathcal{A})$ consisting of the restrictions of $\lambda_p$ 
to $\mathcal{A}$, for each $p \in X$.  If $\mathcal{A}$ is invariant
under complex conjugation, then $\mathcal{H}_1(\mathcal{A}, X)$ is
dense in $\mathcal{H}_1(\mathcal{A})$ with respect to the weak$^*$
topology on $\mathcal{A}^*$.  Otherwise, suppose that $\lambda \in
\mathcal{H}_1(\mathcal{A})$ is not in the weak$^*$ closure of
$\mathcal{H}_1(\mathcal{A}, X)$.  This means that there are finitely
many elements $f_1, \ldots, f_n$ of $\mathcal{A}$ and an $r > 0$ such that
\begin{equation}
\label{max_{1 le j le n} |lambda(f_j) - lambda_p(f_j)| ge r}
        \max_{1 \le j \le n} |\lambda(f_j) - \lambda_p(f_j)| \ge r
\end{equation}
for each $p \in X$.  Put $g_j(x) = f_j(x) - \lambda(f_j) \, {\bf 1}_X(x)$
for $j = 1, \ldots, n$, so that
\begin{equation}
\label{lambda(g_j) = ... = lambda(f_j) - lambda(f_j) = 0}
        \lambda(g_j) = \lambda(f_j) - \lambda(f_j) \, \lambda({\bf 1}_X) 
                     = \lambda(f_j) - \lambda(f_j) = 0
\end{equation}
for each $j$.  Of course, $\lambda_p(f_j) = f_j(p)$ for each $p \in X$
and $j = 1, \ldots, n$, and hence
\begin{equation}
\label{max_{1 le j le n} |g_j(p)| ge r}
        \max_{1 \le j \le n} |g_j(p)| \ge r
\end{equation}
for each $p \in X$, by (\ref{max_{1 le j le n} |lambda(f_j) -
  lambda_p(f_j)| ge r}).  Observe that
\begin{equation}
\label{g(x) = sum_{j = 1}^n |g_j(x)|^2 = sum_{j = 1}^n g_j(x) overline{g_j(x)}}
 g(x) = \sum_{j = 1}^n |g_j(x)|^2 = \sum_{j = 1}^n g_j(x) \, \overline{g_j(x)}
\end{equation}
is an element of $\mathcal{A}$, because $g_j \in \mathcal{A}$ for each
$j = 1, \ldots, n$, and thus $\overline{g_j} \in \mathcal{A}$ for each
$j = 1, \ldots, n$ too.  It is easy to see that $\lambda(g) = 0$,
using (\ref{lambda(g_j) = ... = lambda(f_j) - lambda(f_j) = 0}) and
the fact that $\lambda$ is a homomorphism on $\mathcal{A}$.  However,
(\ref{max_{1 le j le n} |g_j(p)| ge r}) implies that $|g(p)|^2 \ge
r^2$ for each $p \in X$, so that $g$ is invertible in $C_b(X)$, and
hence $g$ is also invertible in $\mathcal{A}$, as in the previous
paragraph.  This is a contradiction, because the invertibility of $g$
in $\mathcal{A}$ implies that $\lambda(g) \ne 0$.  

        In particular, $\mathcal{H}_1(\mathcal{A}, X)$ is dense in
$\mathcal{H}_1(\mathcal{A})$ with respect to the weak$^*$ topology
when $\mathcal{A} = C_b(X)$, in which case
$\mathcal{H}_1(\mathcal{A})$ is the \emph{Stone--{\v C}ech
  compactification} of $X$.  Suppose for the moment that $X$ is
compact, so that $C_b(X)$ is the same as the algebra $C(X)$ of
all complex-valued continuous functions on $X$.  As before,
the mapping from $p \in X$ to $\lambda_p$ as in (\ref{lambda_p(f) = f(p)})
is continuous with respect to the weak$^*$ topology on $C(X)^*$.
If $X$ is compact, then it follows that the set of $\lambda_p$
with $p \in X$ is compact in $C(X)^*$ with respect to the weak$^*$
topology.  This implies that every nonzero complex homomorphism on
$C(X)$ is of the form $\lambda_p$ for some $p \in X$ when $X$ is compact,
because the set of these homomorphisms on $C(X)$ is dense in the set
of all nonzero complex homomorphisms on $C(X)$ with respect to the weak$^*$
topology on $C(X)^*$, as before.

        Let $\mathcal{A}$ be a closed subalgebra of $C_b(X)$ that
contains the constant functions and is invariant under complex
conjugation again, where $X$ is not necessarily compact.  If $\lambda$
is any nonzero complex homomorphism on $C_b(X)$, then the restriction
of $\lambda$ to $\mathcal{A}$ is a nonzero complex homomorphism on
$\mathcal{A}$.  This defines a mapping from $\mathcal{H}_1(C_b(X))$
into $\mathcal{H}_1(\mathcal{A})$, which is continuous with respect to
the weak$^*$ topologies on $C_b(X)^*$ and $\mathcal{A}^*$.  
Because $\mathcal{H}_1(C_b(X))$ is compact with respect to the weak$^*$
topology on $C_b(X)^*$, it follows that the set of complex
homomorphisms on $\mathcal{A}$ that are restrictions of nonzero
complex homomorphisms on $C_b(X)$ to $\mathcal{A}$ is compact with
respect to the weak$^*$ topology on $\mathcal{A}^*$.  Of course,
$\lambda_p \in \mathcal{H}_1(C_b(X))$ for each $p \in X$, and the
set $\mathcal{H}_1(\mathcal{A}, X)$ of the restrictions of the
$\lambda_p$'s to $\mathcal{A}$ is dense in $\mathcal{H}_1(\mathcal{A})$ 
with respect to the weak$^*$ topology on $\mathcal{A}^*$, as before.
Thus the set of restrictions of nonzero complex homomorphisms on
$C_b(X)$ to $\mathcal{A}$ is also dense in $\mathcal{H}_1(\mathcal{A})$ 
with respect to the weak$^*$ topology on $\mathcal{A}^*$.  
This shows that every nonzero complex homomorphism on $\mathcal{A}$
can be obtained from the restriction of a nonzero complex homomorphism
on $C_b(X)$ to $\mathcal{A}$ under these conditions, because the set
of these restrictions is both dense and compact in $\mathcal{H}_1(\mathcal{A})$
with respect to the weak$^*$ topology on $\mathcal{A}^*$.  In particular,
if $X$ is compact, then every nonzero complex homomorphism on $\mathcal{A}$
can be expressed as the restriction of $\lambda_p$ to $\mathcal{A}$ for
some $p \in X$.

        Note that the closed subalgebra $\mathcal{A}$ of $C_b(X)$ given by 
(\ref{{g circ h : g in C_b(Y)}}) contains the constant functions on $X$ and 
is invariant under complex conjugation when $h$ is a continuous mapping from 
$X$ onto a dense subset of a topological space $Y$.

\section{The Bohr compactification}
\label{bohr compactification}

        Let $A$ be a commutative topological group again, and let $\widehat{A}$
be the dual group of continuous homomorphisms from $A$ into ${\bf T}$, as
usual.  To be precise, let $\widehat{A}_d$ be $\widehat{A}$ as a commutative
group equipped with the discrete topology.  Thus the dual 
$B = \widehat{(\widehat{A}_d)}$ of $\widehat{A}_d$ as a discrete group
is a compact commutative topological group in a natural way.
As in Section \ref{compactifications, duality}, the mapping from
$a \in A$ to $\Psi_a$ defined in (\ref{Psi_a(phi) = phi(a), 4}) is
a continuous homomorphism from $A$ onto a dense subgroup of $B$,
which is injective exactly when $\widehat{A}$ separates points on $A$.
This group $B$ is known as the \emph{Bohr compactification} of $A$.

        Of course, if $A$ is compact, then $B$ is the same as the
second dual of $A$, and the mapping from $a \in A$ to $\Psi_a$ is an
isomorphism from $A$ onto $B$.  Conversely, if the mapping from $a \in
A$ to $\Psi_a$ is a homeomorphism from $A$ onto $B$, then $A$ is
compact.  Suppose for the moment that the mapping from $a \in A$ to
$\Psi_a$ is a homeomorphism from $A$ onto its image in $B$, with
respect to the induced topology.  If $A$ is complete in the sense
discussed in Chapter \ref{completeness}, then it would follow that the
image of $A$ in $B$ is complete as well.  This would imply that the
image of $A$ in $B$ is a closed subgroup of $B$, and hence is equal to
$B$, because the image of $A$ is dense in $B$.  As before, this would
mean that $A$ should be compact under these conditions.  In
particular, the completeness property of $A$ used here holds when $A$
is locally compact, and when $A$ is metrizable and complete with
respect to a translation-invariant metric.

        Let $h$ be the usual homomorphism from $A$ into $B$, and consider
the mapping that sends a continuous complex-valued function $g$ on $B$
to $g \circ h \in C_b(A)$.  As in the previous section, this is a
homomorphism from the algebra $C(B)$ of all continuous complex-valued
functions on $B$ into $C_b(A)$ with respect to pointwise addition and
multiplication, and an isometry with respect to the corresponding
supremum norms, because $h(A)$ is dense in $B$.  It follows that the
image $\mathcal{A}$ of $C(B)$ in $C_b(A)$ under the mapping $g \mapsto
g \circ h$ is a closed subalgebra of $C_b(A)$ with respect to the
supremum norm, which is contained in the algebra $\mathcal{AP}(A)$ of
almost periodic functions on $A$, as before.  Because $B =
\widehat{(\widehat{A}_d)}$, the dual of $B$ can be identified with
$\widehat{A}_d$ in the usual way.  Thus characters on $A$ correspond
to characters on $B$ by construction, which implies that characters on
$A$ are contained in $\mathcal{A}$.  This shows that $\mathcal{A} =
\mathcal{AP}(A)$, because finite linear combinations of characters on
$A$ are dense in $\mathcal{AP}(A)$, as in Section \ref{almost
periodic functions}.

        Each element of $B$ determines a nonzero homomorphism from
$C(B)$ into ${\bf C}$, as in the previous section, and conversely
every nonzero complex homomorphism on $C(B)$ is of this form.
It follows that $B$ can also be identified with the set of nonzero
complex homomorphisms on $\mathcal{AP}(A)$, since $C(B)$ is
isomorphic to $\mathcal{AP}(A)$ as a Banach algebra.  One can show
more directly that nonzero complex homomorphisms on $\mathcal{AP}(A)$
determine characters on $\widehat{A}_d$, using the fact that characters
on $A$ are almost periodic.

\chapter{Completeness}
\label{completeness}

\section{Directed systems and nets}
\label{directed systems, nets}

        A nonempty partially-ordered set $(I, \prec)$ is said to be a 
\emph{directed system}\index{directed systems} if for every finite
collection $i_1, \ldots, i_n$ of elements of $I$ there is a $j \in I$
such that $i_k \prec j$ for each $k = 1, \ldots, n$.  A
\emph{net}\index{nets} $\{x_i\}_{i \in I}$ of elements of a set $X$
indexed by $I$ is basically the same as a function defined on $I$ with
values in $X$, which associates to each $i \in I$ an element $x_i$ of
$X$.  If $X$ is a topological space, then a net $\{x_i\}_{i \in I}$ of
elements of $X$ is said to \emph{converge}\index{convergent nets} to
$x \in X$ if for each open set $U$ in $X$ with $x \in U$ there is an
$i(U) \in I$ such that
\begin{equation}
\label{x_j in U}
        x_j \in U
\end{equation}
for every $j \in I$ that satisfies $i(U) \prec j$.  If $I = {\bf Z}_+$
with the standard ordering, then this reduces to ordinary convergence
of sequences.  It is easy to see that the limit of a convergent net in
a Hausdorff topological space is unique.

        Let $X$ be a topological space, let $p$ be an element of $X$,
and let $\mathcal{B}(p)$ be a local base for the topology of $X$ at
$p$.  Thus $\mathcal{B}(p)$ is a collection of open subsets of $X$,
each of which contains $p$ as an element, and if $V$ is any other open
set in $X$ that contains $p$ as an element, then there is an element
$U$ of $\mathcal{B}(p)$ such that $U \subseteq V$.  Also let $\prec$
be the partial ordering on $\mathcal{B}(p)$ corresponding to
reverse-inclusion, so that $U \prec W$ when $U, W \in \mathcal{B}(p)$
satisfy $W \subseteq U$.  If $U_1, \ldots, U_n$ are finitely many elements
of $\mathcal{B}(p)$, then $\bigcap_{j = 1}^n U_j$ is an open set in $X$
that contains $p$, and hence there is an element $W$ of $\mathcal{B}(p)$
such that $W \subseteq \bigcap_{j = 1}^n U_j$.  Equivalently, $W \subseteq
U_j$ for $j = 1, \ldots, n$, which is the same as saying that $U_j \prec W$
for each $j$, so that $\mathcal{B}(p)$ is a directed system with respect
to $\prec$.

        Suppose that for each $U \in \mathcal{B}(p)$, $x(U)$ is an element
of $U$.  Under these conditions, it is easy to see that $\{x(U)\}_{U
  \in \mathcal{B}(p)}$ converges to $p$ as a net of elements of $X$
indexed by $\mathcal{B}(p)$.  If $p$ is an element of the closure of
a set $E \subseteq X$, then one can choose $x(U)$ so that $x(U) \in E
\cap U$ for each $U \in \mathcal{B}(p)$.  Similarly, if there is a
countable local base for the topology of $X$ at $p$, and if $p$ is in
the closure of $E \subseteq X$, then there is a sequence of elements
of $E$ that converges to $p$.  Conversely, if there is a net of elements
of $E$ that converges to $p$ in $X$, then $p$ is contained in the closure
of $E$.

        Let $Y$ be another topological space, and let $f$ be a mapping
from $X$ into $Y$.  If $f$ is continuous at $p \in X$ and $\{x_i\}_{i
  \in I}$ is a net of elements of $X$ that converges to $p$, then it
is easy to see that $\{f(x_i)\}_{i \in I}$ converges to $f(p)$ as a
net of elements of $Y$.  Conversely, if $f$ is not continuous at $p$,
then there is an open set $V$ in $Y$ such that $f(p) \in V$ and
$f(U) \not\subseteq V$ for every open set $U$ in $X$ that contains
$p$ as an element.  Let $\mathcal{B}(p)$ be a local base for the
topology of $X$ at $p$, and for each $U \in \mathcal{B}(p)$ let $x(U)$
be an element of $U$ such that $f(x(U)) \not\in V$.  As before,
$\{x(U)\}_{U \in \mathcal{B}(p)}$ converges to $p$ as a net of elements
of $X$ indexed by $\mathcal{B}(p)$, where $\mathcal{B}(p)$ is ordered
by reverse-inclusion.  Of course, $\{f(x(U))\}_{U \in \mathcal{B}(p)}$
does not converge to $f(p)$ as a net of elements of $Y$, because
$f(x(U)) \not\in V$ for every $U \in \mathcal{B}(p)$.  Similarly,
if there is a countable local base for the topology of $X$ at $p$
and $f$ is not continuous at $p$, then there is a sequence $\{x_j\}_{j
  = 1}^\infty$ of elements of $X$ that converges to $p$ such that
$\{f(x_j)\}_{j = 1}^\infty$ does not converges to $f(p)$ in $Y$.

\section{Cauchy sequences and nets}
\label{cauchy sequences, nets}

        Remember that a sequence $\{x_j\}_{j = 1}^\infty$ of elements
of a metric space $(M, d(x, y))$ is said to be a \emph{Cauchy
  sequence}\index{Cauchy sequences} if for each $\epsilon > 0$
there is an $L(\epsilon) \ge 1$ such that
\begin{equation}
\label{d(x_j, x_l) < epsilon}
        d(x_j, x_l) < \epsilon
\end{equation}
for every $j, l \ge L(\epsilon)$.  In particular, it is well known and
easy to see that convergent sequences in $M$ are Cauchy sequences.  If
every Cauchy sequence in $M$ converges to an element of $M$, then $M$
is said to be complete as a metric space.

        Similarly, a sequence $\{x_j\}_{j = 1}^\infty$ of elements of a 
commutative topological group $A$ is said a \emph{Cauchy
  sequence}\index{Cauchy sequences} if for each open set $U$ in $A$
with $0 \in U$ there is an $L(U) \ge 1$ such that
\begin{equation}
\label{x_j - x_l in U}
        x_j - x_l \in U
\end{equation}
for every $j, l \ge L(U)$.  As before, one can check that convergent
sequences in $A$ are Cauchy sequences.  If $A$ is a commutative
topological group and $d(x, y)$ is a translation-invariant metric on
$A$ that determines the same topology on $A$, then it is easy to see
that a sequence $\{x_j\}_{j = 1}^\infty$ of elements of $A$ is a
Cauchy sequence in $A$ as a commutative topological group if and only
if $\{x_j\}_{j = 1}^\infty$ is a Cauchy sequence with respect to the
metric $d(x, y)$.  If every Cauchy sequence of elements of $A$
converges to an element of $A$, then we may say that $A$ is
\emph{sequentially complete}.\index{sequential completeness}
However, if $A$ does not have a countable local base for its topology
at $0$, then it may be appropriate to consider nets in $A$ as
well.

        A net $\{x_i\}_{i \in I}$ of elements of a commutative topological
group $A$ is said to be a \emph{Cauchy net}\index{Cauchy nets} if for each
open set $U \subseteq A$ with $0 \in A$ there is an $i(U) \in I$ such that
\begin{equation}
\label{x_j - x_l in U, 2}
        x_j - x_l \in U
\end{equation}
for every $j, l \in I$ with $i(U) \prec j$ and $i(U) \prec l$.  As
before, this reduces to the previous definition of a Cauchy sequence
in $A$ when $I = {\bf Z}_+$, and one can check that convergent nets in
$A$ are always Cauchy nets.  If every Cauchy net of elements of $A$
converges to an element of $A$, then we say that $A$ is
\emph{complete}\index{completeness} as a commutative topological
group.

        Let $\mathcal{B}_0$ be a local base for the topology of $A$ 
at $0$, and let $\prec_0$ be the partial ordering on $\mathcal{B}_0$
corresponding to reverse-inclusion, as in the previous section.  Thus
$(\mathcal{B}_0, \prec_0)$ is a directed system, as before.  If
$\{x_i\}_{i \in I}$ is a Cauchy net of elements of $A$ indexed by any
directed system $(I, \prec)$, then for each $U \in \mathcal{B}_0$
there is an $i(U) \in I$ such that (\ref{x_j - x_l in U, 2}) holds for
every $j, l \in I$ with $i(U) \prec j, l$.  Under these conditions,
one can check that $\{x_{i(U)}\}_{U \in \mathcal{B}_0}$ is also a
Cauchy net in $A$ as a net indexed by $\mathcal{B}_0$.  Similarly, if
$\{x_{i(U)}\}_{U \in \mathcal{B}_0}$ converges to an element $x$ of
$A$, then one can verify that $\{x_i\}_{i \in I}$ converges to $x$
too.  It follows that completeness of $A$ can be characterized
equivalently by the convergence of Cauchy nets indexed by
$\mathcal{B}_0$.  In particular, if there is a countable local base
for the topology of $A$ at $0$ and $A$ is sequentially complete, then
$A$ is complete.

\section{Completeness and compactness}
\label{completeness, compactness}

        Suppose that $K$ is a compact subset of a commutative
topological group $A$.  If $\{x_i\}_{i \in I}$ is a Cauchy
net of elements of $A$ such that $x_i \in K$ for every $i \in I$,
then we would like to show that $\{x_i\}_{i \in I}$ converges in $A$
to an element of $K$.  To do this, for each $i \in I$ let $E_i$
be the closure in $A$ of the set of $x_j$ with $j \in I$ and $i \prec j$,
i.e.,
\begin{equation}
\label{E_i = overline{{x_j : j in I, i prec j}}}
        E_i = \overline{\{x_j : j \in I, \, i \prec j\}}.
\end{equation}
Let $i_1, \ldots, i_n$ be finitely many elements of $I$, and let $l$
be an element of $I$ such that $i_r \prec l$ for each $r = 1, \ldots,
n$.  This implies that
\begin{equation}
\label{E_l subseteq bigcap_{r = 1}^n E_{i_r}}
        E_l \subseteq \bigcap_{r = 1}^n E_{i_r},
\end{equation}
and in particular that $\bigcap_{r = 1}^n E_{i_r} \ne \emptyset$.
Thus $E_i$, $i \in I$, is a family of nonempty closed subsets of $K$
with the ``finite intersection property'', and a well-known
reformulation of compactness implies that
\begin{equation}
\label{bigcap_{i in I} E_i ne emptyset}
        \bigcap_{i \in I} E_i \ne \emptyset.
\end{equation}
If $x \in \bigcap_{i \in I} E_i$, then one can check that $\{x_i\}_{i
  \in I}$ converges to $x$ in $A$, as desired.

        Suppose now that $A$ is a locally compact commutative topological
group, and let $\{x_i\}_{i \in I}$ be a Cauchy net of elements of $A$.
Let $U$ be an open set in $A$ such that $0 \in U$ and $\overline{U}$
is a compact set in $A$, and let $i(U)$ be an element of $I$ such that
(\ref{x_j - x_l in U, 2}) holds for every $j, l \in I$ with
$i(U) \prec j, l$.  This implies that
\begin{equation}
\label{x_j in x_{i(U)} + U subseteq x_{i(U)} + overline{U}}
        x_j \in x_{i(U)} + U \subseteq x_{i(U)} + \overline{U}
\end{equation}
for every $j \in I$ with $i(U) \prec j$, where $x_{i(U)} + \overline{U}$
is also a compact set in $A$.  The net consisting of $x_j$ with $j \in I$
such that $i(U) \prec j$ is a Cauchy net of elements of $A$ contained
in the compact set $x_{i(U)} + \overline{U}$, which converges to an
element of $x_{i(U)} + \overline{U}$ by the argument in the previous 
paragraph.  It follows that the original net $\{x_i\}_{i \in I}$
converges to the same element of $x_{i(U)} + \overline{U}$, and hence
that $A$ is complete.

        Let $B$ be a commutative topological group, and let $A$ be a 
subgroup of $B$.  If $x \in B$ is an element of the closure of $A$ in
$B$, then there is a net $\{x_i\}_{i \in I}$ of elements of $A$ that
converges to $x$ in $B$, as in Section \ref{directed systems, nets}.  
As before, $\{x_i\}_{i \in I}$ is a Cauchy net in $B$, and in fact
$\{x_i\}_{i \in I}$ is a Cauchy net in $A$ as a topological group with
the topology induced by the one on $B$.  If $A$ is complete, then
$\{x_i\}_{i \in I}$ also converges to an element $y$ of $A$.  Because
topological groups are Hausdorff, it follows that $x = y$ under these
conditions, and hence that $x \in A$.  This shows that $A$ is a closed
subgroup of $B$ when $A$ is complete as a commutative topological
group with respect to the topology induced by the one on $B$.  If
there is a countable local base for the topology of $B$ at $0$, then
one can take $I = {\bf Z}_+$ with the standard ordering, and it
suffices that $A$ be sequentially complete.  In particular, if $A$
is locally compact with respect to the topology induced by the one
on $B$, then $A$ is complete, as in the preceding paragraph,
and hence $A$ is a closed subset of $B$.

        As another application, suppose that $A$ is a dense subgroup
of a commutative topological group $B$, and that $h$ is a continuous
homomorphism from $A$ into $C$, with respect to the topology on $A$
induced by the one on $B$.  Let $x$ be an element of $B$, and let
$\{x_i\}_{i \in I}$ be a net of elements of $A$ that converges to $x$
in $B$.  As before, $\{x_i\}_{i \in I}$ is a Cauchy net of elements of
$A$ with respect to the topology induced on $A$ by the one on $B$.
Under these conditions, it is easy to see that $\{h(x_i)\}_{i \in I}$
is a Cauchy net in $C$.  If $C$ is complete, then it follows that
$\{h(x_i)\}_{i \in I}$ converges to an element of $C$.  In particular,
this holds when $C$ is locally compact.  In this case, we would like to
define $h$ at $x$ to be the limit of $\{h(x_i)\}_{i \in I}$ in $C$.
To do this, one should check this value of $h(x)$ does not depend on
the choice of net $\{x_i\}_{i \in I}$ of elements of $A$ converging
to $x$.  This is not too difficult, and it is more pleasant when the
nets are indexed by the same directed system.  One might as well use
nets indexed by a local base $\mathcal{B}_0$ for the topology of $B$
at $0$, which determines a local base for the topology of $B$ at any
point, by translation.  It is easy to see that this extension of $h$
is a homomorphism from $B$ into $C$ under these conditions.  One can
also check that this extension is a continuous mapping from $B$ into
$C$, using the fact that $C$ is regular as a topological space.
If there is a countable local base for the topology of $B$ at $0$,
then one can simply use sequences instead of nets, and it suffices
for $C$ to be sequentially complete.

\section{Continuous functions}
\label{continuous functions, 3}

        Let $X$ be a topological space, and let $C(X)$ be the space of
continuous complex-valued functions on $X$.  Thus $C(X)$ is a vector
space over the complex numbers with respect to pointwise addition and
scalar multiplication, and
\begin{equation}
\label{||f||_K = sup_{x in K} |f(x)|, 2}
        \|f\|_K = \sup_{x \in K} |f(x)|
\end{equation}
is a seminorm on $C(X)$ for each nonempty compact subset $K$ of $X$.
As usual, the collection of these seminorms defines a topology on
$C(X)$ that makes $C(X)$ into a locally convex topological vector
space, and a commutative topological group with respect to addition in
particular.  Remember that
\begin{equation}
\label{||f g||_K le ||f||_K ||g||_K, 2}
        \|f \, g\|_K \le \|f\|_K \, \|g\|_K
\end{equation}
for every $f, g \in C(X)$ and nonempty compact set $K \subseteq X$.
This implies that pointwise multiplication of functions on $X$ defines
a continuous mapping from $C(X) \times C(X)$ into $C(X)$, so that
$C(X)$ is a commutative topological algebra.

        Suppose that $\{f_i\}_{i \in I}$ is a Cauchy net of continuous
complex-valued functions on $X$, as a net of elements of $C(X)$ as a 
commutative topological group.  This means that for each nonempty compact
set $K \subseteq X$ and $\epsilon > 0$ there is an $i(K, \epsilon) \in I$
such that
\begin{equation}
\label{||f_j - f_l||_K < epsilon}
        \|f_j - f_l\|_K < \epsilon
\end{equation}
for every $j, l \in I$ that satisfy $i(K, \epsilon) \prec j, l$.
In particular, $\{f_i(x)\}_{i \in I}$ is a Cauchy net of complex numbers
for each $x \in X$, since one can take $K = \{x\}$.  It follows that
$\{f_i(x)\}_{i \in I}$ converges to a complex number $f(x)$ for each
$x \in X$, because ${\bf C}$ is complete.

        Combining this with (\ref{||f_j - f_l||_K < epsilon}), we get that
for each nonempty compact set $K \subseteq X$, $\epsilon > 0$, and $x \in K$,
\begin{equation}
\label{|f_j(x) - f(x)| le epsilon}
        |f_j(x) - f(x)| \le \epsilon
\end{equation}
for every $j \in I$ such that $i(K, \epsilon) \prec j$.  This
basically says that $\{f_i\}_{i \in I}$ converges to $f$ uniformly on
compact subsets of $K$, which implies that the restriction of $f$ to
any nonempty compact set $K \subseteq X$ is continuous, by standard
arguments.  If $f$ is continuous on $X$, then it follows that
$\{f_i\}_{i \in I}$ converges to $f$ in $C(X)$.  In particular,
if $X$ is locally compact, then the continuity of $f$ on compact
subsets of $X$ implies that $f$ is continuous on $X$.

        Alternatively, let $\{p_r\}_{r = 1}^\infty$ be a sequence of
elements of $X$ that converges to another element $p$ of $X$, and let
$K$ be the subset of $X$ consisting of the $p_r$'s with $r \in {\bf
  Z}_+$ and $p$.  It is easy to see that $K$ is compact under these
conditions, so that the restriction of $f$ to $K$ is continuous, as
before.  This implies that $\{f(p_r)\}_{r = 1}^\infty$ converges to
$f(p)$ as a sequence of complex numbers, which is to say that $f$ is
sequentially continuous at every point in $X$.  If there is a
countable local base for the topology of $X$ at each point in $X$,
then it follows that $f$ is continuous on $X$.

        Let $A$ be a commutative topological group, and let $\widehat{A}$
be the dual group of continuous homomorphisms from $A$ into the
multiplicative group ${\bf T}$ of complex numbers with modulus $1$, as
usual.  It is easy to see that $\widehat{A}$ is a closed set in $C(A)$
with respect to the topology described above with $X = A$, and that
$\widehat{A}$ is a topological group with respect to the topology
induced by the one on $C(A)$.  If $\{\phi_i\}_{i \in I}$ is a net of
elements of $\widehat{A}$, then one can also check that $\{\phi_i\}_{i \in I}$
is a Cauchy net in $\widehat{A}$ as a commutative topological group with
respect to multiplication if and only if $\{\phi_i\}_{i \in I}$ is a Cauchy
net in $C(A)$ as a commutative topological group with respect to addition.
More precisely, this uses the fact that
\begin{equation}
\label{|z w^{-1} - 1| = |z - w|}
        |z \, w^{-1} - 1| = |z - w|
\end{equation}
for any two complex numbers $z$, $w$ with modulus equal to $1$.  It
follows that $\widehat{A}$ is complete as a commutative topological
group with respect to multiplication when $C(A)$ is complete as a
commutative topological group with respect to addition.  In
particular, this holds when $A$ is locally compact, and when there is
a countable local base for the topology of $A$ at $0$, as before.  Of
course, the latter condition implies that there is a countable local
base for the topology of $A$ at every point, using translations.

        Now let $V$ be any topological vector space over the real numbers,
and let $V^*$ be the corresponding dual space of continuous linear
functionals on $V$.  Thus $V^*$ may be considered as a linear subspace
of the space of all continuous real-valued functions on $V$, and it is
easy to see that $V^*$ is a closed set with respect to the topology
defined by the supremum seminorms associated to nonempty compact
subsets of $V$.  The topology on $V^*$ determined by the restriction
of the supremum seminorms associated to nonempty compact subsets of
$V$ to $V^*$ is the same as the one induced by the corresponding
topology on $C(V)$.  If $C(V)$ is complete as a commutative topological
group with respect to addition, then it follows that $V^*$ is complete
as a commutative topological group with respect to addition too.
As before, this holds in particular when $V$ is locally compact,
and when there is a countable local base for the topology of $V$
at $0$.

        However, it is well known that $V$ is locally compact if and
only if it is finite-dimensional.  Moreover, if $V$ has
finite-dimension $n$, then $V$ is isomorphic as a vector space to
${\bf R}^n$, and any such isomorphism is a homeomorphism as well.  Of
course, $V^*$ is also isomorphic to ${\bf R}^n$ in this case.  If $V$
is any vector space over the real numbers with a norm, then there is a
natural dual norm on $V^*$, as in Section \ref{bounded linear
  functionals}.  As before, $V^*$ is also complete with respect to the
dual norm, but the dual norm is stronger than the supremum seminorms
on $V^*$ corresponding to compact subsets of $V$, except when $V$ is
locally compact and hence finite-dimensional.  Note that there are
classes of topological vector spaces $V$ defined using inductive
limits, for which there is not a countable local base for the topology
of $V$ at $0$, but for which continuity of linear functionals on $V$
can be characterized in terms of sequential continuity.  This still
implies that the dual space $V^*$ is complete with respect to the
topology determined by the collection of supremum of seminorms
associated to nonempty compact subsets of $V$, for the same reasons
as before.

        Let $V$ be a topological vector space over the real numbers again,
and let $\lambda$ be a continuous linear functional on $V$.  As in
Section \ref{homomorphisms into T},
\begin{equation}
\label{phi(v) = exp (i lambda(v)), 2}
        \phi(v) = \exp (i \, \lambda(v))
\end{equation}
is a continuous homomorphism from $V$ as a commutative topological
group with respect to addition into ${\bf T}$, and every continuous
homomorphism from $V$ into ${\bf T}$ is of this form.  The mapping
from $\lambda$ to $\phi$ defines an isomorphism from $V^*$ as a
commutative group with respect to addition onto the dual group
$\widehat{V}$ of $V$ as a commutative topological group with respect
to addition.  Let us check that this mapping is a homeomorphism 
with respect to the appropriate topologies on $V^*$ and $\widehat{V}$.
More precisely, this means the topology on $V^*$ determined by the
collection of supremum seminorms associated to nonempty compact
subsets of $V$, and the topology induced on $\widehat{V}$ by the
one on $C(A)$ determined by the supremum seminorms associated to
nonempty compact subsets of $V$.

        It is easy to see that $\lambda \mapsto \phi$ defines a
continuous mapping from $V^*$ onto $\widehat{V}$.  Basically, the main
point is that if $\lambda$ and $\lambda'$ are continuous linear
functionals on $V$ that are uniformly close to each other on a compact
set $K \subseteq V$, then $\phi = \exp (i \, \lambda)$ and $\phi' =
\exp (i \, \lambda')$ are also uniformly close on $K$, because of the
continuity of the exponential function.  To show that $\lambda \mapsto
\phi$ is a homeomorphism is a bit more complicated, because $\exp (i \, t)$
is a local homeomorphism from ${\bf R}$ onto ${\bf T}$, and not a
homeomorphism.  In particular, $\phi$ and $\phi'$ are close when
$\lambda$ and $\lambda'$ are close modulo $2 \pi$.  To deal with this,
let $K$ be a nonempty compact set in $V$, and consider
\begin{equation}
\label{widetilde{K} = {t v : v in K, t in {bf R}, and 0 le t le 1}}
        \widetilde{K} = \{t \, v : v \in K, \, t \in {\bf R}, \hbox{ and }
                                                              0 \le t \le 1\}.
\end{equation}
Remember that scalar multiplication on $V$ defines a continuous
mapping from ${\bf R} \times V$ into $V$, because $V$ is a topological
vector space.  If $K \subseteq V$ is compact, then $[0, 1] \times K$
is compact in ${\bf R} \times V$ with respect to the product topology,
which implies that $\widetilde{K}$ is compact as well, since
$\widetilde{K}$ is the image of $[0, 1] \times K$ under scalar
multiplication as a mapping from ${\bf R} \times V$ into $V$.  If
$\phi$ and $\phi'$ are uniformly close on $\widetilde{K}$, then one
can check that $\lambda$ and $\lambda'$ are uniformly close on
$\widetilde{K}$, and hence on $K$, using the fact that $\lambda(0) =
\lambda'(0) = 0$.  This permits one to show that $\lambda \mapsto
\phi$ is a homeomorphism from $V^*$ onto $\widehat{V}$, as desired.

\section{Filters}
\label{filters}

        A nonempty collection $\mathcal{F}$ of nonempty subsets of a
set $X$ is said to be a \emph{filter}\index{filters} on $X$ if it
satifies the following two additional conditions.  First,
\begin{equation}
\label{A cap B in mathcal{F}}
        A \cap B \in \mathcal{F}
\end{equation}
for every $A, B \in \mathcal{F}$.  Second, if $A \in \mathcal{F}$,
$E \subseteq X$, and $A \subseteq E$, then
\begin{equation}
\label{E in mathcal{F}}
        E \in \mathcal{F}.
\end{equation}
Similarly, a nonempty collection $\mathcal{F}^*$ of nonempty subsets
of $X$ is said to be a \emph{pre-filter}\index{pre-filters} on $X$ if for
every $A, B \in \mathcal{F}^*$ there is a $C \in \mathcal{F}^*$ such that
\begin{equation}
\label{C subseteq A cap B}
        C \subseteq A \cap B.
\end{equation}
Thus every filter $\mathcal{F}$ on $X$ is a pre-filter, with $C = A
\cap B$.  Conversely, suppose that $\mathcal{F}^*$ is a pre-filter on
$X$, and consider
\begin{equation}
\label{mathcal{F} = {E subseteq X : A subseteq E for some A in mathcal{F}^*}}
        \mathcal{F} = \{E \subseteq X : A \subseteq E \hbox{ for some }
                                                      A \in \mathcal{F}^*\}.
\end{equation}
It is easy to see that $\mathcal{F}$ is a filter on $X$ under these
conditions, which is the filter generated by $\mathcal{F}^*$.

        Suppose now that $X$ is a topological space.  A filter $\mathcal{F}$
on $X$ is said to \emph{converge}\index{convergent filters} to a point
$p \in X$ if for every open set $U$ in $X$ with $p \in U$, we have that
$U \in \mathcal{F}$.  Note that the limit of a convergent filter on $X$
is unique when $X$ is Hausdorff.  If $p$ is any element of a set $X$,
then the collection $\mathcal{F}(p)$ of all subsets $E$ of $X$ such that
$p \in E$ is a filter on $X$.  If $X$ is a topological space, then this
filter $\mathcal{F}(p)$ converges to $p$ on $X$.  If $X$ is equipped with
the discrete topology and $\mathcal{F}$ is a filter on $X$ that converges
to $p$, then $\mathcal{F} = \mathcal{F}(p)$.  If $X$ is any topological
space and $p \in X$, then the collection of all open subsets of $X$
that contain $p$ as an element is a pre-filter on $X$, and the filter
on $X$ generated by this pre-filter converges to $p$.

        Let $\mathcal{F}$ be a filter on a topological space $X$ that
converges to a point $p \in X$, and suppose that $E \in \mathcal{F}$.
If $U$ is an open set in $X$ that contains $p$ as an element, then
$U \in \mathcal{F}$, and hence $U \cap E \in \mathcal{F}$.  This implies
that $U \cap E \ne \emptyset$, so that $p$ is an element of the closure
$\overline{E}$ of $E$ in $X$.  Conversely, suppose that $p \in \overline{E}$,
so that $U \cap E \ne \emptyset$ for every open set $U$ in $X$ with $p \in U$.
It is easy to see that the collection of subsets of $X$ of the form $U \cap E$
for some open set $U \subseteq X$ that contains $p$ as an element is a
pre-filter on $X$, and that the filter on $X$ generated by this pre-filter
converges to $p$ and contains $E$ as an element.

        Now let $(I, \prec)$ be a directed system, and let $\{x_i\}_{i \in I}$
be a net of elements of a set $X$ indexed by $I$.  Let $i$ be an element of
$I$, and put
\begin{equation}
\label{B_i = {x_l : l in I, i prec l}}
        B_i = \{x_l : l \in I, \, i \prec l\}.
\end{equation}
If $i, j \in I$, then there is a $k \in I$ such that $i, j \prec k$,
because $I$ is a directed system, and hence
\begin{equation}
\label{B_k subseteq B_i cap B_j}
        B_k \subseteq B_i \cap B_j.
\end{equation}
This implies that the collection of subsets of $X$ of the form $B_i$
for some $i \in I$ is a pre-filter on $X$.  If $X$ is a topological
space, then it is easy to see that the filter on $X$ generated by this
pre-filter converges to a point $p \in X$ if and only if $\{x_i\}_{i
  \in I}$ converges to $p$ as a net of elements of $X$.

        Let $\mathcal{F}$ be a filter on a set $X$, and suppose that
$I \subseteq \mathcal{F}$ is a pre-filter on $X$ that generates $\mathcal{F}$.
Let $\prec$ be the partial ordering on $I$ defined by
reverse-inclusion, so that $A \prec B$ when $A, B \in I$ and $B
\subseteq A$.  Observe that $(I, \prec)$ is a directed system, because
$I$ is a pre-filter on $X$.  Let $\{x_B\}_{B \in I}$ be a net of
elements of $X$ indexed by $I$ such that $x_B \in B$ for every $B \in
I$.  If $X$ is a topological space, and if $\mathcal{F}$ converges as
a filter on $X$ to a point $p \in X$, then every net of elements of
$X$ indexed by $I$ of this type also converges to $p$, because
$\mathcal{F}$ is generated by $I$.  If $\mathcal{F}$ does not converge
to $p$, then there is an open set $U$ in $X$ such that $p \in U$ and
$U \not\in \mathcal{F}$.  This implies that $B \not\subseteq U$ for
each $U \in I$, so that for each $B \in I$ there is an $x_B \in B
\backslash U$.  In this case, $\{x_B\}_{B \in I}$ is a net of elements
of $X$ associated to $\mathcal{F}$ as before, and $\{x_B\}_{B \in I}$
does not converge to $p$ in $X$.

        Let $X$ and $Y$ be sets, and let $f$ be a mapping from $X$
into $Y$.  If $\mathcal{F}$ is a filter on $X$, then
\begin{equation}
\label{f_*(mathcal{F}) = {E subseteq Y : f^{-1}(E)}}
        f_*(\mathcal{F}) = \{E \subseteq Y : f^{-1}(E)\}
\end{equation}
is a filter on $Y$.  Equivalently, the collection of subsets of $Y$ of
the form $f(B)$ for some $B \in \mathcal{F}$ is a pre-filter on $Y$,
and $f_*(\mathcal{F})$ is the same as the filter on $Y$ generated by
this pre-filter.  If $X$ and $Y$ are topological spaces, $\mathcal{F}$
is a filter on $X$ that converges to a point $p \in X$, and $f : X \to
Y$ is continuous at $p$, then it is easy to see that
$f_*(\mathcal{F})$ converges to $f(p)$ on $Y$.  Conversely, if
$\mathcal{F}$ is the filter on $X$ generated by the pre-filter of all
open sets in $X$ that contain $p$ as an element, and if
$f_*(\mathcal{F})$ converges to $f(p)$ in $Y$, then $f$ is continuous
at $p$.

        A filter $\mathcal{F}$ on a commutative topological group $A$ 
is said to be a \emph{Cauchy filter}\index{Cauchy filters} if for each
open set $U$ in $A$ with $0 \in U$ there is an $E \in \mathcal{F}$
such that
\begin{equation}
\label{E - E = {x - y : x, y in E} subseteq U}
        E - E = \{x - y : x, y \in E\} \subseteq U.
\end{equation}
If $\mathcal{F}$ converges to a point $p \in X$, then it is easy to
see that $\mathcal{F}$ is a Cauchy filter on $A$.  Suppose that
$\{x_i\}_{i \in I}$ is a net of elements of $A$, and let $\mathcal{F}$
be the filter on $X$ generated by the sets $B_i$ with $i \in I$ as in
(\ref{B_i = {x_l : l in I, i prec l}}).  Observe that $\{x_i\}_{i \in
  I}$ is a Cauchy net in $A$ if and only if $\mathcal{F}$ is a Cauchy
filter.  Now let $\mathcal{F}$ be any filter on $A$, let $I \subseteq
\mathcal{F}$ be a pre-filter that generates $\mathcal{F}$, and let
$\{x_B\}_{B \in I}$ be a net of elements of $A$ indexed by $I$ such
that $x_B \in B$ for each $B \in I$.  If $\mathcal{F}$ is a Cauchy
filter on $A$, then $\{x_B\}_{B \in I}$ is a Cauchy net in $A$.  If
$\mathcal{F}$ is a Cauchy filter on $A$ and $\{x_B\}_{B \in I}$
converges to a point $p \in A$, then one can check that $\mathcal{F}$
converges to $p$ as well.  It follows that $A$ is
complete\index{completeness} in the sense that every Cauchy net of
elements of $A$ converges to an element of $A$ if and only if every
Cauchy filter on $A$ converges to an element of $A$.

        Let $h$ be a continuous homomorphism from $A$ into another
commutative topological group $C$.  If $\mathcal{F}$ is a Cauchy
filter on $A$, then it is easy to see that $h_*(\mathcal{F})$ is a
Cauchy filter on $C$.  In particular, suppose that $A$ is a dense
subgroup of a topological group $B$, and let $x$ be any element of $B$.
Let $\mathcal{F}^*$ be the collection of subsets of $A$ of the form
$U \cap A$, where $U$ is an open set in $B$ that contains $x$ as an
element.  It is easy to see that $\mathcal{F}^*$ is a pre-filter on
$A$, and that the filter $\mathcal{F}$ on $A$ generated by $\mathcal{F}^*$
is a Cauchy filter on $A$.  Of course, $\mathcal{F}^*$ may also be considered
as a pre-filter on $B$, that generates a filter on $B$ that converges to
$x$, as before.  At any rate, $h_*(\mathcal{F})$ is a Cauchy filter on $C$,
which converges to an element of $C$ when $C$ is complete.  This gives 
another way to look at the extension of $h$ to a continuous homomorphism
from $B$ into $C$ when $C$ is complete, as discussed at the end of
Section \ref{completeness, compactness}.

\section{Refinements}
\label{refinements}

        Let $\mathcal{F}$ be a filter on a set $X$.  A filter $\mathcal{F}'$
on $X$ is said to be a \emph{refinement} of $\mathcal{F}$ if
$\mathcal{F} \subseteq \mathcal{F}'$ as collections of subsets of $X$.
In particular, $\mathcal{F}$ may be considered as a refinement of
itself.  If $X$ is a topological space and $\mathcal{F}$ converges
to a point $p \in X$, then every refinement of $\mathcal{F}$ converges
to $p$ too.

        Let $\mathcal{F}$ be a filter on a topological space $X$, and
let $p$ be an element of $X$.  If there is a refinement $\mathcal{F}'$
of $\mathcal{F}$ that converges to $p$, then $p \in \overline{E}$ for
every $E \in \mathcal{F}'$, as in the previous section.  In
particular, $p \in \overline{E}$ for every $E \in \mathcal{F}$.
Conversely, suppose that $p \in \overline{E}$ for every $E \in 
\mathcal{F}$.  Let $\mathcal{F}^*$ be the collection of subsets
of $X$ of the form $U \cap E$, where $U$ is an open set in $X$
that contains $p$ as an element, and $E \in \mathcal{F}$.
It is easy to see that $\mathcal{F}^*$ is a pre-filter on $X$,
and that the filter $\mathcal{F}'$ generated by $\mathcal{F}$
is a refinement of $\mathcal{F}$ that converges to $p$.
This shows that there is a refinement of $\mathcal{F}$ that
converges to $p \in X$ if and only if $p \in \bigcap_{E \in \mathcal{F}}
\overline{E}$.

        Let $X$ be a topological space again, and let $K$ be a subset of
$X$.  Also let $I$ be a nonempty set, and suppose that $E_i$ is a closed
set in $X$ for each $i \in I$.  We say that $\{E_i\}_{i \in I}$
satisfies the \emph{finite intersection property} with respect to $K$ if
\begin{equation}
\label{K cap E_{i_1} cap cdots cap E_{i_n} ne emptyset}
        K \cap E_{i_1} \cap \cdots \cap E_{i_n} \ne \emptyset
\end{equation}
for every finite collection of indices $i_1, \ldots, i_n$ in $I$.  If
$\{E_i\}_{i \in I}$ has the finite intersection property with respect
to $K$ and $K$ is compact, then
\begin{equation}
\label{K cap (bigcap_{i in I} E_i) ne emptyset}
        K \cap \Big(\bigcap_{i \in I} E_i\Big) \ne \emptyset.
\end{equation}
To see this, suppose for the sake of a contradiction that $K \cap
\Big(\bigcap_{i \in I} E_i\Big) = \emptyset$, which is the same as
saying that $K \subseteq \bigcup_{i \in I} (X \backslash E_i)$.  Thus
$\{X \backslash E_i\}_{i \in I}$ is an open covering of $K$ in $X$.
If $K$ is compact, then there are finitely many indices $i_1, \ldots,
i_n \in I$ such that $K \subseteq \bigcup_{l = 1}^n (X \backslash
E_{i_l})$, contradicting (\ref{K cap E_{i_1} cap cdots cap E_{i_n} ne
  emptyset}).  Conversely, if (\ref{K cap (bigcap_{i in I} E_i) ne
  emptyset}) holds for every collection $\{E_i\}_{i \in I}$ of closed
subsets of $X$ with the finite intersection property with respect to
$K$, then $K$ is compact.  Indeed, if $\{U_i\}_{i \in I}$ is an open
covering of $K$ in $X$ for which there is no finite subcovering, then
the collection of closed sets $E_i = X \backslash U_i$ with $i \in I$
has the finite intersection property with respect to $K$, but (\ref{K
  cap (bigcap_{i in I} E_i) ne emptyset}) does not hold.

        Suppose that $K \subseteq X$ is compact, and that $\mathcal{F}$
is a filter on $X$ that contains $K$ as an element.  If $E_1, \ldots, E_n$
are finitely many elements of $\mathcal{F}$, then
\begin{equation}
\label{K cap E_1 cap cdots cap E_n}
        K \cap E_1 \cap \cdots \cap E_n
\end{equation}
is also an element of $\mathcal{F}$, and hence is nonempty.  In particular,
\begin{equation}
\label{K cap overline{E_1} cap cdots cap overline{E_n} ne emptyset}
        K \cap \overline{E_1} \cap \cdots \cap \overline{E_n} \ne \emptyset,
\end{equation}
so that the collection of closed subsets of $X$ of the form $\overline{E}$
for some $E \in \mathcal{F}$ has the finite intersection property with
respect to $K$.  This implies that
\begin{equation}
\label{K cap (bigcap_{E in mathcal{F}} overline{E}) ne emptyset}
  K \cap \Big(\bigcap_{E \in \mathcal{F}} \overline{E}\Big) \ne \emptyset,
\end{equation}
because $K$ is compact, as in the previous paragraph.  It follows that
there is a refinement $\mathcal{F}'$ of $\mathcal{F}$ that converges
to an element of $K$, as discussed earlier.

        Conversely, suppose that $K \subseteq X$ has the property that
every filter $\mathcal{F}$ on $X$ that contains $K$ as an element has
a refinement that converges to an element of $K$, and let us show that
$K$ is compact.  To do this, let $\{E_i\}_{i \in I}$ be an arbitrary
collection of closed subsets of $X$ with the finite intersection
property with respect to $K$, and let us check that (\ref{K cap
  (bigcap_{i in I} E_i) ne emptyset}) holds.  Let $\mathcal{F}^*$
be the collection of subsets of $X$ of the form
\begin{equation}
\label{K cap E_{i_1} cap cdots cap E_{i_n}}
        K \cap E_{i_1} \cap \cdots \cap E_{i_n},
\end{equation}
where $i_1, \ldots, i_n$ are finitely many indices in $I$.  Because
$\{E_i\}_{i \in I}$ has the finite intersection property with respect
to $K$, subsets of $X$ of the form (\ref{K cap E_{i_1} cap cdots cap
  E_{i_n}}) are nonempty, and it is easy to see that $\mathcal{F}^*$
is a pre-filter on $X$.  Let $\mathcal{F}$ be the filter on $X$
generated by $\mathcal{F}^*$, and observe that $K \in \mathcal{F}$, by
construction.  Hence there is a refinement $\mathcal{F}'$ of
$\mathcal{F}$ that converges to an element $p$ of $K$, by hypothesis.
It follows that $p \in \overline{E}$ for every $E \in \mathcal{F}'$
and thus every $E \in \mathcal{F}$, as discussed earlier.  In particular,
$p \in E_i$ for each $i \in I$, because $E_i$ is a closed set in $X$,
so that (\ref{K cap (bigcap_{i in I} E_i) ne emptyset}) holds, as desired.

        Now let $A$ be a commutative topological group, and let $\mathcal{F}$
be a Cauchy filter on $A$.  If $\mathcal{F}'$ is a refinement of
$\mathcal{F}$ that converges to a point $p \in X$, then $\mathcal{F}$
converges to $p$ as well.  To see this, let $W$ be an open set in $A$
that contains $0$, and let $U$ and $V$ be open subsets of $A$ that
contain $0$ and satisfy
\begin{equation}
\label{U + V subseteq W}
        U + V \subseteq W.
\end{equation}
Of course, the existence of $U$ and $V$ follows from the continuity of
addition at $0$ on $A$.  Because $\mathcal{F}$ is a Cauchy filter on $A$,
there is an $E \in \mathcal{F}$ such that
\begin{equation}
\label{E - E subseteq U}
        E - E \subseteq U.
\end{equation}
Similarly, $p + V \in \mathcal{F}'$, because $\mathcal{F}'$ converges
to $p$ on $A$.  In particular, $p + V \in \mathcal{F}$, so that $E
\cap (p + V) \in \mathcal{F}$ too, and hence $E \cap (p + V) \ne
\emptyset$.  If $x \in E \cap (p + V)$, then we get that
\begin{equation}
\label{E subseteq x + U subseteq p + U + V subseteq p + W}
        E \subseteq x + U \subseteq p + U + V \subseteq p + W.
\end{equation}
More precisely, this uses (\ref{E - E subseteq U}) and the fact that
$x \in E$ in the first step, the fact that $x \in p + V$ in the second
step, and (\ref{U + V subseteq W}) in the third step.  It follows that
$p + W \in \mathcal{F}$, so that $\mathcal{F}$ converges to $p$ on $X$,
as desired.

        Let $\mathcal{F}$ be a Cauchy filter on $A$ again, and suppose
that $\mathcal{F}$ contains a compact set $K \subseteq X$ as an
element.  As before, there is a refinement $\mathcal{F}'$ of
$\mathcal{F}$ that converges to an element $p$ of $K$, which implies
that $\mathcal{F}$ converges to $p$ as well.  Suppose now that $A$ is
locally compact, and let $U$ be an open set in $A$ such that $0 \in U$
and $\overline{U}$ is compact.  If $\mathcal{F}$ is a Cauchy filter on $A$,
then there is an $E \in \mathcal{F}$ such that $E - E \subseteq U$,
which implies that
\begin{equation}
\label{E subseteq x + U subseteq x + overline{U}}
        E \subseteq x + U \subseteq x + \overline{U}
\end{equation}
for every $x \in E$.  It follows that $x + \overline{U}$ is a compact
set in $X$ which is also an element of $\mathcal{F}$ for each $x \in
E$, and hence that $\mathcal{F}$ converges to an element of $A$ under
these conditions.

\section{Ultrafilters}
\label{ultrafilters}

        A filter $\mathcal{F}$ on a set $X$ is said to be an
\emph{ultrafilter} if it is maximal with respect to refinements.
More precisely, this means that if $\mathcal{F}'$ is a filter on
$X$ which is a refinement of $\mathcal{F}$, then $\mathcal{F}' =
\mathcal{F}$.  If $p \in X$ and $\mathcal{F}(p)$ is the collection
of subsets $E$ of $X$ such that $p \in E$, then it is easy to see
that $\mathcal{F}(p)$ is an ultrafilter on $X$.  If $\mathcal{F}$
is any filter on $X$, then there is a refinement of $\mathcal{F}$
which is an ultrafilter, by standard arguments based on the axiom
of choice through Zorn's lemma or the Hausdorff maximality principle.

        Suppose that $X$ is a topological space, and that $K \subseteq X$
is compact.  If $\mathcal{F}$ is a filter on $X$ that contains $K$ as
an element, then there is a refinement $\mathcal{F}'$ of $\mathcal{F}$
that converges to an element of $K$, as in the previous section.  If
$\mathcal{F}$ is an ultrafilter on $X$, then it follows that
$\mathcal{F}$ converges to an element of $K$.  Conversely, suppose that
$K \subseteq X$ has the property that every ultrafilter on $X$ that
contains $K$ as an element converges to an element of $K$.  Let $\mathcal{F}$
be any filter on $X$ that contains $K$ as an element, and let $\mathcal{F}'$
be a refinement of $\mathcal{F}$ which is an ultrafilter.  In particular,
$\mathcal{F}'$ contains $K$ as an element, and so $\mathcal{F}'$ converges
to an element of $K$, by hypothesis.  Thus $\mathcal{F}$ has a refinement
that converges to an element of $K$, which implies that $K$ is compact,
as in the previous section.

        Let $\mathcal{F}$ be an filter on a set $X$, and suppose that
$B \subseteq X$ has the property that $B \cap E \ne \emptyset$ for every
$E \in \mathcal{F}$.  Let $\mathcal{F}^*$ be the collection of subsets
of $X$ of the form $B \cap E$ for some $E \in \mathcal{F}$.  It is
easy to see that $\mathcal{F}^*$ is a pre-filter on $X$, and that the
filter $\mathcal{F}'$ generated by $\mathcal{F}^*$ is a refinement of
$\mathcal{F}$.  If $\mathcal{F}$ is an ultrafilter on $X$, then it
follows that $\mathcal{F} = \mathcal{F}'$.  In particular, this
implies that $B \in \mathcal{F}$ under these conditions.  Conversely,
suppose that $\mathcal{F}$ is a filter on a set $X$ with the property
that if $B \subseteq X$ satisfies $B \cap E \ne \emptyset$ for every
$E \in \mathcal{F}$, then $B \in \mathcal{F}$.  Let us check that
$\mathcal{F}$ is an ultrafilter on $X$ in this case.  Let
$\mathcal{F}'$ be a refinement of $\mathcal{F}$, and let $B$ be any
element of $\mathcal{F}'$.  If $E \in \mathcal{F}$, then $E \in
\mathcal{F}'$, which implies that $B \cap E \in \mathcal{F}'$, and
hence that $B \cap E \ne \emptyset$.  By hypothesis, it follows that
$B \in \mathcal{F}$, so that $\mathcal{F} = \mathcal{F}'$, as desired.

        Let $\mathcal{F}$ be an ultrafilter on a set $X$, and let $B$
be any subset of $X$.  If $B \cap E \ne \emptyset$ for every $E \in
\mathcal{F}$, then $B \in \mathcal{F}$, as in the preceding paragraph.
Otherwise, if $B \cap E = \emptyset$ for some $E \in \mathcal{F}$,
then $E \subseteq X \backslash B$, which implies that $X \backslash B
\in \mathcal{F}$, because $\mathcal{F}$ is a filter on $X$.
Conversely, suppose that $\mathcal{F}$ is a filter on $X$ with the
property that for each $B \subseteq X$, either $B \in \mathcal{F}$ or
$X \backslash B \in \mathcal{F}$.  Let us show that this implies that
$\mathcal{F}$ is an ultrafilter on $X$.  Let $B$ be a subset of $X$
such that $B \cap E \ne \emptyset$ for each $E \in \mathcal{F}$.  If
$X \backslash B \in \mathcal{F}$, then we get a contradiction, by applying
the previous condition to $E = X \backslash B$.  This implies that
$B \in \mathcal{F}$ in this situation, and hence that $\mathcal{F}$
is an ultrafilter on $X$, by the criterion in the previous paragraph.

        Now let $A$ be a commutative topological group.  A subset $K$
of $A$ is said to be \emph{totally bounded}\index{totally bounded sets}
in $A$ if for each open set $U$ in $A$ with $0 \in U$, $K$ is contained
in the union of finitely many translates of $U$.  If $K$ is compact,
then it is easy to see that $K$ is totally bounded, by covering $K$
by translates of $U$ and reducing to a finite subcovering by compactness.
If $d(x, y)$ is a translation-invariant metric on $A$ that determines
the same topology on $A$, then one can check that $K \subseteq A$ is
totally bounded as a subset of $A$ as a commutative topological group
if and only if $K$ is totally bounded as a subset of $A$ as a metric
space with respect to the metric $d(x, y)$.

        Suppose that $K \subseteq A$ is totally bounded, and that 
$\mathcal{F}$ is an ultrafilter on $A$ that contains $K$ as an element.
Let us check that $\mathcal{F}$ is a Cauchy filter on $A$ under these
conditions.  Let $U$ be an open set in $A$ with $0 \in U$, and let $V$
be an open set in $A$ such that $0 \in V$ and $V - V \subseteq U$.
Because $K$ is totally bounded in $A$, there are finitely many
elements $x_1, \ldots, x_n$ of $A$ such that
\begin{equation}
\label{K subseteq bigcup_{j = 1}^n (x_j + V)}
        K \subseteq \bigcup_{j = 1}^n (x_j + V).
\end{equation}
If $x_j + V \in \mathcal{F}$ for some $j$, then we get the desired
Cauchy condition, since
\begin{equation}
\label{(x_j + V) - (x_j + V) = V - V subseteq U}
        (x_j + V) - (x_j + V) = V - V \subseteq U.
\end{equation}
Otherwise, $A \backslash (x_j + V) \in \mathcal{F}$ for each $j = 1,
\ldots, n$, because $\mathcal{F}$ is an ultrafilter on $A$.  Of
course, this implies that $\bigcap_{j = 1}^n (A \backslash (x_j + V))
\in \mathcal{F}$, because $\mathcal{F}$ is a filter on $A$.
Equivalently, this means that $A \backslash \Big(\bigcup_{j = 1}^m
(x_j + V)\Big) \in \mathcal{F}$, which contradicts (\ref{K subseteq
  bigcup_{j = 1}^n (x_j + V)}) and the fact that $K \in \mathcal{F}$.

        Suppose now that $A$ is complete, and that $K \subseteq A$ is
both closed and totally bounded.  If $\mathcal{F}$ is an ultrafilter on
$A$ that contains $K$ as an element, then $\mathcal{F}$ is a Cauchy
filter on $A$, as in the preceding paragraph.  This implies that $\mathcal{F}$
converges to an element $p$ of $A$, because $A$ is complete, and that
$p \in K$, because $K$ is closed in $A$.  It follows that $K$ is a
compact set in $A$ under these conditions.

        Let $X$ and $Y$ be sets, and let $f$ be a mapping from $X$ into
$Y$.  If $\mathcal{F}$ is an ultrafilter on $X$, then the corresponding
filter $f_*(\mathcal{F})$ on $Y$ as in (\ref{f_*(mathcal{F}) = {E
    subseteq Y : f^{-1}(E)}}) is an ultrafilter on $Y$.  To see this,
it suffices to check that for each subset $B$ of $Y$, either $B \in
f_*(\mathcal{F})$ or $Y \backslash B \in f_*(\mathcal{F})$.  This is
the same as saying that $f^{-1}(B) \in \mathcal{F}$ or $f^{-1}(Y
\backslash B) \in \mathcal{F}$, because of the way that
$f_*(\mathcal{F})$ is defined.  Of course, $f^{-1}(Y \backslash B) = X
\backslash f^{-1}(B)$, and either $f^{-1}(B) \in \mathcal{F}$ or $X
\backslash f^{-1}(B) \in \mathcal{F}$, as desired, because $\mathcal{F}$
is an ultrafilter on $X$.

        Now let $I$ be a nonempty set, and suppose that $X_i$ is a 
topological space for each $i \in I$.  Also let $X = \prod_{i \in I}
X_i$ be the Cartesian product of the $X_i$'s, equipped with the product
topology.  If $K_i$ is a compact subset of $X_i$ for each $i \in I$,
then Tychonoff's theorem states that $K = \prod_{i \in I} K_i$ is a
compact set in $X$.  There is a nice proof of this using ultrafilters,
as follows.  It suffices to show that if $\mathcal{F}$ is an ultrafilter
on $X$ that contains $K$ as an element, then $\mathcal{F}$ converges to
a point in $K$.  Let $p_i$ be the obvious coordinate projection from $X$
onto $X_i$ for each $i \in I$, so that $(p_i)_*(\mathcal{F})$ is an
ultrafilter on $X_i$ for each $i \in I$, as in the previous paragraph.
It is easy to see that $K_i \in (p_i)_*(\mathcal{F})$ for each $i \in I$,
so that $(p_i)_*(\mathcal{F})$ converges to an element $x_i$ of $K_i$
for each $i \in I$, because $K_i$ is compact.  Using this, one can check that
$\mathcal{F}$ converges to the point $x \in K$ defined by $p_i(x) = x_i$ 
for each $i \in I$, as desired.

\section{Equicontinuity}
\label{equicontinuity, 2}

        Let $X$ be a topological space, and let $C(X)$ be the space of
continuous complex-valued functions o $X$, as usual.  A collection $E$
of complex-valued functions on $X$ is said to be
\emph{equicontinuous}\index{equicontinuity} at a point $p \in X$
if for each $\epsilon > 0$ there is an open set $U \subseteq X$
such that
\begin{equation}
\label{|f(p) - f(q)| < epsilon}
        |f(p) - f(q)| < \epsilon
\end{equation}
for every $q \in U$ and $f \in E$.  Of course, this implies that each
$f \in E$ is continuous at $p$.  Suppose now that $E \subseteq C(X)$
is equicontinuous at every point $p \in X$, and let $K$ be a nonempty
compact subset of $X$.

        Also let $\epsilon > 0$ be given, and for each $p \in K$,
let $U(p)$ be an open set in $X$ such that $p \in U(p)$ and
\begin{equation}
\label{|f(p) - f(q)| < epsilon / 3}
        |f(p) - f(q)| < \epsilon / 3
\end{equation}
for every $q \in U$ and $f \in E$.  Because $K$ is compact, there are
finitely many elements $p_1, \ldots, p_n$ of $K$ such that 
\begin{equation}
\label{K subseteq bigcup_{j = 1}^n U(p_j)}
        K \subseteq \bigcup_{j = 1}^n U(p_j).
\end{equation}
If $f, g \in E$ and $q \in U(p_j)$ for some $j = 1, \ldots, n$, then
it follows that
\begin{eqnarray}
\label{|f(q) - g(q)| le ... < |f(p_j) - g(p_j)| + 2 epsilon / 3}
  \quad |f(q) - g(q)| & \le & |f(q) - f(p_j)| + |f(p_j) - g(p_j)| 
                                                         + |g(p_j) - g(q)| \\
                & < & |f(p_j) - g(p_j)| + 2 \, \epsilon / 3.  \nonumber
\end{eqnarray}
This and (\ref{K subseteq bigcup_{j = 1}^n U(p_j)}) imply that
\begin{equation}
\label{|f(q) - g(q)| < max_{1 le j le n} |f(p_j) - g(p_j)| + 2 epsilon / 3}
 |f(q) - g(q)| < \max_{1 \le j \le n} |f(p_j) - g(p_j)| + 2 \, \epsilon / 3
\end{equation}
for every $f, g \in E$ and $q \in K$.

        Suppose in addition that $E$ is uniformly bounded pointwise
on $X$.  More precisely, this means that for each $p \in X$ there is
a nonnegative real number $C(p)$ such that
\begin{equation}
\label{|f(p)| le C(p)}
        |f(p)| \le C(p)
\end{equation}
for every $f \in E$.  If we consider
\begin{equation}
\label{(f(p_1), ldots, f(p_n))}
        (f(p_1), \ldots, f(p_n))
\end{equation}
as an element of ${\bf C}^n$ for each $f \in C(X)$, then it follows
that the set of (\ref{(f(p_1), ldots, f(p_n))}) with $f \in E$ is a
bounded subset of ${\bf C}^n$.  Since bounded subsets of ${\bf C}^n$
are totally bounded, we get that there are finitely many functions
$f_1, \ldots, f_r$ in $E$ with the property that for each $f \in E$
there is an $l = 1, \ldots, r$ such that
\begin{equation}
\label{max_{1 le j le n} |f(p_j) - f_l(p_j)| < epsilon / 3}
        \max_{1 \le j \le n} |f(p_j) - f_l(p_j)| < \epsilon / 3.
\end{equation}
Combining this with (\ref{|f(q) - g(q)| < max_{1 le j
    le n} |f(p_j) - g(p_j)| + 2 epsilon / 3}), we get that
\begin{equation}
\label{sup_{q in K} |f(q) - f_l(q)| < epsilon}
        \sup_{q \in K} |f(q) - f_l(q)| < \epsilon.
\end{equation}

        As before, $C(X)$ is a topological vector space with respect 
to the topology determined by the collection of supremum seminorms
associated to nonempty compact subsets of $X$.  In particular, $C(X)$
is a commutative topological group with respect to addition and this
topology.  The previous discussion shows that $E \subseteq C(X)$ is
totally bounded in $C(X)$ as a topological group with respect to addition
and this topology when $E$ is equicontinuous at every point in $X$ and 
uniformly bounded pointwise on $X$.

        If $E$ is also a closed set in $C(X)$, then $E$ is compact
with respect to this topology.  This follows from the discussion in
the previous section when $C(X)$ is complete as a commutative
topological group with respect to addition, which we have seen holds
when $X$ is locally compact, and when there is a countable local base
for the topology of $X$ at each point.  Otherwise, one can show
directly that $E \subseteq C(X)$ is complete in the sense that any
Cauchy net of elements of $E$ converges to an element of $E$ when $E$
is equicontinuous at each point in $X$ and closed in $C(X)$.  Indeed,
if a net of elements of $E$ converges pointwise to a function $f$ on $X$,
and if $E$ is equicontinuous at every point in $X$, then it is easy to
see that $f$ is continuous on $X$.  Using this, one can check that any
Cauchy net of element of $E$ converges to an element of $E$ with
respect to the usual topology on $C(X)$ when $E$ is equicontinuous 
at each point in $X$ and closed in $C(X)$, by essentially the same 
arguments as in the other cases of completeness.  As before,
this implies the analogous completeness condition in terms of Cauchy
filters, which can also be checked directly in a similar way.
Alternatively, one can embed $E$ in a Cartesian product, as in Section
\ref{equicontinuity}.

        In the other direction, if $E \subseteq C(X)$ is totally bounded,
then $E$ is uniformly bounded pointwise on $X$, as well as on compact subsets
of $X$.  If, in addition, $X$ is locally compact, then $E$ is equicontinuous
at every point in $X$.  Otherwise, the restrictions of the elements of $E$
to compact subsets of $X$ are equicontinuous at each point.

        Let $C(X, {\bf T})$ be the space of continuous mappings from $X$
into the unit circle ${\bf T}$, which is a closed subset of $C(X)$
with respect to the usual topology.  We have seen that $C(X, {\bf T})$
is also a commutative topological group with respect to pointwise
multiplication of functions and the topology on $C(X, {\bf T})$
induced by the one on $C(X)$.  As in Section \ref{continuous
  functions, 3}, a net of elements of $C(X, {\bf T})$ is a Cauchy net
of with respect to addition if and only if it is a Cauchy net with
respect to multiplication, because of (\ref{|z w^{-1} - 1| = |z -
  w|}).  Similarly, a set $E \subseteq C(X, {\bf T})$ is totally
bounded as a subset of $C(X, {\bf T})$ as a topological group with
respect to multiplication if and only if $E$ is totally bounded as a
subset of $C(X)$ as a topological group with respect to addition.

\backmatter

\newpage

\addcontentsline{toc}{chapter}{Index}

\printindex

\end{document}